%% file: thesis.tex
\documentclass{article}

\title{Euler Characteristics of Teichm\"uller Curves in Genus Two}
  \author{Matt Bainbridge}

\usepackage{amssymb} 
\usepackage{amsmath}
\usepackage{amsxtra}
\usepackage{amsfonts} 
\usepackage{amsthm}
\usepackage[all]{xy}
\usepackage[dvips]{graphicx}
\usepackage{color}
\usepackage{xspace}
\usepackage{ifthen}
\usepackage{array}

\input{preamble}

\numberwithin{equation}{section}

\newtheorem{theorem}{Theorem}[section] 
\newtheorem{prop}[theorem]{Proposition} 
\newtheorem{cor}[theorem]{Corollary}
\newtheorem{lemma}[theorem]{Lemma}

\theoremstyle{definition}
\newtheorem*{definition}{Definition}

\theoremstyle{remark}
\newtheorem*{remark}{Remark}

\begin{document}
\bibliographystyle{halpha} \maketitle
\tableofcontents

\pagebreak

\begin{abstract}
  For any integer $D\equiv 0$ or $1\pmod 4$, there is a Teichm\"uller
  curve $\W$ immersed in the moduli space $\moduli$ of genus two
  Riemann surfaces.  The curve $\W$ is naturally embedded in a Hilbert
  modular surface $\X$, Our main result is that the Euler
  characteristic of $\W$ is proportional to the Euler characteristic
  of $\X$.  More precisely,
  \begin{equation*}
    \chi(\W) = -\frac{9}{2} \chi(\X),
  \end{equation*}
  when $D$ is not square.
  
  When $D\equiv 1\pmod 8$, the curve $\W$ has two connected
  components, $\Wzero$ and $\Wone$.  We also calculate the Euler
  characteristics of these components.  When $D$ is not square, we
  show that
  \begin{equation*}
    \chi(\Wzero) = \chi(\Wone).
  \end{equation*}
  When $D$ is square, we show that
  \begin{align*}
    \chi(W_{d^2}^0)&=-\frac{1}{32} d^2 (d-1)\sum_{r | d}\frac{\mu(r)}{r^2}, \quad \text{and} \\
    \notag \chi(W_{d^2}^1)&=-\frac{1}{32} d^2 (d-3)\sum_{r |
      d}\frac{\mu(r)}{r^2}.
  \end{align*}
  
  The idea of the calculation of $\chi(\W)$ is to use techniques from
  algebraic geometry to compute the fundamental class of the closure
  $\W$ in a compactification of $\X$.  We define a compactification
  $\Y$ of $\X$ which maps to the Deligne-Mumford compactification of
  $\moduli$ by a finite morphism.  We then exhibit $\W$ as the zero
  locus of a meromorphic section of a line bundle over $\Y$, which
  allows us to calculate the fundamental class of $\W$.
  
  To calculate the Euler characteristics of the connected components
  $\We$, we find several relations involving the fundamental class of
  the closure of $\We$ which allow us to solve for $\chi(\We)$.  For
  example, we calculate the self-intersection numbers of the $\We$ in
  terms of $\chi(\We)$.

  We apply these results to calculate the Siegel-Veech constants for
  counting closed billiards paths in certain L-shaped polygons.  We
  then calculate the Lyapunov exponents of the Kontsevich-Zorich cocycle
  for any ergodic, $\SLtwoR$-invariant measure on the moduli space
  $\Omega_1\moduli$ of Abelian differentials in genus two (previously
  calculated in unpublished work of Kontsevich and Zorich).
\end{abstract}

\include{introduction}

\include{abeliansurfaces}

\include{prototypes}

\include{abeliandifferentials}

\include{delignemumford}

\include{localcoordinates}

\include{limits}

\include{geometric}

\include{bundles}

\include{euler}

\include{fundamental}

\include{normalbundles}

\include{euler2}

\include{siegelveech}

\include{lyapunov}

\appendix

\include{normalization}

\bibliography{my}

\end{document}

%% file: preamble.tex
\newcommand{\A}[1][D]{\mathcal{A}_{#1}}
\newcommand{\augteich}{\overline{\teich}}

\newcommand{\bA}[1][d^2]{\widehat{\mathcal{A}}_{#1}} 
\newcommand{\bX}[1][D]{\widehat{X}_{#1}} 
\newcommand{\barmoduli}[1][2]{{\overline{\mathcal M}}_{#1}}
 
\newcommand{\barP}[1][D]{\overline{P}_{#1}}

\newcommand{\barRtwo}[1][d^2]{\overline{R}^2_{#1}}
\newcommand{\barSone}[1][d^2]{\overline{S}^1_{#1}}
\newcommand{\barStwo}[1][d^2]{\overline{S}^2_{#1}}
\newcommand{\barSi}[1][d^2]{\overline{S}^i_{#1}} 
\newcommand{\barW}[1][D]{\overline{W}_{#1}}
\newcommand{\barWe}[1][D]{\overline{W}^\epsilon_{#1}}
\newcommand{\barWone}[1][D]{\overline{W}^1_{#1}}
\newcommand{\barWzero}[1][D]{\overline{W}^0_{#1}}
\newcommand{\barX}[1][D]{\overline{X}_{#1}}

\newcommand{\bdry}{\partial}

\newcommand{\system}[2][{}]{
  T_{
    \ifthenelse{\equal{#2}{1}}{0,1}{
      \ifthenelse{\equal{#2}{2}}{1,0}{
        \ifthenelse{\equal{#2}{3}}{2,0}{
          \ifthenelse{\equal{#2}{4}}{2,1}{
            \ifthenelse{\equal{#2}{5}}{3,0}{1,1}
          }
        }
      }
    }
  }^{#1}
}

\newcommand{\Def}{\mathcal{D}}

\newcommand{\dX}[1][D]{\partial X_{#1}}
\newcommand{\E}[1][D]{E_{#1}}

\newcommand{\F}[1][D]{\mathcal{F}_{#1}}
\newcommand{\GLtwoR}{{\rm GL}_2\reals}
\newcommand{\GLtwoRplus}{{\rm GL}_2^+\reals}
\newcommand{\fhoriz}{\mathcal{F}_h}
\newcommand{\fmod}[1][2]{\mathcal{FM}_{#1}}

\renewcommand{\Im}{\IM}
\newcommand{\isom}{\cong}
\newcommand{\kron}[2]{\left(\frac{#1}{#2}\right)}
\newcommand{\moduli}[1][2]{{\mathcal M}_{#1}}

\newcommand{\omn}[1][g]{\Omega\moduli[#1]({\bf n})}

\newcommand{\ord}[1][D]{\mathcal{O}_{#1}}
\newcommand{\otn}[1][g]{\Omega\teich_{#1}({\bf n})}
\renewcommand{\P}[1][D]{P_{#1}}
\newcommand{\partialbar}{\overline{\partial}}
\newcommand{\Pprot}[1][D]{\mathcal{P}_{#1}}

\newcommand{\PSL}{{\rm PSL}}

\newcommand{\PSLtwoZ}{\PSL_2 \zed}
\newcommand{\Q}[1][D]{Q_{#1}}
\newcommand{\Rone}[1][d^2]{R^1_{#1}}
\newcommand{\Rtwo}[1][d^2]{R^2_{#1}}
\newcommand{\Ri}[1][d^2]{R^i_{#1}} 
\newcommand{\Sone}[1][d^2]{S^1_{#1}}
\newcommand{\Stwo}[1][d^2]{S^2_{#1}}
\newcommand{\Si}[1][d^2]{S^i_{#1}} 
\newcommand{\satsiegelmod}[1][2]{\widehat{\mathcal{A}}_{#1}}
\newcommand{\siegelmod}[1][2]{\mathcal{A}_{#1}} 
\newcommand{\siegelhalf}[1][2]{\mathcal{H}_{#1}}
\newcommand{\semisiegelhalf}[1][2]{\mathcal{H}_{#1}^\dagger}
\newcommand{\SL}{{\mathrm{SL}}}

\newcommand{\SLtwoord}[1][D]{\SL(\ord[#1]\oplus\ord[#1]^\vee)}
\newcommand{\SLtwoR}{\SL_2 \reals} 
\newcommand{\SLtwoZ}{\SL_2 \zed}
\newcommand{\SP}{{\mathrm{Sp}}}
\newcommand{\SPQ}[1]{\SP_{#1} \ratls} 
\newcommand{\SPR}[1]{\SP_{#1} \reals}
\newcommand{\SPZ}[1]{\SP_{#1} \zed}
\newcommand{\struct}{\mathcal{O}}
\newcommand{\Sum}{\sum_{r|d}\frac{\mu(r)}{r^2}}
\newcommand{\teich}{\mathcal{T}}
\newcommand{\teichs}[1][g, n]{\mathcal{T}(\Sigma_{#1})}

\newcommand{\tildemoduli}[1][2]{\widetilde{\mathcal{M}_{#1}}}
\newcommand{\W}[1][D]{W_{#1}}
\newcommand{\We}[1][D]{W^\epsilon_{#1}}
\newcommand{\Wprot}[1][D]{\mathcal{W}_{#1}}
\newcommand{\Wzero}[1][D]{W_{#1}^0}
\newcommand{\Wone}[1][D]{W_{#1}^1} 
\newcommand{\X}[1][D]{X_{#1}} 
\newcommand{\Y}[1][D]{Y_{#1}}
\newcommand{\Yprot}[1][D]{\mathcal{Y}_{#1}}

\newcommand{\cx}{{\mathbb C}} 
\newcommand{\half}{{\mathbb H}} 
\newcommand{\integers}{{\mathbb Z}}
\newcommand{\nats}{{\mathbb N}}
\newcommand{\ratls}{{\mathbb Q}} 
\newcommand{\reals}{{\mathbb R}} 
\newcommand{\proj}{{\mathbb P}}
\newcommand{\zed}{\integers}

\DeclareMathOperator{\Area}{Area}
\DeclareMathOperator{\const}{const}
\DeclareMathOperator{\Hom}{Hom}
\DeclareMathOperator{\IM}{Im}
\DeclareMathOperator{\Ann}{Ann}
\DeclareMathOperator{\Aut}{Aut}
\DeclareMathOperator{\diameter}{diameter}
\DeclareMathOperator{\Div}{Div}
\DeclareMathOperator{\End}{End}

\DeclareMathOperator{\height}{height}
\DeclareMathOperator{\Jac}{Jac}
\DeclareMathOperator{\Ker}{Ker}

\DeclareMathOperator{\Mod}{Mod}

\DeclareMathOperator{\mult}{mult}
\DeclareMathOperator{\Pic}{Pic}

\DeclareMathOperator{\Res}{Res}
\DeclareMathOperator{\sign}{sign}
\DeclareMathOperator{\Stab}{Stab}
\DeclareMathOperator{\tr}{Tr}
\DeclareMathOperator{\Twist}{Tw}
\DeclareMathOperator{\vol}{vol}


%% file: introduction.tex
\section{Introduction}
\label{sec:intro}

Let $\moduli$ be the moduli space of genus two Riemann surfaces.  Given any $D>0$
with $D\equiv 0$ or $1\pmod 4$, let $\W'\subset\moduli$ be the locus of Riemann surfaces $X$ such that:
\begin{itemize}
\item The Jacobian $\Jac(X)$ has real multiplication by $\ord$, the unique real quadratic order of
  discriminant $D$, and
\item There is an Abelian differential $\omega$ on $X$ which is an eigenform for real multiplication and has a
  double zero.
\end{itemize}
Let $\W$ be the normalization of $\W'$, a possibly disconnected curve.

The \emph{Hilbert modular surface},
$$\X = \half\times\half / \SLtwoord,$$
is the moduli space for Abelian surfaces with real multiplication by
$\ord$.  In $\X$, let $\P$ be the Shimura curve consisting of those Abelian surfaces in $\X$ which are
polarized products of elliptic curves.  There is a commutative diagram,
$$\xymatrix{
  \W \ar@{^{(}->}[r] \ar[dr] & \X\setminus\P \ar[d] \\
  & \moduli}
$$
where the top map is an embedding,
and the vertical map is a two-to-one map sending an Abelian surface $A$ to the unique Riemann surface
$X\in\moduli$ such that $\Jac(X)\isom A$.

The curve $\W$ is not a Shimura curve on $\X$; it is, however, a \emph{Teichm\"uller curve}, a curve which is
isometrically immersed in the moduli space $\moduli$ with respect to the Teichm\"uller metric.  In fact, the
curves $\W$ with $D$ nonsquare are all but one of the \emph{primitive} Teichm\"uller curves in $\moduli$,
where a Teichm\"uller curve is said to be primitive if it does not arise from a Teichm\"uller curve of lower
genus by a certain branched covering construction.  The curves $\W$ arise from the study of billiards in
certain $L$-shaped polygons, and the study of these curves has applications to the dynamics of billiards in
these polygons.

\paragraph{Euler characteristics.}

The main object of this paper is to calculate the Euler characteristics of the curves $\W$.  The idea is to
relate $\chi(\W)$ to $\chi(\X)$ and $\chi(\P)$.  It turns out that these Euler characteristics are all
proportional if $D$ is not square.  The following is our main result.

\begin{theorem}
  \label{thm:eulerW}
  If $D$ is not square, then
  \begin{equation}
    \label{eq:WD}
    \chi(\W)= -\frac{9}{2} \chi(\X).
  \end{equation}
\end{theorem}

\begin{remark}
  Zagier conjectured formula \eqref{eq:WD} from numerical evidence.  In \cite{siegel36}, Siegel calculated the
  volume of $\X$ when $D$ is a fundamental discriminant (that is $D$ is not of the form $D=f^2E$ for some
  $E\equiv 0$ or $1\pmod 4$ and $f>1$).  This yields formulas for $\chi(\X)$ for all $D$, and
  for $D\neq1$ a fundamental discriminant, we obtain 
  \begin{equation*}
    \chi(\W[f^2D])= -9 \zeta_{\ratls(\sqrt{D})}(-1) f^3\sum_{r|f}\kron{D}{r}\frac{\mu(r)}{r^2},
  \end{equation*}
  where $\mu$ is the M\"obius function, and $\kron{D}{r}$ is the Kronecker symbol.  When $f=1$, this reduces
  to
  \begin{equation*}
    \chi(\W) = -9\zeta_{\ratls(\sqrt{D})}(-1).
  \end{equation*}

  The Euler characteristics $\chi(\W)$ when $D$ is not square are given by the Fourier coefficients of a
  modular form.  More precisely, there is a function $H(2, D)$ such that for $D\neq 1$ a fundamental
  discriminant and $f\in\nats$,
  $$\sum_{r|f}\chi(\W[r^2 D]) = \frac{3}{4}H(2, f^2 D), \quad \text{and}$$
  $$\sum_{\substack
    {D\equiv 0, 1 \,(4)\\
      D\geq 0
    }}
  H(2, D) q^D$$
  is a modular form studied by Cohen in \cite{cohen75}.  See \S\ref{subsec:hilbertmodular} for more about this
  form.
\end{remark}

When $D=d^2$, Eskin, Masur, and Schmoll calculated $\chi(\W[d^2])$ in \cite{ems}.  Our methods give a
new proof of their formula.  In this case, it is no longer true that $\chi(\W[d^2])$ is proportional to
$\chi(\X[d^2])$.

\begin{theorem}
  For any $d>1$,
  \begin{equation}
    \label{eq:ems}
    \chi(W_{d^2})=-\frac{1}{16} d^2(d-2)\sum_{r|d}\frac{\mu(r)}{r^2}.
  \end{equation}
\end{theorem}

\paragraph{Connected components of $\W$.}

It is known that when $D\neq9$ and $D\equiv 1 \pmod 8$, $\W$ has two connected components, $\Wzero$ and
$\Wone$; otherwise $\W$ is connected \cite{mcmullenspin}.

We will also show that when $D$ is not square, the connected components of $\W$ have the same Euler
characteristic.  It should be true that these components are actually homeomorphic, but we have not proved
this.  It is not true that the connected components of $\W$ are isomorphic as curves over $\cx$.

\begin{theorem}
  \label{thm:chiequal}
  If $D\equiv 1\pmod 8$ is not square, then
  \begin{equation*}
    \chi(\Wzero) = \chi(\Wone).
  \end{equation*}
\end{theorem}

There is a canonical involution $\tau\colon\X\to\X$ which is covered by the involution $(z, w)\mapsto(w, z)$
of $\half\times \half$.  This involution replaces an Abelian surface $A$ together with a choice of real multiplication,

$$\rho\colon\ord\to\End A,$$
with the Galois conjugate real multiplication $\rho'$ obtained by precomposing
$\rho$ with the Galois automorphism of $\ord$.  It is not true that $\tau$ permutes the connected components
of $\W$ (otherwise Theorem~\ref{thm:chiequal} would be trivial).  In fact, $\W$ is not even invariant under $\tau$.

We also calculate $\chi(\We[d^2])$, but it is no longer true in the square discriminant case that the
connected components of $\W[d^2]$ have the same Euler characteristic.

\begin{theorem}
  \label{thm:hlformula}
  If $d>0$ and $d^2\equiv 1 \pmod 8$, then
  \begin{align}
    \label{eq:hlformula}
    \chi(W_{d^2}^0)&=-\frac{1}{32} d^2 (d-1)\sum_{r | d}\frac{\mu(r)}{r^2}, \quad \text{and} \\
    \notag
    \chi(W_{d^2}^1)&=-\frac{1}{32} d^2 (d-3)\sum_{r | d}\frac{\mu(r)}{r^2}.
   \end{align}
\end{theorem}

\begin{remark}
  These formulas were established independently by Leli\`evre and Royer in \cite{lelievreroyer}.
\end{remark}

\paragraph{Siegel-Veech constants.}

Given a positive, nonsquare integer $D\equiv 0$ or $1 \pmod 4$, let
$P(D)$ be the L-shaped polygon as in Figure~\ref{fig:lshaped}, where the side
lengths are given by
$$ a = b = \frac{1 + \sqrt{D}}{2}$$
if $D$ is odd, and
$$a = \frac{\sqrt{D}}{2} \quad\text{and}\quad b = 1 +
\frac{\sqrt{D}}{2}$$
if $D$ is even.  By a standard unfolding
construction, a $L$-shaped polygon determines a genus two Riemann
surface $X$ equipped with an Abelian differential $\omega$, and the
$(X, \omega)$ determined by $P(D)$ lies on $\Omega\W$.

\begin{figure}[htbp]
  \centering
  \input{lshaped.pstex_t}
  \caption{$P(D)$}
  \label{fig:lshaped}
\end{figure}

A billiards path on $P(D)$ is a path which is geodesic on the interior
of $P(D)$ and bounces off the walls as a physical billiards ball would
(angle of incidence equals angle of reflection).  We allow billiards
paths to pass through the corners of angle $\pi/2$ but not through the
one of angle $3\pi/2$.   Closed billiards paths occur in parallel
families which are bounded by two paths which each start and end at
the corner of angle $3\pi/2$.  Let
\begin{multline*}
  N(P(D), L) = \#\{\text{families of simple closed billiards paths on
    $P(D)$} \\ \text{of length at
    most $L$}\}.
\end{multline*}
Note that we do not count paths which go around a single closed path
several times.  Also paths are not oriented, so we do not count a path
going once in each direction.  Since the unfolding lies on a
Teichm\"uller curve, we have by \cite{veech89}
$$N(P(D), L)\sim c(D) \frac{\pi}{4 \Area(P(D))} L^2.$$
for some constants $c(D)$ called Siegel-Veech constants.  Evaluating the constants
$c(D)$ amounts to calculating the Euler characteristics of the $\W$ together with
the volumes of certain neighborhoods of the cusps of $\W$.  In
\S\ref{sec:siegelveech}, we calculate these constants.  For example,
we obtain:
\begin{theorem}
  For small values of $D$, the Siegel-Veech constants $c(D)$ are as given
  in Table~\ref{tab:siegelveech}.
\end{theorem}

\begin{table}[hbt]
  \centering
    \begin{tabular}{| l | c c c c c c c c c c|}
    \hline
    $D$ & $5$ & $8$ & $12$ & $13$ & $17$ & $20$ & $21$ & $24$ & $28$ &
    $29$\\
    \hline
     $c(D)$ & $\displaystyle \frac{\vphantom{\int}25}{\vphantom{\int}3}$ & $\displaystyle \frac{28}{3}$ & $\displaystyle \frac{26}{3}$ &
    $\displaystyle \frac{91}{9}$ & $\displaystyle \frac{221}{24} + \frac{1}{8} \displaystyle \sqrt{17}$ &
    $\displaystyle \frac{31}{3}$ & $\displaystyle \frac{133}{15}$ &
    $\displaystyle \frac{148}{15}$ & $\displaystyle \frac{82}{9}$ &
    $\displaystyle \frac{377}{35}$ \\
    \hline
  \end{tabular}
  \caption{Siegel-Veech constants for $P(D)$}
  \label{tab:siegelveech}
\end{table}

Note that there are other L-shaped polygons, classified in
\cite{mcmullenbild}, whose unfoldings lie on the Teichm\"uller curves
$\W$ and our results apply equally well to those.  The Siegel-Veech
constants only depend on the connected component of $\W$ on which the
unfolding lies.  We also show in \S\ref{sec:siegelveech}:

\begin{theorem}
  If $\W$ has two connected components,
  then the associated Siegel-Veech constants are Galois-conjugate
  elements of $\ratls(\sqrt{D})$.  If $\W$ is connected, then the
  Siegel-Veech constants are rational.
\end{theorem}

\paragraph{The Kontsevich-Zorich cocycle.}

There is a well-known action of $\SLtwoR$ on $\Omega\moduli[g]$, the moduli space of Abelian differentials,
which preserves the subspace $\Omega_1\moduli[g]$ of Abelian differentials of norm one with respect to the
norm,
$$\|\omega\| = \left(\int_X |\omega|^2\right)^{1/2},$$
on $\Omega(X)$.
The diagonal subgroup $A\subset\SLtwoR$ induces a flow $g_t$ on
$\Omega\moduli[g]$, the \emph{Teichm\"uller geodesic flow}. Kontsevich
and Zorich \cite{kontsevich, zorich99} introduced a linear cocycle
over this flow which is closely related to the dynamics of flows on
surfaces.  There is a bundle
$\mathcal{H}_1(\reals)\to\Omega_1\moduli[g]$ whose fiber over an
Abelian differential $(X, \omega)$ is $H_1(X; \reals)/\Aut(X,
\omega)$.  This bundle has a natural flat connection, the
\emph{Gauss-Manin connection}.  By coupling the geodesic flow $g_t$
with the Gauss-Manin connection, we get a flow $\tilde{g}_t$, the
\emph{Kontsevich-Zorich cocycle}, on $\mathcal{H}_1(\reals)$ which is
linear on the fibers and covers the flow $g_t$ on
$\Omega_1\moduli[g]$.

Consider the closure $S\subset\Omega_1\moduli$ of an $\SLtwoR$ orbit.
By results of \cite{mcmullenabel}, $S$ is a suborbifold of
$\Omega_1\moduli$ which comes equipped with a unique $\SLtwoR$-invariant,
ergodic, absolutely continuous, probability measure, the \emph{period
  measure}.  In the case of the Teichm\"uller curves, period measure
is just the normalized hyperbolic area measure (because Haar measure
is unique).

Associated to any ergodic, $g_t$-invariant measure $\mu$ are the $2g$ \emph{Lyapunov exponents},
$$1 = \lambda_1(\mu) > \dots> \lambda_g(\mu) > -\lambda_g(\mu) > \dots
>- \lambda_1(\mu).$$
Kontsevich
gave a formula in \cite{kontsevich} which calculates $\sum\lambda_i(\mu)$ as a certain integral over
moduli space when $\mu$ is $\SLtwoR$ invariant.  The same methods which we use to calculate $\chi(W_D)$ allow
us to evaluate these integrals and to calculate $\lambda_2(\mu)$ for
all ergodic, $\SLtwoR$-invariant measures on
$\Omega_1\moduli$.  Let $\Omega_1\moduli(2)\subset\Omega_1\moduli$
be the locus of forms with a double zero.

\begin{theorem}
  \label{thm:kz}
  If $\mu$ is any finite, ergodic, $\SLtwoR$-invariant measure on
  $\Omega_1\moduli$, then
  $$\lambda_2(\mu)=
  \begin{cases}
    1/3, & \text{if $\mu$ is supported on $\Omega_1\moduli(2)$};\\
    1/2, & \text{if $\mu$ is not supported on $\Omega_1\moduli(2)$.}
  \end{cases}
  $$
\end{theorem}
\begin{remark}
  Theorem~\ref{thm:kz} is an unpublished result of Kontsevich and Zorich which is mentioned in \cite{zorich06}.
\end{remark}

\paragraph{Outline of proof of Theorem \ref{thm:eulerW}.}

The strategy for computing the Euler characteristic of $\W$ is to study the relationship between $\W$ and the
Hilbert modular surface $\X$ on which it lies.  More precisely, we will define a compactification $\Y$ of $\X$
and compute the fundamental class of $\barW$ in $\Y$ by expressing it as the zero locus of a meromorphic
section of a line bundle over $\Y$.

Here is a sketch of the calculation of the Euler characteristic $\chi(\W)$ in the case when $D$ is nonsquare.
When $D$ is a square, the calculation is more complicated because there are some extra curves in
$\Y[d^2]\setminus\X[d^2]$ which have to be taken into account.

\begin{enumerate}
\item Given a nonsquare $D\in \nats$ with $D\equiv 0$ or $1\pmod 4$, let
  $$\ord = \zed[T]/(T^2 + bT+c)$$
  for integers $b, c\in\zed$ such that $b^2-4c = D$.  The discriminant
  $D$ determines $\ord$ up to isomorphism, and $\ord$ is naturally embedded in its field of fractions $K_D$,
  which is isomorphic to $\ratls(\sqrt{D})$.  Let $\ord^\vee\subset K_D$ be the inverse different of $\ord$, and
  let $\SLtwoord\subset{\rm SL}_2 K_D$ be the subgroup which preserves $\ord\oplus\ord^\vee$.
  
  There are two closed two-forms $\omega_i$ on
  $$\X=\half\times\half/\SLtwoord$$
  covered by the forms
  \begin{equation}
    \label{eq:omega}
    \frac{1}{2\pi}\frac{dx_i\wedge dy_i}{y_i^2}.
  \end{equation}
  on $\half\times\half$.  The inverse image of $\W$ in $\half\times\half$ is the union of the graphs of
  countably many holomorphic functions $\half\to\half$.  It follows that the restriction of $\omega_1$ to $\W$
  is $1/2\pi$ times the hyperbolic volume form on $\W$.  Gauss-Bonnet then tells us that
  $$\chi(\W)=-\int_{\W}\omega_1.$$
\item  There is a Shimura curve $\P\subset\X$ which parameterizes Abelian surfaces which are polarized products
  of elliptic curves.  It is covered in $\tilde{X}_D$ by a union of graphs of M\"obius transformations
  $\half\to\half$.  We introduce these curves in \S\ref{subsec:hilbertmodular} and show that
  \begin{equation}
    \label{eq:chiPD}
    \chi(\P) = -\frac{5}{2}\chi(\X).
  \end{equation}
\item The bundle $\Omega\moduli$ over $\moduli$ whose fiber over a
  Riemann surface $X$ is the space $\Omega(X)\setminus\{0\}$ of
  nonzero Abelian differentials on $X$ extends to a bundle
  $\Omega\barmoduli$ over the Deligne-Mumford compactification which
  we discuss in \S\ref{sec:moduli}.  In \S\ref{sec:limitsofeigenforms}
  and \S\ref{sec:geometric}, we construct a compactification $\Y$ of
  $\X$, a complex projective orbifold which is obtained by taking the
  normalization of the closure of the image of an embedding
  $\X\to\proj\Omega\barmoduli$.  We call this the \emph{geometric
    compactification} of $\X$.  The complement $\dX=\Y\setminus\X$
  consists of finitely many rational curves, each labeled by a
  discrete invariant which we call a $\Y$-prototype.  We discuss these
  invariants in \S\ref{sec:prototypes}.  For each $\Y$-prototype $P$,
  there is one rational curve $C_P\subset\dX$.
\item Associated to each point $p\in\Y$, there is a stable Riemann surface $X\in\barmoduli$ together with an
  action of $\ord$ on $\Jac(X)$ by real multiplication.  This determines a splitting of the space $\Omega(X)$
  of stable Abelian differentials on $X$ into two eigenspaces for real multiplication,
  $$\Omega(X) = \Omega^1(X)\oplus \Omega^2(X).$$
  (See \S\ref{sec:moduli} for information about stable Riemann
  surfaces and stable Abelian differentials.)  A choice of an embedding $\ord\to\reals$ determines one of
  these two eigenspaces, so the point $p\in\Y$ determines an eigenform in $\Omega(X)$
  up to scale if we fix such an embedding.   There is a holomorphic line bundle $\Omega\Y$ whose fiber over
  $p\in\Y$ is the eigenspace $\Omega^1(X)$.

  In \S\ref{sec:bundles}, we study a holomorphic foliation $\A$ of
  $\Y$ along whose leaves the absolute periods of eigenforms in $\Omega\Y$ are locally constant.  Each
  rational curve $C_P\subset\dX$ is a leaf of $\A$, and $\A$ has isolated singularities at the intersection of
  two curves $C_P$ and $C_Q$.  A leaf of $\A$ can be locally parameterized away from $\W$ and $\P$ by the
  \emph{relative periods} of $(X, \omega)$, the integrals along paths joining the two zeros $\omega$.
  Following \cite{mcmullenhilbert}, we use this parameterization to define a quadratic differential on each
  leaf of $\A$.  These quadratic differentials can be pieced together to define a section $q$ of the line
  bundle
  $$\mathcal{L}=(\Omega\Y)^{-2}\otimes (T^*\A)^2$$
  over $\Y$.
  This section $q$ vanishes along $\barW$, has a simple pole along $\barP$, and
  is elsewhere finite and nonzero.  This gives a relation between the fundamental classes of $\barW$, $\barP$,
  and the Chern class of $\mathcal{L}$, resulting in the formula,
  \begin{equation}
    \label{eq:barWclass}
    [\overline{W}_D] = [\barP] + c_1(\mathcal{L}).
  \end{equation}
  in $H^2(\Y;\ratls)$.  We derive this formula in \S\ref{sec:modular}.  This is the main idea of the
  calculation of $\chi(\W)$, and the reader who only wants the gist of the argument might find it useful to skip directly
  to this section.
\item In \S\ref{sec:bundles}, we calculate $c_1(\mathcal{L})$.  We define a Hermitian metric $h$ on
  $\Omega\Y$ which is singular along $\dX$.  We show that $h$ is a \emph{good metric} in the sense of Mumford
  \cite{mumford77}.  This implies that the Chern form $c_1(\Omega\Y, h)=\omega_1/2$ is a closed current on $\Y$,
  and
  \begin{equation}
    \label{eq:c1omegaY}
    c_1(\Omega\Y) = \frac{1}{2}[\omega_1]
  \end{equation}
  in $H^2(\Y; \ratls)$.

  The canonical involution $\tau$ of $\X$ extends to an involution of $\Y$, which we continue to
  call $\tau$.  We also show that, 
  $$\tau^*(\Omega\Y)^2 = T^*\A,$$
  which implies,
  \begin{equation}
    \label{eq:c1tstarA}
    c_1(T^*\A) = [\omega_2].
  \end{equation}
  Putting together \eqref{eq:c1omegaY} and \eqref{eq:c1tstarA}, we obtain
  \begin{equation}
    \label{eq:c1L}
    c_1(\mathcal{L}) = -[\omega_1] + 2[\omega_2].
  \end{equation}
\item Combining \eqref{eq:barWclass} and \eqref{eq:c1L}, we obtain
  \begin{equation}
    \label{eq:barWclass2}
    [\barW] = [\barP] - [\omega_1] + 2[\omega_2].
  \end{equation}
  By pairing $-[\omega_1]$ with both sides and applying the Gauss-Bonnet Theorem, we obtain
  \begin{equation}
    \label{eq:chirelation}
    \chi(\W) = \chi(\P) - 2 \chi(\X).
  \end{equation}
  Putting \eqref{eq:chiPD} in \eqref{eq:chirelation}, we obtain \eqref{eq:WD}.
\end{enumerate}

The proof of Theorem~\ref{thm:eulerW} uses very little about the compactification $\Y$.  Equation
\eqref{eq:barWclass2} is true, when interpreted as an equation of cohomology classes of closed currents on
$\X$; however, we are not allowed to pair $\omega_1$ with this equation because $\omega_1$ is not compactly
supported.  Working in the compactification $\Y$ allows us to justify this paring.  The proof of
Theorem~\ref{thm:chiequal} uses the compactification $\Y$ in a more essential way.

When $D=d^2$, there are some extra curves $\Sone$ and $\Stwo$ in $\Y[d^2]\setminus\X[d^2]$, which
complicates the calculation of $\chi(\W[d^2])$.  One difference from the square case is that the section $q$
of $\mathcal{L}$ vanishes on $\barStwo$ as well as $\barP[d^2]$.  Also, \eqref{eq:c1tstarA} becomes 
\begin{equation*}
  c_1(T^*\A) = [\omega_2] -[\barStwo]. 
\end{equation*}
Instead of \eqref{eq:chirelation}, we get
\begin{equation*}
  \chi(\W[d^2]) = \chi(\P[d^2]) - 2 \chi(\X[d^2]) -\chi(\Stwo).
\end{equation*}
In this case, $\chi(\P[d^2])$ and $\chi(\Stwo)$ are not proportional to $\chi(\X[d^2])$, so neither is
$\chi(\W[d^2])$.

\paragraph{Outline of proof of Theorem~\ref{thm:chiequal}.}

The idea of the proof of Theorem~\ref{thm:chiequal} is to find as many equations as we can involving the
fundamental classes $[\barWzero]$ and $[\barWone]$.  Once we have enough equations, we are able to solve for the
pairing $[\omega_1]\cdot[\barWe]=-\chi(\We)$.

\begin{enumerate}
\item In \S\ref{sec:limitsofeigenforms} and \S\ref{sec:geometric}, we will prove the following properties of
  $\Y$:
  \begin{enumerate}
  \item The closures $\barW$ and $\barP$ are disjoint suborbifolds of $\Y$.  \label{item:a}
  \item Each rational curve $C_P\subset\Y\setminus\X$ meets $\barW$ and $\barP$ transversely and meets each in
    the same number of points. \label{item:b}
  \item The canonical involution $\tau$ of $\Y$ has the property that for each rational curve $C_P\subset\dX$,
    $$\barWzero\cdot C_P = \barWone\cdot\tau(C_P).$$ \label{item:c}
  \item The rational cohomology of $\Y$ is an orthogonal direct sum:
    $$H^2(\Y;\ratls)\isom B \oplus \langle [\omega_1], [\omega_2] \rangle \oplus J,$$
    where $B$ is the subspace generated by the rational curves $C\subset \dX$, and $J$ is the orthogonal
    complement of the other two terms.  $J$ contains all of $H^{2,0}(\Y;\ratls)$ and $H^{0, 2}(\Y;\ratls)$ and most of
    $H^{1,1}(\Y;\ratls)$.
  \end{enumerate}
  \item There is a well-known formula relating the first Chern class of the bundle $\Omega\barmoduli$
    over $\barmoduli$ to the fundamental classes of divisors in $\barmoduli\setminus\moduli$.  Pulling back
    this formula by the natural map $\Y\to\barmoduli$, we obtain a formula for the fundamental class of $[\barP]$:
  \begin{equation}
    \label{eq:barPclass}
    [\barP] = \frac{5}{2}[\omega_1] + \frac{5}{2}[\omega_2]+ \pi_B[\barP],
  \end{equation}
  where $\pi_B$ is the orthogonal projection of $H^2(\Y;\ratls)$ to $B$.  Combining this with
  \eqref{eq:barWclass2}, we obtain
  \begin{equation*}
    [\barW] = \frac{3}{2}[\omega_1] + \frac{9}{2}[\omega_2] + \pi_B[\barW].
  \end{equation*}
  It follows from \eqref{eq:barWclass2} that $\pi_B[\barP] = \pi_B[\barW]$.  We derive these formulas in
  \S\ref{sec:fundamental}.
\item In \S\ref{sec:normalbundles}, we study the normal bundles of $\barW$ and $\barP$ in $\Y$.  The
  restriction of $\Omega\Y$ to $\barWe$ is a line bundle $\Omega\barWe$.  The degree of normal bundle $N(\barWe)$ can be
  calculated in terms of $\Omega\barWe$.  We obtain the following formula for the self intersection number of $\barWe$:
  \begin{equation}
    \label{eq:normalWe}
    [\barWe]^2 = \deg N(\barWe) =  -\frac{2}{3}\deg\Omega\barWe = \frac{1}{3}\chi(\We).
  \end{equation}
  The proof of the middle equality involves defining a three-to-one map from a tubular neighborhood $\barWe$ to $(\Omega\barWe)^{-2}$
  via an operation called ``collapsing a saddle connection,'' which we introduce in \S
  \ref{subsec:flatgeometry}.  This operation replaces an Abelian differential together with  a
  saddle connection joining a pair of simple zeros with an Abelian differential with a double zero.
  
  Since $\barWe$ is transverse to the foliation $\A$ of $\Y$, we also have the relation,
  $$N(\barWe)\isom T\A,$$
  which together with \eqref{eq:c1tstarA} implies
  \begin{equation}
    \label{eq:selfintersection}
    [\barWe]^2=\int_{\We}\omega_2.
  \end{equation}
  From \eqref{eq:normalWe} and \eqref{eq:selfintersection}, we obtain
  \begin{equation}
    \label{eq:omegairelation}
    \int_{\We}\omega_2 = \frac{1}{3}\int_{\We}\omega_1.
  \end{equation}
\item We show in \S\ref{sec:eulerWe} that the fundamental classes  $[\barWe]$ are given by
  \begin{equation}
    \label{eq:barWeclass}
    [\barWe]= \frac{3}{4}[\omega_1] + \frac{9}{4}[\omega_2] + \pi_B[\We]\pm j,\\
  \end{equation}
  for some $j\in J$.  The plan is to leave the coefficient of $[\omega_1]$ as an unknown $a$ and then to use
  our knowledge of $\Y$ to write down equations involving the class $[\barWe]$ that allow us to solve for $a$.
  The following equations follow from Properties (b) and (c) of $\Y$:
  \begin{gather*}
    [\barWzero]\cdot[\barP]=0,\\
    [\barWzero]\cdot[\barWzero]=[\barWone]\cdot[\barWone],
  \end{gather*}
  These together with \eqref{eq:omegairelation} give us enough equations to solve for $a$.
  Theorem~\ref{thm:chiequal} follows from \eqref{eq:barWeclass} by pairing $-[\omega_1]$ with both sides as in
  the proof of Theorem~\ref{thm:eulerW}.
\end{enumerate}

\paragraph{Notes and references.}

The orbits of the $\SLtwoR$ action on $\Omega_1\moduli[g]$ project to
immersions $\half\to\moduli[g]$ of the hyperbolic plane into
$\moduli[g]$ which are isometric and totally geodesic with respect to
the Teichm\"uller metric on $\moduli[g]$.  A totally geodesic
immersion $\half\to\moduli[g]$ is called a \emph{Teichm\"uller disk}.
It sometimes happens that a Teichm\"uller disk covers an algebraic
curve $C$ in $\moduli[g]$.  In that case, the normalization of $C$ is
called a \emph{Teichm\"uller curve}.  The curves $\W$ are examples of
Teichm\"uller curves in genus two.  This action of $\SLtwoR$ is also
closely related to the study of billiards in rational angled polygons,
as well as the study of interval exchange maps.  For information
about Teichm\"uller curves, the action of $\SLtwoR$ on
$\Omega\moduli[g]$, and its relation with rational billiards and
interval exchange maps, see for example \cite{kms}, \cite{masur82},
\cite{masur02}, \cite{veech86}, \cite{veech89}, \cite{veech90},
\cite{veech92}, \cite{kontsevichzorich}, \cite{zorich99},
\cite{moeller}, or \cite{bm}.

By analogy with the work of Ratner \cite{ratner} on the dynamics of actions of groups generated by unipotent
elements on homogeneous spaces, it is conjectured that the closure of every Teichm\"uller disk in $\moduli[g]$
is an algebraic suborbifold of $\moduli[g]$.  McMullen's work establishes this conjecture in genus two for
Teichm\"uller disks generated by Abelian differentials.  McMullen analyzed the dynamics of $\SLtwoR$ on
$\Omega_1\moduli$ in the series of papers, \cite{mcmullenbild}, \cite{mcmulleninfinite}, \cite{mcmullenabel},
\cite{mcmullenspin}, \cite{mcmullendecagon}, and \cite{mcmullentor}.  In these papers, he classified the
closures of $\SLtwoR$ orbits in $\Omega_1\moduli$, and classified Teichm\"uller curves in genus two.  He also
introduced the curves $\W$ in these papers and showed that they are all of the Teichm\"uller curves
which are generated by an Abelian differential with a double zero.  According to \cite{mcmullenbild}, the
Teichm\"uller curve $\W$ is primitive exactly when $D$ is nonsquare.  The form generated by identifying
opposite sides of a regular decagon has two simple zeros and generates a primitive Teichm\"uller curve
$D_{10}$ in $\moduli$.  In \cite{mcmullentor}, McMullen showed that $D_{10}$ and the $\W$ for $D$ nonsquare
are in fact all of the primitive Teichm\"uller curves in $\moduli$ generated by an Abelian differential.

Teichm\"uller curves in genus two were also studied by Calta in \cite{calta}, using the
Kenyon-Smillie invariant introduced in \cite{ks}.

Associated to every known $\SLtwoR$ orbit closure in $\Omega_1\moduli[g]$ is a canonical finite, ergodic measure
$\mu$.  The volume of $\mu$ is an  interesting quantity, which according to 
Veech \cite{veech98}  and Eskin-Masur \cite{eskinmasur} gives information about the dynamics of the geodesic
flow of a generic Abelian differential $(X, \omega)\in S$ with respect to the canonical flat metric on $X$
determined by $\omega$.  Volumes of orbit closures and their applications to dynamics are studied in
\cite{eo}, \cite{eop}, \cite{ems}, \cite{emz}, \cite{hubertlelievre}, \cite{lelievreroyer}, and \cite{lelievre}.

A closed $\SLtwoR$-orbit $S$ lies over a Teichm\"uller curve $C$.  In this case, the measure $\mu$ descends to a multiple of
the hyperbolic area measure on $C$.  These means that calculating the volume of $\mu$ is equivalent to
calculating $\chi(C)$.  In genus two, the Euler characteristics of the Teichm\"uller curves $\W[d^2]$ (as well as the volume
of a related $\SLtwoR$-invariant measure on $\X[d^2]$) were studied in \cite{ems}.  They established \eqref{eq:ems} by
counting certain special Abelian differentials, called square-tiled surfaces.  The connected components of
$\W$ were classified in \cite{hubertlelievre} when $D$ is the square of a prime, and in \cite{mcmullenspin} for
arbitrary $D$.  Theorem~\ref{thm:hlformula} calculating $\chi(\We[d^2])$ was established in
\cite{hubertlelievre} when $d$ is prime and was conjectured for arbitrary $d$.  In \cite{lelievreroyer} Leli\`evre and Royer
established Theorem~\ref{thm:hlformula} independently by counting square-tiled surfaces using the theory of
quasimodular forms.

With the calculation of $\chi(\W)$, we now know the Euler characteristics of all of the Teichm\"uller curves
in $\moduli$ which are generated by Abelian differentials with a double zero.  The number of cusps of $\W$ was
calculated by McMullen in \cite{mcmullenspin}.  We can also calculate the number of elliptic points on $\W$
using known formulas for the numbers of elliptic points on $\X$ and $\P$.  Thus we can calculate the genus of
$\W$ for any $D$.  At this point we don't know the number of elliptic points on the components $\We$, so we
can't calculate the genus of these components, but we would conjecture that the two components have the same
number of elliptic points and genera if $D$ is not square.  The curve $\W$ is defined over $\ratls$, so the
Galois automorphism $\sigma$ of $\ratls(\sqrt{D})$ acts on the $\overline{\ratls}$ points of $\W$.  An alternative
approach to showing that $\chi(\Wzero) = \chi(\Wone)$, as well as showing that the components have the same
number of elliptic points and genera, would be to show that $\sigma$ permutes the two components of $\W$.

\paragraph{Acknowledgments.}

This paper is an expanded version of my Ph.D. thesis at Harvard.  I would
like to thank my advisor, Curt McMullen, for his invaluable help with
every stage of this work.  I would also like to thank Sabin Cautis,
Alex Eskin, Florian Herzig, and Howie Masur for useful conversations,
and I would like to thank Joe Harris and Yum-Tong Siu for reading an
earlier version of this paper.


%% file: lshaped.pstex_t
\begin{picture}(0,0)%
\includegraphics{lshaped.pstex}%
\end{picture}%
\setlength{\unitlength}{3947sp}%
\begingroup\makeatletter\ifx\SetFigFont\undefined%
\gdef\SetFigFont#1#2#3#4#5{%
  \reset@font\fontsize{#1}{#2pt}%
  \fontfamily{#3}\fontseries{#4}\fontshape{#5}%
  \selectfont}%
\fi\endgroup%
\begin{picture}(2526,2431)(5701,-4910)
\put(5701,-3811){\makebox(0,0)[lb]{\smash{{\SetFigFont{12}{14.4}{\rmdefault}{\mddefault}{\updefault}{\color[rgb]{0,0,0}$a$}%
}}}}
\put(6751,-4861){\makebox(0,0)[lb]{\smash{{\SetFigFont{12}{14.4}{\rmdefault}{\mddefault}{\updefault}{\color[rgb]{0,0,0}$b$}%
}}}}
\put(6301,-2611){\makebox(0,0)[lb]{\smash{{\SetFigFont{12}{14.4}{\rmdefault}{\mddefault}{\updefault}{\color[rgb]{0,0,0}$1$}%
}}}}
\put(7951,-4111){\makebox(0,0)[lb]{\smash{{\SetFigFont{12}{14.4}{\rmdefault}{\mddefault}{\updefault}{\color[rgb]{0,0,0}$1$}%
}}}}
\end{picture}%

%% file: abeliansurfaces.tex
\section{Abelian surfaces and real multiplication}
\label{sec:abelianvarieties}

We discuss in this section the necessary preliminaries involving Abelian surfaces.  In
\S\ref{subsec:abeliansurfaces}, we discuss the Siegel modular varieties which parameterize Abelian varieties.
In \S\ref{subsec:realmultiplication}, we introduce real multiplication, and in \S\ref{subsec:hilbertmodular},
we discuss  Hilbert modular surfaces.  

\subsection{Abelian surfaces}
\label{subsec:abeliansurfaces}

\paragraph{Abelian varieties.}

A \emph{complex torus} is a quotient $A=V/\Lambda$, where $\Lambda$ is a lattice in a finite dimensional
complex vector space $V$.  A \emph{principal polarization} on $A$ is a Hermitian
form $H$ on $V$ such that $\Im H$ takes integral values on $\Lambda\times \Lambda$, and the pairing
$$\Im H\colon \Lambda\times \Lambda \to \zed$$
is unimodular.  A \emph{principally polarized Abelian variety}
is a complex torus equip\-ped with a principal polarization.

\paragraph{Siegel modular varieties.}

The Siegel upper half space is
$$\siegelhalf[g] = \{Z\in M_g(\cx) : Z^t = Z, \Im Z > 0\},$$
an open subset in the $g(g+1)/2$ dimensional space of symmetric matrices in $M_g(\cx)$.  The group $\SPR{2g}$
acts on $\siegelhalf[g]$ as a group of biholomorphic transformations by
$$
\begin{pmatrix}
  \alpha & \beta \\
  \gamma & \delta
\end{pmatrix}
\cdot Z = (\alpha Z + \beta)(\gamma Z + \delta)^{-1}.
$$

Equip $\Lambda=\zed^{2g}$ with the usual symplectic form defined by the matrix,
\begin{equation}
  \label{eq:sympform}
  \begin{pmatrix}
    \phantom{-}0 & I \\
    -I & 0
  \end{pmatrix}.
\end{equation}
Following \cite{bl}, we define for each $Z\in\siegelhalf[g]$ a principally polarized Abelian surface $X_Z$
together with a symplectic isomorphism $\Lambda\to H_1(X_Z; \zed)$.  For $Z\in \siegelhalf[g]$, let
$$\Lambda_Z = (Z, I)\Lambda,$$
a lattice in $\cx^g$.  Let $X_Z = \cx^g/\Lambda_Z$, and give $X_Z$ the polarization coming from the Hermitian
form on $\cx^g$ defined by the matrix,
$$(\Im Z)^{-1}.$$
It can be shown that two given points $Z, Z' \in \siegelhalf[g]$, the polarized Abelian surfaces $X_Z$ and
$X_{Z'}$ are isomorphic if and only if $Z$ and $Z'$ are equivalent under the action of $\SPZ{2g}$.

The Siegel modular variety is  $\siegelmod[g] = \siegelhalf[g]/\SPZ{2g}$.  From the previous paragraph, we obtain:
\begin{theorem}
  The normal analytic space $\siegelmod[g]$ is the moduli space of principally polarized Abelian varieties of
  dimension $g$.
\end{theorem}

\paragraph{Satake compactification.}

Following \cite{vandergeer88}, we briefly describe the Satake compactification $\satsiegelmod[g]$ of
$\siegelmod[g]$, introduced by Satake in \cite{satake56a}.

The transformation,
$$Z\mapsto (1 + i Z)(1 - i Z)^{-1},$$
maps $\siegelhalf[g]$ biholomorphically onto the bounded symmetric domain,
$$\mathcal{D}_g = \{Z\in M_g(\cx) : Z = Z^t, 1_g - Z \overline{Z} >0 \}.$$
For $r\leq n$, we map $\mathcal{D}_r$ into $\overline{\mathcal{D}}_g$ by
$$Z\mapsto 
\begin{pmatrix}
  Z & 0 \\
  0 & 1_{n-r}
\end{pmatrix},
$$
and let
$$\mathcal{D}^*_g = \bigcup_{0\leq r\leq g} \bigcup_{h\in\SPQ{2g}} h \mathcal{D}_r.$$
Satake defined a natural
topology on $\mathcal{D}^*_g$ for which the action of $\SPQ{2g}$ extends continuously to an action on
$\mathcal{D}^*_g$.  With this topology, the Satake showed that the quotient $\satsiegelmod[g] =
\mathcal{D}^*_g / \SPZ{2g}$ is a compact, normal complex analytic space which contains one copy of $\siegelmod[r]$ for
each $r\leq g$.  Baily showed in \cite{baily58} that $\satsiegelmod[g]$ has the structure of a normal projective variety.

\subsection{Real multiplication}
\label{subsec:realmultiplication}

\paragraph{Quadratic orders.}

Let $K$ be a quadratic field or $\ratls\oplus\ratls$.  A \emph{quadratic order} is a subring $\mathcal{O}$ of
$K$ such that $1\in\mathcal{O}$ and $\mathcal{O}\otimes\ratls
= K$.  Any quadratic order is isomorphic to one of
$$\ord = \zed[T]/(T^2 + b T + c),$$
where $b,c\in\zed$ and $b^2-4c = D$.  The isomorphism class of $\ord$ only depends on $D$, so this defines a
unique quadratic order for every nonzero integer $D\equiv 0$ or $1\pmod 4$.  This integer $D$ is called the
\emph{discriminant} of $\ord$, and $D$ is a \emph{fundamental discriminant} if $D$ is not of the form $f^2 E$
for some integers $f$ and $E$ with $E\equiv 0$ or $1\pmod 4$ and $f>1$.

If $D$ is not square, then $\ord$ is a subring of $\ratls(\sqrt{D})$, and $\ord$ is the ring of integers of
$\ratls(\sqrt{D})$ if and only if $D$ is a fundamental discriminant; otherwise, $\ord$ is a subring of the
ring of integers.  If $D$ is square, then $\ord$ is a subring of $\ratls\oplus\ratls$.  In that case,
\begin{equation}
  \label{eq:odsquared}
  \ord[d^2]=\{(x, y)\in\zed\times\zed : x\equiv y\pmod d\}.
\end{equation}

We will regard $\ratls\oplus\ratls$ as an extension of $\ratls$ by the diagonal map
$\ratls\to\ratls\oplus\ratls$.  The Galois automorphism of $\ratls\oplus\ratls$ is $(x, y)'=(y, x)$, and we
can use this to define norm and trace on $\ratls\oplus\ratls$ as for a field.

For the rest of this paper, fix two nonzero homomorphisms $\iota_i\colon\ord\to\reals$.  If $D$ is not square,
then $\iota_i$ is an embedding, and if $D$ is square, then $\iota_i$ is induced by one of the two projections
$\zed\oplus\zed\to\zed$.  We will often abuse terminology and call $\iota_i$ an embedding even when $D$ is
square.  We will also use the notation $\lambda^{(i)} = \iota_i(\lambda)$.

\paragraph{Inverse different.}

Given a quadratic order $\ord$, the inverse different is the fractional ideal,
$$\ord^\vee = \{x\in K_D : \tr(xy)\in\zed,  \forall y \in\ord\}.$$
More concretely,
$$\ord^\vee = \frac{1}{\sqrt{D}}\ord,$$
where if $D=d^2$, we interpret $\sqrt{D}$ to be $(d, -d)$.

We equip $\ord\oplus\ord^\vee$ with the unimodular symplectic pairing,
$$\langle (x_1, y_1), (x_2, y_2)\rangle = \tr(x_1 y_2 - x_2 y_1).$$

\paragraph{Real multiplication.}

Consider an Abelian surface $A=V/\Lambda$.  \emph{Real multiplication by $\ord$} on $A$ is a monomorphism
$\rho\colon\ord\to \End(A)$ (where $\End(A)$ is the ring of holomorphic endomorphisms of $A$) with the
following properties:
\begin{itemize}
\item For each $\lambda\in \ord$, the lift $\tilde{\rho}(\lambda)\colon V \to V$ is self-adjoint with respect
  to the Hermitian form on $V$ given by the polarization of $A$.
\item $\rho$ is \emph{proper} in the sense that it doesn't extend to a monomorphism
  $$\rho'\colon\ord[E]\to\End(A)$$ for some $\ord[E]\supset\ord$ -- that is, for some $\ord[E]$ with $E=f^2 D$ for
  some $f>1$.
\end{itemize}

Let $\Omega(A)\isom V^*$ be the space of holomorphic one-forms on $A$.  Since the real multiplication is
self-adjoint, there is an eigenspace decomposition
$$\Omega(A) = \Omega^1(A)\oplus \Omega^2(A).$$
We order the eigenspaces so that $\rho(\lambda)\cdot\omega = \lambda^{(i)}\omega$ for $\omega\in\Omega^i(A)$.

\subsection{Hilbert modular surfaces}
\label{subsec:hilbertmodular}

For any quadratic discriminant $D$, let $K_D = \ord\otimes\ratls$, which is $\ratls(\sqrt{D})$ if $D$ is not
square and is $\ratls\oplus\ratls$ if $D$ is square.  Define the group,
$$\SLtwoord = \left\{
  \begin{pmatrix}
    a & b \\
    c & d
  \end{pmatrix}
  \in {\rm SL}_2 K_D : a\in\ord, b\in(\ord^\vee)^{-1}, c\in\ord^\vee, d\in\ord \right\}.
$$
$\SLtwoord$ has two embeddings in $\SLtwoR$ induced by the two embeddings $\iota_i\colon K_D\to\reals$.

\begin{definition}
  The \emph{Hilbert modular surface} of discriminant $D$ is the quotient,
  $$\X=\half\times\half/\SLtwoord,$$
  where $\SLtwoord$ acts on $\half \times \half$ by 
  $$
  \begin{pmatrix}
    a & b \\
    c & d
  \end{pmatrix}
  \cdot(z_i, z_2) = \left(\frac{a^{(1)} z_1 + b^{(1)}}{c^{(1)} z_1 + d^{(1)}}, \frac{a^{(2)} z_2 +
      b^{(2)}}{c^{(2)} z_2 + d^{(2)}}\right).
  $$
\end{definition}

There is an isomorphism $\SLtwoord\to{\rm SL}_2 \ord$ defined by
$$
\begin{pmatrix}
  a & b \\
  c & d
\end{pmatrix}
\mapsto
\begin{pmatrix}
  a & \frac{1}{\sqrt{D}} b \\
  \sqrt{D} c & d
\end{pmatrix},
$$
and the map $T\colon \half\times\half\to\half\times(-\half)$ induces an isomorphism,
$$\X\isom\half\times(-\half)/{\rm SL}_2\ord.$$

\begin{theorem}
  \label{thm:hilbertmoduli}
  The Hilbert modular surface $\X$ is the moduli space of all pairs $(A, \rho)$, where $A$ is a principally
  polarized Abelian surface, and $\rho\colon\ord\to\End(A)$ is a choice of real multiplication on $A$.
\end{theorem}

\begin{proof}[Sketch of proof (following \cite{mcmullenhilbert})]
  
  Given $\tau=(\tau_1, \tau_2)\in\half\times\half$, let
  $$\phi_\tau\colon \ord\oplus\ord^\vee$$
  be the embedding,
  $$\phi_\tau(x, y) = (x^{(1)} + y^{(1)}\tau_1, x^{(2)} + y^{(2)}\tau_2).$$
  Let $A_\tau =
  \cx^2/\phi_\tau(\ord\oplus\ord^\vee)$ with the principal polarization induced by the symplectic pairing on
  $\ord\oplus\ord^\vee$.  This polarization is also given by the Hermitian form,
  \begin{equation}
    \label{eq:polarization}
    H_\tau(z, w) = \frac{1}{\Im\tau_1} z_1 \bar{w}_1 + \frac{1}{\Im\tau_2} z_2 \bar{w}_2.  
  \end{equation}
    Define real multiplication on $A_\tau$ by
  $$\lambda\cdot(z_1, z_2) = (\lambda^{(1)} z_1, \lambda^{(2)} z_2).$$
  We thus get a map from
  $\half\times\half$ to the set of all triples $(A, \rho, \phi)$, where $(A, \rho)$ is a principally polarized
  Abelian surface with real multiplication by $\ord$, and $\phi$ is a choice of an $\ord$-linear, symplectic
  isomorphism $\phi\colon\ord\oplus\ord^\vee\to H_1(A; \zed)$.

  Given
  $$g=
  \begin{pmatrix}
    a & b \\
    c & d
  \end{pmatrix}
  \in {\rm SL}_2 K_D,
  $$
  let
  $$g^* = 
  \begin{pmatrix}
    a & -b \\
    -c & d
  \end{pmatrix},$$
  and define
  \begin{equation}
    \label{eq:chigtau}
    \chi(g, \tau) =
    \begin{pmatrix}
      (c^{(1)} \tau_1 + d^{(1)})^{-1} & 0 \\
      0 & (c^{(2)}\tau_2 + d^{(2)})^{-1}
    \end{pmatrix}.
  \end{equation}
  We have the following commutative diagram.
  $$\xymatrix{
    \ord\oplus\ord^\vee \ar[r]^-{\phi_\tau}\ar[d]_{g^*} & \cx^2 \ar[d]^{\chi(g, \tau)} \\
    \ord\oplus\ord^\vee \ar[r]_-{\phi_{g\cdot\tau}} & \cx^2
  }
  $$
  Thus $\chi(g, \tau)$ induces an isomorphism between $A_\tau$ and $A_{g\cdot\tau}$ which preserves the
  polarizations and commutes with the action of real multiplication.  We then get a map from $\X$ to the set of
  all principally polarized Abelian surfaces with a choice of real multiplication, which can be shown to be
  a bijection.
\end{proof}

Replacing each pair $(A, \rho)$ with $(A, \rho')$, where $\rho'$ is the composition of $\rho$ with the Galois
automorphism of $\ord$, induces an involution $\tau$ of $\X$.  The lift of $\tau$ to the universal cover
$\half\times\half$ of $\X$ is given by $\tilde{\tau}(z_1, z_2) = (z_2, z_1)$.

There is a natural map $j\colon\X\to\siegelmod$ which forgets the choice of real multiplication.  This map $j$
is generically two to one and is equivariant with respect to $\tau$.

\paragraph{Baily-Borel compactification.}

We can regard the boundary $\bdry \half\subset \proj^1(\cx)$ as $\proj^1(\reals) \isom \reals\cup\{\infty\}$.
Given a real quadratic field $K\subset\reals$, define an embedding $\proj^1(K)\to \proj(\reals)^2$ by
$$[x : y] \mapsto ([x: y], [x': y']).$$
When $D$ is not square, define via this embedding
$$(\half\times\half)_D = (\half\times\half)\cup \proj^1(K_D).$$
When $D$ is square, define
$$(\half\times\half)_D = (\half\cup\proj^1(\ratls))\times(\half\cup\proj^1(\ratls)).$$

We give $\half\cup\proj^1(\ratls)$ the usual topology where if $r\in\proj^1(\ratls)$, then a basis of open
neighborhoods of $r$ is given by sets of the form $U\cup\{r\}$, where $U\subset\half$ is an open horoball
resting on $r$.  We then give $(\half\times\half)_{d^2}$ the product topology.  When $D$ is not square, there
is a similar natural topology on $(\half\times\half)_D$, which is described in \cite{vandergeer88}.

The action of $\SLtwoord$ on $\half\times\half$ extends continuously to $(\half\times\half)_D$.  The quotient,
$$\bX = (\half\times\half)_D/\SLtwoord,$$
is compact and Hausdorff.  The space $\bX$ is the Baily-Borel compactification of $\X$.

\begin{theorem}[{\cite{baily58}}]
  The compactification $\bX$ is a normal, projective, algebraic variety.  
\end{theorem}

When $D$ is not square, $\bX\setminus \X$ consists of finitely many points, which we will call the cusps of
$\X$.  When $D=d^2$, the image of $\proj^1(\ratls)\times\proj^1(\ratls)$ in $\bX[d^2]$ also consists of
finitely many points which we will call the cusps of $\bX[d^2]$.   In $\bX[d^2]$, let
$$\Rone = \pi\left(\bigcup_{r\in\proj^1(\ratls)} \half\times\{r\}\right),$$
and let
$$\Rtwo = \pi\left(\bigcup_{r\in\proj^1(\ratls)} \{r\}\times\half\right),$$
where $\pi\colon(\half\times\half)_D\to\bX$ is the natural quotient map.  Then we have a disjoint union,
$$\bX[d^2] = \X[d^2]\amalg\Rone\amalg\Rtwo\amalg C,$$
where $C$ is the set of cusps of $\bX[d^2]$.

For $d\in\nats$, define
$$\Gamma_1(d)=\left\{
  \begin{pmatrix}
    a & b\\
    c & e
  \end{pmatrix}
  \in \SLtwoZ : a\equiv e\equiv 1 \pmod d,\text{ and } c\equiv 0 \pmod d \right\}.$$

\begin{prop}
  The curve $\Ri$ is irreducible, and
  $$\Ri\isom\half/\Gamma_1(d).$$
\end{prop}

\begin{proof}
  By \eqref{eq:odsquared}, we can regard ${\rm SL}_2\ord[d^2]$ as
  \begin{equation}
    \label{eq:sl2ord}
    \{(A, B)\in (\SLtwoZ)^2 : A\equiv B\pmod d\}.
  \end{equation}
  Since $\SLtwoZ$ acts transitively on $\proj^1(\ratls)$, we have
  $$\Rone = \half\times\{\infty\} / \Stab_{\half\times\{\infty\}}.$$
  The matrices $A$ in the pair $(A, B)\in \Stab_{\half\times\{\infty\}}$ are exactly the matrices which are
  congruent to an upper triangular matrix mod $d$.  Thus $\Rone$ is as claimed.
\end{proof}

When $D$ is not square, the cusps of $\bX$ are complicated singularities; however, when $D$ is square, they
are just orbifold singularities.

\begin{prop}
  The Baily-Borel compactification $\bX[d^2]$ is a compact, complex orbifold with singularities at the elliptic points
  of $\X[d^2]$, the elliptic points of $\Ri$, and possibly the cusps of $\bX[d^2]$.
\end{prop}

\begin{proof}
  Consider the principal congruence subgroup,
  $$\Gamma(d) = \{A\in\SLtwoZ : A\equiv I \pmod d\}.$$
  The product,
  $$\Gamma(d)\times\Gamma(d)\subset{\rm SL}_2\ord[d^2],$$
  is a finite index, normal subgroup.  We have,
  $$(\half\times\half)_{d^2} / (\Gamma(d)\times\Gamma(d)) = \overline{\half/\Gamma(d)}\times
  \overline{\half/\Gamma(d)},$$
  which is a manifold, so
  $$\bX[d^2] = (\overline{\half/\Gamma(d)}\times \overline{\half/\Gamma(d)})/G,$$
  where $G = {\rm SL}_2\ord[d^2]/(\Gamma(d)\times\Gamma(d))$, a finite group.  Thus $\bX[d^2]$ is compact
  orbifold with singularities at the fixed points of $G$.
\end{proof}

The involution $(z_1, z_2)\mapsto(z_2, z_1)$ of $\half\times\half$ extends continuously to
$(\half\times\half)_D$, so the involution $\tau$ extends to an involution $\tau$ of $\bX$, which preserves the
cusps of $\bX$.  If $D=d^2$, then we also have $\tau(\Rone) = \Rtwo$.

\begin{theorem}
  The natural map $j\colon\X\to\siegelmod$ extends to a finite morphism $j\colon\bX\to\satsiegelmod$ which
  sends the cusps of $\bX$ to the point $\siegelmod[0]\in\siegelmod$ and sends the curves $\Ri[D]$ to $\siegelmod[1]$.
\end{theorem}

This  allows us to define $\bX$ alternatively in terms of the image of $\X$ in $\satsiegelmod$. Recall
the notion of the normalization of a variety in a finite algebraic extension of its function field discussed
in \S\ref{sec:normal}.

\begin{cor}
  The Baily-Borel compactification $\bX$ is the normalization of the closure of $j(\X)$ in $\satsiegelmod$ in
  the field $K(\X)$, the function field of $\X$.
\end{cor}

\paragraph{Line bundles on $\X$.}

We now discuss some bundles on $\X$, which will be used in \S\ref{sec:bundles}.  For $i=1$ or $2$, define an action of
$\SLtwoord$ on $(\half\times\half)\times\cx$ by
$$
\begin{pmatrix}
  a & b \\
  c & d
\end{pmatrix}
\cdot (z_1, z_2, w) =  \left(\frac{a^{(1)} z_1 + b^{(1)}}{c^{(1)} z_1 + d^{(1)}}, \frac{a^{(2)} z_2 +
    b^{(2)}}{c^{(2)} z_2 + d^{(2)}}, (c^{(i)} z_i + d^{(i)})^2 w\right).
$$
The quotient is a line bundle over $\X$, which we call $L_i$.

We define a Hermitian metric $\tilde{h}_i$ on $(\half\times\half)\times\cx$ by defining on the fiber over $(z_1,z_2)$,
$$\tilde{h}_i(w, w) = y_i^2 |w|^2.$$
The metric $\tilde{h}_i$ is $\SLtwoord$-invariant, so it descends to a Hermitian metric $h_i$ on $L_i$.  The
Chern form of $\tilde{h}_i$ is,
\begin{align}
  c_1(L_i, h_i) &= -\frac{i}{\pi} \partial\partialbar \log(y_i) \notag\\
  &= \frac{1}{2\pi}\frac{dx_i\wedge dy_i}{y_i^2}, \label{eq:omegailift}
\end{align}
so the Chern form of $h_i$ is $c_1(L_i, h_i) = \omega_i$, where $\omega_i$ is the 2-form on $\X$ covered by
the form \eqref{eq:omegailift} on $\half\times\half$.

Define line bundles,
$$\Omega^i\X = \{(A, \omega) : A\in\X \text{ and } \omega\in\Omega^i(A)\},$$
and
$$Q^i\X = (\Omega^i\X)^{\otimes 2},$$
for $i=1$ or $2$.  We will mostly use these bundles when $i=1$, so we will write $\Omega\X$ or $Q\X$ for
$\Omega^1\X$ or $Q^1\X$.

For $A\in\X$, the polarization defines a Hermitian metric on $\Omega(A)$, which restricts to a Hermitian
metric on $\Omega^i(A)$.  Put this Hermitian metric on each $\Omega^i(A)$ to define a Hermitian metric
$h_\Omega^i$ on $\Omega^i\X$.  Let $h_Q^i$ be the induced metric on $Q^i\X$.

\begin{prop}
  \label{prop:chHQ}
  There is an isomorphism $L_i\to Q^i\X$ which preserves the Hermitian metrics on these bundles.  Thus,
  $$c_1(Q^i\X, h_Q^i) = \omega_i.$$
\end{prop}

\begin{proof}
  Define a map $f\colon (\half\times\half)\times\cx\to Q^i\X$ by
  $$(\tau, w)\mapsto (A_\tau, w\, dz_i^2),$$
  where $A_\tau = \cx^2/\phi_\tau(\ord\oplus\ord^\vee)$ is as in the proof of Theorem~\ref{thm:hilbertmoduli}.
  With $\chi(g, \tau)$ as in \eqref{eq:chigtau}, we have
  $$\chi(g, \tau)_*dz_i^2 = (c^{(i)} \tau_i + d^{(i)})^2 dz_i^2.$$
  This means that $f$ descends to an isomorphism $L_i\to Q^i\X$ as claimed.  The fact that this isomorphism
  preserves the Hermitian metrics follows directly from the formula \eqref{eq:polarization} for the
  polarization.
\end{proof}

The product foliation of $\half \times\half$ with leaves of the form $\{c\}\times\half$ is invariant under
$\SLtwoord$.  Let $\A$ be the induced foliation of $\X$, and let
$$T^*\A\to\X$$
be the line bundle whose fiber
over a point $p$ is the cotangent bundle to the leaf of $\A$ through $p$.  We give the leaves of $\A$ their
hyperbolic metric, and give $T^*\A$ the induced Hermitian metric $h_{\A}$.  The following is easy to check.

\begin{prop}
  \label{prop:c1Tstar}
  There is an isomorphism $L_2\to T^*\A$ preserving the Hermitian metrics.  Thus,
  $$c_1(T^*\A, h_{\A}) = \omega_2.$$
\end{prop}

The following relation between $\Omega\X$ and $T^*\A$ will be used in \S\ref{sec:modular}.

\begin{prop}
  \label{prop:periodsconstant}
  For any $p\in\X$, there is a neighborhood $U$ of $p$ in the leaf $L$ of $\A$ through $p$ and a section
  $\omega$ of $\Omega\X$ over $U$ such that the periods of the forms $\omega(z)$ are constant over $U$.
\end{prop}

\begin{proof}
  Let
  $$(\tau_1, \tau_2)\in \tilde{X}_D = \half \times \half$$
  lie over $p$ in the universal cover of $\X$.
  Define a section $s$ of $\Omega\tilde{X}_D$ by
  $$s(\tau_1, \tau_2) = dz_1\in \Omega^1(A_{(\tau_1, \tau_2)}),$$
  with $A_\tau$ as in the proof of Theorem~\ref{thm:hilbertmoduli}.  By the definition of $A_\tau$, the
  periods of $s(\tau_1, \tau_2)$ only depend on $\tau_1$, so are constant along $\A$. 
\end{proof}

\paragraph{Euler characteristic of $\X$.}

We now discuss the Euler characteristic $\chi(\X)$, which by the generalized Gauss-Bonnet theorem is given by
the volume,
$$\chi(\X) = \int_{\X} \omega_1\wedge\omega_2,$$
where $\omega_1$ are the 2-forms on $\X$ defined by
\eqref{eq:omegailift}.  The volume of $\X$ was calculated by Siegel when $D$ is a fundamental discriminant.

For any number field $K$, the \emph{Dedekind zeta-function} of $K$ is defined by
$$\zeta_K(s) = \sum_{\mathfrak{a}} \frac{1}{N^K_\ratls(\mathfrak{a})^s},$$
where the sum is over all nonzero ideals $\mathfrak{a}\subset\mathcal{O}$ with $\mathcal{O}$ the ring of integers in
$K$.  If $K=\ratls$, then $\zeta_K$ is just the Riemann zeta-function.  The definition of $\zeta_K$ also makes
sense if $K=\ratls\oplus\ratls$.  In that case,
$$\zeta_{\ratls\oplus\ratls}(s) = \zeta_\ratls(s)^2.$$
$\zeta_K$ can be analytically continued to a meromorphic function on $\cx$ which has a simple pole at $s=1$.  

\begin{theorem}[\cite{siegel36}]
  \label{thm:siegeleuler}
  When $D$ is a fundamental discriminant, and $D\neq1$,
  \begin{equation*}
    \chi(\X) = 2 \zeta_{K_D}(-1).
  \end{equation*}
\end{theorem}

We will also want to know $\chi(\X)$ when $D$ is not a fundamental discriminant.  For integers $a$ and $b$
with $b>0$, let
$$\kron{a}{b}$$
be the Kronecker symbol, defined in \cite[p.~82]{miyake}.

\begin{theorem}
  \label{thm:chiX}
  If $D$ is a fundamental discriminant, and $f\in\nats$, then
  \begin{equation*}
    \chi(\X[f^2D]) =
    \begin{cases}
      1/36 & \text{if $D=f=1$;} \\
      1/6 & \text{if $D=1$ and $f=2$.}
    \end{cases}
  \end{equation*}
  Otherwise,
  \begin{equation}
    \label{eq:chiX}
    \chi(\X[f^2D]) = 2 f^3 \zeta_{K_D}(-1) \sum_{r|f}\kron{D}{r}\frac{\mu(r)}{r^2}.
  \end{equation}
\end{theorem}

\begin{remark}
  When $D=1$, \eqref{eq:chiX} reduces to
  \begin{equation*}
    \chi(\X[d^2]) = \frac{1}{72}d^3\sum_{r|d}\frac{\mu(r)}{r^2}.
  \end{equation*}
\end{remark}

We will derive Theorem~\ref{thm:chiX} from Theorem~\ref{thm:siegeleuler} in a sequence of lemmas.

\begin{lemma}
  \label{lem:index1}
  We have,
  \begin{equation*}
    [\PSLtwoZ\times\PSLtwoZ : {\rm PSL}_2\ord[4]] = 6,
  \end{equation*}
  and
  \begin{equation*}
    [\PSLtwoZ\times\PSLtwoZ : {\rm PSL}_2\ord[d^2]] = \frac{1}{2}d^3\sum_{r|d}\frac{\mu(r)}{r^2}
  \end{equation*}
  if $d>2$.
\end{lemma}

\begin{proof}
  We will regard ${\rm SL}_2\ord[d^2]$ as the subgroup of $\SLtwoZ\times\SLtwoZ$ given in \eqref{eq:sl2ord}.
  The natural map,
  $$\SLtwoZ\to{\rm SL}_2(\zed/d),$$
  is surjective by \cite[Theorem~4.2.1]{miyake} and has kernel the principal
  congruence subgroup $\Gamma(d)$, so we have the exact sequence,
  $$0\to\Gamma(d)\times\Gamma(d)\to{\rm SL}_2\ord[d^2] \to {\rm SL}_2(\zed/d)\to 0.$$
  Thus,
  $$[{\rm SL}_2\ord[d^2] : \Gamma(d)\times \Gamma(d)] = |{\rm SL}_2(\zed/d)|,$$
  and
  \begin{align}
    \label{eq:index}
    [\SLtwoZ\times\SLtwoZ : {\rm SL}_2\ord[d^2]] &= \frac{[\SLtwoZ\times\SLtwoZ : \Gamma(d)\times
      \Gamma(d)]}{[{\rm SL}_2\ord[d^2] : \Gamma(d)\times\Gamma(d)]} \\
    \notag
    &= |{\rm SL}_2(\zed/d)|.    
  \end{align}
  By \cite[Theorem~4.2.4]{miyake},
  $$|{\rm SL}_2(\zed/d)| = d^3\sum_{r|d}\frac{\mu(r)}{r^2}.$$

  The kernel of $\SLtwoZ\times\SLtwoZ \to \PSLtwoZ\times\PSLtwoZ$ has order $4$, and the kernel of ${\rm
    SL}_2\ord[d^2]\to{\rm PSL}_2\ord[d^2]$ has order $4$ if $d=2$, and order $2$ if $d>2$.  The desired
  formulas follow from this and \eqref{eq:index}.
\end{proof}

\begin{lemma}
  \label{lem:index2}
  For any fundamental discriminant $D\neq 1$, we have,
  \begin{align*}
    [{\rm PSL}_2\ord : {\rm PSL}_2\ord[f^2D]] &= f^3 \prod_{p|f}\left(1-\kron{D}{p}p^{-2}\right) \\
    &= f^3\sum_{r|f}\kron{D}{r}\frac{\mu(r)}{r^2}.
  \end{align*}
\end{lemma}

\begin{proof}
  Since the map  ${\rm SL_2} \ord[f^2 D] \to {\rm PSL_2} \ord[f^2 D]$ has kernel $\pm I$, we have
  $$[{\rm PSL}_2\ord : {\rm PSL}_2\ord[f^2D]] = [{\rm SL}_2\ord : {\rm SL}_2\ord[f^2D]].$$

  The natural map,
  $${\rm SL}_2\ord / {\rm SL}_2\ord[f^2D] \to \prod_{\substack{
      {p|f} \\
      p \text{ prime}
    }}
  {\rm SL}_2(\ord/p) / {\rm SL}_2(\zed/p),$$
  is a bijection.  Thus, we need only to show that
  $$[{\rm SL}_2(\ord/p) : {\rm SL}_2\zed] = p^3 - \kron{D}{p} p$$
  for every prime $p|f$.

  First, suppose that $(D/p) =-1.$  This means that $p$ remains prime in $\ord$, so
  $\ord/p\isom\mathbb{F}_{p^2}$, where we write $\mathbb{F}_q$ for the unique finite field of order $q$.  It
  is easy to show that
  $$|{\rm SL}_2 \mathbb{F}_q| = q^3 - q,$$
  for any $q$.  Thus
  $$[{\rm SL}_2\mathbb{F}_{p^2} : {\rm SL}_2\mathbb{F}_{p}] = p^3+p$$
  as desired.

  Now suppose that $(D/p) = 1$.  This means that $p$ splits in $\ord$, so $\ord/p\isom
  \mathbb{F}_p\oplus\mathbb{F}_p$.
  Therefore,
  \begin{align*}
    [{\rm SL}_2(\ord/p) : {\rm SL}_2(\zed/p)] &= [{\rm SL}_2 \mathbb{F}_p \times {\rm SL}_2 \mathbb{F}_p : {\rm
      SL}_2 \mathbb{F}_p] \\
    &=|{\rm SL}_2 \mathbb{F}_p| \\
    &= p^3 - p.
  \end{align*}

  Finally, suppose that $(D/p)=0.$  This means that $p$ ramifies in $\ord$, so
  $\ord/p\isom\mathbb{F}_p(\epsilon)/(\epsilon^2$).  We have
  $$
  \begin{pmatrix}
    a + a'\epsilon & b + b'\epsilon \\
    c + c'\epsilon & d + d'\epsilon
  \end{pmatrix}
  \in{\rm SL}_2 (\mathbb{F}_p(\epsilon)/(\epsilon^2))
  $$
  if and only if
  \begin{gather}
    a d - b c=1, \quad\text{and} \label{eq:homer1}\\
    d a' + a d' - c b' - b c' =0. \label{eq:homer2}
  \end{gather}
  For any $(a, b, c, d)$ satisfying \eqref{eq:homer1}, there are $p^3$ solutions to \eqref{eq:homer2}.  Thus,
  $$[{\rm SL}_2 (\mathbb{F}_p(\epsilon)/(\epsilon^2)) : {\rm SL}_2\mathbb{F}_p] = p^3.$$
\end{proof}

\begin{proof}[Proof of Theorem~\ref{thm:chiX}]
  Since the natural map $\X[f^2D]\to \X[D]$ is an orbifold covering map, we have
  $$\chi(\X[f^2 D]) = \chi(\X[D]) [{\rm PSL}_2\ord[D]: {\rm PSL}_2\ord[f^2D]],$$
  so the claim follows from Lemma~\ref{lem:index1}, Lemma~\ref{lem:index2}, and Theorem~\ref{thm:siegeleuler}
  together with
  $$\chi(\X[1]) = \frac{1}{36}$$
  because $\X[1]\isom \half/\SLtwoZ \times\half/\SLtwoZ$, and
  $$\chi(\half/\SLtwoZ) = -\frac{1}{6}.$$
\end{proof}

\paragraph{Modular forms.}

For any fundamental discriminant $D>0$, define
$$H(2, f^2 D) = -12 \zeta_{K_D}(-1) \sum_{r|f}\mu(r)\kron{D}{r}r \sigma_3\left(\frac{f}{r}\right),$$
where
$$\sigma_m(n) = \sum_{d|n}d^m,$$
and we also adopt the convention that
$$\sigma_m(0)= \frac{1}{2}\zeta_\ratls(-m),$$
and $\sigma_m(n)=0$ if $n<0$.
Also define
$$H(2, 0) = \zeta_\ratls(-3)=\frac{1}{120}.$$
The function $\mathcal{H}\colon\half\to\cx,$
\begin{align*}
  \mathcal{H}(\tau)&=\sum_{\substack{D\equiv 0, 1\,(4)\\ D\geq0}} H(2, D)q^D, \quad (\text{with $q(\tau)=e^{2\pi i\tau}$}) \\
  &= -\frac{1}{120} - \frac{1}{12}q-\frac{7}{12}q^4-\frac{2}{5}q^5-q^8-\frac{25}{12}q^9-2
  q^{12}-2q^{13}-\frac{55}{12}q^{16} - \dots,
\end{align*}
was shown by Cohen in \cite{cohen75}, to be a modular form of weight $5/2$ for the
group $\Gamma_0(4)$.  See the discussion in Chapter IX of \cite{vandergeer88} for more about this form and its
relation to $\X$.

It is an elementary argument using M\"obius inversion to show that for any fundamental discriminant $D$ and $f\in\nats,$
\begin{equation}
  \label{eq:h2dsum}
  H(2, f^2D) = -12 \sum_{s|f}\left(\zeta_{K_D}(-1)s^3\sum_{r|s}\kron{D}{r}\frac{\mu(r)}{r^2}\right).
\end{equation}
From this and Theorem~\ref{thm:chiX}, we obtain:

\begin{theorem}
  \label{thm:sumchiX}
  If $D\neq1$ is a fundamental discriminant, then
  \begin{equation*}
    \sum_{r|f} \chi(\X[r^2 D]) = -\frac{1}{6} H(2, f^2D).
  \end{equation*}
\end{theorem}

Siegel and Cohen gave the following simple formula for $H(2, D)$.

\begin{theorem}[{\cite{siegel69, cohen75}}]
  \label{thm:siegelformula}
  If $D$ is not square, then
  \begin{equation*}
    H(2, D) = -\frac{1}{5}\sum_{e\equiv D\,(2)} \sigma_1\left(\frac{D-e^2}{4}\right).
  \end{equation*}
  When $D$ is square,
  \begin{equation*}
     H(2, D) = -\frac{1}{5}\sum_{e\equiv D\,(2)}
    \sigma_1\left(\frac{D-e^2}{4}\right)-\frac{D}{10}.
  \end{equation*}
\end{theorem}

We will obtain an alternative proof of these formulas in \S\ref{sec:fundamental}.

\paragraph{The product locus.}

In $\X$, define the \emph{product locus} $\P$ to be the set of all $A\in \X$ such that $A$ is a polarized
product of elliptic curves.  The locus $\P$ was studied in \cite{mcmullenspin} and \cite{mcmullenhilbert}.
The locus $\P$ is also a union of \emph{modular curves}, which were studied in \cite{hirzebruch} and \cite{vandergeer88}.

Let $\tilde{P}_D\subset\half\times\half$ be the inverse image of $\P$ in $\X$.  The locus $P_D$ is a linear
subvariety of $\X$ in the following sense.

\begin{prop}[\cite{mcmullenspin}]
  $\tilde{P}_D$ is a countable union of graphs of M\"obius transformations $\half\to\half$.
\end{prop}

The goal of the rest of this section is to calculate $\chi(\P)$.  We will work with an auxiliary covering space
$\Q$ of $\P$ defined as follows.

Let $\Omega Q$ be the space of elliptic curves equipped with Abelian differentials; that is,
\begin{equation}
  \label{eq:omegaQ}
  \Omega Q = \Omega\siegelmod[1] \times \Omega\siegelmod[1] = \GLtwoRplus/\SLtwoZ \times \GLtwoRplus/\SLtwoZ, 
\end{equation}
where $\GLtwoRplus$ is the subgroup of $\GLtwoR$ consisting of matrices of positive determinant.  In $\Omega
Q$, let $\Omega\Q$ be the locus of pairs $((E_1, \omega_1), (E_2, \omega_2))$ such that the product
\begin{equation}
  \label{eq:defofsum}
  (E_1, \omega_1)\oplus (E_2, \omega_2) := (E_1\times E_2, \omega_1 + \omega_2)
\end{equation}
is an eigenform for real multiplication by $\ord$.  Let
$$\Q = \proj\Omega\Q.$$

We now describe the connected components of $\Omega\Q$ following \cite{mcmullenspin}.  A \emph{prototype} for
real multiplication by $\ord$ is a triple of integers $(e, l, m)$ such that
$$D=e^2+4l^2m, \quad l,m>0, \quad \text{and } \gcd(e, l)=1.$$
To each such prototype, we associate a
\emph{prototypical eigenform} $\mathcal{Q}(e,l,m)$ as follows.  Let $\lambda$ be the unique positive root of
$\lambda^2-e\lambda-l^2m=0$, and consider the lattices,
\begin{align}
  \label{eq:lattices}
  \Lambda_1 &= \zed(\lambda, 0) \oplus \zed(0, \lambda) \\
  \notag
  \Lambda_2 &= \zed(lm, 0)\oplus \zed(0, l).
\end{align}
Let $E_i=\cx/\Lambda_i$, and equip $E_i$ with the form $\omega_i$ covered by $dz$ on $\cx$.  Let
$$\mathcal{Q}(e,l,m) = (E_1, \omega_1)\oplus (E_2, \omega_2).$$

From the description of $\Omega Q$ in \eqref{eq:omegaQ}, we obtain an action of $\GLtwoRplus$ on $\Omega Q$
induced by the diagonal action of $\GLtwoRplus$ on $\GLtwoRplus \times \GLtwoRplus$.

\begin{theorem}[{\cite[Theorem~2.1]{mcmullenspin}}]
  \label{thm:Qcomponents}
  The locus $\Omega\Q$ is invariant under the action of $\GLtwoRplus$.  Each component is a $\GLtwoRplus$
  orbit which contains exactly one prototypical eigenform $\mathcal{Q}(e,l,m)$ and is isomorphic to
  $\GLtwoRplus/\Gamma_0(m).$
\end{theorem}

\begin{remark}
  Here $\Gamma_0(m)$ is the congruence subgroup of $\SLtwoZ$ defined by
  $$\Gamma_0(m) = \left\{
  \begin{pmatrix}
    a & b \\
    c & d
  \end{pmatrix}
  \in \SLtwoZ : c\equiv 0 \pmod m \right\}.
$$
\end{remark}

\begin{theorem}
  \label{thm:chiQ}
  For any nonsquare discriminant $D$, we have
  \begin{equation}
    \label{eq:chiQ}
    \chi(\Q[D]) = -5\chi(\X).
  \end{equation}
  If $D=d^2$ with $d\geq 2$, then we have
  \begin{equation}
    \label{eq:chiQ2}
    \chi(\Q[d^2]) = -\frac{1}{72}d^2(5d-6)\sum_{r|d}\frac{\mu(r)}{r^2}.
  \end{equation}
\end{theorem}

\begin{proof}
  We first claim that
  \begin{equation}
    \label{eq:sumQ}
    \sum_{s|f}\chi(\Q[s^2D]) = -\frac{1}{6} \sum_{\substack{
        e\equiv f^2D \: (2) \\
        -f\sqrt{D}<e<f\sqrt{D}
      }}
    \sigma_1\left(\frac{f^2D-e^2}{4}\right)
  \end{equation}
  for any fundamental discriminant $D$ and $f\in\nats$.  By Theorem~\ref{thm:Qcomponents}, there is one
  component of $\Q$ for each prototype $(e,l,m)$ which is isomorphic to $\half/\Gamma_0(m)$.
  By \cite[Theorem~4.2.5]{miyake},
  $$\chi(\half/\Gamma_0(m)) = \psi(m),$$ where
  $$\psi(m) = -\frac{1}{6}m\prod_{\substack{
      p|m \\
      p \text{ prime}
    }}
  \left(1+\frac{1}{p}\right).$$
  It is elementary to show that for any $n\in\nats$,
  \begin{equation*}
    -\frac{1}{6}\sigma_1(n) = \sum_{\substack{
        n = l^2m\\
        l,m\in\nats
      }}
    \psi(m).
  \end{equation*}

  By Theorem~\ref{thm:Qcomponents},
  \begin{align*}
    \sum_{r|f}\chi(Q_{f^2D}) &= \sum_{\substack{
        f^2D = e^2+4l^2m \\
        l,m>0
      }}
    \psi(m)\\
    &= -\frac{1}{6}\sum_{\substack{
         e\equiv f^2D \: (2) \\
        -f\sqrt{D}<e<f\sqrt{D}
      }}
    \sigma_1\left(\frac{f^2D-e^2}{4}\right),
  \end{align*}
  which proves \eqref{eq:sumQ}.

  If $D\neq 1$ is a fundamental discriminant, then for any $f\in\nats$, we have
  \begin{equation*}
    \sum_{r|f}\chi(\Q[r^2D]) = \frac{5}{6} H(2, f^2 D) = -5\sum_{r|f}\chi(\X[r^2D])
  \end{equation*}
  by Theorem~\ref{thm:siegelformula} and \eqref{eq:sumQ}, which proves \eqref{eq:chiQ}.

  Now for any $d\in\nats$, we have
  \begin{equation}
    \label{eq:meow1}
    \frac{1}{72}+\sum_{\substack{
        r|d \\
        r>1}}
    \chi(\Q[r^2]) = \frac{5}{6}H(2, d^2) + \frac{d^2}{12}
  \end{equation}
  by Theorem~\ref{thm:siegelformula} and \eqref{eq:sumQ}, and using the convention that $\sigma_1(0) =
  -1/24$.  From \eqref{eq:h2dsum}, we get for any $d\in\nats$,
  \begin{equation}
    \label{eq:meow2}
    H(2, d^2) = \sum_{s|d}\left(-\frac{1}{12}s^3\sum_{r|s}\frac{\mu(r)}{r^2}\right).
  \end{equation}
  It follows from M\"obius inversion that for any $d\in\nats$,
  \begin{equation}
    \label{eq:meow3}
    d^2 = \sum_{s|d}\left(s^2\sum_{r|s}\frac{\mu(r)}{r^2}\right).
  \end{equation}
  Define
  \begin{equation*}
    f(d) =
    \begin{cases}
      1/72 & \text{if $d=1$} \\
      \chi(\Q[d^2]) & \text{if $d>1$}.
    \end{cases}
  \end{equation*}
  Combining \eqref{eq:meow1}, \eqref{eq:meow2}, and \eqref{eq:meow3}, we obtain
  \begin{equation*}
    \sum_{s|d}f(s) = \sum_{s|d}\left(\left(-\frac{5}{72}s^3+\frac{1}{12}s^2\right)\sum_{r|s}\frac{\mu(r)}{r^2}\right),
  \end{equation*}
  which implies \eqref{eq:chiQ2} by M\"obius inversion.
\end{proof}

There is an involution of $\Omega \Q$ which interchanges the order of the two factors in \eqref{eq:omegaQ}.
This involution preserves $\Omega\Q$ and so restricts to an involution $\tau$ of $\Omega\Q$.  The involution
$\tau$ preserves the fibers of the bundle $\Omega\Q\to\Q$, so it induces an involution $\overline{\tau}$ of
$\Q$.

\begin{lemma}
  \label{lem:taufixes}
  The involution $\overline{\tau}$ fixes pointwise the component of $\Q$ containing the prototypical
  eigenform $\mathcal{Q}(e,l,m)$ if and only if
  \begin{equation}
    \label{eq:taufixesQ}
    \overline{\tau} \mathcal{Q}(e,l,m) =  \mathcal{Q}(e,l,m).
  \end{equation}
  The only such component on which $\tau$ is the identity is the one containing $\mathcal{Q}(0,1,1)$.
\end{lemma}

\begin{proof}
  First, note that \eqref{eq:taufixesQ} holds if and only if $\Lambda_1=\Lambda_2$ in \eqref{eq:lattices},
  which is true if and only if $m=1$ and $\lambda=l$.  The only prototype for which this holds is $(0,1,1)$,
  so \eqref{eq:taufixesQ} is true if and only if $(e,l,m)=(0,1,1)$.

  Clearly, \eqref{eq:taufixesQ} must hold for $\overline{\tau}$ to be the identity on the component containing
  $\mathcal{Q}(e,l,m)$.  Conversely, if \eqref{eq:taufixesQ} does hold, then $(e,l,m)=(0,1,1)$, and
  \begin{equation*}
    \tau(g\cdot \mathcal{Q}(0,1,1)) = g\cdot \mathcal{Q}(0,1,1)
  \end{equation*}
  for any $g\in\GLtwoRplus$ because $\tau$ commutes with the action of $\GLtwoRplus$ on $\Omega\Q$.  Thus
  $\tau$ fixes every point of the $\GLtwoRplus$-orbit of $\mathcal{Q}(0,1,1)$, which is exactly the component
  of $\Omega\Q[4]$ containing $\mathcal{Q}(0,1,1)$.
\end{proof}

There is natural map $\pi\colon\Omega\Q \to \Omega\P$ defined by
\begin{equation}
  \label{eq:defofpi}
  \pi((E_1, \omega_1), (E_2,\omega_2)) = (E_1, \omega_1)\oplus(E_2, \omega_2)
\end{equation}
which descends to a map $\overline{\pi}\colon\Q\to\P$.
 
\begin{lemma}
  \label{lem:tauquotient}
  The map $\overline{\pi}\colon\Q\to\P$ factors through to an isomorphism of orbifolds,
  $\tilde{\pi}\colon\Q/\overline{\tau}\to\P$.
\end{lemma}

\begin{proof}
  The map $\overline{\pi}$ factors because the right hand side of \eqref{eq:defofpi} does not depend on the
  order of the $E_i$.  It is clearly onto.  That $\overline{\pi}$ is one-to-one follows from the following
  fact about Abelian surfaces:  a principally polarized Abelian surface can have at most one representation as
  a polarized product of elliptic curves.
\end{proof}

\begin{theorem}
  \label{thm:chiP}
  If $D$ is not square, then
  \begin{equation*}
    \chi(\P) = -\frac{5}{2}\chi(\X).
  \end{equation*}
  If $D=d^2$, then
  \begin{equation*}
    \chi(\P[d^2])=-\frac{1}{144} d^2(5 d -6)\sum_{r|d}\frac{\mu(r)}{r^2}
  \end{equation*}
  when $d>2$, and
  \begin{equation*}
    \chi(\P[4]) = -\frac{1}{6}.
  \end{equation*}
\end{theorem}

\begin{proof}
  By Theorem~\ref{thm:Qcomponents}, $\Q[4]$ consists of a single component, which contains
  $\mathcal{Q}(0,1,1)$.  Thus by Lemmas~\ref{lem:taufixes} and \ref{lem:tauquotient}, $\P[4]\isom\Q[4]$, and
  $\Q$ is a twofold cover of $\P$ when $D>4$.  The claim then follows from Theorem~\ref{thm:chiQ}.
\end{proof}


%% file: prototypes.tex
\section{Prototypes}
\label{sec:prototypes}

\subsection{$P_D$, $W_D$, and $Y_D$-prototypes}
\label{subsec:prototypes}

Curves in $\Y\setminus\X$ will be classified by certain discrete numerical
invariants which we call $\Y$-prototypes.  They are almost (but not quite) the same as the \emph{splitting
  prototypes} introduced by McMullen in \cite{mcmullenspin} to classify the cusps of $\W$.  

\begin{definition}
  A \emph{$\Y$-prototype} of discriminant $D$ is a quadruple, $(a, b, c, \bar{q})$ with $a, b,  
  c\in\integers$ and $\bar{q}\in\integers/\gcd(a, b, c)$ which satisfies the following six properties:
  \begin{enumerate}
  \item $b^2-4 a c=D$
  \item $a>0$
  \item $c \leq 0$
  \item $\gcd(a, b, c, \bar{q})=1$
  \item $a+b+c\leq 0 $
  \item $a+b+c$ and $c$ are not both zero.
  \end{enumerate}
\end{definition}

Let $\Yprot$ be the set of $\Y$-prototypes.  Define a map $\lambda\colon\Yprot\to K_D$ by 
associating to each prototype $P=(a, b, c, \bar{q})$ of discriminant $P$ the unique algebraic number
$\lambda(P)\in K_D$ such that $a \lambda(P)^2 + b\lambda(P) + c=0$ and $\lambda^{(1)}>0$.
This makes sense because if $c<0$, then the two roots of $ax^2+bx+c=0$ have opposite signs, and if $c=0$, then
the condition that $a+b<0$ implies that the nonzero root of $a \lambda^2 + b\lambda =0$ is positive. It is easy to check that
the last condition $a + b + c\leq0$ is equivalent to $\lambda^{(1)} \geq 1$.

We say that a prototype is \emph{terminal} if $a+b+c=0$; it is \emph{initial} if $a-b+c=0$; and it is
degenerate if $c=0$. Terminal, initial, and degenerate prototypes only arise if $D$ is square.  We define
the involution,
\begin{equation}
  \label{eq:involutionone}
  (a, b, c, \bar{q})\mapsto (-c, -b, -a, \bar{q}),
\end{equation}
on the set of terminal prototypes.  We consider two terminal prototypes to be the same if they are related
by this involution.  Similarly, we consider two degenerate prototypes to be the same if they are related by
the involution,
\begin{equation}
  \label{eq:involutiontwo}
  (a, b, 0, \bar{q}) \mapsto (-b - a, b, 0, \bar{q})
\end{equation}
on the set of degenerate prototypes.

\paragraph{Operations on $\Y$-prototypes.}
Given a nonterminal prototype $P$, define the \emph{next prototype} $P^+$ by
\begin{align*}
  P^+&=
  \begin{cases}
    (a, 2a+b, a+b+c, \bar{q}), &\text{if } 4a + 2b + c \leq 0;\\
    (-a-b-c, -2a-b, -a, \bar{q}), &\text{if } 4a + 2b + c \geq 0.
  \end{cases}\\
  \intertext{Given a nondegenerate prototype $P$, define the \emph{previous prototype} $P^-$ by}
  P^-&=
  \begin{cases}
    (a, -2a+b, a-b+c, \bar{q}), &\text{if } a-b+c\leq 0;\\
    (-c, -b+2c, -a+b-c, \bar{q}), &\text{if } a-b+c \geq 0.
  \end{cases}
\end{align*}
Its easy to check that $P^+$ and $P^-$ are actually prototypes of the same discriminant and that
$$(P^+)^-=(P^-)^+=P$$
when these operations are defined.

On the level of $\lambda(P)$,  we have,
\begin{align*}
  \lambda(P^+)&=
  \begin{cases}
    \lambda(P)-1, & \text{if }
    \lambda(P)^{(1)}-1\geq1;\\
    1/(\lambda(P)-1), & \text{if } \lambda(P)^{(1)}-1<1;
  \end{cases}\\
  \intertext{and} \lambda(P^-)&=
  \begin{cases}
    \lambda(P)+1, & \text{if }    N^{K_D}_\ratls(\lambda(P)+1)\leq 0;\\
    (\lambda(P)+1)/\lambda(P), & \text{if } N^{K_D}_\ratls(\lambda(P)+1)>0.
  \end{cases}
\end{align*}

Define an involution $t$ on the set of prototypes of discriminant $D$ by
$$
t(a, b, c, \bar{q})=
\begin{cases}
  (a, -b, c, \bar{q}), & \text{if } a - b + c \leq
  0;\\
  (-c, b, -a, \bar{q}), & \text{if } a - b + c \geq 0.
\end{cases}
$$

Define the \emph{multiplicity} of a nondegenerate prototype $P$ by
$$\mult(P)=\frac{\gcd(a, c)}{\gcd(a, b, c)}.$$

\paragraph{$\P$ and $\W$-prototypes.} In addition to the $\Y$-prototypes, we will define similar objects which
we call $\P$ and $\W$-prototypes.  These will classify the cusps of the curves $\P$ and $\W$.

\begin{definition}
   A \emph{$\P$-prototype} of discriminant $D$ is a quadruple, $(a, b, c, \bar{q})$ with $a, b,  
  c\in\integers$ and $\bar{q}\in\integers/\gcd(a,  c)$ which satisfies the following five properties:
  \begin{enumerate}
  \item $b^2-4 a c=D$
  \item $a>0$
  \item $c < 0$
  \item $\gcd(a, b, c, \bar{q})=1$
  \item $a+b+c\leq 0 $
  \end{enumerate}
  Note that the only difference between $\Y$ and $\P$-prototypes are that for $\P$-prototypes, $\bar{q}$
  ranges in a potentially larger group, and $c$ is not allowed to be zero.
  
   A \emph{$\W$-prototype} of discriminant $D$ is a quadruple, $(a, b, c, \bar{q})$ with $a, b,  
  c\in\integers$ and $\bar{q}\in\integers/\gcd(a,  c)$ which satisfies the same properties as above except
  that the last one is replaced with
  $$a+b+c<0.$$
\end{definition}

Let $\Pprot$ and $\Wprot$ be the sets of $\P$ and $\W$-prototypes.  Just as for the $\Y$-prototypes, we have
natural maps $\lambda\colon\Wprot\to K_D$ and $\lambda\colon \Pprot\to K_D$.  By reducing $\bar{q}$
modulo $\gcd(a, b, c)$, there are natural maps $\Pprot\to\Yprot$ and $\Wprot\to \Yprot$.

We will see in \S\ref{subsec:curvesinY} that these maps encode intersections of cusps of $\W$ and $\P$ with
curves in $\Y \setminus \X$.

\paragraph{Relation with splitting prototypes.}

In \cite{mcmullenspin}, McMullen defined a splitting prototype to be a quadruple of integers $(a, b, c, e)$.
Splitting prototypes correspond bijectively to our $\W$-prototypes by sending the splitting prototype $(a, b,
c, e)$ to the $\W$-prototype $(c, e, -b, \bar{a})$.

\subsection{Quasi-invertible $\mathcal{O}_D$-modules}
\label{subsec:quasiinvertiblemodules}

In this section, we will study a class of $\ord$-modules which arises naturally in the study of our
compactification of the Hilbert modular surface.  We will also study a class of bases of these modules whose
combinatorics is closely related to the geometry of the compactification.  This material will only be used in
\S\ref{sec:limitsofeigenforms}.

\begin{definition}
  An $\ord$-module $M$ is \emph{quasi-invertible} if $M\isom \zed\oplus\zed$ as an Abelian group, and $M$ contains
  some element $x$ such that $\Ann(x) =0$.
\end{definition}

\paragraph{Examples.}

If $D$ is a fundamental discriminant, then every quasi-invertible $\ord$-module is actually invertible.  Since
$\ord$ is a Dedekind domain, the set of isomorphism classes of such modules forms a group, the ideal class
group of $\ord$.

If $D$ is not square, then a quasi-invertible $\ord$-module is just an invertible $\ord[E]$-module over some order
$\ord[E]$ containing $\ord$. A quasi-invertible $\ord$-module is invertible if and only if it is primitive in
the sense that it is not also a module over any $\ord[E]$ containing $\ord$.

In this paper, quasi-invertible $\ord$-modules will arise in the following way:

\begin{prop}
  \label{prop:lagrangianquasiinvertible}
  If $M\subset \ord\oplus\ord^\vee$ is a Lagrangian $\ord$-submodule which has $\zed$-rank two, then $M$ is quasi-invertible.
\end{prop}

\begin{proof}
  This is trivial if $D$ is not square because then $\ord\oplus\ord^\vee$ is torsion-free, so suppose $D=d^2$.
  Let $M \subset \ord[d^2]\oplus\ord[d^2]^\vee$ be a rank two Lagrangian submodule which is not
  quasi-invertible, and let $\{u_1, u_2\}$ be a basis of $M$.  We must then have either $(d, 0)\cdot u_i=0$
  for $i=1, 2$, or $(0, d)\cdot u_i =0$ for $i=1, 2$, or else there would be some $x\in M$ such that
  $\Ann(x)=0$.  For concreteness, suppose we are in the second case.  Then we must have
  $$u_i=((a_i, 0), (b_i, 0))$$
  for some $a_1, b_i\in\ratls$.

  Since $M$ is Lagrangian, we know that,
  \begin{align*}
    0 = \langle u_1, u_2 \rangle &= \tr_\ratls^{\ratls\oplus\ratls} 
      \begin{vmatrix}
        (a_1, 0) & (a_2, 0) \\
        (b_1, 0) & (b_2, 0)
      \end{vmatrix}\\
      &=a_1 b_2 - a_2 b_2,
  \end{align*}
  so $u_1 = r u_2$ for some $r\in \ratls$.  This contradicts the fact that the $u_i$ are a basis.
\end{proof}

\paragraph{Embeddings in $K_D$.}

A \emph{lattice} in $K_D$ is a rank two Abelian subgroup.  Given a lattice $M\subset K_D$, the set
$$\{ x\in K_D : xM\subset M\}$$
is called the \emph{coefficient ring} of $M$.  It is an order in $K_D$.

\begin{prop}
  \label{prop:quasi-invertibleembedding}
  An $\ord$-module $M$ is quasi-invertible if and only if it is isomorphic to some lattice $M\subset K_D$
  whose coefficient ring contains $\ord$.

  If $M$ is quasi-invertible, then the embedding $M\to K_D$ is unique up to multiplication by a
  non-zero-divisor in $K_D$.
\end{prop}

\begin{proof}
  If $D$ is not square, then it is clear that a lattice $M\subset K_D$ whose coefficient ring contains
  $\ord$ is a quasi-invertible $\ord$-module because $M$ is automatically torsion-free.  If $D=d^2$, then $M$
  could have  $\Ann(x)\neq 0$ for every $x\in M$ only if $M \subset \ratls \oplus \{0\}$ or $M\subset
  \{0\}\oplus\ratls$.  This cannot happen because $M$ has rank two.

  Conversely, suppose $M$ is a quasi-invertible $\ord$-module with $x\in M$ such that $\Ann(x)=0$.  Since $M$
  is torsion-free as an Abelian group, the natural map $M\to M\otimes\ratls$ is injective.  The tensor product
  $M\otimes\ratls$ is a $K_D$-module and a two-dimensional vector space over $\ratls$, so the map $K\to
  M\otimes \ratls$ defined by $r \mapsto r\cdot x$ is an isomorphism.  The inverse of this map embeds $M$ as a
  lattice whose coefficient ring contains $\ord$.  Since the image of $x$ determined the embedding, it is
  unique up to constant multiple.
\end{proof}

It is easy to identify the coefficient ring of a lattice in $K_D$.  Let $M_\lambda$ be the lattice
generated by $\{1, \lambda\}$.  For $\lambda\in K_D$, let $\phi_\lambda(t)$ be the minimal polynomial of
$\lambda$, the unique polynomial $\phi_\lambda(t)=at^2+bt+c$ such that $a, b, c\in \zed$, $a>0$,
$\phi_\lambda(\lambda)=0$, and $\gcd(a, b, c)=1$.

\begin{prop}[{\cite[p.136]{bs}}]
  \label{prop:coefficientring}
  If $\lambda\in K_D\setminus\ratls$ with $\phi_\gamma(t)=at^2+bt+c$, then the coefficient ring of
  $M_\lambda$ is the order $\zed[a\lambda]$, which is isomorphic to $\ord$, where $D=b^2-4ac$.
\end{prop}

\paragraph{Admissible bases.}

Since a quasi-invertible $\ord$-module $M$ can be embedded in $K_D$ uniquely up to constant multiple, for
any $u$, $v\in M$ with $\Ann(v)\neq 0$, the ratio $u/v$ is a well-defined element of $K_D$.

\begin{definition}
  A basis $\{u, v\}$ of a quasi-invertible $\ord$-module $M$ is
  \emph{admissible} if $\Ann(v)=0$ and
  $N_\ratls^{K_D}<0$.
\end{definition}

We will consider two admissible bases of $M$ to be equivalent if they are the same as subsets of $M/\pm 1$.  The following
proposition gives a complete classification of quasi-invertible $\ord$-modules together with an admissible basis.

\begin{prop}
  \label{prop:admissibleclassification}
  Every pair $(M, \{u, v\})$, where $M$ is a quasi-invertible $\ord$-module together with an admissible basis $\{u, v\}$,
  is equivalent to one of the form $(M_\lambda, \{1,\lambda\})$, where $\lambda^{(1)}\geq1$ and
  $\psi(\lambda)=0$, where
  $$\psi(t)=at^2 + bt + c$$
  for some integers $(a, b, c)$ such that:
  \begin{enumerate}
  \item $b^2-4 a c=D$
  \item $a>0$
  \item $c < 0$
  \item $a+b+c\leq 0$
  \item $\ord = \zed[a \lambda]$
  \end{enumerate}

  Pairs $(M_\lambda, \{1, \lambda\})$ and $(M_\mu, \{1, \mu\})$ of this form are equivalent if and only if $\lambda=\mu^{-1}$.
\end{prop}

\begin{proof}
  If $\lambda>1$, then by Proposition~\ref{prop:quasi-invertibleembedding}, there is a unique embedding $M\to K_D$ which -- after
  possibly switching $u$ and $v$ and changing their signs -- sends $u$ to 1 and $v$ to some $\lambda\in K$ with
  $\lambda^{(1)}>1$.  The coefficient ring of the image $M_\lambda$ is $\ord[E]$, where $D=s^2E$ for some
  quadratic discriminant $E$.  Let
  $\phi_\lambda(t)= a't^2 + b't+c'$ be the minimal polynomial of $\lambda$.  By
  Proposition~\ref{prop:coefficientring}, we have $b'^2-4a'c'=E$ and $\ord[E]=\zed[a'\lambda]$.  We have
  $c'<0$ because $N_\ratls^{K_D}<0$.  Since $\lambda^{(1)}>1$, we know that $a'+b'+c'\leq0$.  The
  integers $(a, b, c) = s(a', b', c')$ then have the required properties.  

  If $\lambda^{(1)}=1$, then there are two embeddings $M\to K_D$ as above, one sending $u$ to $1$ and $v$ to
  $\lambda$; the other sending $u$ to $\lambda^{-1}$ and $v$ to $1$.  As in the above paragraph, both
  embeddings yield a presentation of $(M, \{u, v\})$ in the desired form.
\end{proof}

We can define a map from the set of nondegenerate $\Y$-prototypes to the set of isomorphism classes of
admissible bases of quasi-invertible $\ord$-modules by sending $P$ to $(M_{\lambda(P)}, \{1, \lambda(P)\})$.
This proposition implies that this map is onto and that two prototypes are sent to the same admissible basis if and
only if their integers $a$, $b$, and $c$ are the same.

\paragraph{Characterization of admissible bases.}

We now describe a useful characterization of admissible bases.  Given two quasi-invertible $\ord$-modules $M$
and $N$, a perfect pairing is a bilinear map
$$h\colon M\times N\to\zed$$
such that
$$h(\lambda\cdot x, y) = h(x, \lambda\cdot y)$$
for each $x\in M$, $y\in N$, and $\lambda\in\ord$, and such that the induced map $M\to\Hom(N, \zed)$ is an
isomorphism.

Given such a perfect pairing and a basis $\{u, v\}$ of $M$, let $\{u^*, v^*\}$ be a dual basis of $N$, and
define
$$\sign(u,v;h) = \sign(\iota_M(u)\iota_M(v) \iota_N(u^*) \iota_N(v^*)),$$
where $\iota_M$ and $\iota_N$ are $\iota_1$-linear, nonzero maps $M\to\reals$ and $N\to\reals$.  Note
that $\sign(u, v;h)$ is invariant under changing the sign of $u$ or $v$ as well as independent of the choice
of $\iota_M$ and $\iota_N$.

\begin{theorem}
  \label{thm:criterionforadmissible}
  Given a quasi-invertible $\ord$-module $M$, a basis $\{u, v\}$ is an admissible basis if and only if
  $\sign(u, v; h)=1$ for some (or any) perfect pairing $h\colon M\times N\to\zed$.
\end{theorem}

\begin{proof}
  Assume that $M$ is embedded in $K_D$.
  If $u^{(1)}=0$ or $v^{(1)}=0$, then $\sign(u, v;h)=0$, and the basis is not admissible, so we are
  done.  Since $M$ is quasi-invertible, it can't happen that $u^{(2)}=v^{(2)}=0$, so assume $u^{(2)}\neq 0$.
  
  Assume without loss of generality that $u=1$.  By possibly changing the sign of $v$, we can suppose $v^{(1)}>0$.  Then
  $v^{(2)}<0$ if and only if $\{u,v\}$ is admissible.  Since $N$ and the perfect pairing $h$ are unique up to
  isomorphism, we can suppose without loss of generality that $N=M^\vee$, the inverse different of $M$, and
  that $h$ is the trace pairing $M\times M^\vee\to\zed$.

  Define
  \begin{align*}
    u^* &= \frac{-v'}{v-v'}\\
    v^* &= \frac{1}{v-v'}.
  \end{align*}
  It is easy to check that $\{u^*, v^*\}$ is dual to $\{u, v\}$ with respect to the trace pairing.  If $\{u,
  v\}$ is not admissible, then $v^{(2)}<0$, in which case $(u^*)^{(1)}, (v^*)^{(1)}>0$, so $\sign(u,v;h)>0$.
  If $\{u, v\}$ is not admissible, then $v^{(2)}\geq 0$.  If $v^{(2)}>0$, then $(u^*)^{(1)}$ and $(v^*)^{(1)}$
  have opposite signs, so $\sign(u, v;h)=-1$.  If $v^{(2)}=0$, then $(u^*)^{(1)}=0$, so $\sign(u,v;h)=0$.
\end{proof}

\paragraph{Admissible triples.}

We are also interested in special triples of elements of quasi-invertible $\ord$-modules which will arise
naturally in the study of the compactification of the Hilbert modular surface.

\begin{definition}
  An \emph{admissible triple} in a quasi-invertible $\ord$-module $M$ is an unordered triple of elements $\{a,
  b, c\}\subset M/\pm 1$ such that:
  \begin{itemize}
  \item $\pm a \pm b \pm c =0$ for some choice of signs. 
  \item Some pair of elements of the triple form an admissible basis of $M$.
  \end{itemize}
\end{definition}

\begin{prop}
  \label{prop:tripleclassification}
  Every pair $(M, T)$, where $M$ is a quasi-invertible $\ord$-module with an admissible triple $T\subset M$,
  is equivalent to one of the form
  $$(M_\lambda, \{1,\lambda, \lambda-1\}),$$
  where $\lambda^{(1)}\geq1$ and
  $$a\lambda^2 + b\lambda + c =0$$ for some integers $(a, b, c)$ such that:
  \begin{enumerate}
  \item $b^2-4 a c=D$
  \item $a>0$
  \item $c \leq 0$
  \item $a+b+c\leq 0$
  \item $a+b+c$ and $c$ are not both zero.
  \item $\ord = \zed[a \lambda]$
  \end{enumerate}

  Two such pairs $(M_\lambda, \{1, \lambda, \lambda-1\})$ and $(M_\mu,\{1, \mu, \mu-1\})$ are equivalent if
  and only if $\lambda^{(1)}=1$ and $\mu=\lambda^{-1}$ or if $\lambda^{(2)}=0$ and $\mu = \lambda/(\lambda-1)$.
\end{prop}

\begin{proof}
  Let $\{u, v, w\}\subset M$ be an admissible triple and assume that $M$ is embedded in $K_D$.  This triple
  contains at least one admissible basis, so assume without loss of generality that $\{u, v\}$ is an
  admissible basis.  There are now three cases to consider, depending on whether $\Ann(w)=0$, $w^{(1)}=0$,
  or $w^{(2)}=0$.

  First suppose that $\Ann(w)=0$.  Then $\{u, v, w\}$ contains exactly two admissible bases because
  $$N^{K_D}_\ratls\left(\frac{u}{v}\right)N^{K_D}_\ratls\left(\frac{v}{w}\right)N^{K_D}_\ratls\left(\frac{w}{u}\right)=1.$$
  Suppose without loss of generality that $u$ is contained in both bases.  Dividing everything by $u$, we can
  assume $u=1$.  By possibly switching $v$ and $w$ and changing their signs, we can put the triple uniquely in
  the form $\{1, \lambda, \lambda-1\}$ with $\lambda^{(1)}>0$.  Let $(a, b, c)$ be as in Proposition~\ref{prop:coefficientring}
  applied to $M_\lambda$.   Since $\lambda^{(1)}\lambda^{(2)}<0$ and
  $(\lambda^{(1)}-1)(\lambda^{(2)}-1)<0$, we must actually have $\lambda^{(1)}> 1$, and it follows that $a+b+c<0$.  The other
  necessary properties of $a$, $b$, and $c$ follow from Proposition~\ref{prop:coefficientring}.

  Now assume $w^{(1)}=0$.  In this case, $K_D=\ratls\oplus\ratls$, and there are two ways two put the triple
  in the form, $\{1, \lambda, \lambda-1\}$, where $\lambda=(1, s)$  with $s<0$: either divide $\{u, v, w\}$ by
  $u$, or divide by $v$ and swap the first two elements.  With $(a, b, c)$ as in
  Proposition~\ref{prop:coefficientring} applied to $M_\lambda$, it follows that $a+b+c=0$ and $c<0$.  The other necessary
  properties are clear, and the two $\lambda$ which arise in this way are $u/v$ and $v/u$, so they are
  reciprocal. 

  Now assume $w^{(2)}=0$.  Dividing the triple by $u$ and switching $v$ and $w$ puts it in the form $\{1,
  \lambda, \lambda-1\}$ with $\lambda=(s,0)$.  Since $\lambda-1 = (s-1, -1)$ has negative norm, we must have
  $s>1$.  It follows that if $(a, b, c)$ are as in Proposition~\ref{prop:coefficientring} applied to $M_\lambda$, we must have
  $c=0$ and $a+b<0$.  Alternatively, we could have divided the triple by $v$ and rearranged the elements to
  put it in the form $\{1, \mu, \mu-1\}$ with $\mu=(t, 0)$.  Then $t>1$ by the same argument, and $t=s/(s-1)$.

  The converse statement that triples of the given form are admissible is not hard and will be left to the reader.
\end{proof}

Just as for admissible pairs, we can define a map from the set of $\Y$ prototypes to the set of isomorphism classes
of quasi-invertible $\ord$-modules with admissible triples, sending $(a, b, c, \bar{q})$ to $(M_\lambda, \{1,
\lambda, \lambda-1\})$ as above.  Since we identified two terminal or degenerate prototypes if they are
related by the involutions \eqref{eq:involutionone} and \eqref{eq:involutiontwo}, two prototypes have the same
image if and only if they have the same $\bar{q}$.


%% file: abeliandifferentials.tex
\section{Abelian differentials}
\label{sec:abelian}

In this section, we recall known material about Abelian differentials (holomorphic one-forms) on compact
Riemann surfaces.  In \S\ref{subsec:flatgeometry}, we discuss the flat geometry associated to an Abelian
differential.  In \S\ref{subsec:abelianmoduli}, we discuss moduli spaces of Abelian differentials and the
action of $\SLtwoR$ on these spaces.  In \S\ref{subsec:action}, we discuss McMullen's results on the dynamics
of this action in genus two.

\subsection{Flat geometry of Abelian differentials}\label{subsec:flatgeometry}

A Riemann surface with a nonzero Abelian differential (holomorphic one-form) carries a
canonical flat geometry, which is closely related to the study of billiards in rational angled polygons, as
well as the study of flows on moduli space.  In this paper, this geometry will be useful because it will allow
us to deform Abelian differentials using some concrete cut-and-paste operations which we will describe in this
section.

\paragraph{Translation surfaces.}
Let $X$ be a Riemann surface with a nonzero Abelian differential $\omega$, and let $Z(\omega)$ be the
discrete set of zeros of $\omega$.  Every point of $X$ which is not a zero of $\omega$ has a neighborhood $U$
and a conformal map $\phi\colon U\to\cx$ such that $\omega|_U=\phi^*(dz)$. The conformal map $\phi$ can be
defined explicitly by 
$$\phi(z)=\int_{z_0}^z \omega$$
for a choice of base point $z_0$. These coordinates $\phi$ are unique up to translation by a constant. 

We can also put $\omega$ in a standard form in the neighborhood of a zero or a pole at $z_0$.  There is a
neighborhood $U$ of $z_0$ and conformal map 
$$\phi_0\colon (U, z_0)\to(\cx, 0)$$ 
such that
\begin{equation}
  \label{eq:conecoordinates}
  \omega|_U=
  \begin{cases}
    \phi_0^*(z^n dz), &\text{if $z_0$ is not a simple pole;}\\
    \phi_0^*(a\, dz/z), &\text{if $z_0$ is a simple pole with nonzero residue $a$.}
  \end{cases}
\end{equation}

Since the local coordinates $\phi$ away from the zeros of $\omega$ are unique up to translation by a constant,
any translation invariant geometric structure on $\cx$ is inherited by $X$.  In particular, $X$ inherits a
flat metric and an oriented foliation $\fhoriz$ coming from the foliation of the plane by horizontal lines.
More generally, for slope $s\in\proj^1(\reals)$, the surface $X$ has an orientable foliation $\mathcal{F}_s$
coming from the foliation of $\cx$ by lines of slope $s$.

In terms of the Abelian differential $\omega$, the flat metric is just $|\omega|$.  A vector $v$ is tangent
to $\fhoriz$ in the positive direction if $\omega(v)>0$.

\begin{figure}[htbp]
  \centering
  \includegraphics{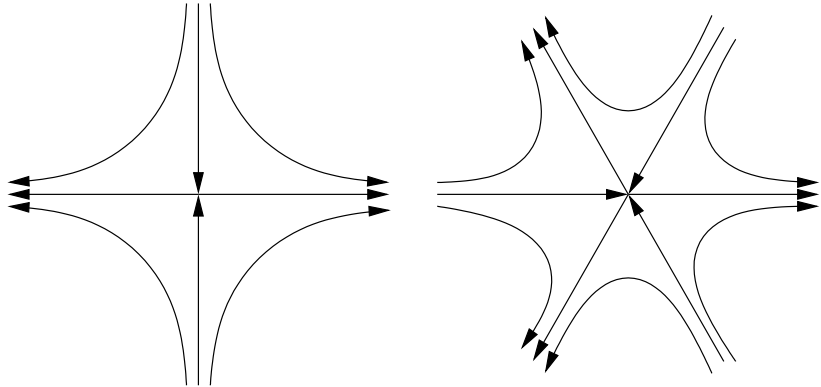}
  \caption{Zeros of order one and two.}
  \label{fig:zeros}
\end{figure}

The flat metric has singularities at the zeros.  From the coordinates \eqref{eq:conecoordinates}, we see that a zero of
order $n$ has a neighborhood isometric to a cone with cone angle $2 \pi(n+1)$, and the foliation $\fhoriz$ has
$2n+2$ leaves meeting at the zero.  Zeros of orders one and two together with the foliation
$\fhoriz$ are pictured in Figure \ref{fig:zeros}.

To summarize, we have associated to an Abelian differential $(X, \omega)$, a flat metric on $X\setminus
Z(\omega)$ together with a horizontal foliation $\fhoriz$ which is parallel with respect to a metric such that
the points of $Z(\omega)$ are cone singularities of the metric.  A surface with such a structure is sometimes
called a \emph{translation surface}.  A translation surface is equivalent to a Riemann surface with a nonzero
Abelian differential.

A geodesic on $(X, \omega)$ is called \emph{straight} if it does not pass through any zeros of $\omega$.  A
straight geodesic which joins two zeros is called a \emph{saddle connection}.  Any straight, closed geodesic
is contained in a \emph{cylinder} on $X$ by taking nearby parallel geodesics.  If the genus of $X$ is greater
than one, any cylinder can be extended until either end contains a zero of $\omega$ and each boundary
component is a finite union of saddle connections.  Such a cylinder is called a \emph{maximal cylinder}.

\paragraph{Plumbing a cylinder.}

Cylinders arise from simple poles: by \eqref{eq:conecoordinates}, a simple pole of an Abelian differential has
a neighborhood which is isometric to a half-infinite cylinder.
Given a meromorphic Abelian differential $(X, \omega)$ with simple poles at $p$ and $q$ such that
\begin{equation}
  \label{eq:residuecondition}
  \Res_p\omega=-\Res_q\omega,
\end{equation}
there is a simple surgery operation which allows us to replace the poles at
$p$ and $q$ with a cylinders.  Cut $X$ along two closed geodesics, one in each of the two half-infinite
cylinders around $p$ and $q$, and then glue $X$ along the resulting boundary components by an isometry.  The
condition \eqref{eq:residuecondition} means exactly that this gluing map is locally a
translation; therefore, we get a new Abelian differential with two fewer poles.  Call this operation
\emph{plumbing a cylinder}.  It depends on two parameters: the height of the resulting cylinder and the amount
of twisting of the gluing map.

We can reverse this operation to replace a cylinder with two simple poles.  Just cut a cylinder along a closed
geodesic, and then glue two half-infinite cylinders to the resulting boundary components Call this operation
\emph{unplumbing a cylinder}.  These operations will be used in \S\ref{sec:localcoordinates} to define
polar coordinates around the boundary of the Deligne-Mumford compactification of moduli space. 

\paragraph{Connected Sums.}

Two Abelian differentials $(X_1, \omega_1)$ and $(X_2, \omega_2)$ can be combined into
one by taking a \emph{connected sum}.  This operation is studied in detail in \cite{mcmullenabel}.

Let $I$ be a line segment in the complex plane.  Suppose $I$ is embedded in each of the surfaces $(X_i,
\omega_i)$ by an isometric embedding $\epsilon_i\colon I\to X_i$ preserving the slope of $I$ as well as its
length.  Now cut each of the $X_i$ along $I$, and glue the surfaces together by gluing one boundary component
along $I$ on one of the surfaces to the opposite boundary component on the other surface.  The gluing maps are
just translations in the flat structures on the $(X_i, \omega_i)$, so the resulting surface has a flat
structure as well away from the ends of $I$.  By the Riemann removable singularity theorem, the conformal
structure and Abelian differential can actually be extended to the ends of $I$.  Thus we obtain a new Abelian
differential $(X, \omega)$.  If the segments $\epsilon_i(I)\subset X_i$ are disjoint from the zeros of
$\omega$, then the $\omega$ has simple zeros.

We call the resulting Abelian differential,
$$(X, \omega)=(X_1, \omega_1)\#_I(X_2, \omega_2),$$
the \emph{connected sum} of $(X_1, \omega_1)$ and $(X_2, \omega_2)$ along $I$.

This construction can be modified in the obvious way to perform a self connected sum of an Abelian
differential with itself, given two parallel embeddings of $I$ in that Abelian differential.

An Abelian differential $(X, \omega)$ resulting from a connected sum has two simple zeros $p$ and $q$ and two oriented
saddle connections $I_1$ and $I_2$ going from $p$ to $q$ such that
\begin{equation}
  \label{eq:sameperiods}
  \int_{I_1}\omega=\int_{I_2}\omega.
\end{equation}
Given any Abelian differential $(X, \omega)$ with a pair of oriented, embedded saddle connections $I_1$ and $I_2$ both
beginning and ending at the same zeros and satisfying \eqref{eq:sameperiods}, we can reverse the connected sum
operation.  To do this, cut $X$ along $I_1$ and $I_2$ and then reglue to get a new Abelian differential $(X', \omega')$.
Equation \eqref{eq:sameperiods} implies that this gluing can be done by a translation.  This operation is
called \emph{splitting along $I_1$ and $I_2$} and is inverse to the connected sum operation.

In this paper, there are two main cases where we will use these constructions.  First, suppose $(X_1, \omega_1)$
and $(X_2, \omega_2)$ are both genus one, and a segment $I\subset\cx$ is embedded in each by embeddings
$\epsilon_i$ as above.  Then we can form the connected sum along $I$, and the resulting Abelian differential
has genus two with two simple zeros.  Conversely, splitting a genus two Abelian differential along a pair of
saddle connections $I_1$ and $I_2$ such that $I_1\cup I_2$ separates the surface yields a pair of genus one
Abelian differentials.

Second, given an genus one Abelian differential with two embedded segments $I_i$ which are parallel and of the same length, we
can form a self connected sum along these segments.  This again yields a genus two Abelian differential with
two simple zeros.  Conversely, splitting along a pair of saddle connections on a genus two Abelian
differential satisfying \eqref{eq:sameperiods} which don't separate the surface yields a single genus one Abelian differential.

\paragraph{Splitting a double zero.}
There is also a cut-and-paste operation which replaces a zero of an Abelian differential with two zeros of
lower order.  This operation is explained in detail in \cite{kz03} and \cite{emz}.  We will describe this
operation in the case of a double zero, the only case we need.

\begin{figure}[htb]
  \centering
  \input{splittingzero.pstex_t}
  \caption{Splitting a double zero.}
  \label{fig:splittingadoublezero}
\end{figure}
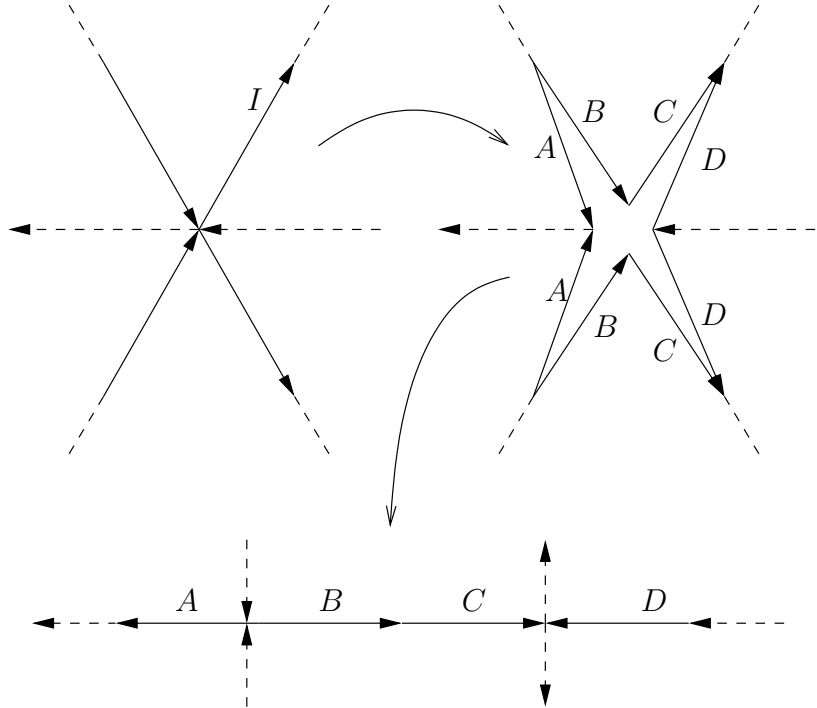

Let $(X, \omega)$ be an Abelian differential with a double zero at $p$, and choose a straight geodesic segment
$I$ starting at $p$.  The segment $I$ is a leaf of the foliation $\mathcal{F}_\theta$ of some slope $\theta$.
On $X$, draw an ``X'' composed of the segment $I$ and three other geodesic segments of the same length and
slope as in Figure \ref{fig:splittingadoublezero}.  (In this figure, leaves of $\mathcal{F}_\theta$ are
represented by dotted or solid lines, with the orientation indicated by an arrow, and the ``X'' is represented
by solid lines.  Two consecutive segments meeting at the zero form a $180^\circ$ angle.)  Assume that $I$
was chosen to be short enough so that each arm of the ``X'' is an embedded straight geodesic on $X$ not
meeting any of the other arms.

Now cut the surface along the ``X'' to obtain a surface with eight geodesic boundary components, and then glue
the boundary components pairwise so that the components marked with the same letter in
Figure~\ref{fig:splittingadoublezero} are glued together. The result is a new Abelian differential with simple
zeros at the two points where four leaves of $\mathcal{F}_\theta$ meet in the figure.  We will denote this new
Abelian differential by
$$(X, \omega) \#_I.$$

This operation can also be reversed in the obvious way.  Given a saddle connection of length $l$ joining two
distinct zeros, draw two more segments of length $l/2$ emanating from the zeros with the same slopes as in the bottom of
Figure \ref{fig:splittingadoublezero}.  If these segments can be drawn without intersecting each other or any
other zero, then the process of splitting a double zero can be reversed: cut along these three segments and
then reglue, following Figure \ref{fig:splittingadoublezero} in reverse.  This operation is called
\emph{collapsing a saddle connection}.

This operation, together with the connected sum operation, will be used in \S \ref{sec:normalbundles} to
parameterize tubular neighborhoods of $\barW$ and $\barP$ in $\Y$.

\subsection{Moduli of Abelian differentials}
\label{subsec:abelianmoduli}

In this subsection, we introduce the moduli space of Abelian differentials and discuss its geometry.

\paragraph{Teichm\"uller space.}

Let $\Sigma_g$ be a connected, closed, oriented, topological surface of genus $g$, and let $\Sigma_{g, n}$ be
a genus $g$ topological surface with $n$ marked points. A marked Riemann surface is a Riemann surface $X$
together with a homeomorphism $\Sigma_{g, n}\to X$.  Two marked Riemann $(f, X)$ and $(g, Y)$ marked by
$\Sigma_{g, n}$ are considered to be equivalent if $g \circ f^{-1}$ is homotopic to a conformal isomorphism by
a homotopy fixing the marked points.

Let the \emph{Teichm\"uller space} $\teich(\Sigma_{g, n})$ be the space of all Riemann surfaces marked by $\Sigma_{g,
  n}$ up to equivalence.  It has a topology induced by the well-known Teichm\"uller metric and is homeomorphic
to  $\cx^N$, where
$$N= \begin{cases}
    1 & \text{if $g=1$ and $n=0$ or $1$;}\\
  n -1 &\text{if $g=1$ and $n>1$;}\\
  n-3 &\text{if $g=0$ and $n>2$;}\\
  3g-3 +n & \text{if } g>1.
\end{cases}$$
We will use the abbreviation $\teich_{g, n}$ or $\teich_g$ when we don't need to emphasize the surface $\Sigma_g$.

Bers gave $\teich_{g, n}$ a complex structure by defining an embedding $B\colon \teich_{g, n}\to \cx^N$ which
is a homeomorphism onto its image, a bounded domain in $\cx^N$.

\paragraph{Moduli space.}

The modular group $\Mod(\Sigma_{g})$ is the group of all self homeomorphisms of $\Sigma_g$ up to isotopy.
Similarly, $\Mod(\Sigma_{g, n})$ is the group of all self homeomorphisms which preserve the marked points, up
to isotopy preserving the marked points.

An element $\gamma\in\Mod(\Sigma_{g, n})$ acts on $\teich(\Sigma_{g, n})$ by replacing a marking $$f\colon
\Sigma_{g, n}\to X$$ with $f\circ \gamma^{-1}$.  This defines a biholomorphic action of $\Mod(\Sigma_{g, n})$
on $\teich(\Sigma_{g, n})$ which is properly discontinuous.  Let $Z\subset \Sigma_{g, n}$ be the set of marked
points.  The stabilizer of a point $(f, X)$ is isomorphic to the group $\Aut(X, f(Z))$ of conformal
automorphisms of $X$ preserving the marked points.  In genus two, every Riemann surface has an order-two
automorphism $J$, the hyperelliptic involution.  This yields an element of order two in $\Mod(\Sigma_2)$ which
acts trivially on $\teich(\Sigma_2)$.

The \emph{moduli space} of genus $g$ Riemann surfaces with $n$ marked points is the quotient
$$\moduli[g, n]=\teichs/\Mod(\Sigma_{g, n}).$$
$\moduli[g, n]$ is a complex orbifold.

\paragraph{Bundles of Abelian differentials.}

Given a Riemann surface $X$, let $\Omega(X)$ be the space of Abelian differentials on $X$, a rank $g$ complex
vector space.  Let $\Omega\teich_g$ be the space of all pairs $(X, \omega)$ with $X\in\teich_g$ and
$\omega\in\Omega(X)$ a nonzero Abelian differential.  We can give $\Omega\teich_g$ the structure of a
trivial holomorphic punctured vector bundle as follows.

Over $\teich_g$, there is the universal curve
$\mathcal{C}\teich_g$, defined by Bers in \cite{bers73a}.  It is a complex manifold and comes with a
proper map $\pi\colon\mathcal{C}\teich_g\to\teich_g$ whose fiber over a Riemann surface $X$ is
isomorphic to $X$ itself.  The cotangent bundle to the fibers of $\pi$ is a line bundle
$\mathcal{L}\to\mathcal{C}\teich_g$.  The push-forward $\pi_*\mathcal{O}(\mathcal{L})$ of the sheaf of
sections of $\mathcal{L}$ is a sheaf on $\teich_g$.  The following theorem follows from \cite{bers61a}.

\begin{theorem}
  The sheaf $\pi_*\mathcal{O}(\mathcal{L})$ is the sheaf of sections of a trivial bundle over $\teich_g$ whose
  fiber over a Riemann surface $X$ is $\Omega(X)$.
\end{theorem}

The action of the mapping class group $\Mod_g$ on $\teich_g$ extends to an action on
$\Omega\teich_g$.  The quotient is a rank $g$ orbifold vector bundle $\Omega\moduli[g]$
over $\moduli[g]$ whose fiber over a Riemann surface $X$ is the quotient $\Omega(X)/\Aut(X)$.  This bundle
is sometimes called the \emph{Hodge bundle}.

In general, when $S$ is any sort of space of Riemann surfaces, $\Omega S$ will denote the natural bundle of
 Abelian differentials over $S$.

\paragraph{Strata.}

The bundles $\Omega\teich_g$ have a natural stratification in terms of the types of zeros of
the Abelian differentials.  Each nonzero Abelian differential on a nonsingular Riemann surface has $2g-2$ zeros,
counting multiplicity.  Given a sequence of integers $\mathbf{n}=(n_i)_{i=1}^r$ such that $\sum n_i=2g-2$, let
$\otn$ be the locus of all Abelian differentials which have $r$ zeros whose
multiplicities are given by the $n_i$.  This locus is a locally closed subset of $\Omega\teich_g$, and by Veech
\cite{veech90} it is actually a complex submanifold.  The quotient
$$\omn = \otn / \Mod$$
is then a complex suborbifold of $\Omega\moduli[g]$.

\paragraph{Period coordinates.}

These strata have natural coordinates defined in terms
of their periods by Veech \cite{veech90} and Masur \cite{masur82} which give the stratum $\otn$ the structure
of an affine manifold.

There is a fiber bundle of homology groups
$$\mathcal{H}_1\to\otn,$$
whose fiber over an Abelian differential $(X, \omega)$
is $H_1(X, Z(\omega); \zed),$ where $Z(\omega)$ is the set of zeros of $\omega$.
Given two Abelian differentials $(X_1, \omega_1)$, and $(X_2,
\omega_2)\in\otn$ with $(X_2, \omega_2)$ sufficiently close to $(X_1, \omega_1)$
there is a natural isomorphism
$$H_1(X_2, Z(\omega_2))\to H_1(X_1, Z(\omega_1)).$$
This defines a
flat connection on $\mathcal{H}_1$, the \emph{Gauss-Manin connection}.

Consider $(X, \omega)\in \otn$ and a small neighborhood $U$ of $(X, \omega)$.
Any $(X', \omega')\in U$ defines an element of $H^1(X', Z(\omega'); \cx)$ via the periods of $\omega'$.  Composing this with
the isomorphism
$$H^1(X', Z(\omega'); \cx)\to H^1(X, Z(\omega); \cx)$$
from the Gauss-Manin connection, we get a map
$$\phi\colon U\to H^1(X, Z(\omega); \cx),$$
the \emph{period coordinates}.  Veech
\cite{veech90} showed that these are in fact biholomorphic coordinate charts.

A choice of basis of $H_1(X, Z(\omega);\zed)$ defines an isomorphism
$$H^1(X, Z(\omega); \cx)\to\cx^{2g+n-1},$$
where $n$ is the number of zeros of $\omega$.  We can suppose that the basis is of the form,
\begin{equation}
  \label{eq:basisform}
  \{u_1, \ldots, u_{2g}, v_1, \ldots, v_{n-1}\},
\end{equation}
where $\{u_1,\ldots,u_{2g}\}$ is a symplectic basis of $H_1(X)$.  Two such bases are related by a matrix in
the group,
$$G=\left\{
  \begin{pmatrix}
    A & B \\
    0 & C
  \end{pmatrix}
  : A\in {\rm SP}_{2g}\zed,\quad B\in M_{2g, n}(\zed), \quad C\in{\rm GL}_n\zed \right\}.$$

Changing the basis by a matrix in $G$ changes
the period coordinates by the transpose of this matrix.  This gives $\otn$ the structure of an $(G,
\cx^{2g+n-1})$-manifold, and $\omn$ inherits the structure of a $(G,
\cx^{2g+n-1})$-orbifold.

\paragraph{Measures.}

Since the action of ${\rm GL}_n(\zed)$ on $\cx$ preserves Lebesgue measure, we can pull back this measure by
the period coordinates to get the \emph{period measure} $\mu({\bf n})$ on $\omn$.

The flat metric defined by an Abelian differential has area given by,
$$\Area(\omega)=\frac{i}{2}\int_X \omega\wedge\overline{\omega}.$$
Let $\Omega_{\leq 1}\moduli[g]({\bf n})$ be
the locus of Abelian differentials with $\Area(\omega)\leq1$, and $\Omega_1\moduli[g]({\bf n})$ the locus of
Abelian differentials with $\Area(\omega)=1$ (we will use the prefixes $\Omega_{\leq1}$ and $\Omega_1$ to
denote the analogous subsets of any space of Abelian differentials).  We can define a measure $\mu'({\bf n})$
on $\Omega_1\moduli[g]({\bf n})$ by restricting $\mu({\bf n})$ to $\Omega_{\leq\ 1}\moduli[g]({\bf n})$ and
then projecting by the natural map $\Omega_{\leq 1}\moduli[g]({\bf n})\to \Omega_1\moduli[g]({\bf n})$.
\label{page:measures}

The following theorem was proved by Veech and Masur: 
\begin{theorem}[{\cite{masur82},\cite{veech90}}]
  The measures $\mu'({\bf n})$  have finite total volume.
\end{theorem}

\paragraph{Action of $\SLtwoR$.}

Given a nonzero Abelian differential $(X, \omega)\in\Omega\teich_g$, we can as in \S \ref{subsec:flatgeometry}
choose an atlas of coordinate charts $\phi_i\colon U_i\to\cx$, which cover the complement of the zeros, such
that $\phi_i$ pulls back the form $dz$ on $\cx$ to the form $\omega$ on $X$.  These coordinate charts differ
by translations on their overlap.  Now, given an element $A\in\SLtwoR$, define new coordinate charts
$\phi_i'=A\circ\phi_i$ by composing with the usual action of $A$ in the complex plane.  The new coordinates
still differ by translations, and so define a new translation surface. This defines a new Abelian differential
$A\cdot(X, \omega)\in \Omega\teich_g$ which has a different complex structure than $X$ unless $A$ happens to
be a rotation in ${\rm SO}_2\reals$.  If $(X', \omega')=A\cdot(X, \omega)$, then there is a natural
real-affine map $h_A\colon (X, \omega)\to (X', \omega')$ which takes zeros to zeros.  The periods satisfy the
relation,
\begin{equation}
  \label{eq:periodrelation}
  \omega'((h_A)_*(\gamma))=a\cdot\omega(\gamma),
\end{equation}
for any $\gamma\in H_1(X, Z(\omega))$ with $Z(\omega)$ the set of zeros of $\omega$.

This construction defines an action of $\SLtwoR$ on $\Omega\teich_g$.  This action commutes with the action of
$\Mod$ and so defines an action of $\SLtwoR$ on $\Omega\moduli[g]$ which preserves to locus
$\Omega_1\moduli[g]$ of Abelian differentials with area one.

Equation \eqref{eq:periodrelation} can be used to show that the measure $\mu({\bf n})$ is $\SLtwoR$ invariant,
as is $\mu'({\bf n})$.

This action is closely related to the Teichm\"uller geodesic flow: the projection
$\Omega\moduli[g]\to\moduli[g]$ sends $\SLtwoR$ orbits to copies of the hyperbolic plane in $\moduli[g]$ which
are isometrically embedded with respect to the Teichm\"uller metric, and the restriction of the action to the
one parameter subgroup of diagonal matrices is the Teichm\"uller geodesic flow.

The strata $\omn$ and $\Omega_1\moduli[g]({\bf n})$ are invariant under the action of $\SLtwoR$.  The measures
$\mu$ and $\mu'$ defined in \S\ref{subsec:teichmuller} are invariant measures by \eqref{eq:periodrelation}.

The natural projection $\Omega_1\moduli[g]\to\proj\Omega\moduli[g]$ sends $\SLtwoR$ orbits to immersed copies
of the hyperbolic plane, giving a foliation $\fmod[g]$ of $\proj\Omega\moduli[g]$ by immersed hyperbolic
planes.  It is possible for a leaf of $\fmod[g]$ to be a closed subset of $\proj\Omega\moduli[g]$.  In that
case it is called a \emph{Teichm\"uller curve} because it is an algebraic curve whose projection to
$\moduli[g]$ is isometrically immersed with respect to the Teichm\"uller metric. It is an important unsolved
problem to classify Teichm\"uller curves in $\proj\Omega\moduli[g]$ or more generally to classify orbit
closures or invariant measures.

\subsection{${\rm SL}_2{\mathbb R}$ orbits in genus two.}
\label{subsec:action}

In genus two, McMullen has completely classified the orbit closures for the action of $\SLtwoR$ on $\Omega_1\moduli$
and the ergodic, invariant measures in the series of papers \cite{mcmullenbild, mcmullenabel, mcmullenspin, mcmullentor}.
In this subsection, we will discuss the subsets of $\Omega_1\moduli$ which appear in this classification and
which will be the focus of this paper.

A $\SLtwoR$-invariant subset $S$ of $\Omega_1\moduli$ corresponds to a subset $\proj S$ of $\proj\Omega\moduli$ which is
\emph{saturated} in the sense that leaves of the foliation $\fmod$ are either contained in $\proj S$ or
disjoint from $\proj S$.  Invariant measures on $\Omega_1\moduli$ correspond by a disintegration construction
to holonomy invariant measures on $\fmod$.  In this paper, we will adopt this point of view and focus on
leaves of this foliation rather than on $\SLtwoR$ orbits in $\Omega_1\moduli$.

\paragraph{Eigenform loci.}
\label{subsec:eigenformloci}

A form $(X, [\omega])\in\proj\Omega\moduli$ is an \emph{eigenform for real multiplication by $\ord$} if the
Jacobian of $X$ has real multiplication by $\ord$ with $\omega$ a nonzero eigenform.  Let the \emph{eigenform
  locus} $\E$ be the locus of all such pairs.  The following proposition follows from \cite[Corollary
5.7]{mcmullenabel}. 

\begin{prop}
  $\E$ is a closed, saturated subset of $\proj\Omega\moduli$.
\end{prop}

An embedding $\iota_i\colon \ord\to\reals$ determines a map
$$j_i\colon\E\to\X$$
by sending an eigenform $(X, \omega)\in\E$ to the pair $(\Jac(X), \rho)$, where
$$\rho\colon\ord\to\End\Jac(X)$$
is chosen so that
$$\omega(\rho(\lambda)\cdot\gamma)=\lambda^{(i)}\omega(\gamma)$$
for each $\gamma\in H_1(X;\zed)$.  Recall
that in \S\ref{sec:abelianvarieties}, we introduced the locus $\P$ of Abelian varieties in $\X$ which are
polarized products of elliptic curves.  The following follows from Proposition~\ref{prop:jacobians}.

\begin{prop}
  \label{prop:jisomorphism}
  The map $j_i$ is an isomorphism of $\E$ onto $\X\setminus\P$. 
\end{prop}

In this paper, we will implicitly identify $\E$ with $\X\setminus\P$ by the isomorphism $j_1$.

Let $\F$ be the foliation of $\E$ by Riemann surfaces induced by the foliation $\mathcal{F}$ of
$\proj\Omega\moduli$.  The foliation $\F$ of $\E$ extends to a foliation -- which we will continue to call
$\F$ -- of $\X$ defined by adding the connected components of $\P$ as leaves of $\F$.  McMullen
\cite{mcmullenhilbert} proved that this is actually a foliation of $\X$.

\paragraph{Elliptic differentials.}

There is a useful alternative characterization of eigenforms $(X, \omega)$ for
real multiplication by $\ord[d^2]$.

A branched cover $f\colon X\to E$ from a Riemann surface to an elliptic curve is said to be primitive if it
does not factor through an isogeny of elliptic curves $g\colon E'\to E$ of degree greater than one.
Equivalently, $f$ is primitive if the map on homology $f_*\colon H_1(X;\zed)\to H_1(E;\zed)$ is surjective.

An Abelian differential $(X, \omega)$ is called an \emph{elliptic differential} if it is the pullback of a
nonzero Abelian differential on an elliptic curve by a primitive branched cover.  The
\emph{degree} of an elliptic differential is the degree of the cover.

\begin{prop}[{\cite[Theorem~4.10]{mcmullenabel}}]
  \label{prop:ellipticdifferentials}
  The locus of degree $d$ elliptic differentials in $\Omega\moduli$ is exactly $\Omega\E[d^2]$.
\end{prop}

There is a special class of elliptic differentials called \emph{square-tiled surfaces}.  A square-tiled
surface is an Abelian differential which is pulled back from the square elliptic curve $(\cx/\zed[i], dz)$ by
some (not necessarily primitive) cover branched only over $0$.  Equivalently, an Abelian differential is
square-tiled if and only if its absolute and relative periods all lie in the Gaussian integers $\zed[i]$.
All square-tiled surfaces  lie on Teichm\"uller curves.

\paragraph{Weierstrass curves.}

A genus two Abelian differential is a \emph{Weierstrass form} if it is an
eigenform for real multiplication by some quadratic order $\ord$ and if it has a double zero.  Its
\emph{discriminant} is the discriminant $D$ of the order $\ord$.  The Weierstrass forms are parameterized by
the \emph{Weierstrass form bundle} $\Omega\W$, which is a line bundle over the \emph{Weierstrass curve},
$$\W:=\proj\Omega\W\subset\X.$$

\begin{theorem}[{\cite{mcmullenspin}}]
  $\W$ is a union of Teichm\"uller curves which is nonempty if and only if $D\geq 5$.  If $D=9$ or if
  $D\not\equiv 1 \mod 8$, then $\W$ is connected.  Otherwise $\W$ has two connected components.
\end{theorem}

\begin{remark}
  When $D=p^2$ for some prime $p$, this statement was first proved in \cite{hubertlelievre}.
\end{remark}

When $\W$ has two connected components, they are denoted by $\Wzero$ and $\Wone$.  They are distinguished by a
topological invariant called the \emph{spin invariant}.  This is  easy to describe when $D=d^2$.  In that
case, an Abelian differential $(X, \omega)\in\Omega\W[d^2]$ is branched over an elliptic curve $(E, \nu)$ by a
$d$-fold branched cover $f\colon X\to E$.  Let $p$ be the unique zero of $\omega$, one of the six Weierstrass
points of $X$.  Of the six Weierstrass points, $N$ of them have the same image in $E$ as $p$, with
either $N=1$ or $N=3$.  If $N=1$, then $(X, \omega)$ lies in $\Omega\Wzero$,
and if $N=3$, then $(X, \omega)$ lies in $\Omega\Wone$.

Let $\widetilde{W}_D$ be the inverse image of $\W$ in the universal cover $\half\times\half$ of $\X$.

\begin{prop}[{\cite{mcmullenbild}}]
  $\widetilde{W}_D$ is a countable union of graphs of transcendental holomorphic
  maps $\half\to\half$.
\end{prop}

Equivalently, $\W$ is transverse to the absolute period foliation $\A$ of $\X$ introduced in
\S\ref{subsec:hilbertmodular}.

\paragraph{Period coordinates for $\Omega\E$.}

Let $\E(1,1) = \E\setminus\W$, and let $\Omega\E(1, 1)$ be the bundle of nonzero eigenforms.  We can define
period coordinates on $\Omega\E(1, 1)$ in the same way as we defined period coordinates on the strata $\omn$.
This material will only be used in \S\ref{sec:lyapunov}.

Let
$$\widetilde{E}_D(1, 1)=\half\times\half\setminus\widetilde{W}_D,$$
the inverse image of $\E(1, 1)$ in the universal cover of $\X$, and let $\Omega\widetilde{E}_D(1, 1)$ be the (trivial)
bundle of $\iota_1$-eigenforms.
Given  an $(X,
\omega)\in\Omega\widetilde{E}_D(1, 1)$,  the homology group $H_1(X, Z(\omega);\zed)$ contains $H_1(X;\zed)$ as a
subgroup which is equipped with an isomorphism
$$\ord\oplus\ord^\vee\to H_1(X;\zed).$$
The condition that
$\omega$ is an eigenform is equivalent to the period map
$$P_\omega\colon H_1(X; \zed)\to\cx$$
being
$\ord$-linear, with $\ord$ acting on $\cx$ via the embedding $\iota_1$.  Let
$$H_{\ord}^1(X, Z(\omega); \cx)$$
be the subspace of $H^1(X, Z(\omega);\cx)$ consisting of all linear maps
$$H_1(X, Z(\omega);\zed)\to\cx$$
that are $\ord$-linear on $H_1(X;\zed)$.

Now consider $(X, \omega)\in \Omega\widetilde{E}_D(1,1)$.  For a sufficiently small neighborhood $U$ of $(X,
\omega)$, any $(X', \omega')\in U$ defines an element of $H^1_{\ord}(X, Z(\omega); \cx)$ by composing the
period map $P_{\omega'}$ with the isomorphism coming from the Gauss-Manin connection,
$$H_1(X, Z(\omega); \zed)\to H_1(X', Z(\omega'); \zed),$$
which is $\ord$-linear on $H_1(X;\zed)$.
This give biholomorphic period coordinates 
$$\phi\colon U\to H^1_{\ord}(X, Z(\omega); \cx).$$

Consider a triple $(\alpha, \beta, \gamma)\subset H_1(X, Z(\omega);\zed)$ such that:
\begin{itemize}
\item  $(\alpha, \beta)$ is a basis of $H_1(X;\zed)$ over $\ord$.
\item $H_1(X;\zed)\oplus\langle\gamma\rangle = H_1(X, Z(\omega);\zed)$.
\end{itemize}
Such a triple determines an isomorphism
$$H_1(X, Z(\omega);\zed)\to\cx^3.$$
Two such triples are related by a matrix in the group,
$$G = \left\{
  \begin{pmatrix}
    A & B \\
    0 & C
  \end{pmatrix}
  \in{\rm GL}_3 K(\ord): A\in\SLtwoord, \: B\in \ord\oplus\ord^\vee, \: C=\pm 1\right\}.
$$

Changing the triple by a matrix in $G$ changes
the period coordinates by the transpose of this matrix.  This gives $\Omega\widetilde{E}_D(1,1)$ the structure of an $(G,
\cx^3)$-manifold, and $\Omega E_D(1,1)$ inherits the structure of a $(G,
\cx^3)$-orbifold.

These period coordinates are compatible with the action of $\SLtwoR$ in the sense that
$\phi\colon U\to\cx^3$ commutes with the two $\SLtwoR$ actions, where $\SLtwoR$ acts on $\cx$ coordinate-wise by
identifying each $\cx$ factor with $\reals^2$.

By pulling back Lebesgue measure on $\cx^3$ via these charts, we define a measure $\mu_D$ on $\E(1,1)$.  We can use this
measure to define a measure $\mu'_{D}$ on $\Omega_1\E(1,1)$, using the same trick we used to define the
measures $\mu'({\bf n})$ on p.~\pageref{page:measures}.

\begin{theorem}[\cite{mcmullenabel}]
  \label{thm:volfinite}
  The measures $\mu'_{D}$ are finite, ergodic, $\SLtwoR$-invariant measures.  
\end{theorem}

\paragraph{Cusps of $\W$.} 

We describe here McMullen's classification of cusps of $\W$ from \cite{mcmullenspin},
which will play a crucial role in our calculation of $\chi(\We)$.

Given an Abelian differential $(X, \omega)\in \Omega_1 \moduli[g]$ and a slope $s\in\proj^1(\reals)$, the
foliation $\mathcal{F}_s$ of $X$ is \emph{completely periodic} if every leaf of the foliation is either a
saddle connection joining zeros of $\omega$, or a closed loop on $X$.  A completely periodic foliation divides
$X$ into finitely many maximal cylinders $C_i$ foliated by closed leaves of $\mathcal{F}_s$.  The complement
of $\bigcup C_i$ is a union of saddle connections on $X$ called the spine of $(X, \omega)$.

Two completely periodic foliations $\mathcal{F}_{s_i}$ of $(X_i, \omega_i)\in\Omega_1\moduli[g]$ are
\emph{equivalent} if there is some $A\in\SLtwoR$ such that $A\cdot(X_1, \omega_1)=(X_2, \omega_2)$ and $A\cdot
s_1=s_2$.

Now restrict to the case of Weierstrass forms $(X, \omega)\in\W$.  If $D$ is not square, then a completely
periodic foliation decomposes $X$ into two cylinders; if $D$ is square, then a completely periodic foliation
decomposes $X$ into either one or two cylinders (see \cite[Theorem 4.3]{mcmullenspin}).

A Weierstrass form $(X, \omega)\in\Omega_1 \W$ together with a completely periodic foliation $\mathcal{F}_s$ determines a cusp
of $\W$.  Let $N\subset\SLtwoR$ be the upper-triangular subgroup consisting of all
$$N_t=
\begin{pmatrix}
  1 & t\\
  0 & 1
\end{pmatrix}.
$$
Let $g\in\SLtwoR$ be some matrix taking $s$ to $0$ (so that the horizontal foliation of $g\cdot(X, \omega)$ is
completely periodic).  The map
$$t \mapsto g^{-1} N_t g$$
defines a path on $\Omega_1\W$, which happens to cover a closed horocycle $h$ on
$\W$. Replacing $(X, \omega, s)$ with an equivalent completely periodic foliation gives a horocycle homotopic
to $h$.  Homotopy classes of closed horocycles on a Riemann surfaces $S$ correspond to cusps of $S$, so this
construction associates a cusp of $\W$ to every equivalence class of completely periodic foliations.  This
correspondence is in fact a bijection; see, for example, \cite[Theorem 4.1]{mcmullenspin}.

\begin{figure}[htbp]
  \centering
  \includegraphics{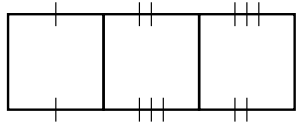}
  \caption{One cylinder cusp of $\W[9]$}
  \label{fig:onecylinder}
\end{figure}

We say that a cusp of $\W$ is a \emph{one-cylinder cusp} or a \emph{two-cylinder cusp} if the associated
completely periodic foliation has one cylinder or two cylinders respectively.  One cylinder cusps only arise
on $\W$ if $D$ is square.  An example of a foliation
associated to a one-cylinder cusp is given by gluing the edges of a three-by-one rectangle as in Figure~\ref{fig:onecylinder}
and taking the horizontal foliation.

McMullen \cite{mcmullenspin} classified the cusps of $\W$ by identifying them with splitting
prototypes.  We will describe his classification here using the equivalent $\W$-prototypes from
\S\ref{sec:prototypes}.

\begin{figure}[htbp]
  \centering
  \input{cusp.pstex_t}
  \caption{Cusp of $\W$}
  \label{fig:cusp}
\end{figure}
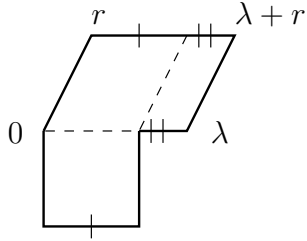

\label{page:cusp}
To the prototype $P=(a, b, c, \bar{q})$, we associate the surface $(X_P, \omega_P)$ formed by gluing a square in
$\cx$ with unit length sides to the parallelogram with sides $0$, $\lambda=\lambda(P)$, $r$, and $\lambda+r$, where
$\lambda$ is the unique positive root of
$$a\lambda^2+b\lambda+c=0,$$
and
$$r=-\frac{q}{c}\lambda - i \frac{a}{c}\lambda,$$
and then gluing opposite sides of the resulting polygon (see Figure~\ref{fig:cusp}).  This $(X_P, \omega_P)$
is an eigenform for real multiplication by $\ord$ where $D$ is the discriminant of $P$.
The horizontal foliation of $(X_P, \omega_P)$ is completely periodic and so determines a cusp $w_P$ of $\W$.
\begin{theorem}[{\cite[Theorem 4.1]{mcmullenspin}}]
  \label{thm:cuspstoprots}
  The map $P\mapsto w_P$ described above determines a bijection between the set of $\W$-prototypes and the set
  of two-cylinder cusps of $\W$.
\end{theorem}

We will also need to know which connected component of $\W$ contains a given  cusp $w_P$.  Given an order
$\ord$, define the \emph{conductor} of $\ord$ to be the integer $f$ such that $D=f^2 E$ with $E$ a fundamental
discriminant.

\begin{theorem}[{\cite[Theorem 5.3]{mcmullenspin}}]
  \label{thm:cuspsbyspin}
  The cusp $w_P$ of $\W$ associated to the $\W$-prototype $P=(a, b, c, \bar{q})$ is contained in the component
  $\W^{\epsilon(P)}$, where
  \begin{equation}
    \label{eq:epsilonofP}
    \epsilon(P)\equiv\frac{b-f}{2} + (a+1)(q+c+qc)\pmod 2,
  \end{equation}
  and $f$ is the conductor of $\ord$.
\end{theorem}

\paragraph{Other Teichm\"uller curves.}
\label{subsec:otherteichmullercurves}

When the discriminant $D$ is square, there is an infinite family of Teichm\"uller curves on $\X[d^2]$
parameterizing Abelian differentials with two simple zeros.  Given $(X, \omega)\in \Omega\X[d^2]$,  let $f\colon
(X, \omega)\to(E, \nu)$ its associated degree $d$ torus cover. Define $\Omega\W[d^2][n]$ to be the locus of
all differentials such that the two branch points of $f$ in $E$ differ by torsion of degree exactly $n$ in the
group law on $E$, and define $\W[d^2][n]\subset \X[d^2]$ be its projectivization.  The locus $\W[d^2][n]$ is a
union of Teichm\"uller curves.  Note that $\W[d^2][0]$ contains $\W[d^2]$ but is in general larger because $X$
can have distinct zeros which are branched over the same point of $E$.

There is one more example of a genus two Teichm\"uller curve which comes from an infinite family discovered by
Veech.  For even $n$, consider the Abelian differential obtained by gluing opposite sides of the regular
$n$-gon.  Veech \cite{veech92} showed that the $\SLtwoR$ orbit of this differential lies on a Teichm\"uller
curve $D_n$.  The Abelian differential coming from the decagon is an genus two Abelian differential with two
simple zeros.  Its orbit is a Teichm\"uller $D_{10}$ curve lying on $\X[5]$.

\paragraph{Classification of $\SLtwoR$ orbit closures in genus two.}

It happens that every $\SLtwoR$ orbit-closure on $\Omega_1\moduli$ or ergodic, invariant measure is one of the ones just described.

\begin{theorem}[{\cite{mcmullenabel},    \cite{mcmullentor}}]
  Every closure of an $\SLtwoR$ orbit in $\Omega_1\moduli$ is either all of $\Omega_1\moduli$, or one of the
  following submanifolds:  $\Omega_1\moduli(2)$, $\Omega_1\E$, $\Omega_1 D_{10}$, a connected component of
  $\Omega_1\W$, or a component of $\Omega_1\W[d^2][n]$.

  Furthermore each of these orbit closure
  carries a unique ergodic, absolutely continuous, $\SLtwoR$-invariant probability measure. These are all
  of the ergodic, $\SLtwoR$-invariant probability measures on $\Omega_1\moduli$. 
\end{theorem}

The only gap remaining in this classification is to describe the connected
components of $\W[d^2][n]$.


%% file: splittingzero.pstex_t
\begin{picture}(0,0)%
\includegraphics{splittingzero.pstex}%
\end{picture}%
\setlength{\unitlength}{3947sp}%
\begingroup\makeatletter\ifx\SetFigFont\undefined%
\gdef\SetFigFont#1#2#3#4#5{%
  \reset@font\fontsize{#1}{#2pt}%
  \fontfamily{#3}\fontseries{#4}\fontshape{#5}%
  \selectfont}%
\fi\endgroup%
\begin{picture}(5124,4449)(3589,-6973)
\put(7651,-4786){\makebox(0,0)[lb]{\smash{{\SetFigFont{12}{14.4}{\rmdefault}{\mddefault}{\updefault}{\color[rgb]{0,0,0}$C$}%
}}}}
\put(4651,-6361){\makebox(0,0)[lb]{\smash{{\SetFigFont{12}{14.4}{\rmdefault}{\mddefault}{\updefault}{\color[rgb]{0,0,0}$A$}%
}}}}
\put(5551,-6361){\makebox(0,0)[lb]{\smash{{\SetFigFont{12}{14.4}{\rmdefault}{\mddefault}{\updefault}{\color[rgb]{0,0,0}$B$}%
}}}}
\put(6451,-6361){\makebox(0,0)[lb]{\smash{{\SetFigFont{12}{14.4}{\rmdefault}{\mddefault}{\updefault}{\color[rgb]{0,0,0}$C$}%
}}}}
\put(7576,-6361){\makebox(0,0)[lb]{\smash{{\SetFigFont{12}{14.4}{\rmdefault}{\mddefault}{\updefault}{\color[rgb]{0,0,0}$D$}%
}}}}
\put(6976,-4411){\makebox(0,0)[lb]{\smash{{\SetFigFont{12}{14.4}{\rmdefault}{\mddefault}{\updefault}{\color[rgb]{0,0,0}$A$}%
}}}}
\put(7201,-3286){\makebox(0,0)[lb]{\smash{{\SetFigFont{12}{14.4}{\rmdefault}{\mddefault}{\updefault}{\color[rgb]{0,0,0}$B$}%
}}}}
\put(7951,-3586){\makebox(0,0)[lb]{\smash{{\SetFigFont{12}{14.4}{\rmdefault}{\mddefault}{\updefault}{\color[rgb]{0,0,0}$D$}%
}}}}
\put(7951,-4561){\makebox(0,0)[lb]{\smash{{\SetFigFont{12}{14.4}{\rmdefault}{\mddefault}{\updefault}{\color[rgb]{0,0,0}$D$}%
}}}}
\put(7276,-4636){\makebox(0,0)[lb]{\smash{{\SetFigFont{12}{14.4}{\rmdefault}{\mddefault}{\updefault}{\color[rgb]{0,0,0}$B$}%
}}}}
\put(5101,-3211){\makebox(0,0)[lb]{\smash{{\SetFigFont{12}{14.4}{\rmdefault}{\mddefault}{\updefault}{\color[rgb]{0,0,0}$I$}%
}}}}
\put(6901,-3511){\makebox(0,0)[lb]{\smash{{\SetFigFont{12}{14.4}{\rmdefault}{\mddefault}{\updefault}{\color[rgb]{0,0,0}$A$}%
}}}}
\put(7651,-3286){\makebox(0,0)[lb]{\smash{{\SetFigFont{12}{14.4}{\rmdefault}{\mddefault}{\updefault}{\color[rgb]{0,0,0}$C$}%
}}}}
\end{picture}%

%% file: cusp.pstex_t
\begin{picture}(0,0)%
\includegraphics{cusp.pstex}%
\end{picture}%
\setlength{\unitlength}{3947sp}%
\begingroup\makeatletter\ifx\SetFigFont\undefined%
\gdef\SetFigFont#1#2#3#4#5{%
  \reset@font\fontsize{#1}{#2pt}%
  \fontfamily{#3}\fontseries{#4}\fontshape{#5}%
  \selectfont}%
\fi\endgroup%
\begin{picture}(1447,1494)(4576,-5398)
\put(6001,-4036){\makebox(0,0)[lb]{\smash{{\SetFigFont{12}{14.4}{\rmdefault}{\mddefault}{\updefault}{\color[rgb]{0,0,0}$\lambda + r$}%
}}}}
\put(4576,-4786){\makebox(0,0)[lb]{\smash{{\SetFigFont{12}{14.4}{\rmdefault}{\mddefault}{\updefault}{\color[rgb]{0,0,0}$0$}%
}}}}
\put(5851,-4786){\makebox(0,0)[lb]{\smash{{\SetFigFont{12}{14.4}{\rmdefault}{\mddefault}{\updefault}{\color[rgb]{0,0,0}$\lambda$}%
}}}}
\put(5101,-4036){\makebox(0,0)[lb]{\smash{{\SetFigFont{12}{14.4}{\rmdefault}{\mddefault}{\updefault}{\color[rgb]{0,0,0}$r$}%
}}}}
\end{picture}%

%% file: delignemumford.tex
\section{Deligne-Mumford compactification of  moduli space}
\label{sec:moduli}

\subsection{Stable Riemann surfaces}
\label{subsec:rswithnodes}

A \emph{nodal Riemann surface} is a connected, compact, one-dimensional, complex analytic space with only
nodes as singularities (a node is a transverse crossing of two nonsingular branches).  Equivalently, a nodal
Riemann surface can be regarded as a finite type Riemann surface with finitely many cusps which have been
identified pairwise to form nodes.  A connected component of a nodal Riemann surface $X$ with its nodes
removed is called a \emph{part} of $X$, and the closure of a part of $X$ is an \emph{irreducible component} of
$X$.  In this paper the \emph{genus} of a nodal Riemann surface $X$ will mean its arithmetic genus,
$$g = 1- \chi(\mathcal{O}_X),$$
where $\mathcal{O}_X$ is the structure sheaf of $X$.  In topological terms,
the arithmetic genus of a nodal Riemann surface $X$ is the genus of the nonsingular Riemann
surface obtained by replacing each node of $X$ with an annulus.

A \emph{stable Riemann surface} is a connected nodal Riemann surface for which each part has nonabelian fundamental
group (or equivalently negative Euler characteristic).

A nodal Riemann surface $X$ has a \emph{normalization} $\widetilde{X}\to X$ defined by separating the two branches
passing through each node of $X$.  

\paragraph{Stable Abelian differentials.}

A \emph{stable Abelian differential} on a stable Riemann surface $X$ is a
holomorphic 1-form on $X$ minus its nodes such that:

\begin{itemize}
\item Its restriction to each part of $X$ has at worst simple poles at the cusps.
\item At two cusps which have been identified to form a node, the differential has opposite residues.
\end{itemize}

These properties can be conveniently rephrased using the normalization.  If $\omega$ is a meromorphic Abelian
differential on $X$, and $f\colon \tilde{X}\to X$ is the normalization, then $\omega$ is stable if and only
if $f^*\omega$ has at worst simple poles, and for every $q\in X$,
$$\sum_{f(p)=q} \Res_p(f^*\omega)=0.$$

A stable Riemann surface $X$ has a \emph{dualizing sheaf} $\omega_X$.   A stable Abelian
differential on a stable Riemann surface is just a global section of $\omega_X$.
This is discussed in \cite{harrismorrison} and \cite{hartshorne63}.   

The stable Abelian differentials on a genus $g$ singular Riemann surface $X$ form a complex vector space
which we will write as $\Omega(X)$.  The complex dimension of $\Omega(X)$ is $g$:
\begin{align*}
  h^0(X, \omega_X) &= h^1(X, \struct_X)\\
  &= 1-\chi(\struct_X)\\
  &=g.
\end{align*}
This fact also follows easily from the Riemann-Roch Theorem.  It is also proved in
\cite[Theorem IV.2]{serre59}.

We will often use the term ``Abelian differential'' as shorthand for ``Riemann surface together with an
Abelian differential''. We will call a node of a stable Abelian differential $(X, \omega)$ a \emph{polar node}
if $\omega$ has a pole there and a \emph{holomorphic node} otherwise.

\paragraph{Jacobians.}

The classical notion of the Jacobian variety of a nonsingular Riemann surface can be extended to singular
Riemann surfaces.

If $X$ is a possibly singular Riemann surface, let $X_0$ be the set of nonsingular points of $X$.  Then there
is a natural map $H_1(X_0;\zed)\to \Omega^*(X)$, given by integrating forms over homology classes. 
The \emph{Jacobian variety} of  $X$ is the variety,
$$\Jac(X)=\Omega^*(X)/H_1(X_0; \zed).$$

\begin{figure}[htbp]
  \centering
  \input{node.pstex_t}
  \caption{Curves around a node.}
  \label{fig:noderelation}
\end{figure}
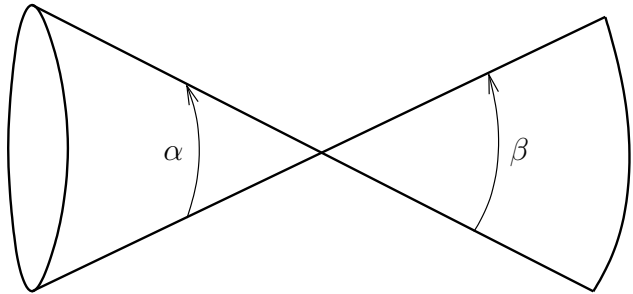

In the case of a stable Riemann surface $X$, it is easy to describe the kernel of the map
$H_1(X_0;\zed)\to\Omega^*(X)$.  It is the subgroup of $H_1(X_0;\zed)$ generated by the relations $\alpha-\beta$ when $\alpha$ and $\beta$
are homology classes generated by curves going around the same node of $X$ in the ``same direction'' on
opposite parts of $X$ as in Figure \ref{fig:noderelation}.

There is an exact sequence which relates the Jacobian of a nodal Riemann surface to the Jacobian of its
normalization.  If $\widetilde{X}\to X$ is the normalization of $X$, then holomorphic Abelian differentials on
$\widetilde{X}$ restrict to stable Abelian differentials on $X$, so there is a natural map
$\Jac(X)\to\Jac(\widetilde{X})$.  We get an exact sequence:
$$0\to(\cx^*)^n\to\Jac(X)\to\Jac(\widetilde{X})\to 0,$$
where $n$ is the difference of the genera of $X$ and $\widetilde{X}$.  This realizes $\Jac(X)$ as a
\emph{semi-Abelian variety}.

For example, if $X$ has genus two and one separating node, then $\widetilde{X}$ is the disjoint union of two
elliptic curves.  The Jacobian of $X$ is then the product of these two elliptic curves.

If $X$ has genus two and one nonseparating node, then the normalization of $X$ is an elliptic curve $E$.  The
Jacobian of $X$ is then an extension of  $E$ by $\cx^*$.  The original stable Riemann surface $X$ can be recovered
from $\Jac(X)$. 

If $X$ has genus two and two (or three) nonseparating nodes, then $\widetilde{X}$ is $\proj^1$ (or $\proj^1
\cup \proj^1$ respectively), so
$\Jac(X)\isom\cx^*\times\cx^*$.  This example shows that two distinct nodal Riemann surfaces may have the same
Jacobian.

The Jacobian of $X$ can also be identified with $\Pic^0(X)$, the group of all line bundles on $X$ which have
degree zero on each irreducible component of $X$ \cite[p.250]{harrismorrison}.

\subsection{Deligne-Mumford compactification}
\label{subsec:teichmuller}

\paragraph{Marked stable Riemann surfaces.}

Given a topological stable  surface $X$ of genus $g$, a \emph{collapse} of $\Sigma_g$ onto
$X$ is a continuous surjection $f\colon\Sigma_g\to X$ with the following properties.
\begin{itemize}
\item The inverse image of each node of $X$ is a Jordan curve on $\Sigma_g$.
\item Each component of $\Sigma\setminus f^{-1}(N)$, where $N$ is the set of nodes of $X$, maps
  homeomorphically, preserving the orientation, onto a part of $X$.
\end{itemize}
If $X$ is nonsingular, a collapse is just a homeomorphism $\Sigma_g\to X$.  These maps were introduced by
Bers; in his terminology, a collapse is called a strong deformation.

A \emph{marked stable Riemann surface} is a stable Riemann surface $X$, together with a collapse
$\Sigma_g\to X$.  Two markings $f_i\colon \Sigma_g\to X_i$ are \emph{equivalent} if there is a conformal
isomorphism $g\colon X_1\to X_2$ such that the following diagram commutes up to homotopy.

$$\xymatrix{
  \Sigma_g \ar[r]^{f_1} \ar[dr]_{f_2} & X_1 \ar[d]^{g} \\
  & X_2}$$

We will sometimes denote by $[f\colon\Sigma_g\to X]$ the class of all marked surfaces equivalent to $f\colon \Sigma_g\to
X$.

\paragraph{Augmented Teichm\"uller space.}

The Teichm\"uller space $\teich(\Sigma_g)$ is contained in the \emph{Augmented Teichm\"uller space}, $\augteich(\Sigma_g)$,
the set of all marked stable Riemann surfaces up to equivalence.  Let
$\bdry\teich(\Sigma_g)=\augteich(\Sigma_g)\setminus\teich(\Sigma_g)$.

We give $\augteich(\Sigma_g)$ a topology as follows.  Given a closed curve $\gamma$ on $\Sigma_g$, define a
function $l_\gamma\colon\augteich(\Sigma_g)\to\reals\cup\{\infty\}$:  if $\gamma$ is not homotopic
to a curve on $X$ disjoint from the nodes, let $l_\gamma(X)=\infty$; if $\gamma$ is homotopic to a node of
$X$, let $l_\gamma(X)=0$; otherwise let $l_\gamma(X)$ be the length of the unique geodesic homotopic to
$\gamma$ in the Poincar\`e metric on $X$ minus its nodes.  Give $\augteich(\Sigma_g)$ the smallest topology
such that $l_\gamma$ is continuous for every closed curve $\gamma$ on $\Sigma_g$.  The induced subspace
topology on $\teich(\Sigma_g)$ agrees with that defined by the Teichm\"uller metric on $\teich(\Sigma_g)$.

Abikoff \cite{abikoff77} showed that this topology is equivalent to other natural topologies on
$\augteich(\Sigma_g)$, such as those obtained by looking at quasiconformal maps or quasi-isometries defined
outside a neighborhood of the nodes.  

Another useful topology which is equivalent to this one is the \emph{conformal topology}
\label{page:conformaltopology} which is defined as
follows.  Let $[f\colon\Sigma_g\to X]\in\bdry\teich(\Sigma_g)$, and let $U\subset X$ be a
neighborhood of the nodes of $X$.  Let $V_U$ be the set of all $[g\colon\Sigma_g\to Y]\in\augteich(\Sigma_g)$
for which $f$ and $g$ can be adjusted by homotopies so that $g\circ f^{-1}|_{X\setminus \overline{U}}$ is
conformal.  The set of all $V_U$ as $U$ runs over all neighborhoods of the nodes of $X$ define a
neighborhood basis of $X$.  Together with the open sets of $\teich(\Sigma_g)$, this defines a topology on
$\augteich(\Sigma_g)$.  It is well-known that this topology is equivalent to the one defined above; however, a
proof does not exist in the literature.  

A \emph{curve system} on $\Sigma_g$ is a collection of simple closed curves on $\Sigma_g$, none of which are
isotopic to any other or to a point.  For each curve system $S$, there is a subspace $\teich(\Sigma_g,
S)\subset \bdry\teich(\Sigma_g)$ consisting of marked stable Riemann surfaces $[f\colon \Sigma_g\to X]$ which
collapse each curve of a curve system homotopic to $S$ to a point.  For each connected component $\Sigma_g^i$
of $\Sigma_g\setminus S$, let $S_i$ be the closed surface with marked points obtained by collapsing the boundary
components of $\Sigma_g^i$ to points and regarding the images of the boundary components as marked points.
There is a natural isomorphism,
$$\teich(\Sigma_g, S) \isom \prod_i \teich(S_i),$$

\paragraph{Deligne-Mumford compactification.}

The action of $\Mod(\Sigma_g)$ of $\teich(\Sigma_g)$ extends to an action on $\augteich(\Sigma_g)$.  The
quotient,
$$\barmoduli[g]=\augteich(\Sigma_g)/\Mod(\Sigma_g),$$
is the Deligne-Mumford compactification of moduli space.
It is a compact orbifold whose points naturally parameterize stable Riemann surfaces of genus $g$.

Given a curve system $S\subset
\Sigma_g$, let $\moduli[g](S)\subset\barmoduli[g]$ be the stratum of stable Riemann surfaces homeomorphic to the
topological stable surface $\Sigma_g/ S$.  This is a locally closed subset of $\barmoduli[g]$, and $\barmoduli[g]$
is the disjoint union of all of the $\moduli[g](S)$ as $S$ ranges over all isotopy classes of curve systems
on $\Sigma_g$ up to the action of the modular group.

\paragraph{Dehn space.}

Given a curve system $S$ on $\Sigma_g$, let $\Twist(S)$ be the group generated by Dehn twists
around the curves of $\Sigma$.  It is an Abelian group isomorphic to $\zed^h$, where $h$ is the number of
curves in $S$.   Define the \emph{Dehn space} $\Def(\Sigma_g, S)$  to be
$$\Def(\Sigma_g, S) = (\teich(\Sigma_g)\cup \teich(\Sigma_g, S))/\Twist(S).$$
We will sometimes use the notation $\Def_g(S) = \Def(\Sigma_g, S)$ when we don't need to emphasize the surface
$\Sigma_g$.  Bers called $\Def(\Sigma_g, S)$  the strong deformation space.

Let $\Stab(S)$ be the subgroup of $\Mod(\Sigma_g)$ which maps $S$ to an isotopic curve system.
$\Twist(S)$ is a normal subgroup of $\Stab(S)$, and the quotient,
$$\Mod(\Sigma_g, S) = \Stab(S)/\Twist(S),$$
is the mapping class group of the topological stable surface obtained by collapsing each curve in $S$ to a
point.

The group $\Mod(\Sigma_g, S)$ acts on $\Def(\Sigma_g, S)$, and the natural map
$$\pi\colon\Def(\Sigma_g, S)\to\barmoduli[g]$$
is invariant under this action.  Bers showed:
\begin{prop}[{\cite[p. 1221]{bers74}}]\label{prop:bersorbifold}
  If $X\in\Def(\Sigma_g, S)$, and $G\subset\Mod(\Sigma_g, S)$ is the stabilizer of $X$, then there is a neighborhood
  $U\subset\Def(\Sigma_g, S)$ of $X$, stable under the action of $G$ such that $\pi$ factors through to a map,
  $\bar{\pi}\colon U/G\to\barmoduli[g]$, which is homeomorphic onto its image.
\end{prop}

\paragraph{Complex structure.}

We now describe the complex structure on $\barmoduli[g]$, following Bers's approach \cite{bers81}.

In \cite{bers61}, Bers defined an embedding $B\colon\teich(\Sigma_g)\to\cx^n$, where
$$n=
\begin{cases}
  1 & \text{if } g=1;\\
  3g-3 & \text{if } g>1.
\end{cases}
$$
It is a homeomorphism onto a bounded domain, and $\teich(\Sigma_g)$ inherits the complex structure on
$\cx^n$.  Ahlfors \cite{ahlfors60} showed that this is the unique complex structure on $\teich(\Sigma_g)$
for which the periods of Abelian differentials vary holomorphically.

The modular group $\Mod(\Sigma_g)$ acts biholomorphically on $\teich(\Sigma_g)$ with this complex structure,
so $\moduli[g]$ inherits the structure of a complex orbifold.

Bers gives $\barmoduli[g]$ a complex structure by first giving $\Def(\Sigma_g, S)$ a complex structure for
each curve system $S$.  First define a sheaf of rings $\struct$ on $\Def(\Sigma_g, S)$ to be the sheaf
consisting of all continuous functions $\Def(\Sigma_g, S)\to\cx$ which are holomorphic on
$$\Def_0(\Sigma_g, S)=\teich(\Sigma_g)/\Twist(S)\subset \Def(\Sigma_g, S)$$
with respect to the complex
structure on $\Def_0(\Sigma_g, S)$ induced by the Bers embedding of Teichm\"uller space.  Given a domain,
$U\subset\Def(\Sigma_g, S)$, we consider a function $f\colon U\to\cx$ to be holomorphic if it is in $\struct$.
Bers showed that this defines an integrable complex structure on $\Def(\Sigma_g, S)$ by giving a biholomorphic
isomorphism of $\Def(\Sigma_g, S)$ with a bounded domain in $\cx^n$ parameterizing a certain family of
Kleinian groups:

\begin{prop}[{\cite{bers81}}]
  \label{prop:modulispaceorbifold}
  Let $S=\bigcup S_i$ be a curve system on $\Sigma_g$.  There is a biholomorphic map
  $$B\colon \Def(\Sigma_g, S)\to \cx^n$$
  onto a bounded domain in $\cx^n$.  A curve $S_i$ is homotopic to a node on a marked stable
  Riemann surface $X\in\Def(\Sigma_g, S)$ if and only if 
  $$z_i(B(X))=0,$$
  where the $z_i$ are the coordinates on $\cx^n$.
\end{prop}
Similar embeddings of $\Def(\Sigma_g, S)$ are constructed in \cite{marden86} and \cite{kra90}.

Using Proposition \ref{prop:bersorbifold}, we give $\barmoduli[g]$ the unique complex structure which makes
all of the maps $\Def(\Sigma_g, S)\to\barmoduli[g]$ holomorphic.

Baily showed in \cite{baily60, baily62} that $\moduli[g]$ also has the structure of a quasi-projective
variety.  Wolpert \cite{wolpert85} showed using the Weil-Petersson metric on $\moduli[g]$ that $\barmoduli[g]$
has the structure of a projective variety.  Alternatively, Deligne-Mumford \cite{delignemumford,mumford71}
constructed $\barmoduli[g]$ algebraically and showed that it is a coarse moduli space for stable genus $g$
curves, and  Knudsen-Mumford \cite{knudsenmumford}  showed that $\barmoduli[g]$ is a projective variety.

\paragraph{Jacobians in genus two.}

By associating a Riemann surface to its Jacobian, there is a natural morphism
$\Jac\colon \moduli[g]\to\siegelmod[g]$.

\begin{theorem}[{\cite{namikawa73}}]
  The morphism $\Jac$ extends to a morphism
  $$\overline{\Jac}\colon\barmoduli[g]\to\satsiegelmod[g],$$
  which sends a Riemann surface $X$ to the Jacobian of the normalization of $X$.
\end{theorem}

Determining the image of $\Jac$ map is in general very difficult; however, this is simple in genus two.  Let
$\tildemoduli$ be the subvariety of $\barmoduli$ consisting of nonsingular Riemann surfaces together with
pairs of elliptic curves joined at a single node.  The following proposition is well known.

\begin{prop}
  \label{prop:jacobians}
  The image of $\Jac\colon \moduli\to\siegelmod$ is exactly those Abelian varieties which are not polarized products
  of elliptic curves. Furthermore, $\Jac$ extends to an isomorphism $\widetilde{\Jac}\colon \tildemoduli\to\siegelmod$.
\end{prop}

\begin{proof}[Sketch of proof]
  The image of $\overline{\Jac}$ contains an open set because these varieties have the same dimension, and
  $\Jac$ restricted to $\moduli$ is injective by the Torelli theorem.  Any map between complete irreducible
  varieties of the same dimension whose image contains an open set is onto, so $\overline{\Jac}$ is onto.

  It's easy to check that the Jacobian of a genus two Riemann surface $X$ can't be the polarized product of two elliptic
  curves (for example this would give a degree one map of $X$ to an elliptic curve), so the image of $\moduli$
  is exactly the complement of the locus of polarized products.

  The map $\widetilde{\Jac}$ is also injective (since the Jacobian of a Riemann surface $X$ formed from two
  elliptic curves joined at a node is $X$ itself), so it is a bijection.  It follows that $\widetilde{\Jac}$
  is an isomorphism because it is a bijection between two normal varieties.  
\end{proof}

\paragraph{Abelian differentials over $\barmoduli[g]$.}

The bundle $\Omega\teich(\Sigma_g)$ extends to a trivial bundle $\Omega\augteich(\Sigma_g)$ over
$\augteich(\Sigma_g)$.  The quotient,
$$\Omega\barmoduli[g] = \Omega\augteich(\Sigma_g) / \Mod(\Sigma_g),$$
is the moduli space of stable Abelian differentials of genus $g$.  It is an orbifold vector bundle over
$\barmoduli[g]$ whose fiber over $X$ is $(\Omega(X)\setminus \{0\})/\Aut(X)$.

Similarly, define the bundle $\Omega\Def(\Sigma_g, S)$ to be the restriction of
$$\Omega\augteich(\Sigma_g, S)/\Twist(S)$$
to $\Def(\Sigma_g, S)$.  This is a trivial bundle whose fiber over a marked, stable Abelian differential $(f,
X)$ is $\Omega(X)$.  To see that it is trivial, choose a Lagrangian subspace
$\Lambda\subset H_1(\Sigma_g; \zed)$ which is orthogonal to the homology class of each of the curves in
$S$.  Then $\Lambda$ is fixed pointwise by the action of $\Twist(S)$, and
$$\Omega\Def(\Sigma_g, S)\isom \Def(\Sigma_g, S)\times \Hom(\Lambda, \cx).$$

\subsection{Degenerating Abelian differentials}
\label{subsec:polarcoordinates}

We now study in more detail the geometry of Abelian differentials as they approach the boundary of moduli
space.  For example, we will see that Abelian differentials close to a stable Abelian differential $(X, \omega)$ with a
polar node develop very long cylinders which are pinched off to form a node in the limit.  We will then use
our understanding of this geometry to study the behavior of Jacobians of Riemann surfaces near the boundary of moduli
space.  Finally, we will give two applications concerning real multiplication which will be used in \S\ref{sec:limitsofeigenforms}.

\paragraph{Long cylinders.}

The conformal topology on $\augteich(\Sigma_g)$ defined on p.~\pageref{page:conformaltopology} is compatible with the topology on the
trivial bundle $\Omega\augteich(\Sigma_g)$ in the following sense.  Consider $[f\colon\Sigma_g\to
X]\in\bdry\teich(\Sigma_g)$ with $\omega\in\Omega(X)$.  Let $U\subset X$ be a neighborhood of the nodes of $X$,
and let $K>1$.  Let $V_{U,K}$ be the set of all $([g\colon\Sigma_g\to Y], \eta)\in\Omega\augteich(\Sigma_g)$ for
which $f$ and $g$ can be adjusted by homotopies so that $g\circ f^{-1}|_{X\setminus \overline{U}}$ is
conformal, and
$$K^{-1}<\left|\frac{(g\circ f^{-1})^*\eta}{\omega}\right| < K$$
on $X\setminus \overline{U}$.  It is an unpublished result that the $V_{U, K}$ are a neighborhood basis of $([f\colon\Sigma_g\to
X], \omega)$ in $\Omega\augteich(\Sigma_g)$.

We can use this description of the topology on $\Omega\augteich(\Sigma_g)$ to prove the following theorem
which gives a more precise description of the shape of an Abelian differential close to a stable Abelian
differential $(X, \omega)$ with $X\in\bdry\moduli$.

\begin{theorem}
  \label{thm:longcylinders}
  Let $(X, \omega)\in\Omega\augteich(\Sigma_g)$ with $X\in\augteich(\Sigma_g, S)$, and let $C\subset X$ be a
  compact subset disjoint from the nodes.  For any $\epsilon, h>0$ and $K>1$, there is a neighborhood $U$ of
  $(X, \omega)$ in $\Omega\augteich(\Sigma_g)$ such that each $(Y, \eta)\in U$ has the following properties.
  \begin{itemize}
  \item The collapse $Y \to X$ induced by the markings is homotopic to a collapse $f\colon Y\to X$ which is
    conformal on $C$ and satisfies
    \begin{equation}
      \label{eq:quasitranslation}
      K^{-1} < \left|\frac{f^*\omega}{\eta}\right| < K
    \end{equation}
    on $C$.
  \item Each curve in $S$ which represents a polar node of $(X, \omega)$ is homotopic to a cylinder on $Y$ of height at least $h$.
  \item When $g=2$, for each curve $\gamma$ in $S$ which represents a polar node of $X$, if
    \begin{equation}
      \label{eq:periodzero}
      \int_\gamma \eta=0.
    \end{equation}
    then there are saddle connections $I_1$ and $I_2$ on $(Y, \eta)$ joining two distinct zeros of $\eta$ of length
    less than $\epsilon$ such that $\gamma$ is homotopic to $I_1 \cup I_2$ and $(Y, \eta)$ can be split along
    $I_1\cup I_2$.
  \end{itemize}
\end{theorem}

\begin{proof}
  The first statement follows directly from the definition of the conformal topology on
  $\Omega\augteich(\Sigma_g)$ above.

  Now let $p$ be a polar node of $(X, \omega)$.  By possibly enlarging the compact set $C$, we can suppose
  that $C$ contains a cylinder $D$ of height $2h$ around $p$.  If $f\colon Y\to X$ is a collapse satisfying
  \eqref{eq:quasitranslation} with $K<2$, then the two boundary components of $f^{-1}(D)$ are at least
  distance $h$ apart.  Furthermore, if $K$ is sufficiently small, then the boundary components of $f^{-1}(D)$
  both have winding number zero, so $f^{-1}(D)$ contains no zeros.  Therefore $f^{-1}(D)$ contains a cylinder
  of height at least $h$, which proves the second statement.
  
  Now suppose that $p$ is a holomorphic node of $(X, \omega)$ and suppose that \eqref{eq:periodzero} holds.
  If $\epsilon$ is sufficiently small, then there is a neighborhood $U$ of $p$ isomorphic to two
  disks of radius $\epsilon/\pi$ joined at a point.  By possibly enlarging the compact subset $C\subset X$, we
  have $\bdry U \subset C$.  If $([f\colon\Sigma_2\to Y], \eta)\in V_{X\setminus C, K}$ with $K<2$, then
  $f^{-1}(U)$ is bounded by two curves $\beta_i$ of circumference less than $2\epsilon$.  Furthermore, if $K$
  is sufficiently small, then the derivatives of $f$ will be close to the identity, and we can assume that the
  curves $\beta_i$ have winding number one and positive curvature.  Since the $\beta_i$ have winding number
  one, it follows from the Gauss-Bonnet Theorem that $\eta$ has two zeros in $f^{-1}(U)$, counting
  multiplicity.  Since the $\beta_i$ have positive curvature, there is a shortest geodesic $\delta\in
  f^{-1}(U)$ generating $\pi_1(f^{-1}(U))$.  By \eqref{eq:periodzero}, we have
  \begin{equation}
    \label{eq:periodzero2}
    \int_\delta \eta =0,
  \end{equation}
  which implies that $\delta$ is not a straight geodesic, so it must be a union of saddle connections joining
  distinct zeros.  As there are at most two zeros in $f^{-1}(U)$, the curve $\delta$ must be a union of two
  saddle connections $I_1$ and $I_2$ going from a zero $q$ to a zero $r$ of $\eta$.  By
  \eqref{eq:periodzero2}, we have
  \begin{equation}
    \label{eq:periodsame}
    \int_{I_1}\eta = \int_{I_2}\eta.
  \end{equation}
  This means that we can split along $I_1\cup I_2$ by the discussion of splitting along a union of saddle
  connections in \S\ref{subsec:flatgeometry}.  Finally, since $\delta$ is the shortest curve in its homotopy
  class in $f^{-1}(U)$ and the length of $\beta_i$ is less than $2\epsilon$, the length of $I_i$ must be less
  than $\epsilon$ as claimed.   
\end{proof}

\paragraph{Semi-Abelian varieties.}

We now introduce a topology on the space of semi-Abelian varieties and show that the map which associates to a
marked stable Riemann surface its Jacobian is continuous.

A \emph{semi complex torus} is a quotient $A=V/\Lambda$, where $V$ is a complex vector space containing a
discrete subgroup $\Lambda$, such that there is an exact sequence,
$$0\xrightarrow{}(\cx^*)^n\xrightarrow{} A \xrightarrow{\phi}B\xrightarrow{} 0,$$
where $B=V'/\Lambda'$ is a compact complex torus.  A principal polarization $h'$ on $B$ induces a degenerate Hermitian
form $h$ on $V$ by pulling back $h'$ by the lift $\tilde{\phi}\colon V \to V'$.  Such a Hermitian form on $V$ is a
\emph{principal polarization} on $A$.

Equip $\zed^{2g}$ with the usual symplectic form defined by \eqref{eq:sympform}.  A \emph{marked semi-Abelian
  variety} is a triple $(A, \Gamma, \phi)$, where
\begin{itemize}
\item $A$ is a semi-Abelian variety,
\item $\Gamma\subset\zed^{2g}$ is a subgroup of $\zed^{2g}$ on which the symplectic form vanishes and which is
  saturated in the sense that if $x\in\zed^{2g}$ and $mx\in\Gamma$ for some $m\in\zed$, then $x\in\Gamma$, and
\item $\phi$ is a symplectic isomorphism $\delta\colon \Gamma^\perp\to H_1(A; \zed)$.
\end{itemize}

Let $\semisiegelhalf[g]$ be the space of all marked semi-Abelian varieties.  The subspace of
$\semisiegelhalf[g]$ consisting of marked Abelian varieties is naturally isomorphic to the Siegel upper half
plane $\siegelhalf[g]$ by the discussion is \S\ref{subsec:abeliansurfaces}.

We give $\semisiegelhalf[g]$ a topology as follows.  If
$$(A_n = V_n/\Lambda_n, \Gamma_n, \phi_n)$$
is a sequence in $\semisiegelhalf[g]$, then we say that this sequence converges to $(A=V/\Lambda, \Gamma,
\phi)$ if eventually $\Gamma_n\subset\Gamma$ and there are linear isomorphisms $\psi_n\colon V_n\to V$ with
the following properties:
\begin{itemize}
\item For all $\alpha\in \Gamma^\perp$,
  $$\lim_{n\to\infty} \psi_n\circ\phi_n(\alpha) = \phi(\alpha),$$
\item For all $\alpha\in \zed^{2g}\setminus\Gamma^\perp$, the sequence $\psi_n\circ\phi_n(\alpha)$ eventually
  leaves every compact subset of $V$.
\item If $h_n$ and $h$ are the polarizations on $A_n$ and $A$, then
  $$\lim_n(\phi_n)_*(h_n) = h$$
  as Hermitian forms on $V$,
\end{itemize}
This defines a Hausdorff topology on $\semisiegelhalf[g]$.

Over $\semisiegelhalf[g]$ is the trivial rank-two bundle $\Omega\semisiegelhalf[h]$ consisting of pairs $(A,
\omega)$, where $A\in\semisiegelhalf[g]$ and $\omega\in\Omega(A)$, the space of holomorphic one-forms on $A$.

A choice of a symplectic isomorphism $H_1(\Sigma_2; \zed)\to \zed^4$ defines natural maps
$\Jac\colon\augteich(\Sigma_2) \to \semisiegelhalf[2]$ and $\Omega\Jac\colon\Omega\augteich(\Sigma_2) \to
\Omega\semisiegelhalf[2]$, sending a marked stable Riemann surface to its Jacobian with the induced marking.

\begin{theorem}
  \label{thm:jaccontinuous}
  The maps $\Jac$ and $\Omega\Jac$ defined above are continuous.
\end{theorem}

We will sketch the proof of this theorem below.

Given any $X\in\teich(\Sigma_2, S)$, we can define a norm $\|\cdot\|_X$ on $H_1(\Sigma_2;\reals)$, the
\emph{Hodge norm} as follows.  Let $V_S\subset H_1(\Sigma_2;\reals)$ be the subspace generated by the curves
in $S$.  If $\gamma\notin V_S^\perp\subset H_1(\Sigma_2; \reals)$, set $\|\gamma\|_X = \infty$.  On
$V_S^\perp$, let $\|\cdot\|_X$ be the norm induced by the Hermitian metric on $\Omega(X)^*$ coming
from the polarization via the embedding $V_S^\perp\to\Omega(X)^*$.  Alternatively, for $\gamma\in
V_S^\perp$,
\begin{equation*}
  \|\gamma\|_X = \sup_{\substack{
      \omega\in\Omega(X)\\
      \|\omega\| = 1}}
  |\omega(\gamma)|,
\end{equation*}
where
$$\|\omega\| = \left(\int_X |\omega|^2\right)^{1/2}.$$

\begin{theorem}
  \label{thm:jacestimates}
  Let $V_{S, \zed}\subset H_1(\Sigma_2; \zed)$ be the subgroup generated by the curves in $S$.  If $\{X_n\}$
  is a sequence in $\augteich(\Sigma_2)$ converging to $X$ in $\teich(\Sigma_2, S)$, then
  \begin{itemize}
  \item There exists a $C>0$ such that $\|\gamma\|_{X_n}>C$ for all $\gamma\in V_{S, \zed}^\perp\setminus V_{S;\zed}$.
  \item For all $D>0$, there exists $N>0$ such that $\|\gamma\|_{X_n}>D$ for all $\gamma\in
    H_1(\Sigma_2;\zed)\setminus V_{S, \zed}^\perp$ if $n>N$.
  \end{itemize}
\end{theorem}

Define for any curve system $S$ on $\Sigma_2$,
$$\augteich(\Sigma_2, S) = \teich(\Sigma_2)\cup \bigcup_{T\subset S} \teich(\Sigma_2, T),$$
where the union is
over all curve systems $T\subset S$.  The proof of these theorems will rely on the following lemma.

\begin{lemma}
  \label{lem:idoneoussection}
  Given any linear map $R\colon V_S\to\reals$, there is a unique section $Z\mapsto\omega_X$ of
  $\Omega\augteich(\Sigma_2, S)$ over $\augteich(\Sigma_2, S)$ such that:
  \begin{itemize}
  \item $\omega_X(\gamma)= R(\gamma)$ for all $\gamma\in V_S$, and
  \item $\Im \omega_X(\gamma) = 0 $ for all $\gamma\in V_S^\perp$.
  \end{itemize}
\end{lemma}

\begin{proof}
  Let $X\in \augteich(\Sigma_2, S)$.  We have an exact sequence,
  \begin{equation*}
    0\xrightarrow{} \Omega(\tilde{X})\xrightarrow{} \Omega(X)\xrightarrow{} \Hom(V_S, \cx)\xrightarrow{}0,
  \end{equation*}
  where $\tilde{X}$ is the normalization of $X$.  Choose some $\eta\in\Omega(X)$ such that
  $\omega(\gamma) = R(\gamma)$ for each $\gamma\in V_S$.   The form $\eta$ defines a map,
  $$S\colon V_S^\perp/V_S\to\reals,$$
  by $S(\gamma) = \Im \eta(\gamma)$.  We have isomorphisms,
  $$H_1(\tilde{X}; \reals)\to V_S^\perp/V_S,$$
  and
  $$\Omega(\tilde{X})\to \Hom(H_1(\tilde{X}; \reals), \reals),$$
  defined by $\omega \mapsto (\gamma \mapsto \Im \omega(\gamma))$.  Let $\nu\in \Omega(\tilde{X})$ induce the
  map $S$ above.  Then $\omega_X=\eta-\nu$ is the desired form, giving the desired section  over
  $\augteich(\Sigma_2, S)$.
\end{proof}

\begin{proof}[Sketch of Proof of Theorem~\ref{thm:jaccontinuous}]
  Suppose we have a sequence $\{X_n\}$ in $\augteich(\Sigma_2)$ converging to some $X\in \teich(\Sigma_2,S)$.
  To prove that $\Jac$ is continuous, we need to show:
  \begin{itemize}
  \item For every section $Y\mapsto\omega_Y$ of $\Omega\augteich(\Sigma_2, S)$ and $\gamma\in V_S^\perp$, we
    have
    $$\omega_{X_n}(\gamma)\to \omega_X(\gamma),$$
    and
  \item For every $\gamma\in H_1(\Sigma_2;\reals)\setminus V_S^\perp$, there is a section $Y\mapsto \eta_Y$ of
    $\Omega\augteich(\Sigma_2, S)$ such that $|\eta_{X_n}(\gamma)|\to\infty$.
  \end{itemize}
  We should also show that the polarizations converge, but we will omit the proof of this.
  
  If $X\mapsto \omega_X$ is a section of $\Omega\augteich(\Sigma_2, S)$, then $\omega_{X_n}$ converges
  uniformly to $\omega_X$ away from the nodes of $X$ as in Theorem~\ref{thm:longcylinders}.  If $\gamma\in
  V_S^\perp$, then $\gamma$ is represented by a curve on $X$ which is disjoint from the nodes, so
  $\omega_{X_n}(\gamma)\to\omega(\gamma)$ by uniform convergence.

  Suppose $\gamma\in  H_1(\Sigma_2;\reals)\setminus V_S^\perp$.  By Lemma~\ref{lem:idoneoussection}, there is a
  section $X\mapsto \eta_X$ of $\Omega\augteich(\Sigma_2, S)$ such that
  $$\sign \eta(\alpha) = \sign (\alpha\cdot\gamma)$$
  for every $\alpha\in H_1(\Sigma_2)$ representing a curve of $S$.  If $n$ is large, then each $\alpha\in
  H_1(\Sigma_2)$ representing a curve of $S$ is homologous to a very tall cylinder $C_\alpha$ on $(X_n,
  \eta_{X_n})$, and this cylinder contributes positively to $\Im\eta_{X_n}(\gamma)$.  If $n$ is sufficiently
  large, then the contributions from these cylinders will be much greater then that from the rest of $(X_n,
  \eta_{X_n})$, and we will have $\Im\eta_{X_n}(\gamma)\to \infty$.
\end{proof}

\begin{proof}[Proof of Theorem~\ref{thm:jacestimates}]
  Let $\{\eta^i\}$ be a basis of $\Omega(\tilde{X})$ such that $\|\eta^i\|=1$.  There is some $C'>0$ such that
  \begin{equation*}
    \sup_i |\eta^i(\gamma)|>C'
  \end{equation*}
  for all $\gamma\in V_{S,\zed}^\perp/V_{S,\zed}$.  Let
  $$C=\frac{1}{2}\sup_i C_i.$$
  For each $i$, there is a sequence $\eta^i_n\in\Omega(X_n)$ such that
  \begin{equation}
    \label{eq:etaconverges}
    \eta^i_n\to\eta^i.
  \end{equation}
  By Theorem~\ref{thm:longcylinders}, we have
  $$\|\eta^i_n\|\to\|\eta^i\|,$$
  so we can normalize $\eta^i_n$ so that $\|\eta^i_n\|=1$ for all $n$, and \eqref{eq:etaconverges} still
  holds.  By the uniform convergence statement of Theorem~\ref{thm:longcylinders}, we have for each $i$ and
  $\gamma\in V_{S,\zed}^\perp$,
  $$|\eta^i_n(\gamma)|>\frac{1}{2}|\eta^i(\gamma)|$$
  for $n$ sufficiently large.  The first claim of the Theorem follows.

  Let $\alpha_i\in H_1(\Sigma_2; \zed)$ be homology classes representing each nonseparating curve $S_i$ of
  $S$.  Choose a section $Y\mapsto\omega_Y$ of $\Omega\augteich(\Sigma_2; S)$ as in
  Lemma~\ref{lem:idoneoussection} such that
  \begin{equation*}    
    \omega(\alpha_i) = \reals \setminus \{0\}
  \end{equation*}
  for each $i$.  We claim that for every $C>0$ and $\gamma\in H_1(\Sigma_2; \zed)\setminus V_{S,\zed}^\perp$ such that
  \begin{equation}
    \label{eq:signs2}
    \sign \alpha\cdot\gamma = \sign \omega_X(\alpha),
  \end{equation}
  there is some $N$ for which we have
  \begin{equation*}
    \frac{|\omega_{X_n}(\gamma)|}{\|\omega_{X_n}\|}>C
  \end{equation*}
  for all $n>N$.  The second claim of the Theorem follows because we
  can choose finitely many such sections so that for all $\gamma\in
  H_1(\Sigma_2; \zed)\setminus V_{S,\zed}^\perp$, \eqref{eq:signs2}
  holds for some section.
  
  Let $M\subset X$ be a compact subset which is a deformation retract of the complement of the nodes.  By
  Theorem~\ref{thm:longcylinders}, there are collapses $g_n\colon X_n\to X$ such that
  \begin{equation*}
    K_n^{-1}<\left|\frac{(g_n)^*\omega_X}{\omega_{X_n}}\right| < K_n
  \end{equation*}
  on $M_n = g_n^{-1}(M)$ with $K_n\geq 1$ and $K_n \to 1$.  There are also cylinders $C_{i,n}$ in $X_n\setminus
  M_n$ which represent the homology classes $\alpha_i$ such that if $H_{i, n}=\height C_{i, n}$, then $\lim
  H_{i,n} = \infty$.

  The class $\gamma$ is represented on each $X_n$ by a simple closed curve $\gamma_n$.  On $X_n$, each
  cylinder $C_{i, n}$ contributes
  \begin{equation*}
    |\alpha_i\cdot\gamma| H_{i, n}
  \end{equation*}
  to $\Im\int_\gamma\omega_{X_n}$.  There are also segments of $\gamma_n$ which pass through $M_n$ to get from
  one cylinder to another.  Let $\beta$ be such a segment.  Since we chose $\omega_{X_n}$ so that its periods
  over $V_{S,\zed}^\perp$ are all real, $\Im\int_{\beta}\omega_{X_n}$ doesn't depend on the path $\beta$ takes.
  If $n$ is large, then this integral will be uniformly bounded over each such path $\beta$, say by
  $$\left|\Im\int_\beta \omega_{X_n}\right| < c = \diameter(M) + 1.$$
  
  We have
  \begin{equation*}
    \|\omega_{X_n}\|=\left(\sum_i H_{i, n} + \Area(M)\right)^{1/2} \leq \left(r \max_i H_{i, n}+ \Area(M)\right)^{1/2},
  \end{equation*}
  where $r$ is the number of curves in $S$.  Also, we have
  $$\Im \int_{\gamma_n} \omega \geq\max_i H_{i, n} - c.$$
  Thus
  \begin{equation*}
    \|\gamma\|_{X_n} \geq
    \frac{\left|\int_{\gamma_n}\omega_{X_n}\right|}{\|\omega_{X_n}\|}\geq\frac{\Im\int_{\gamma_n}\omega_{X_n}}{\|\omega_{X_n}\|} \geq \frac{\max_i H_{i,n}- c}{\sqrt{r\max_iH_{i, n} + \Area(M)}} \to \infty 
  \end{equation*}
  as $n\to\infty$.
\end{proof}

\paragraph{Real multiplication near $\bdry \moduli$.}

We now use the two previous theorems to study eigenforms for real multiplication near $\bdry\moduli$.  As in
the previous paragraph, we will write $V_{S,\zed}$ for the subgroup of $H_1(\Sigma_2; \zed)$ generated by the
curves in a curve system $S$ on $\Sigma_2$.

\begin{theorem}
  \label{cor:rmnearboundary}
  Let $X\in\teich(\Sigma_2, S)$ be a stable Riemann surface whose Jacobian has real multiplication by $\ord$.
  There is a neighborhood $U$ of $X$ in $\augteich(\Sigma_2)$ such that for each stable Riemann surface $Y\in U$
  which has real multiplication by $\ord$, the real multiplication preserves the two subspaces $V_{S, \zed}$ and
  $V_{S,\zed}^\perp$ of $H_1(\Sigma_2;\zed)$.
\end{theorem}

\begin{proof}
  If $\gamma\in V_{S,\zed}$, then we have $\|\gamma\|_X=0$, so $\|\gamma\|_{X_n}\to 0$ as $n\to\infty$.
  Let $\rho_n\colon\ord\to\End\Jac(X_n)$ be a choice of real multiplication, and let $\lambda\in\ord$ be a
  generator of $\ord$ over $\zed$.  If the claim of the Theorem is
  false, then taking a subsequence, we can assume that $\rho_n(\lambda)\cdot\gamma\notin V_{S,\zed}$ for all
  $n$.  We have
  \begin{equation*}
    \|\rho_n(\lambda)\cdot\gamma\|_{X_n} \leq \sup_i |\lambda^{(i)}| \cdot \|\gamma\|_{X_n}.
  \end{equation*}
  This is a contradiction because the right hand side of this equation goes to zero while the left hand side
  is bounded below by Theorem~\ref{thm:jacestimates}.

  The assertion that $V_{S,\zed}^\perp$ is preserved is proved in the same way.
\end{proof}

Recall that we have the subspace $\Omega\E\subset\Omega\moduli$ of eigenforms for real multiplication by
$\ord$.  We can define real multiplication on semi-Abelian varieties just as we did for Abelian varieties in
\S\ref{subsec:realmultiplication}; however, we will not require real multiplication on a semi-Abelian variety
to be proper.  Let $\overline{\Omega E}_D\subset\Omega\barmoduli$ be the closure of $\Omega\E$.

\begin{theorem}
   \label{cor:rmclosed}
  Every stable Abelian differential $(X, \omega)\in \overline{\Omega E}_D$ is a (not necessarily proper)
  eigenform for real multiplication by $\ord$.
\end{theorem}

\begin{proof}
  Suppose $(X_n, \omega_n^1)$ is a sequence in $\Omega\teich(\Sigma_2)$ is a sequence of $\iota_1$-eigenforms
  for real multiplication by $\ord$ converging to some $(X, \omega^1)$ with $X\in\teich(\Sigma_2, S)$.  Let
  $\omega^2_n\in \Omega(X_n)$ be a sequence of $\iota_2$-eigenforms.  Taking a subsequence and suitably
  normalizing the $\omega^2_n$, we can assume that
  $\omega^2_n\to\omega^2$ for some nonzero $\omega^2\in\Omega(X)$.

  We claim that the $\omega^i$ span $\Omega(X)$.  Suppose not.  Then $\omega^1 = c \omega^2$ for some $c\in\cx$.
  Choose some $\alpha\in V_{S, \zed}^\perp$ such that $\omega^1(\alpha)\neq0$, and let
  $\lambda\in\ord\setminus\zed$ be a generator of $\ord$ over $\zed$.  Then for $n$ large,
  \begin{equation*}
    \lambda^{(1)}\omega^1(\alpha) \sim \omega^1_n(\rho_n(\lambda)\cdot\alpha) \sim c
    \omega^2_n(\rho_n(\lambda)\cdot\alpha) = c \lambda^{(2)}\omega^2_n(\alpha)\sim \lambda^{(2)}\omega^1(\alpha),
  \end{equation*}
  a contradiction because $\lambda^{(1)}\neq \lambda^{(2)}$.

  Now via the dual bases to $\{\omega_n^i\}$ and $\{\omega^i\}$, we can identify $\Omega(X_n)^*$ and
  $\Omega(X)^*$ with $\cx^2$.  We then have natural embeddings
  \begin{gather*}\phi_n\colon H_1(\Sigma_2; \zed)\to\cx^2, \quad \text{and} \\
    \phi\colon V_{S, \zed}^\perp \to \cx^2
  \end{gather*}
  such that
  \begin{align*}
    \Jac(X_n) & \isom \cx^2/\Im\phi_n,\quad \text{and} \\
    \Jac(X) & \isom \cx^2/\Im\phi,
  \end{align*}
  and $\phi_n\to\phi$ as $n\to\infty$.  The real multiplication on
  $\Jac(X_n)$ lifts to
  $$\tilde{\rho}(z_1, z_2) = (\lambda^{(1)} z_1, \lambda^{(2)} z_2)$$
  on $\cx^2$.  Since $\phi_n\to\phi$, the map $\tilde{\rho}$ preserves $\Im \phi$, so defines real
  multiplication on $\Jac(X)$ with $\omega^1$ an eigenform as desired.
\end{proof}


%% file: node.pstex_t
\begin{picture}(0,0)%
\includegraphics{node.pstex}%
\end{picture}%
\setlength{\unitlength}{3947sp}%
\begingroup\makeatletter\ifx\SetFigFont\undefined%
\gdef\SetFigFont#1#2#3#4#5{%
  \reset@font\fontsize{#1}{#2pt}%
  \fontfamily{#3}\fontseries{#4}\fontshape{#5}%
  \selectfont}%
\fi\endgroup%
\begin{picture}(3946,1844)(4629,-5558)
\put(5626,-4711){\makebox(0,0)[lb]{\smash{{\SetFigFont{12}{14.4}{\rmdefault}{\mddefault}{\updefault}{\color[rgb]{0,0,0}$\alpha$}%
}}}}
\put(7801,-4711){\makebox(0,0)[lb]{\smash{{\SetFigFont{12}{14.4}{\rmdefault}{\mddefault}{\updefault}{\color[rgb]{0,0,0}$\beta$}%
}}}}
\end{picture}%

%% file: localcoordinates.tex
\section{Local coordinates for strata in $\Omega\overline{{\cal M}}_2$}
\label{sec:localcoordinates}

While the stratification of $\Omega\moduli$ is very simple, consisting of only two strata, the natural
stratification of $\Omega\barmoduli$ is much more complicated.  With the appropriate definition of a stratum,
there are seventeen different strata in $\Omega\barmoduli$.  It will be necessary for our study of the
compactification of $\X$ to understand them.  In particular, we will need to give local
coordinates on $\Omega\barmoduli$ around most of these strata.  In \S\ref{subsec:strata}, we will define the
stratification of $\Omega\barmoduli$, give notation for the various strata, and list the seventeen strata in
$\Omega\barmoduli$.  In the rest of this section, we will then discuss
each stratum in turn.

\subsection{Strata in $\Omega\overline{{\cal M}}_2$}
\label{subsec:strata}

Given two stable Abelian differentials $(X_1, \omega_1)$ and $(X_2, \omega_2)\in \Omega\barmoduli$, say that
they are in the same \emph{stratum} of $\Omega\barmoduli$ if there is a homeomorphism $f\colon X_1\to X_2$
with the following properties:
\begin{itemize}
\item $f$ takes zeros of $\omega_1$ to zeros of $\omega_2$ of the same order.
\item $f$ takes polar nodes of $(X_1, \omega_1)$ to polar nodes of $(X_2, \omega_2)$.
\item $f$ takes irreducible components of $X_1$ on which $\omega_1$ vanishes to irreducible components of
  $X_2$ on which $\omega_2$ vanishes.
\end{itemize}

\paragraph{Notation for strata.}

We now introduce notation for the various strata in $\Omega\barmoduli$.  Let $S\subset
\Sigma_2$ be a curve system, and define
$\Omega\teich_2^0(S)$
to be the locus of all $(X, \omega)\in \Omega\augteich_2$ such that:
\begin{itemize}
\item $X\in \teich_2(S)$;
\item The form $\omega$ has poles at the nodes corresponding to nonseparating curves in $S$.
\end{itemize}
If $T\in S$ is a nonseparating curve, then define
$\Omega\teich_2^0(S, T)$
to be the locus of all $(X, \omega)\in \Omega\augteich_2$ such that:
\begin{itemize}
\item $X\in \teich_2(S)$;
\item The form $\omega$ has poles at the nodes corresponding to nonseparating curves in $S\setminus T$;
\item The form $\omega$ is holomorphic at the node corresponding to the curve $T$.
\end{itemize}

If ${\bf n}$ is either $(2)$ or $(1, 1)$, then let
$$\Omega\teich_2^0(S; {\bf n})\subset \Omega\teich_2^0(S)$$
be the locus of $(X, \omega)\in \Omega\teich_2^0(S)$ where in addition, $\omega$ has zeros whose orders are
given by ${\bf n}$.

Let $\Omega\moduli^0(S, T)$ and $\Omega\moduli^0(S)$ and $\Omega\moduli^0(S; {\bf n})$ be the images in
$\Omega\barmoduli(S)$ of the corresponding spaces defined above.

Let $\Omega\Def_2'(S, T)$ be the locus of Abelian differentials $(X, \omega)\in \Omega\Def_2(S)$ such that
either $\int_T\omega=0$ or $T$ represents a holomorphic node on $(X, \omega)$.  This is a nonsingular
hypersurface in $\Omega\Def_2(S)$.

Let $\Omega\Def^0_2(S)$ be the locus of Abelian differentials $(X,
\omega)\in \Omega\Def_2(S)$ such that each nonseparating curve of
$S\setminus T$ is represented by either a polar node of $(X, \omega)$
or a cylinder on $(X, \omega)$.  Given a curve $T\in S$, let
$$\Omega\Def_2^0(S, T)\subset \Omega\Def_2'(S, T)$$
be the locus of $(X, \omega)$ with the following properties:
\begin{itemize}
\item Each nonseparating curve of $S\setminus T$ is represented by either a polar node of $(X, \omega)$ or a
  cylinder on $(X, \omega)$.
\item A curve $\gamma$ which  is either $T$ or a separating curve in $S$ is either represented by a holomorphic node of $(X,
  \omega)$ or is homotopic to a union $I\cup J(I)$ of a saddle connection $I$ with its image under the
  hyperelliptic involution $J$.
\end{itemize}
By Theorem~\ref{thm:longcylinders}, $\Omega\Def^0_2(S, T)$ is a neighborhood of
$\Omega\teich_2^0(S, T)$ in $\Omega\Def_0'(S, T)$.

Also, let
$$\Omega\Def^0_2(S, {\bf n})\subset\Omega\Def^0_2(S)$$
be those $(X, \omega)\in\Omega\Def_2^0(S, T)$ such that $\omega$ has zeros whose orders are given by ${\bf n}$.

\begin{figure}[htbp]
  \centering
  \input{curvesystems.pstex_t}
  \caption{Six curve systems on $\Sigma_2$.}
  \label{fig:curvesystems}
\end{figure}
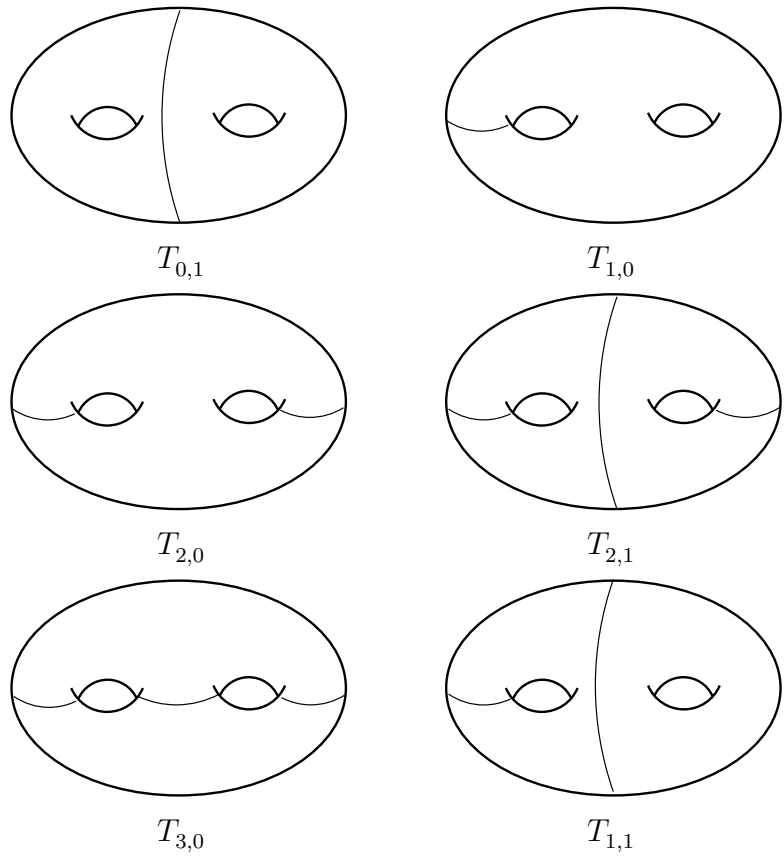

Up to the action of the mapping class group, there are six curve systems in $\Sigma_2$, shown in
Figure~\ref{fig:curvesystems}.  We will denote by $T_{i,j}$ a curve system with $i$ separating curves and $j$
nonseparating curves.  Such a curve system is unique up to the action of the mapping class group.  We will
write $T_{i,j}^k$  with $k=1, \ldots n$ for the $n$ individual curves of $T_{i,j}$, and we will always order the
curves so that if $T_{i,j}$ contains a separating curve, then it is the last curve $T_{i,j}^n$.

\paragraph{The seventeen strata.}

Here is a list of the seventeen strata in $\Omega\barmoduli$.

\begin{enumerate}
\item $\Omega\moduli(1,1)$ is the stratum of nonsingular Abelian differentials with two simple zeros.
\item $\Omega\moduli(2)$ is the stratum of nonsingular Abelian differentials with a double zero.
\item $\Omega\moduli^0(\system{1})$ is the stratum of stable Abelian differentials consisting of two nonzero
  genus one Abelian differentials joined at a node.
\item $\Omega\moduli^0(\system{2}; 1,1)$ is the stratum of stable Abelian differentials with one nonseparating polar
  node and two simple zeros.
\item $\Omega\moduli^0(\system{2}; 2)$ is the stratum of stable Abelian differentials with one nonseparating polar
  node and one double zero.
\item $\Omega\moduli^0(\system{2}, \system[1]{2})$ is the stratum of stable Abelian differentials with one
  nonseparating holomorphic node.  These can be regarded as genus one Abelian differentials with two points joined to
  a node, and they have no zeros.
\item $\Omega\moduli^0(\system{3}, \system[2]{3})$ is the stratum of stable Abelian differentials with one
  nonseparating holomorphic node and one nonseparating polar node.  These can be regarded as the differential
  $a\, dz/z$ on $\cx^*$ with the points $0$ and $\infty$ identified to form a polar node, and with two points in
  $\cx^*$ identified to form a holomorphic node.
\item $\Omega\moduli^0(\system{3}; 1, 1)$ is the stratum of stable Abelian differentials with two nonseparating
  polar nodes and two simple zeros.
\item $\Omega\moduli^0(\system{3}; 2)$ is the stratum of stable Abelian differentials with two nonseparating
  polar nodes and a double zero.
\item $\Omega\moduli^0(\system{4})$ is the stratum of stable Abelian differentials with two nonseparating polar
  nodes and one separating holomorphic node. These can be regarded as two infinite cylinders with the ends of
  the cylinders identified to form polar nodes and one point from each cylinder identified to form a
  holomorphic node.
\item $\Omega\moduli^0(\system{5})$ is the stratum of stable Abelian differentials $(X, \omega)$ with three nonseparating
  polar nodes.  These have two simple zeros, one in each irreducible component of $X$.
\item $\Omega\moduli^0(\system{5}, \system[3]{5})$ is the stratum of stable Abelian differentials with two nonseparating
  polar nodes and one nonseparating holomorphic node.  These differentials have no zeros.
\item $\Omega\moduli^0(\system{6})$ is the stratum of stable Abelian differentials with one nonseparating polar
  node and one separating holomorphic node.
\item[14-17.] There are four other strata consisting of stable Abelian differentials $(X, \omega)$ which
  vanish on some irreducible component of $X$.  There is one such stratum corresponding to the curve system
  $\system{1}$, one corresponding to the curve system $\system{4}$, and two corresponding to the curve system
  $\system{6}$. 
\end{enumerate}

We will now go down this list and discuss each of these strata.  We will not discuss the stratum
$\Omega\moduli(1,1)$ because we don't need to know anything more about it, and we will not discuss the strata
13-17 on this list because they do not arise in this paper.

\subsection{The strata $\Omega{\cal M}_2(2)$ and $\Omega{\cal M}_2(T_{0,1})$}
\label{subsec:stratum1}

The strata $\Omega\moduli(2)$ and $\Omega\moduli(\system{1})$ are suborbifolds of $\Omega\barmoduli$.  We will
see this explicitly by extending the period coordinates on these strata to period coordinates on a
neighborhood of these strata.

\paragraph{Period coordinates for $\Omega\moduli(2)$.}

The map ${\rm Per}\colon\Omega\teich_2\to H^1(\Sigma_2; \cx)$ is a submersion by the period coordinates from
\S\ref{subsec:abelianmoduli}, so the fibers of this map are the leaves of a foliation $\mathcal{A}$ of
$\Omega\teich_2$ by Riemann surfaces along whose leaves the absolute periods of Abelian differentials
are constant.  This foliation meets the stratum $\Omega\teich_2(2)$ transversely because the
restriction of ${\rm Per}$ to this stratum is locally biholomorphic.

We can locally parameterize the leaf $L$ of $\mathcal{A}$ through some $(X, \omega)\in\Omega\teich_2(2)$
using the operation of splitting a double zero from \S\ref{subsec:flatgeometry} as follows.  Choose
$\epsilon>0$ such that every geodesic segment starting at the zero $p$ of $(X, \omega)$ of length less than $\epsilon^3$
is embedded in $X$.  Since $p$ is a cone point of cone angle $6\pi$, we can parameterize short
geodesic segments starting at $p$ by the $\epsilon$-disk $\Delta_\epsilon\subset \cx$ by associating continuously to a point
$z\in\Delta_\epsilon$ a geodesic segment $I(z)\subset X$ such that
\begin{equation}
   \int_{I(z)}\omega = z^3.
\end{equation}
The map $\phi\colon \Delta_\epsilon\to L$ defined by
$$\phi(z) = (X, \omega) \#_{I(z)}$$
is holomorphic and satisfies $\phi(z) = \phi(-z)$ because the segments
$I(z)$ and $I(-z)$ differ by an angle of $3\pi$, and so these segments determine the same ``X'' in $(X,
\omega)$ (see Figure~\ref{fig:splittingadoublezero}), along which we perform surgery to form $\phi(z)$.  Furthermore, $\phi(z)=\phi(w)$ if
and only if $z=\pm w$, so   $\phi$ lifts to a map $\tilde{\phi}$ defined by
$$\tilde{\phi}(z) = (X, \omega) \#_{I(z^{1/2})},$$
which maps $\Delta_\epsilon$ conformally onto a neighborhood of $(X, \omega)$ in $L$. 

Now let $U\subset\Omega\teich_2(2)$ be a small neighborhood of $(X, \omega)$.  We can parameterize the
above construction to define continuously for each $(Y, \eta)\in U$ and $z\in\Delta_\epsilon$ (possibly making
$\epsilon$ smaller) a geodesic segment $I_{(Y, \eta)}(z)$ on $Y$ starting at the zero of $\eta$ such that
$$\int_{I_{(Y, \eta)}(z)}\eta = z^3.$$
Define $\Phi\colon U\times \Delta_\epsilon\to \Omega\teich_2$ by
$$\Phi((Y, \eta), z) = (Y, \eta) \#_{I_{(Y, \eta)}(z^{1/2})}.$$

$\Phi$ gives a conformal isomorphism of $\{(Y, \eta)\}\times\Delta_\epsilon$ onto a neighborhood of $(Y,
\eta)$ in the leaf of $\mathcal{A}$ through $(Y, \eta)$.  Since $\mathcal{A}$ is transverse to
$\Omega\teich_2(2)$, it follows that $\Phi$ is biholomorphic onto its image in $\Omega\teich_2(2)$ if $U$ and
$\epsilon$ are sufficiently small.  By considering the inverse of $\Phi$, we obtain the following local
coordinates around points in $\Omega\teich_2(2)$.

\begin{prop}
  \label{prop:parameterizedsplitting}
  Every $(X, \omega)\in\Omega\teich(2)$ has a neighborhood $V$ such that:
  \begin{itemize}
  \item Each $(Y, \eta)\in V$ has either a shortest saddle connecting distinct zeros or a double zero.
  \item The multi-valued function,
    \begin{equation*}
      \label{eq:yrelation}
      (Y, \eta) \mapsto \left(\int_I\eta\right)^{2/3},
    \end{equation*}
    has a single-valued  branch $y$.
  \item The product ${\rm Per}\times y\colon V\to H^1(\Sigma_2; \cx)\times\cx$ is biholomorphic onto its
    image.
  \end{itemize}
\end{prop}

\paragraph{Period coordinates for $\Omega\moduli(\system{1})$.}

The stratum $\Omega\moduli(\system{1})$ is a suborbifold of $\Omega\barmoduli$ isomorphic to
the symmetric square of $\Omega\moduli[1,1]$.  We  give local coordinates for this stratum just as for
$\Omega\moduli(2)$ using the connected sum construction.  As above, let ${\rm Per}\colon
\Omega\Def_2(\system{1})\to H^1(\Sigma_2; \cx)$ be the natural period map.

Each $(X, \omega)\in\Omega\teich_2(\system{1})$ is a one point connected sum of two
genus one differentials:
$$(X, \omega) = (X_1, \omega_1) \# (X_2, \omega_2).$$
Given a small open set $U\subset \Omega\teich_2(\system{1})$ and a sufficiently small $\epsilon>0$, we
define a holomorphic map $\Phi\colon U\times \Delta_\epsilon\to\Omega\Def_2(\system{1})$ by
$$\Phi((X, \omega), z) = (X_1, \omega_1) \#_{I(z^{1/2})} (X_2, \omega_2),$$
where $I(z)\subset\cx$ is the segment joining $0$ to $\epsilon$.

If $(Y, \eta)\in \Omega\Def_2(\system{1})$ has a unique shortest saddle connection $I$, then define
\begin{equation*}
  y(Y, \eta) = \left(\int_I \eta \right)^2.
\end{equation*}
By taking an inverse of $\Phi$, we obtain local coordinates around $\Omega\teich_2(\system{1})$ as above.

\begin{prop}
   \label{prop:parameterizedsum}
    Every $(X, \omega)\in\Omega\teich_2(\system{1})$ has a neighborhood $V\subset\Omega\Def_2(\system{1})$ such that:
  \begin{itemize}
  \item Each $(Y, \eta)\in V$ has either a shortest saddle connecting distinct zeros or a separating node.
  \item The function ${\rm Per}\times y\colon V\to H^1(\Sigma_2;\cx)\times\cx$ is biholomorphic onto its
    image.
  \end{itemize}
\end{prop}

\subsection{The strata $\Omega{\cal M}_2^0(T_{1,0}; 1, 1)$ and $\Omega{\cal M}_2^0(T_{1,0}; 2)$}
\label{subsec:stratum2}

The locus $\Omega\moduli^0(\system{2})$ --  the union of the strata $\Omega\moduli^0(\system{2}; 1,
1)$ and $\Omega\moduli^0(\system{2}; 2)$ -- is a four dimensional orbifold.  It can be naturally embedded in the
bundle $\Omega'\moduli[1,2]$ consisting of meromorphic Abelian differentials $(X, \omega)$ on an elliptic
curve $X$ with simple poles at two marked points $p$ and $q$ such that
$$\Res_p\omega=-\Res_q\omega.$$
Its image is the complement of the sub-line-bundle where these residues are zero.  The stratum
$\Omega\moduli^0(\system{2};1,1)$ is an open, dense suborbifold of $\Omega\moduli^0(\system{2})$, and
$\Omega\moduli^0(\system{2};2)$ is a closed, three-dimensional suborbifold.

\paragraph{Local coordinates.}

We now give local coordinates on neighborhoods of these strata.  Choose a symplectic basis $\{\alpha_i,
\beta_i\}_{i=1}^2$ of $H_1(\Sigma_2; \zed)$ such that $\alpha_1$ represents the single curve in $\system{2}$.

Define functions on $\Omega\Def^0_2(\system{2})$ as follows: given $(Y, \eta)\in
\Omega\Def^0_2(\system{2})$, define 
\begin{align*}
  v &= \eta(\alpha_1) & w &= \eta(\alpha_2) \\
  x &= \eta(\beta_2) & z &= e^{2\pi i \eta(\beta_1)/\eta(\alpha_1)},
\end{align*}
where we consider $z(Y, \eta)$ to be $0$ if $Y$ has a node.  These are well-defined functions because
$\alpha_1$, $\alpha_2$, and $\beta_2$ give well-defined homology classes on $Y$ via the marking, and $\beta_1$
gives a homology class on $Y$ which is well-defined up to adding a multiple of $\alpha_1$.  These are all
holomorphic functions on $U$.  This follows from the fact that these functions are all holomorphic on
$\Omega\Def^0_2(\system{2})\setminus\Omega\teich^0_2(\system{2})$ because periods vary holomorphically, together
with the fact that these functions are all continuous on $U$, which follows from
Theorem~\ref{thm:longcylinders}.  Functions which are continuous and holomorphic on an open dense subset of
their domain are holomorphic, so these functions $v$, $w$, $x$, and $z$ are all holomorphic.

Now, let $(X, \omega)\in \Omega\teich^0_2(\system{2})$, and let $U\subset\Omega\Def^0_2(\system{2})$ be an open
neighborhood of $(X, \omega)$ with a continuous choice of a saddle connection $I_{(Y, \eta)}$ on each $(Y,
\eta)\in U$ joining distinct zeros.  If $(X, \omega)$ has two simple zeros, then we can also continuously
orient these saddle connections, but this is not possible if $(X, \omega)$ has a double zero.  For $(Y,
\eta)\in U$, define
$$y =
\begin{cases}
  \int_{I_{(Y, \eta)}}\eta & \text{if $(X, \omega)$ has two simple zeros;}\\
  \left(\int_{I_{(Y, \eta)}}\eta\right)^{2/3} & \text{if $(X, \omega)$ has a double zero,}
\end{cases}
$$

If $U$ is sufficiently small, then the functions $(v, w, x, y)$ are a system of holomorphic local coordinates
on $\Omega\teich^0_2(\system{2})\cap U$ by (a simple generalization of) the period coordinates from
\cite{veech90} together with the argument of Proposition~\ref{prop:parameterizedsplitting}.

Recall that in \S\ref{subsec:flatgeometry}, we defined the unplumbing operation which replaced an Abelian
differential with a choice of a cylinder with a new meromorphic Abelian differential with two simple poles.
Since these poles have opposite residues, we can join them to a node and obtain a stable Abelian differential.

Define the \emph{unplumbing map} to be the function,
$$\psi\colon\Omega\Def^0_2(\system{2})\to\Omega\teich^0_2(\system{2}),$$
defined by sending $(Y, \eta)$ to the
Abelian differential obtained by  unplumbing the cylinder on $Y$ which is homotopic to the single curve in the
curve system $\system{2}$.  The function $\psi$ is continuous and is holomorphic on
$$\Omega\Def^0_2(\system{2})\setminus \Omega\teich^0_2(\system{2}),$$
as can easily be seen using period
coordinates, so $\psi$ is holomorphic on all of $\Omega\Def^0_2(\system{2})$.  Let $\Psi=\psi\times z$.

\begin{prop}
  \label{prop:coordinates1}
  The map $\Psi$ is a biholomorphic map of $\Omega\Def^0_2(\system{2})$ onto its image in
  $\Omega\teich^0_2(\system{2})\times\cx$.
  
  It follows that if $U$ is sufficiently small, the functions $(v, w, x, y, z)$ are a holomorphic system of
  local coordinates on the neighborhood $U$ of $(X, \omega)$.
\end{prop}

\begin{proof}
  We can show that $\Psi$ is biholomorphic onto its image by constructing a holomorphic inverse.
  Given $(Y, \eta)\in \Omega\teich^0_2(\system{2})$ and $u\in\cx^*$, if $|u|$ is sufficiently small, then we can
  uniquely plumb a cylinder into the node of $(Y, \eta)$ to form a new nonsingular Abelian differential $(Y',
  \eta')$ such that $z(Y', \eta')=u$.  The point is that the plumbing construction is determined by two
  quantities: the height of the resulting cylinder, and the twist of the gluing map; the height is determined
  by $|u|$, and the twist is determined by $\arg u$.

  Let $V\subset\Omega\teich^0_2(\system{2})\times\cx$ be the open neighborhood of $\Omega\teich^0_2(\system{2})\times\{0\}$ where
  this operation is defined, and let $\Phi\colon V\to\Omega\Def^0_2(\system{2})$ be the plumbing map.  The map
  $\Phi$ is continuous and is holomorphic on $V \setminus (\Omega\teich^0_2(\system{2})\times\{0\})$ by period
  coordinates, so it is everywhere holomorphic.  The maps $\Phi$ and $\Psi$ are inverse to each other, so both
  are biholomorphic.

  The last statement then follows from the fact that $(v, w, x, y)$ are a system of local coordinates on
  $\Omega\teich^0_2(\system{2})$ on a neighborhood of $(X, \omega)$.
\end{proof}

\paragraph{Cylinder covers.}

The stratum $\moduli(\system{2})\subset\barmoduli$ of stable Riemann surfaces with one nonseparating node is
naturally isomorphic to the moduli space $\moduli[1,2]$ of elliptic curves with two marked points, and we will
implicitly identify them.  In $\moduli(\system{2})$, let $\moduli(\system{2})(d)$ be the locus of all $X$ such
that the two marked points on $X$ differ by exactly $d$-torsion in the group law on the underlying elliptic
curve. This locus $\moduli(\system{2})(d)$ is isomorphic to modular curve $\half/\Gamma_1(d)$.

Given a marked elliptic curve $X$ in $\moduli[1,2]$, call a meromorphic function $f\colon X\to\proj^1$ a
\emph{cylinder cover} if the poles and zeros of $f$ are located at the two marked points of $X$.  By
Riemann-Hurwitz, such an $f$ has two other branch points counting multiplicity.  We say that $f$ is primitive
if $f$ is not of the form $g^r$ for some meromorphic function $g$ and $r>1$.  A meromorphic Abelian
differential $(X, \omega)$ in $\Omega' \moduli[1,2]$ is a \emph{cylinder covering differential} if
$\omega=df/f$ for some primitive cylinder cover.  The \emph{degree} of $\omega$ is the degree of the primitive
cover.  In terms of stable Riemann surfaces, a cylinder cover is a morphism $f\colon X\to C$, where
$X\in\moduli(\system{2})$, the curve $C$ is $\proj^1$ with $0$ and $\infty$ identified to form a node, and the
inverse image of the node of $C$ is the node of $X$.

Let $\Omega\moduli^0(\system{2})(d)$ be the locus of
degree $d$ cylinder covering differentials in $\Omega\moduli^0(\system{2})$.

\begin{prop}
  \label{prop:cylindricalclassification}
  A stable Riemann surface $X\in\moduli(\system{2})$ has a degree $d$ cylinder covering differential if and
  only if $X\in \moduli(\system{2})(d)$, in which case it unique up to constant multiple.
\end{prop}

\begin{proof}
  To prove uniqueness, assume that $f$ and $g$ are two degree $d$ cylinder covers. This means that $f$ and $g$
  have the same zeros and poles, so $f=cg$ for some $c\in\cx$.  Thus $df/f=dg/g$.
  
  Recall that if $E$ is an elliptic curve and $\Div^0(E)$ is the group of divisors of degree $0$ on $E$, then
  there is a natural map $\mu\colon \Div^0(E)\to E,$ defined via the group law on $E$. Abel's Theorem says
  that a divisor $D\in\Div^0(E)$ is the divisor of a meromorphic function on $E$ if and only if $\mu(D)=0$
  (for example, see \cite[p.  235]{griffithsharris}).

  Now let $X\in\moduli(\system{2})$, which we'll regard as an elliptic curve with two marked points $p$ and
  $q$.  Choose an identity point $0$ on $X$, fixing the group law.  If $X$ has a degree $d$ cylinder cover
  $f\colon X\to \proj^1$, with its zeros at $p$ and its poles at $q$, then $(f)=d(p-q)$ (where $(f)$ denotes
  the divisor associated to $f$), and by Abel's Theorem, $d(p-q)=0$ in the group law on $E$; that is, $p$ and
  $q$ differ by $d$-torsion.
  
  We need to show that there is no $e<d$ with $e | d$ and $e(p-q)=0$.  Suppose there is such an $e$.  Then by
  Abel's Theorem again, there would be a meromorphic function $h$ on $X$ with $(h)=e(p-q)$, which implies that
  $h$ is a degree $e$ cylinder cover.  By the uniqueness statement, $f=ch^{d/e}$ for some $c\in\cx$,
  contradicting primitivity of $f$.
  
  The proof of the converse of this proposition is also a straightforward application of Abel's Theorem and
  will be left to the reader.
\end{proof}

\begin{cor}
  The suborbifold $\proj\Omega\moduli^0(\system{2})(d)$ of $\proj\Omega\moduli^0(\system{2})$ is isomorphic to
  $\half/\Gamma_1(d).$
\end{cor}

\subsection{The stratum $\Omega{\cal M}_2^0(T_{1,0}, T_{1,0}^1)$}
\label{subsec:stratum3}

We now give local coordinates around the stratum $\Omega\moduli^0(\system{2}, \system[1]{2})$.  We will actually
only give local coordinates on a hypersurface containing this stratum, but this will be sufficient for our
purposes.

Let $\{\alpha_i,
\beta_i\}_{i=1}^2$ be a symplectic basis of $H_1(\Sigma_2; \zed)$ as in the previous section.  Given $(Y,
\eta)\in\Omega\Def^0_2( \system{2}, \system[1]{2})$, define
\begin{align*}
  w &= \eta(\alpha_2) & x &= \eta(\beta_2) \\
  y & = \eta(\beta_1) & z &= \left(\int_I\eta\right)^2,
\end{align*}
where $I$ is one of the two saddle connections on $(Y, \eta)$ such that $I \cup J(I)$ is homotopic to the
curve $\system[1]{2}$ in $\system{2}$.  These are holomorphic functions on $\Omega\Def^0_2( \system{2},
\system[1]{2})$.  If $(Y, \eta)\in \Omega\teich^0_2( \system{2}, \system[1]{2})$ is regarded as a genus one
differential with two marked points, then $y(Y, \eta)$ is the integral of $\eta$ along a path joining the
marked points.  The functions $(w, x, y)$ are a system of holomorphic local coordinates on
$\Omega\teich^0_2( \system{2}, \system[1]{2})$.

\begin{prop}
  \label{prop:coordinates2}
  For any $(X, \omega)\in\Omega\teich^0_2( \system{2}, \system[1]{2})$, the functions $(w, x, y, z)$ define a
  system of holomorphic local coordinates on a sufficiently small neighborhood of $(X, \omega)$ in
  $\Omega\Def^0_2( \system{2}, \system[1]{2})$.
\end{prop}

\begin{proof}[Sketch of proof]
  The proof of this proposition is completely analogous to the proof of Proposition~\ref{prop:coordinates1}.
  Let
  $$\psi\colon\Omega\Def^0_2( \system{2}, \system[1]{2})\to\Omega\teich^0_2( \system{2}, \system[1]{2})$$
  be the map which sends an Abelian differential to the one obtained by splitting along the union of saddle
  connections $I\cup J(I)$ homotopic to the curve $\system[1]{2}$.  Let $\Psi = \psi\times z$.  This is a
  holomorphic map,
  $$\Psi\colon\Omega\Def^0_2( \system{2}, \system[1]{2})\to\Omega\teich^0_2( \system{2}, \system[1]{2})\times\cx$$

  Given any $u\in\cx$, let $I(u)\subset\cx$ be the segment connecting $0$ to $u$.  For any $(Y, \eta)\in
  \Omega\teich^0_2( \system{2}, \system[1]{2})$, regarded as a genus one differential with two marked points, if
  $|u|$ is sufficiently small, then there are two parallel embeddings $$\epsilon_i\colon I(u^{1/2})\to(Y, \eta),$$ each
  sending $0$ to one of the marked points of $(Y, \eta)$.  Let
  $$V\subset\Omega\teich^0_2( \system{2}, \system[1]{2})\times\cx$$
  be the neighborhood of
  $\Omega\teich^0_2( \system{2}, \system[1]{2})\times\{0\}$ where these two embeddings $\epsilon_i$
  exist, and define
  $$\Phi\colon V \to \Omega\Def^0_2( \system{2}, \system[1]{2})$$
  to be the map which replaces $(Y, \eta)$
  with the self connected sum
  $$(Y, \eta) \#_{I(z^{1/2})}.$$

  $\Phi$ is a holomorphic map, and since $\Phi$ and $\Psi$ are inverse
  to each other, both are biholomorphic.  It then follows from the fact that $(w, x, y)$ give  local coordinates
  on $\Omega\teich^0_2( \system{2}, \system[1]{2})$ around $(X, \omega)$ that $(w, x, y, z)$ give local coordinates on
  $\Omega\Def^0_2( \system{2}, \system[1]{2})$ around $(X, \omega)$.
\end{proof}

\paragraph{The $d$-torsion locus.}

Let
$$\proj\Omega\moduli^0(\system{2}, \system[1]{2})(d) \subset \proj\Omega\moduli^0(\system{2}, \system[1]{2})$$
be the locus of genus one differentials such that the two marked points differ by torsion of degree exactly
$d$.  One can show using period coordinates that $\proj\Omega\moduli^0(\system{2}, \system[1]{2})(d)$ is a
suborbifold of $\proj\Omega\moduli^0(\system{2}, \system[1]{2})$.  We have the isomorphism,
$$\proj\Omega\moduli^0(\system{2}, \system[1]{2})(d) \isom \moduli(\system{2})(d)\isom\half/\Gamma_1(d).$$

\subsection{The stratum $\Omega{\cal M}_2^0(T_{2,0}, T_{2,0}^2)$}
\label{subsec:stratum4}

We now define local coordinates on a hypersurface containing the stratum $\Omega\moduli^0(\system{3},
\system[2]{3})$.  Choose a symplectic basis $\{\alpha_i, \beta_i\}_{i=1}^2$ of $H_1(\Sigma_2; \zed)$ such that
$\alpha_i$ represents the curve $\system[i]{3}$ in $\system{3}$.  Given $(Y, \eta)\in \Omega\Def^0_2(
\system{3}, \system[2]{3})$, define
\begin{align*}
  w &= \eta(\alpha_1) & x &= \eta(\beta_2) \\
  y &= e^{2 \pi i \eta(\beta_1)/\eta(\alpha_1)} & z &= \left( \int_I \eta\right)^2,
\end{align*}
where $I$ is a saddle connection joining distinct zeros on $(Y, \eta)$ such that $I\cup J(I)$ is homotopic to
the curve $\system[2]{3}$.  These are holomorphic functions on $\Omega\Def^0_2( \system{3},
\system[2]{3})$, and the functions $(w, x)$ are holomorphic local coordinates on $\Omega\teich^0_2(\system{3},
\system[2]{3})$.

\begin{prop}
  \label{prop:coordinates3}
  Given $(X, \omega)\in \Omega\teich^0_2(\system{3}, \system[2]{3})$, the functions $(w, x, y, z)$ give
  a system of holomorphic local coordinates on a sufficiently small neighborhood $U$ of $(X, \omega)$ in $\Omega\Def^0_2(
  \system{3}, \system[2]{3})$.
\end{prop}

The proof is a straightforward combination of the proofs of Propositions~\ref{prop:coordinates1} and
\ref{prop:coordinates2}, so we will omit it as well of the proofs of further similar descriptions of local
coordinates.

\paragraph{Points in $\proj\Omega\moduli^0(\system{3}, \system[2]{3})$.}

Every Abelian differential in $\proj\Omega\moduli^0(\system{3}, \system[2]{3})$ can be obtained by the following
construction: given $r\in\cx/\zed$, let $f_r$ be the Abelian differential obtained from the cylinder
$(\cx/\zed, dz)$ by identifying the ends  of the cylinder to form a polar node and identifying the
points $0$ and $r$ to form a holomorphic node.  Two Abelian differentials $f_{r}$ and $f_{r'}$ are the same if
and only if $r$ and $r'$ are related by the group $G$ generated by the transformations $z\mapsto z+1$ and
$z\mapsto -z$.  This gives an isomorphism of $\proj\Omega\moduli^0(\system{3}, \system[2]{3})$ with
$$(\cx\setminus\zed)/G\isom(\cx\setminus\{0, 1\})/(z\mapsto z^{-1}).$$

In $\proj\Omega\moduli^0(\system{3}, \system[2]{3})$, let
$$\proj\Omega\moduli^0(\system{3}, \system[2]{3})(d) = \bigcup_{q\in\Lambda} f_{q/d},$$
where
$$\Lambda = \{q\in \zed : \gcd(q, d)=1\}.$$
We can regard $\proj\Omega\moduli^0(\system{3}, \system[2]{3})(d)$ as the 
locus of stable Abelian differentials $(X, [\omega])\in \proj\Omega\moduli^0(\system{3}, \system[2]{3})$ such that
-- if we regard $(X, \omega)$ as $(\cx^*, dz/z)$ with two marked points -- the two marked points differ by
exactly $d$-torsion in $\cx^*$.

$\proj\Omega\moduli^0(\system{3}, \system[2]{3})$ contains $N$ points, where
\begin{equation*}
  N =
  \begin{cases}
    1 & \text{if $d=2$};\\
    \frac{1}{2}\phi(d) & \text{if $d>2$.}
  \end{cases}
\end{equation*}

\subsection{The strata $\Omega{\cal M}_2^0(T_{2,0}; 1, 1)$, $\Omega{\cal M}_2^0(T_{2,0}; 2)$, and
  $\Omega{\cal M}_2^0(T_{2,1})$}
\label{subsec:stratum5}

We now give local coordinates around the strata $\Omega\moduli^0(\system{3}; 1, 1)$, $\Omega\moduli^0(\system{3};
2)$, and $\Omega\moduli^0(\system{4})$. In each of these cases, let $\{\alpha_i, \beta_i\}_{i=1}^2$ be a
symplectic basis of $H_1(\Sigma_2; \zed)$ such that the $\alpha_i$ represent the nonseparating curves
of the curve system.  Fix an Abelian differential $(X, \omega)$ in $\Omega\teich^0_2( \system{3}; 1,
1)$, $\Omega\teich^0_2( \system{3}; 2)$, or $\Omega\teich^0_2( \system{4})$, and choose a small neighborhood $U$
of $(X, \omega)$ in $\Omega\Def^0_2( \system{3}; 1, 1)$, $\Omega\Def^0_2( \system{3})$, or $\Omega\Def^0_2(
\system{4})$ respectively.  If $U$ is sufficiently small, then we can continuously choose for each $(Y,
\eta)\in U$ having distinct zeros a saddle connection $I_{(Y, \eta)}$ connecting those zeros.  If $(X,
\omega)$ has two distinct zeros, then we can also continuously choose an orientation for these saddle
connections, but this is not possible if it has a double zero or a separating node.  Given $(Y, \eta)\in U$, define
\begin{align*}
  v &= \eta(\alpha_1) & w &= \eta(\alpha_2) \\
  y &= e^{2 \pi i \eta(\beta_1)/\eta(\alpha_1)} & z &= e^{2 \pi i \eta(\beta_2)/\eta(\alpha_2)} 
\end{align*}
\begin{equation*}
  x =
  \begin{cases}
    \int_I\eta & \text{if $(X, \omega)\in \Omega\teich^0_2( \system{3}; 1, 1)$;}\\
    \left(\int_I\eta\right)^{2/3} & \text{if $(X, \omega)\in \Omega\teich^0_2( \system{3}; 2)$;}\\
    \left(\int_I\eta\right)^{2} & \text{if $(X, \omega)\in \Omega\teich^0_2( \system{4})$.}
  \end{cases}
\end{equation*}
These are holomorphic functions on $U$, and $(v, w, x)$ give a holomorphic system of local coordinates on
$\Omega\teich^0_2( \system{3})$ in the first two cases.  In the last case where $(X,
\omega)\in\Omega\teich^0_2(\system{4})$, the $(v, w, x)$ give a system of local coordinates near $(X, \omega)$
on the subspace of $\Omega\Def^0_2(\system{4})$ where the nonseparating curves of $\system{4}$ remain nodes.

\begin{prop}
  The functions $(v, w, x, y, z)$ give a system of holomorphic local coordinates on a 
  neighborhood $U$ of $(X, \omega)$ in $\Omega\Def^0_2( \system{3}; 1, 1)$, $\Omega\Def^0_2(
  \system{3}; 2)$, or $\Omega\Def^0_2( \system{4})$.
\end{prop}

\paragraph{Curves in $\proj\Omega\barmoduli$.}

We now discuss some curves in $\proj\Omega\barmoduli$ which will play a prominent role in the rest of this
paper.  Given $\lambda\in\cx^*$, let $C_\lambda$ be the locus of points $(X, [\omega])\in\proj\Omega\barmoduli$ having two nonseparating
nodes and possibly a separating node such that the ratio of the residues of $\omega$ at the two nonseparating
nodes of $X$ is $\pm \lambda^{\pm1}$.  Two such curves $C_\lambda$ and $C_{\lambda'}$ coincide if and only
if $\lambda'=\pm\lambda^{\pm 1}$.  Also, let $c_\lambda\in\proj\Omega\barmoduli$ be the Abelian differential
with three nonseparating nodes and residues $1$, $\lambda$, and $\lambda-1$ as in Figure~\ref{fig:clambda}.

\begin{figure}[htbp]
  \centering
  \input{clambda.pstex_t}
  \caption{The point $c_\lambda\in\proj\Omega\barmoduli$.}
  \label{fig:clambda}
\end{figure}

\begin{prop}
  \label{prop:Clambda}
  The curve $C_\lambda\subset\proj\Omega\barmoduli$ is a twice-punctured sphere.  If $\lambda=\pm 1$, then it
  is an orbifold locus of order two; otherwise it contains no orbifold points of $\proj\Omega\barmoduli$.  The
  punctures of $C_\lambda$ are the two points $c_\lambda$ and $c_{\lambda+1}$.  The curve $C_\lambda$ contains
  exactly one point of $\proj \Omega\moduli^0(\system{4})$.
\end{prop}

\begin{figure}[htbp]
  \centering
  \input{threecurves.pstex_t}
  \caption{Three curves determining cusps of $C_\lambda^0$.}
  \label{fig:threecurves}
\end{figure}
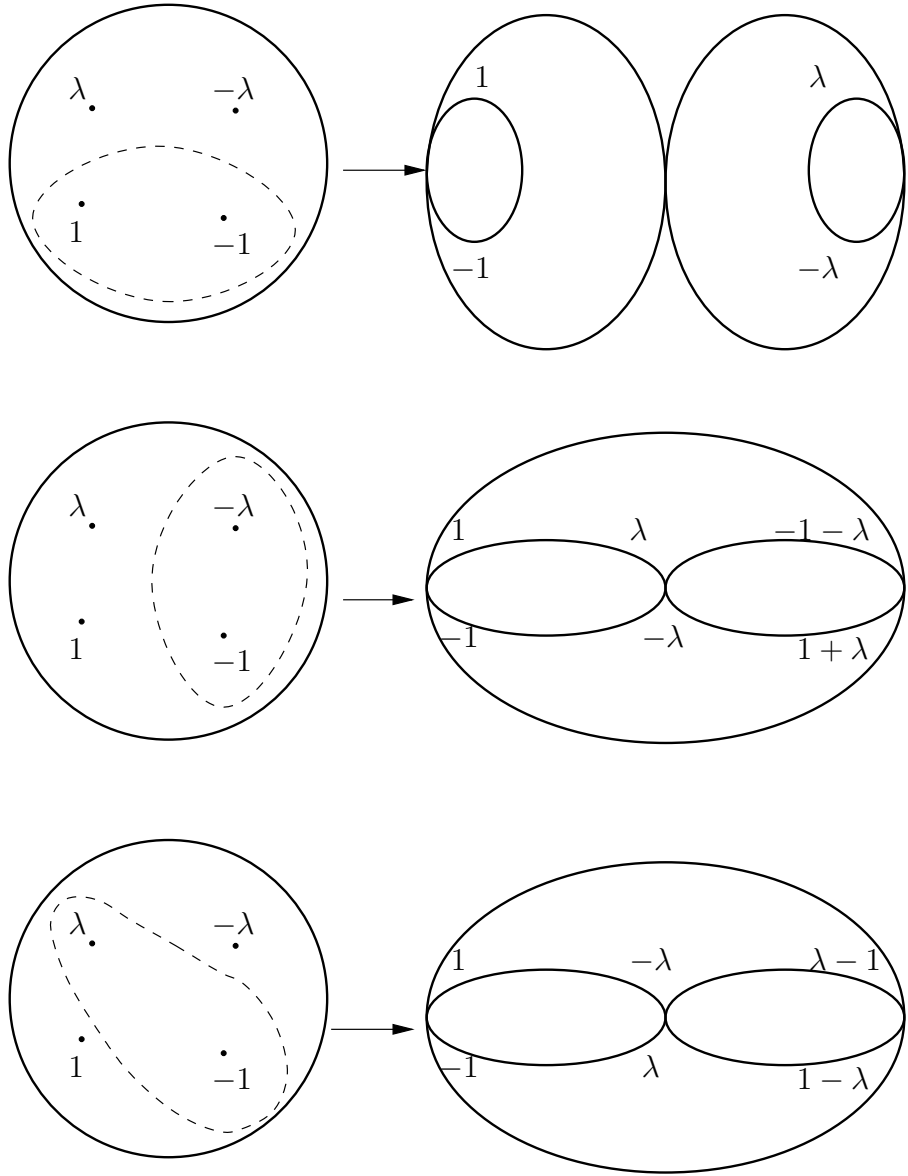

\begin{proof}
  Let $C_\lambda^0\subset C_\lambda$ be the subset of  differentials with no separating
  node.  We claim that $C_\lambda^0$ is a thrice-punctured sphere.  The natural map
  $C_\lambda^0\to\moduli(\system{3})$ is an isomorphism.  Let $\moduli[0,4]^{\rm num}$ be the moduli space of
  four \emph{numbered} points on $\proj^1$.  There is a natural map $\pi\colon\moduli[0,4]^{\rm num}\to\moduli(\system{3})$ which
  sends the point $(p_1, p_2, p_3, p_4)\in\moduli[0,4]^{\rm num}$ to the stable Riemann surface obtained by
  identifying $(p_1, p_2)$ and $(p_3, p_4)$ to nodes.  There are four ways to renumber the points $p_i$ to get
  the same stable Riemann surface, so a fiber of $\pi$ potentially contains four points; however, the
  automorphism group of the $p_i$ contains a Klein four-group whose action on the $p_i$ realizes any
  renumbering of the $p_i$ which preserves the decomposition into the two sets $\{p_1, p_2\}$ and $\{p_3,
  p_4\}$.  Therefore $\pi$ is an isomorphism.  There is thus an isomorphism $\tilde{\pi}\colon\moduli[0,4]^{\rm
    num}\to C_\lambda^0$, defined by sending the point $Q=(p_1, p_2, p_3, p_4)\in\moduli[0,4]^{\rm num}$ to $(X, \omega)$, where
  $X=\pi(Q)$, and $\omega$ is induced by the unique meromorphic Abelian differential $\eta$ on
  $\proj^1$ having a simple pole at each $p_i$ with residue given by
  \begin{align*}
    \Res_{p_1}\eta &= 1 & \Res_{p_2}\eta &= -1 \\
    \Res_{p_3}\eta &= \lambda & \Res_{p_4} \eta &= -\lambda.
  \end{align*}  
  It is well known that $\moduli[0,4]^{\rm num}$ is a thrice-punctured sphere, so $C_\lambda^0$ is as well.
  
  Cusps of $\moduli[0,4]^{\rm num}$ correspond to isotopy classes of simple closed curves in
  $\proj^1\setminus\{p_i\}_{i=1}^4$ up to the action of the modular group $\Mod_{0,4}^{\rm num}$ of
  self-homeomorphisms of $\proj^1$ fixing each of the $p_i$ up to isotopy fixing the $p_i$.  There are three
  such isotopy classes of curve.  These are shown in Figure~\ref{fig:threecurves} with the corresponding
  limiting Abelian differential in $\proj\Omega\barmoduli$.  Thus we see that $C_\lambda$ contains one point
  in $\proj\Omega\moduli^0(\system{4})$, and the two cusps of $C_\lambda$ are as claimed.

  To obtain the statement about orbifold points, note that for each $(X, [\omega])\in C_\lambda$, the
  automorphism group of $X$ is the Klein four group.  If $\lambda=\pm 1$, then all of these automorphisms
  stabilize the projective class of $\omega$, while if $\lambda\neq\pm 1$, only the hyperelliptic involution
  stabilizes $[\omega]$.  
\end{proof}

\begin{prop}
  \label{prop:doublezerosinclambda}
  If $\lambda \neq \pm1$, then $\proj\Omega\moduli^0(\system{3};2)$ meets $C_\lambda$ in exactly one point.
  Otherwise they are disjoint.
\end{prop}

\begin{proof}
  Define an isomorphism $\pi\colon\cx\setminus\{0,1\}\to C_\lambda^0$ by $p(w) = (X_w, \omega_w)$, where $X_w$
  is the stable Riemann surface obtained by identifying the points $\{0,1\}$ and $\{w, \infty\}$ to nodes, and
  $\omega_w$ is induced by the meromorphic Abelian differential,
  \begin{align*}
    \omega_w' &= \left(\frac{1}{z} - \frac{1}{z-1} + \frac{\lambda}{z-w}\right)dz \\
    &= \frac{\lambda z^2 + (-1-\lambda)z+w}{z(z-1)(z-w)}dz,
  \end{align*}
  on $\proj^1$.

  The differential $\omega_w'$ has a double zero if and only if
  $$\lambda z^2 +(-1-\lambda)z + w$$
  has a double zero in $\cx\setminus\{0, 1, w\}$, which happens if and only if the discriminant,
  $$\Delta_w = (1+\lambda)^2 -4 \lambda w,$$
  vanishes for some $w\in\cx\setminus\{0,1\}$.  This happens if and only if $\lambda\neq\pm 1$, in which
  case, there is a unique $w$ such that $\Delta_w=0$.  Thus there is a unique Abelian differential in
  $C_\lambda$ with a double zero if $\lambda\neq\pm 1$, and there is no such Abelian differential if $\lambda=\pm1$.
\end{proof}

For $\lambda\neq 0$ or $\pm 1$, let $w_\lambda\in C_\lambda$ be the unique point representing a stable Abelian
differential with a double zero.  Also, for $\lambda\neq 0$, let $p_\lambda\in C_\lambda$ be the unique point
representing a stable Abelian differential with a nonseparating node.  

\subsection{The stratum $\Omega{\cal M}_2^0(T_{3,0})$}
\label{subsec:stratum6}

We now give local coordinates around points in $\Omega\moduli^0(\system{5})$, the stratum of stable Abelian
differentials with three nonseparating polar nodes.  There is a unique stable Riemann surface
$X\in\moduli(\system{5})$ consisting of two thrice-punctured spheres with the cusps joined together to form
three separating nodes.  The automorphism group of $X$ is $S_3\times\zed/2$, with $S_3$ acting on each
irreducible component of $X$ by permuting the nodes and $\zed/2$ acting by permuting the irreducible
components of $X$.  This means that the stratum $\Omega\moduli^0(\system{5})$ is isomorphic to $\cx^2/(S_3\times\zed/2)$.

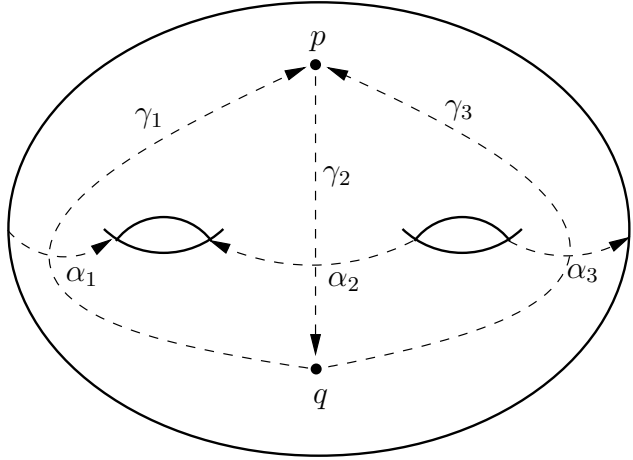
\begin{figure}[htbp]
  \centering
  \input{homologyclasses.pstex_t}
  \caption{Homology classes in $H_1(\Sigma_2, \{p, q\};\zed)$.}
  \label{fig:homologyclasses}
\end{figure}

Choose two points $p, q \in \Sigma_2 \setminus \system{5}$.  Let $\alpha_i\in H_1(\Sigma_2 \setminus \{p,q\}; \zed)$ be a
homology class representing the curve $\system[i]{5}\in\system{5}$.  Let $\gamma_i\in H_1(\Sigma_2,
\{p,q\};\zed)$ be homology classes such that $\alpha_i\cdot\gamma_j=\delta_{ij}$ (see Figure~\ref{fig:homologyclasses}).

Any marked stable Abelian differential $(f, (Y, \eta))\in \Omega\Def^0_2( \system{5})$ has three cylinders
$C_i$ homotopic to the curves $\system[i]{5}$ of $\system{5}$ and two simple zeros, one in each component of
the complement of the cylinders.  We can change $f$ to a marking $f_1$ by an isotopy so that $f_1$ takes each
curve $\system[i]{5}$ into the cylinder $C_i$ and takes the points $p$ and $q$ to zeros of $\eta$.  Two such
markings $f_1$ and $f_2$ which are isotopic to $f$ are isotopic to each other by an isotopy sending for all
time the points $p$ and $q$ to zeros of $\eta$ and sending the curve $\system[i]{5}$ into the cylinder $C_i$.
This means that we can regard the classes $\alpha_i$ as homology classes in $H_1(Y\setminus Z(\eta); \zed)$
and the $\gamma_i$ as homology classes in $H_1(Y, Z(\eta); \zed)$, well-defined up to adding a multiple of
$\alpha_i$ (because the marking $f$ was only defined up to Dehn twist around the curves of $\system{5}$).

Given $(Y, \eta)\in \Omega\Def^0_2( \system{5})$, define
\begin{align*}
  v &= \eta(\alpha_1) & w &= \eta(\alpha_2) \\
  x &= e^{2\pi i \eta(\gamma_1)/\eta(\alpha_1)} &  y &= e^{2\pi i \eta(\gamma_2)/\eta(\alpha_2)} \\
   z &= e^{2\pi i \eta(\gamma_3)/\eta(\alpha_3)}
\end{align*}

These are well-defined holomorphic functions on $\Omega\Def^0_2( \system{5})$, and the coordinates $(v, w)$ define a
map $\Omega\teich^0_2(\system{5})\to\cx^2$ which is biholomorphic onto its image.

\begin{prop}
  The functions $(v, w, x, y, z)$ define a map $\Omega\Def^0_2( \system{5})\to\cx^5$ which is
  biholomorphic onto its image.
\end{prop}

\subsection{The stratum $\Omega{\cal M}_2^0(T_{3,0}, T_{3,0}^3)$}
\label{subsec:stratum7}

Finally, consider the stratum $\Omega\moduli^0(\system{5}, \system[3]{5})$.  An Abelian differential in this
stratum can be regarded as a pair of infinite cylinders with an end of each cylinder identified with an end of
the other cylinder to form two polar nodes, and with a point from each cylinder identified to form a
holomorphic node.  The stratum $\Omega\moduli^0(\system{5}, \system[3]{5})$ is isomorphic to $\cx^*/\pm 1$, with
the isomorphism sending $(X, \omega)$ to the residue of $\omega$ at one of the polar nodes.  The whole stratum
is an orbifold locus of order two in $\Omega\barmoduli$, because each $(X, \omega)$ in this stratum has an
involution which switches the irreducible components of $X$ and preserves $\omega$.  

Choose a symplectic basis $\{\alpha_i, \beta_i\}_{i=1}^2$ of $H_1(\Sigma_2;\zed)$ such that $\alpha_i$
represents the curve $\system[i]{5}$, and let $\alpha_3\in H_1(\Sigma_2; \zed)$ be a class representing
$\system[3]{5}$.  Given a marked stable Abelian differential $(Y, \eta)\in \Omega\Def^0_2(\system{5},
\system[3]{5})$, the class $\beta_i$ determines a class in $H_1(Y;\zed)$ which is well defined up to adding
multiples of $\alpha_i$ and $\alpha_3$.  This means that
$$e^{2\pi i \eta(\beta_i)/\eta(\alpha_i)}$$
is a well defined complex number because $\eta(\alpha_3)=0$.

Given $(Y, \eta)\in \Omega\Def^0_2(\system{5}, \system[3]{5})$, define
\begin{align*}
  w &= \eta(\alpha_1) & x &= e^{2\pi i \eta(\beta_1)/\eta(\alpha_1)} \\
  y &= e^{2\pi i \eta(\beta_2)/\eta(\alpha_2)} & z &= \left(\int_I \eta\right)^2,
\end{align*}
where $I$ is a saddle connection on $(Y, \eta)$ such that $I\cup J(I)$ is homotopic to the curve
$\system[3]{5}$.  The coordinate $w$ defines an isomorphism $\Omega\teich^0_2(\system{5},
\system[3]{5})\to \cx$.

\begin{prop}
  The coordinates $(w, x, y, z)$ give a biholomorphic isomorphism of $\Omega\Def^0_2(\system{5},
  \system[3]{5})$ onto its image in $\cx^4$.
\end{prop}


%% file: curvesystems.pstex_t
\begin{picture}(0,0)%
\includegraphics{curvesystems.pstex}%
\end{picture}%
\setlength{\unitlength}{3947sp}%
\begingroup\makeatletter\ifx\SetFigFont\undefined%
\gdef\SetFigFont#1#2#3#4#5{%
  \reset@font\fontsize{#1}{#2pt}%
  \fontfamily{#3}\fontseries{#4}\fontshape{#5}%
  \selectfont}%
\fi\endgroup%
\begin{picture}(4860,5307)(4471,-8369)
\put(8101,-4711){\makebox(0,0)[lb]{\smash{{\SetFigFont{12}{14.4}{\rmdefault}{\mddefault}{\updefault}{\color[rgb]{0,0,0}$\system{2}$}%
}}}}
\put(5401,-4711){\makebox(0,0)[lb]{\smash{{\SetFigFont{12}{14.4}{\rmdefault}{\mddefault}{\updefault}{\color[rgb]{0,0,0}$\system{1}$}%
}}}}
\put(5401,-6511){\makebox(0,0)[lb]{\smash{{\SetFigFont{12}{14.4}{\rmdefault}{\mddefault}{\updefault}{\color[rgb]{0,0,0}$\system{3}$}%
}}}}
\put(5401,-8311){\makebox(0,0)[lb]{\smash{{\SetFigFont{12}{14.4}{\rmdefault}{\mddefault}{\updefault}{\color[rgb]{0,0,0}$\system{5}$}%
}}}}
\put(8101,-6511){\makebox(0,0)[lb]{\smash{{\SetFigFont{12}{14.4}{\rmdefault}{\mddefault}{\updefault}{\color[rgb]{0,0,0}$\system{4}$}%
}}}}
\put(8101,-8311){\makebox(0,0)[lb]{\smash{{\SetFigFont{12}{14.4}{\rmdefault}{\mddefault}{\updefault}{\color[rgb]{0,0,0}$\system{6}$}%
}}}}
\end{picture}%

%% file: clambda.pstex_t
\begin{picture}(0,0)%
\includegraphics{clambda.pstex}%
\end{picture}%
\setlength{\unitlength}{3947sp}%
\begingroup\makeatletter\ifx\SetFigFont\undefined%
\gdef\SetFigFont#1#2#3#4#5{%
  \reset@font\fontsize{#1}{#2pt}%
  \fontfamily{#3}\fontseries{#4}\fontshape{#5}%
  \selectfont}%
\fi\endgroup%
\begin{picture}(3030,1978)(4186,-4575)
\put(5551,-3961){\makebox(0,0)[lb]{\smash{{\SetFigFont{12}{14.4}{\rmdefault}{\mddefault}{\updefault}{\color[rgb]{0,0,0}$\lambda$}%
}}}}
\put(4351,-3286){\makebox(0,0)[lb]{\smash{{\SetFigFont{12}{14.4}{\rmdefault}{\mddefault}{\updefault}{\color[rgb]{0,0,0}$1$}%
}}}}
\put(5476,-3286){\makebox(0,0)[lb]{\smash{{\SetFigFont{12}{14.4}{\rmdefault}{\mddefault}{\updefault}{\color[rgb]{0,0,0}$-\lambda$}%
}}}}
\put(6601,-3286){\makebox(0,0)[lb]{\smash{{\SetFigFont{12}{14.4}{\rmdefault}{\mddefault}{\updefault}{\color[rgb]{0,0,0}$\lambda-1$}%
}}}}
\put(4276,-3961){\makebox(0,0)[lb]{\smash{{\SetFigFont{12}{14.4}{\rmdefault}{\mddefault}{\updefault}{\color[rgb]{0,0,0}$-1$}%
}}}}
\put(6526,-4036){\makebox(0,0)[lb]{\smash{{\SetFigFont{12}{14.4}{\rmdefault}{\mddefault}{\updefault}{\color[rgb]{0,0,0}$1-\lambda$}%
}}}}
\end{picture}%

%% file: threecurves.pstex_t
\begin{picture}(0,0)%
\includegraphics{threecurves.pstex}%
\end{picture}%
\setlength{\unitlength}{3947sp}%
\begingroup\makeatletter\ifx\SetFigFont\undefined%
\gdef\SetFigFont#1#2#3#4#5{%
  \reset@font\fontsize{#1}{#2pt}%
  \fontfamily{#3}\fontseries{#4}\fontshape{#5}%
  \selectfont}%
\fi\endgroup%
\begin{picture}(5649,7369)(3742,-8475)
\put(8551,-4486){\makebox(0,0)[lb]{\smash{{\SetFigFont{12}{14.4}{\rmdefault}{\mddefault}{\updefault}{\color[rgb]{0,0,0}$-1-\lambda$}%
}}}}
\put(4126,-1711){\makebox(0,0)[lb]{\smash{{\SetFigFont{12}{14.4}{\rmdefault}{\mddefault}{\updefault}{\color[rgb]{0,0,0}$\lambda$}%
}}}}
\put(5026,-1711){\makebox(0,0)[lb]{\smash{{\SetFigFont{12}{14.4}{\rmdefault}{\mddefault}{\updefault}{\color[rgb]{0,0,0}$-\lambda$}%
}}}}
\put(4126,-2611){\makebox(0,0)[lb]{\smash{{\SetFigFont{12}{14.4}{\rmdefault}{\mddefault}{\updefault}{\color[rgb]{0,0,0}$1$}%
}}}}
\put(5026,-2686){\makebox(0,0)[lb]{\smash{{\SetFigFont{12}{14.4}{\rmdefault}{\mddefault}{\updefault}{\color[rgb]{0,0,0}$-1$}%
}}}}
\put(4126,-4336){\makebox(0,0)[lb]{\smash{{\SetFigFont{12}{14.4}{\rmdefault}{\mddefault}{\updefault}{\color[rgb]{0,0,0}$\lambda$}%
}}}}
\put(5026,-4336){\makebox(0,0)[lb]{\smash{{\SetFigFont{12}{14.4}{\rmdefault}{\mddefault}{\updefault}{\color[rgb]{0,0,0}$-\lambda$}%
}}}}
\put(4126,-5236){\makebox(0,0)[lb]{\smash{{\SetFigFont{12}{14.4}{\rmdefault}{\mddefault}{\updefault}{\color[rgb]{0,0,0}$1$}%
}}}}
\put(5026,-5311){\makebox(0,0)[lb]{\smash{{\SetFigFont{12}{14.4}{\rmdefault}{\mddefault}{\updefault}{\color[rgb]{0,0,0}$-1$}%
}}}}
\put(4126,-6961){\makebox(0,0)[lb]{\smash{{\SetFigFont{12}{14.4}{\rmdefault}{\mddefault}{\updefault}{\color[rgb]{0,0,0}$\lambda$}%
}}}}
\put(5026,-6961){\makebox(0,0)[lb]{\smash{{\SetFigFont{12}{14.4}{\rmdefault}{\mddefault}{\updefault}{\color[rgb]{0,0,0}$-\lambda$}%
}}}}
\put(4126,-7861){\makebox(0,0)[lb]{\smash{{\SetFigFont{12}{14.4}{\rmdefault}{\mddefault}{\updefault}{\color[rgb]{0,0,0}$1$}%
}}}}
\put(5026,-7936){\makebox(0,0)[lb]{\smash{{\SetFigFont{12}{14.4}{\rmdefault}{\mddefault}{\updefault}{\color[rgb]{0,0,0}$-1$}%
}}}}
\put(6526,-7186){\makebox(0,0)[lb]{\smash{{\SetFigFont{12}{14.4}{\rmdefault}{\mddefault}{\updefault}{\color[rgb]{0,0,0}$1$}%
}}}}
\put(7651,-7186){\makebox(0,0)[lb]{\smash{{\SetFigFont{12}{14.4}{\rmdefault}{\mddefault}{\updefault}{\color[rgb]{0,0,0}$-\lambda$}%
}}}}
\put(8776,-7186){\makebox(0,0)[lb]{\smash{{\SetFigFont{12}{14.4}{\rmdefault}{\mddefault}{\updefault}{\color[rgb]{0,0,0}$\lambda-1$}%
}}}}
\put(6451,-7861){\makebox(0,0)[lb]{\smash{{\SetFigFont{12}{14.4}{\rmdefault}{\mddefault}{\updefault}{\color[rgb]{0,0,0}$-1$}%
}}}}
\put(8701,-7936){\makebox(0,0)[lb]{\smash{{\SetFigFont{12}{14.4}{\rmdefault}{\mddefault}{\updefault}{\color[rgb]{0,0,0}$1-\lambda$}%
}}}}
\put(7726,-7861){\makebox(0,0)[lb]{\smash{{\SetFigFont{12}{14.4}{\rmdefault}{\mddefault}{\updefault}{\color[rgb]{0,0,0}$\lambda$}%
}}}}
\put(6526,-4486){\makebox(0,0)[lb]{\smash{{\SetFigFont{12}{14.4}{\rmdefault}{\mddefault}{\updefault}{\color[rgb]{0,0,0}$1$}%
}}}}
\put(7651,-4486){\makebox(0,0)[lb]{\smash{{\SetFigFont{12}{14.4}{\rmdefault}{\mddefault}{\updefault}{\color[rgb]{0,0,0}$\lambda$}%
}}}}
\put(6451,-5161){\makebox(0,0)[lb]{\smash{{\SetFigFont{12}{14.4}{\rmdefault}{\mddefault}{\updefault}{\color[rgb]{0,0,0}$-1$}%
}}}}
\put(8701,-5236){\makebox(0,0)[lb]{\smash{{\SetFigFont{12}{14.4}{\rmdefault}{\mddefault}{\updefault}{\color[rgb]{0,0,0}$1+\lambda$}%
}}}}
\put(7726,-5161){\makebox(0,0)[lb]{\smash{{\SetFigFont{12}{14.4}{\rmdefault}{\mddefault}{\updefault}{\color[rgb]{0,0,0}$-\lambda$}%
}}}}
\put(6676,-1636){\makebox(0,0)[lb]{\smash{{\SetFigFont{12}{14.4}{\rmdefault}{\mddefault}{\updefault}{\color[rgb]{0,0,0}$1$}%
}}}}
\put(6526,-2836){\makebox(0,0)[lb]{\smash{{\SetFigFont{12}{14.4}{\rmdefault}{\mddefault}{\updefault}{\color[rgb]{0,0,0}$-1$}%
}}}}
\put(8776,-1636){\makebox(0,0)[lb]{\smash{{\SetFigFont{12}{14.4}{\rmdefault}{\mddefault}{\updefault}{\color[rgb]{0,0,0}$\lambda$}%
}}}}
\put(8701,-2836){\makebox(0,0)[lb]{\smash{{\SetFigFont{12}{14.4}{\rmdefault}{\mddefault}{\updefault}{\color[rgb]{0,0,0}$-\lambda$}%
}}}}
\end{picture}%

%% file: homologyclasses.pstex_t
\begin{picture}(0,0)%
\includegraphics{homologyclasses.pstex}%
\end{picture}%
\setlength{\unitlength}{3947sp}%
\begingroup\makeatletter\ifx\SetFigFont\undefined%
\gdef\SetFigFont#1#2#3#4#5{%
  \reset@font\fontsize{#1}{#2pt}%
  \fontfamily{#3}\fontseries{#4}\fontshape{#5}%
  \selectfont}%
\fi\endgroup%
\begin{picture}(3932,2880)(5236,-6601)
\put(8759,-5467){\makebox(0,0)[lb]{\smash{{\SetFigFont{12}{14.4}{\rmdefault}{\mddefault}{\updefault}{\color[rgb]{0,0,0}$\alpha_3$}%
}}}}
\put(6046,-4487){\makebox(0,0)[lb]{\smash{{\SetFigFont{12}{14.4}{\rmdefault}{\mddefault}{\updefault}{\color[rgb]{0,0,0}$\gamma_1$}%
}}}}
\put(7226,-4854){\makebox(0,0)[lb]{\smash{{\SetFigFont{12}{14.4}{\rmdefault}{\mddefault}{\updefault}{\color[rgb]{0,0,0}$\gamma_2$}%
}}}}
\put(7993,-4441){\makebox(0,0)[lb]{\smash{{\SetFigFont{12}{14.4}{\rmdefault}{\mddefault}{\updefault}{\color[rgb]{0,0,0}$\gamma_3$}%
}}}}
\put(7146,-4007){\makebox(0,0)[lb]{\smash{{\SetFigFont{12}{14.4}{\rmdefault}{\mddefault}{\updefault}{\color[rgb]{0,0,0}$p$}%
}}}}
\put(5612,-5487){\makebox(0,0)[lb]{\smash{{\SetFigFont{12}{14.4}{\rmdefault}{\mddefault}{\updefault}{\color[rgb]{0,0,0}$\alpha_1$}%
}}}}
\put(7259,-5521){\makebox(0,0)[lb]{\smash{{\SetFigFont{12}{14.4}{\rmdefault}{\mddefault}{\updefault}{\color[rgb]{0,0,0}$\alpha_2$}%
}}}}
\put(7166,-6255){\makebox(0,0)[lb]{\smash{{\SetFigFont{12}{14.4}{\rmdefault}{\mddefault}{\updefault}{\color[rgb]{0,0,0}$q$}%
}}}}
\end{picture}%

%% file: limits.tex
\section{Limits of eigenforms}
\label{sec:limitsofeigenforms}

\subsection{Introduction}
\label{subsec:limitsintro}

We now begin the study of the compactification of $\X$.  Recall that we have the locus
$\E\subset\proj\Omega\barmoduli$ of eigenforms for real multiplication by $\ord$, and there is the isomorphism
$j_1\colon \E\to\X\setminus\P$ from Proposition~\ref{prop:jisomorphism}.  The inverse of $j_1$ extends to an
embedding $k_1\colon\X\to\proj\Omega\barmoduli$.  In this section we will study the closure of $\X$ in
$\proj\Omega\barmoduli$, which we denote by $\barX$.  Our goal is to classify exactly which stable Abelian
differentials lie in $\barX\setminus\X$ and to understand the local structure of $\barX$ and its strata around
these points.  More precisely, for each $(X, [\omega])\in\barX$ we will:
\begin{itemize}
\item Choose a neighborhood of $(X, [\omega])$ in $\proj\Omega\barmoduli$ of the form $U/\Aut(X, [\omega])$,
  where $U$ is a neighborhood of $(X, [\omega])$ in an appropriate Dehn space;
\item Give local coordinates on $U$ as in \S\ref{sec:localcoordinates};
\item Give explicit equations for the inverse image  $\pi^{-1}(\barX)$ in $U$.
\end{itemize}
This procedure will be slightly modified when $(X, [\omega])$ has a nonseparating holomorphic node. In that
case, we will only give local coordinates on a hypersurface $V$ in $U$ which we will show contains the inverse
image $\pi^{-1}(\barX)$ in $U$.

We will see that $\barX$ in general has non-normal singularities along curves in $\barX\setminus\X$.  To get a
less singular compactification, we will pass to the normalization $\Y$ of $\barX$, which we will study in
\S\ref{sec:geometric} using the results of this section.

\paragraph{Stratification of $\barX$.}

The stratification of $\Omega\barmoduli$ which we discussed in \S\ref{sec:localcoordinates} induces a
stratification of $\barX$.  Given a stratum $\Omega\moduli^0(S, T; {\bf n})$, let $\X(S, T, {\bf n})$ denote the
stratum,
$$\X(S, T, {\bf n}) = \barX\cap\proj\Omega\moduli^0(S, T, {\bf n}),$$
with the analogous notation if $T$ or
${\bf n}$ is omitted.  Here is a list of all of these strata with a summary of what we will prove about each
one:

\begin{enumerate}
\item $\X(1,1)$ is the stratum of nonsingular eigenforms with two simple zeros.  We called this locus $\E(1, 1)$
  in \S\ref{subsec:action}.  It is an open, dense subset of $\X$.
\item $\X(2)$ is the stratum of nonsingular eigenforms with a double zero, otherwise known as $\W$.
\item $\X(\system{1})$ is the stratum of eigenforms consisting of two genus one differentials joined at a node.
  This is the curve $\P$, which we studied in \S\ref{subsec:hilbertmodular}. 
\item[4-5.] $\X(\system{2})$ is the locus of limiting eigenforms with one nonseparating polar node, containing
  the stratum $\X(\system{2}; 1,1)$ as an open dense set and the stratum $\X(\system{2}; 2)$ as a finite
  subset.  We will show in \S\ref{subsec:emptystrata} that $\X(\system{2})$ is empty unless $D$ is square.  We
  will show in \S\ref{subsec:onenonpolarnode} that $\barX[d^2]$ is  an orbifold around
  $\X[d^2](\system{2})$, and $\X[d^2](\system{2})$ is a suborbifold isomorphic to $\half/\Gamma_1(d)$.  We
  will also show that
  $$\X[d^2](\system{2}) = \proj\Omega\moduli^0(\system{2})(d)$$
  as subsets of $\proj\Omega\barmoduli$, where $\proj\Omega\moduli^0(\system{2}, \system[1]{2})(d)$ is the locus
  of degree $d$ cylinder covering differentials discussed in \S\ref{subsec:stratum2}.  We will show that
  $\X[d^2](\system{2}; 2)$ consists of the intersection points of $\barW[d^2]$ with $\X[d^2](\system{2})$ and
  that these intersections are transverse if $d>3$.
\item[6.] $\X(\system{2},\system[1]{2})$ is the stratum of limiting eigenforms with one nonseparating
  holomorphic node.  We will show in \S\ref{subsec:emptystrata} that this stratum is empty unless $D$ is
  square.  We will show in \S\ref{subsec:onenonholnode} that $\barX[d^2]$ is  an orbifold around
  $\X[d^2](\system{2}, \system[1]{2})$, and $\X[d^2](\system{2}, \system[1]{2})$ is a suborbifold isomorphic
  to $\half/\Gamma_1(d)$.  We will also show that
  $$\X[d^2](\system{2}, \system[1]{2}) = \proj\Omega\moduli^0(\system{2}, \system[1]{2})(d)$$
  as subsets of $\proj\Omega\barmoduli$, where $\proj\Omega\moduli^0(\system{2}, \system[1]{2})(d)$ is the
  $d$-torsion locus introduced in \S\ref{subsec:stratum3}.
\item[7.] $\X(\system{3}, \system[2]{3})$ is the stratum of limiting eigenforms with one nonseparating holomorphic
  node and one nonseparating polar node.  We will show in \S\ref{subsec:emptystrata} that this stratum is
  empty unless $D$ is square.  In \S\ref{subsec:onenonpolandonenonhol} we will show that when $D=d^2$, this stratum is
  the finite set,
  $$\X[d^2](\system{3}, \system[2]{3}) = \proj\Omega\moduli^0(\system{3}, \system[2]{3})(d),$$
  introduced in \S\ref{subsec:stratum4}.
  We will also show that $\barX[d^2]$ is nonsingular at these points and that each of these points is a
  transverse intersection point of the closures of the strata $\X[d^2](\system{2})$ and $\X[d^2](\system{2}, \system[1]{2})$. 
\item[8-10.] The strata $\X(\system{3}; 1, 1)$, $\X(\system{3}; 2)$, and $\X(\system{4})$ together consist of the
  limiting eigenforms with two nonseparating polar nodes, possibly with a separating node.  We will show in
  \S\ref{subsec:twononseppolarnodes} that the union of these strata is the union of the curves $C_{\lambda(P)}$,
  where $P$ is a nondegenerate $\Y$-prototype.  We will see that $\barX$ is in general singular along these
  curves and  can even have several branches passing through them.  We will assign a $\Y$-prototypes to each of
  these branches, and we will give explicit equations for these branches in local coordinates.
\item[11-12.] The strata $\X(\system{5})$ and $\X(\system{5}, \system[3]{5})$ together consist of the limiting
  eigenforms with three nonseparating nodes.  We will show in \S\ref{subsec:threenonsepnodes} that
  $\X(\system{5})$ is the union of the points $c_{\lambda(P)}$ over all nonterminal $\Y$-prototypes $P$, and
  $\X(\system{5}, \system[3]{5})$ is the union of the points $c_{\lambda(P)}$ over all terminal
  $\Y$-prototypes $P$.  We will see that $\barX$ is in general singular at these points.  We will assign
  $\Y$-prototypes to the branches of $\barX$ through each $c_{\lambda(P)}$, and we will give equations for
  these branches in local coordinates.
\item[13-17.] We will show in \S\ref{subsec:emptystrata} that the stratum $\X(\system{6})$ as well as the four
  strata consisting of stable Abelian differentials $(X, \omega)$ where $\omega$ vanishes on some irreducible
  component of $X$ are all empty.
\end{enumerate}

\paragraph{Bundles.}

Over any of the bundles of projective spaces $\proj\Omega\barmoduli$, $\proj\Omega\teich_2$, or
$\proj\Omega\Def_2(S)$, there is a canonical line bundle $\mathcal{O}(-1)$, whose fiber over a projective class of
Abelian differentials, $(X, [\omega])$ is the space of constant multiples of $[\omega]$.  Define a Hermitian
metric on each of these bundles by defining on the fiber over $(X, \omega)$,
$$h(\eta, \eta) = \int_X |\eta|^2.$$
This metric is singular over the $(X, [\omega])$ which have infinite area.
For use in \S\ref{sec:bundles}, we will  give  sections for $\mathcal{O}(-1)$ around
points in $\barX$ and calculate the norms of these sections.

\subsection{Empty strata}
\label{subsec:emptystrata}
   
\begin{prop}
  \label{prop:noonenode}
  If $D$ is not square, and $(X, [\omega])\in\barX\setminus\X$, then $\Jac(X)\isom(\cx^*)^2$, and the period map ${\rm
    Per}_\omega\colon H_1(\Jac(X);\zed)\to\cx$ is injective.
\end{prop}

\begin{proof}
  If $(X, [\omega])\in \barX$, then $(X, \omega)\in\overline{\Omega X}_D$ is a nonzero eigenform for real
  multiplication by $\ord$ by Theorem~\ref{cor:rmclosed}, so $M=H_1(\Jac(X);\zed)$ is a torsion-free
  $\ord$-module.  If $D$ is not square, this implies that the $\zed$-rank of $M$ is even because $M\otimes
  \ratls$ is a $K_D$-vector space.  This means that either $\Jac(X)$ is compact, or
  $\Jac(X)\isom(\cx^*)^2$.  Since for any  $X\in\barX\setminus\X$, the Jacobian $\Jac(X)$ is noncompact, we
  must have $\Jac(X)\isom(\cx^*)^2$.

  Since $\omega$ is an eigenform, ${\rm Per}_\omega\colon M\to\cx$ is $\ord$-linear, so $K=\Ker({\rm
    Per}_\omega)$ is an $\ord$-submodule of $M$.  If $K\neq 0$, then $K$ must be a submodule of $M$ of
  $\zed$-rank two, so $M/K$ would be a finite Abelian group.  If $K\neq M$, then this contradicts the fact that ${\rm
    Per}_\omega$ embeds $M/K$ into $\cx$; if $K=M$, then this contradicts the fact that $\omega$ is nonzero.
\end{proof}

\begin{cor}
  If $D$ is not square, then all of the strata in $\barX\setminus\X$ besides $\X(\system{3}; 1, 1)$,
  $\X(\system{3}; 2)$, $\X(\system{4})$, and $\X(\system{5})$ are empty.
\end{cor}

\begin{prop}
  The stratum $\X[d^2](\system{6})$ as well as the four strata of $\barX[d^2]$ consisting of stable Abelian differentials $(X,
  \omega)$ where $\omega$ vanishes on some irreducible component of $X$ are all empty.
\end{prop}

\begin{proof}
  We know that each $(X, \omega)\in\X[d^2]$ is a degree $d$ branched cover of a genus one Abelian differential
  by Proposition~\ref{prop:ellipticdifferentials}.  For any stable Abelian differential $(X,
  \omega)\in\barX[d^2]$, we can take a limit of these branched covers to get a branched covering of degree at most $d$ from $(X,
  \omega)$ to a genus one  stable Abelian differential $(Y, \eta)$.  If $(X, \omega)\in
  \X[d^2](\system{6})$, then it consists of a one point connected sum of a genus one nonsingular Abelian
  differential (having finite area) with an infinite cylinder (having infinite area).  This is impossible
  because the same $(Y, \eta)$ cannot be finitely covered by both an Abelian differential with infinite area
  and one with finite area.

  The proofs that the other four strata are empty are all similar, so for concreteness we will show that there
  is no stable Abelian differential $(X,\omega)\in\barX[d^2]$ which is the one point connected sum,
  $$(X_1, \omega_1) \# (X_2, 0),$$
  where $X_i$ are elliptic curves with $\omega_1$ nonzero.  Suppose $(X,
  \omega) = \lim (X_n, \omega_n)$ with $(X_n, \omega_n)\in \X[d^2]$.  There is a sequence of separating curves
  $\gamma_n\subset X_n$ which are pinched to the node of $X$ in the limit.  If $n$ is large, then $\gamma_n$ is
  homotopic to a union of short saddle connections $I_1\cup I_2$ by Theorem~\ref{thm:longcylinders}, and we
  can split along these saddle connections to obtain a pair of genus one forms $(E_n^i, \omega^i_n)$. We have
  a degree $d$ branched covering $f_n\colon (X_n, \omega_n)\to (F_n, \nu_n)$ over some genus one form $(F_n,
  \nu_n)$.  Taking a subsequence, $(F_n, \nu_n)$ converges to some genus one form $(F, \nu)$ which is
  nonsingular because otherwise $E_n^1$ would converge to a Riemann surface with a node which would contradict the
  assumption that $X$ has only one node.  Each of the genus one forms $(E_n^i, \omega_n^i)$ is then an
  unbranched cover of $(F_n, \nu_n)$ of degree less than $d$.  Thus we have
  $$\frac{1}{d} \leq \frac{\Area(\omega_n^1)}{\Area(\omega_n^2)} \leq d,$$
  which is a contradiction because exactly one of the sequences $\omega_n^i$ must converge to zero.
\end{proof}

\subsection{The strata $X_D(T_{1,0}; 1, 1)$ and $X_D(T_{1,0}, 2)$}
\label{subsec:onenonpolarnode} 

In this section, we will study the strata $\X(\system{2}; 1, 1)$ and $\X(\system{2}, 2)$, whose union is the
locus $\X(\system{2})$, consisting of Abelian differentials with one polar node.  Since these are empty when
$D$ is not square, we will restrict to the case $D=d^2$.  The goal of the rest of this section is to prove the
following description of these strata.

\begin{theorem}
  \label{thm:X100structure}
  The locus $\X[d^2](\system{2})$ is exactly the locus $\proj\Omega\moduli^0(\system{2})(d)$ of degree
  $d$ cylinder covering differentials discussed in \S\ref{subsec:stratum2}.  The variety $\barX[d^2]$ is
   an orbifold around $\X[d^2](\system{2})$ with $\X[d^2](\system{2})$ a suborbifold isomorphic to
  $\half/\Gamma_1(d)$.  The stratum $\X[d^2](\system{2};2)\subset\X[d^2](\system{2})$ is a finite subset equal
  to
  $$\barW[d^2]\cap\X[d^2](\system{2}).$$
  These intersections are transverse if $d>3$.
\end{theorem}

\paragraph{Local coordinates.}

Let $(X, [\omega])\in \proj\Omega\moduli^0(\system{2})$, and choose some marking $f\colon\Sigma_2\to X$, so that
we can regard $(X, [\omega])$ as a point in $\proj\Omega\teich^0_2( \system{2})$.  Take a symplectic basis
$\{\alpha_i, \beta_i\}$ of $H_1(\Sigma_2; \zed)$ as in \S\ref{subsec:stratum2}.  On a small neighborhood $U$
of $(X, [\omega])$ in $\proj\Omega\Def^0_2( \system{2})$, we can normalize each projective class $(Y,
[\eta])\in U$ so that $\eta(\alpha_1)=1$.  This allows us to identify $U$ with the hypersurface $U_1\subset
\Omega\Def^0_2( \system{2})$ consisting of all $(Y, \eta)$ with $\eta(\alpha_1)=1$.  The restriction of
the local coordinates $(v, w, x, y, z)$ from \S\ref{subsec:stratum2} to $U_1$ then gives a system of local
coordinates $(w, x, y, z)$ on $U$.

Assume that $U$ is small enough that the natural map $U/\Aut(X, [\omega])\to\proj\Omega\barmoduli$ is an
isomorphism onto its image, so that we can regard $U$ as an orbifold coordinate chart around $(X, [\omega])$.
By the definition of $\proj\Omega\Def^0_2(\system{2})$, for each $(Y, \eta)\in U$, the curve $\system[1]{2}$
is represented by a maximal cylinder $C$ on $(Y, \eta)$.  Using Theorem~\ref{thm:longcylinders}, by possibly
shrinking $U$, we can assume that:
\begin{itemize}
\item
  $|\eta(\alpha_2)|<\frac{1}{d+1}\height(C)$
\item
  $|\eta(\beta_2)|<\frac{1}{d+1}\height(C)$
\item
  $\IM\eta(\beta_1)>0$
\end{itemize}

\paragraph{$\X[d^2](\system{2})$ in local coordinates.}

We now turn to the question of when our $(X, [\omega])$ lies in the union of strata $\X[d^2](\system{2})$ and
what these strata look like in coordinates around $(X, [\omega])$.  We will continue to work in the fixed
neighborhood $U$ of  some $(X, [\omega])\in\proj\Omega\moduli^0(\system{2})$ which we chose above.

For any $(Y, \eta)$ in $U$, we have a marking $f\colon\Sigma_2\to Y$, which is defined up to Dehn twist around
$\system{2}$.  Choosing a particular marking, we regard the symplectic basis $\{\alpha_i, \beta_i\}_{i=1}^2$
of $H_1(\Sigma_2; \zed)$ as also a symplectic basis of $H_1(Y; \zed)$.
\begin{lemma}
  \label{lem:matrixofgstar}
  Let $(Y, \eta)\in U\cap\X[d^2]$, and let $g\colon (Y, \eta) \to (E, \nu)$ be a primitive, degree $d$,
  branched cover of a genus one differential $E$.  Then there is a basis $\{a, b\}$ of $H_1(E;\zed)$ and
  integers $p, q$, and $r$ such that $\gcd(d, p, q)=1,$ and
  \begin{align*}
    g_*(\alpha_1)&=da & g_*(\beta_1) &= ra +  b\\
    g_*(\alpha_2) &=p a & g_*(\beta_2) &=q a
  \end{align*}
\end{lemma}

\begin{proof}
  There is some primitive $a\in H_1(E; \zed)$ such that $g_*(\alpha_1)=n a$ for some $n\in\nats$.  Since
  $\eta(\alpha_1)=1$, we must have $\nu(a)>0$.  It follows that we can choose some $b\in H_1(E; \zed)$ such that
  $\IM \nu(b)>0$, and such that $\{a, b\}$ form a basis of $H_1(E; \zed)$.  The cohomology class $b$ is represented by
  some closed geodesic $B$ on $E$.
  
  We claim that
  $$\IM \nu(b)> \frac{1}{d+1}\height(C),$$
  where $C$ is the maximal cylinder on $Y$ homotopic to the single  curve of $\system{2}$.  
  To see this, consider $S=g^{-1}(A)\cap C$.  The set $S$ consists of $r$ parallel closed curves $\{S_i\}_{i=1}^r$ on
  $C$.  Since each point of $B$ has $d$ preimages under $g$, we must have
  $r\leq d$.  Two consecutive curves $S_i$ and $S_{i+1}$ bound a cylinder $C_i$ which has the same height as  the
  cylinder obtained by cutting $E$ along $A$, so $\height(C_i)=\IM\nu(b)$.  There are
  also two more cylinders, $C_0$ and $C_{r}$ bounded by $S_1$ and $S_r$ respectively and the boundary curves
  of $C$.  These cylinders have height less than $\IM\nu(b)$ since they are mapped injectively into $E$ by
  $g$.  Together, these cylinders fill out $C$, so
  $$\height(C)=\sum_{i=0}^{r}\height(C_i)<(r+1)\IM \nu(b)\leq(d+1)\IM\nu(b),$$
  which proves the claim.
  
  It follows from this claim that $g_*(\alpha_2)=p a$ for some $p\in\zed$.  To see this, assume that
  $g_*(\alpha_2)=p a + s b$ with $s\neq 0$, then
  \begin{alignat*}{2}
    \frac{1}{d+1}\height(C)&>|\eta(\alpha_2)| & & \quad\text{(by the definition of $U$)}\\
    & =|\nu(g_*(\alpha_2))|\\
    & \geq|s|\IM\nu(b) & & \quad\text{(because $\nu(a)$ is real)}\\
    &>\frac{1}{d+1}\height(C),
  \end{alignat*}
  a contradiction.  Similarly $g_*(\beta_2)=q a$ for some $q\in\zed$.
  
  The map $g_*\colon H_1(Y; \zed)\to H_1(E; \zed)$ is onto because $g\colon Y\to E$ is primitive. Since we've shown that
  each basis element of $H_1(Y; \zed)$ except $\beta_1$ maps to a multiple of $a$, we must have $g_*(\beta_1)=r a
  \pm b$ for some $r\in\zed$, and because $\IM\eta(\beta_1)>0$ and $\IM \nu(b)>0$, we actually have
  $g_*(\beta_1)=x a + b$.  It follows that
  $$g_*(\alpha_1)\cdot g_*(\beta_1)=(na)\cdot(ra+b)=n,$$
  and similarly $g_*(\alpha_2)\cdot g_*(\beta_2)=0$.
  It follows easily from the fact that $g$ is degree $d$ that
  $$g_*(\alpha_1)\cdot g_*(\beta_1) + g_*(\alpha_2)\cdot g_*(\beta_2)=d.$$
  so $n=d$.  Because $g$ is primitive, we must have $\gcd(d, p, q)=1$.
\end{proof}

Let $\psi\colon U\to\proj\Omega\teich^0_2( \system{2})$ be the unplumbing map defined in
\S\ref{subsec:stratum2}, which sends $(Y, \eta)\in U$ to the Abelian differential obtained by unplumbing the
cylinder on $(Y, \eta)$ homotopic to the curve $\system[1]{2}$.  In our local coordinates, $\psi$ is the
projection $\psi(w, x, y, z)=(w, x, y, 0)$.  In $\proj\Omega\teich^0_2( \system{2})$, let $\proj\Omega\teich^0_2(
\system{2})(d)$ be the inverse image of $\proj\Omega\moduli^0(\system{2})(d)$, the locus of degree $d$ cylinder
covering differentials, under the natural projection.

\begin{theorem}
  \label{thm:ellipticiffcylindrical}
  An Abelian differential $(Y, \eta)\in U\setminus\proj\Omega\teich^0_2( \system{2})$ is in $\X[d^2]$ if
  and only if $\psi(Y, \eta)\in\proj\Omega\teich^0_2(\system{2})(d)$.
\end{theorem}

\begin{proof}
  Let $(Y, \eta)\in U$, and let $C$ be the maximal closed cylinder on $Y$ containing a closed geodesic $A$
  homotopic to the curve of $\system{2}$.  Let $(Z, \nu)$ be the stable Abelian differential obtained by
  cutting $Y$ along $A$ to obtain a surface $(Y', \eta')$ with geodesic boundaries $A_1$ and $A_2$, and then
  gluing in an infinite cylinder to each resulting boundary component of $Y'$.  The surface $Y'$ naturally
  lies in both $Y$ and $Z$, and $(Z, \nu)$ is the unplumbing of $(Y, \eta)$ along $C$.
  
  We need to show that $(Y, \eta)$ is a degree $d$ elliptic differential if and only if $(Z, \nu)$ is a
  degree $d$ cylinder covering differential.
  
  First, assume that $(Z, \nu)$ is a degree $d$ cylinder covering differential.  Let $g\colon (Z, \nu) \to
  (D,\xi)$ be a map to a cylinder $(D, \xi)=(\cx/\zed, dz)$ realizing $Z$ as a degree $d$ cylinder covering
  differential.  The map $g$ sends the horizontal geodesics $A_1$ and $A_2$ to horizontal geodesics $B_1$ and
  $B_2$ on $C$, and by the Open Mapping Theorem, $B_1$ and $B_2$ must be distinct (or else the point $p\in Y'$
  for which $g(p)$ has largest imaginary part would lie in the interior of $Y'$, which contradicts the Open
  Mapping Theorem).  The geodesics $B_1$ and $B_2$ bound a subcylinder $D'$ of $D$, and $g$ maps $Y'$ onto
  $D'$ with degree $d$.  The induced map in homology, $g_*\colon H_1(Y'; \zed)\to H_1(D'; \zed)$ is onto
  because $g_*\colon H_1(Z;\zed)\to H_1(D;\zed)$ is onto by primitivity of $g$.  Now, we can recover $Y$ by
  gluing together $A_1$ and $A_2$ by an appropriate gluing map.  By gluing $D'$ along $B_1$ and $B_2$ in a way
  compatible with the map $g$, we get an elliptic curve $E$ with Abelian differential $\zeta$ and a primitive
  degree $d$ map $h\colon (Y, \eta)\to (E, \zeta),$ realizing $(Y, \eta)$ as a elliptic differential.
  
  Now assume that $(Y, \eta)$ is a degree $d$ elliptic differential.  Let $g\colon (Y, \eta)\to (E, \zeta)$
  be a map realizing $(Y, \eta)$ as an elliptic differential. Take a basis $\{a, b\}$ of $H_1(E;\zed)$ as in Lemma
  \ref{lem:matrixofgstar}. Let $B=g(A)$, a horizontal closed geodesic on $E$.  The geodesics $A$ and $B$ are
  homologous to $\alpha_1$, and $a$ respectively.  Since $g_*(\alpha_1)=d a$ by Lemma \ref{lem:matrixofgstar},
  the restriction $g|_A\colon A\to B$ is degree $d$, so $g^{-1}(B)=A$.

  Cut $E$ along $B$, and call the resulting cylinder with boundary $(D',
  \xi')$.  Let $(D, \xi)$ be the infinite cylinder obtained by gluing half-infinite cylinders to each boundary
  component of $(D', \xi')$.  Since $g^{-1}(B)=A$, there is a natural degree $d$ map $h'\colon (Y',
  \eta')\to (D', \xi')$ and thus a degree $d$ map $h\colon (Z, \nu)\to(D, \xi)$.
  
  To show that $(Z, \nu)$ is a degree $d$ cylinder covering differential, it just remains to show that $h$ is primitive,
  and it is enough to show that $h'_*\colon H_1(Y';\zed)\to H_1(D';\zed)$ is onto.  This follows easily from Lemma
  \ref{lem:matrixofgstar}.  Since the homology classes $\alpha_i$ and $ \beta_2$ can be represented by closed
  curves disjoint from $A$, they naturally define classes on $Y'$.  Then the fact that $\gcd(d, p, q)=1$ from
  Lemma \ref{lem:matrixofgstar}, implies that $h'_*$ is onto.
\end{proof}

\begin{cor}
  \label{cor:X100coordinates}
  We have
  \begin{equation}
    \label{eq:X100equation}
    \X[d^2](\system{2}) = \proj\Omega\moduli^0(\system{2})(d).
  \end{equation}
  Given $(X, [\omega])\in\X[d^2](\system{2})$ contained in a sufficiently small neighborhood
  $U\subset\proj\Omega\moduli^0(\system{2})$ with coordinates on $U$ as defined above, $\barX[d^2]\cap U$ is cut
  out by the equations,
  \begin{align}
    w &= \omega(\alpha_2)/\omega(\alpha_1) \label{eq:X100coordinates}\\
    x &= \omega(\beta_2)/\omega(\alpha_1),\notag
  \end{align}
  and $\X[d^2](\system{2})$ is cut out by the additional equation $z=0$.  If $\omega$ has a double zero, then
  $\barW[d^2]$ is cut out by the equations \eqref{eq:X100coordinates} together with the additional equation
  $y=0$.

  In these coordinates, the foliation $\A[d^2]$ of $\X[d^2]$ has leaves given by $z={\rm const}$.
\end{cor}

\begin{proof}
  It follows immediately from Theorem~\ref{thm:ellipticiffcylindrical} that $(w,x,y,z)\in U$ is in $\X[d^2]$
  if and only if $(w,x,y,0)\in\proj\Omega\teich^0_2( \system{2})(d)$ and $z\neq 0$, and \eqref{eq:X100equation}
  follows by taking the closure.
  
  If $(w, x, y, 0)\in\proj\Omega\teich^0_2( \system{2})(d)$, then $(w, x, y',
  0)\in\proj\Omega\teich^0_2( \system{2})(d)$ for any $y'$ because being a cylinder covering differential
  only depends on absolute periods.  Furthermore, we must have $w, x\in\frac{1}{d}\zed$ by
  Lemma~\ref{lem:matrixofgstar}, so the coordinates $w$ and $x$ are locally constant along
  $\proj\Omega\teich^0_2( \system{2})(d)$.  It follows that if $(X, [\omega])=(w, x, y,
  0)\in\proj\Omega\teich^0_2( \system{2})(d)$, then around this point, $\proj\Omega\teich^0_2(
  \system{2})(d)$ is cut out by the equations \eqref{eq:X100coordinates} together with $z=0$, and
  $\X[d^2]$ is cut out by the equations \eqref{eq:X100coordinates}.

  The equations for $\X[d^2](\system{2})$ and $\barW[d^2]$ are clear from the definitions of the coordinates.
  The leaves of the foliation $\A[d^2]$ are determined by the condition that the absolute periods are locally
  constant along the leaves, which means in these coordinates that $v$, $w$, $x$, and $z$ are all constant.
  Since we're in $U_1$, we have $v=1$, and $w$ and $x$ are locally constant by \eqref{eq:X100coordinates}.
  Thus the foliation is given by $z={\rm const}$.
\end{proof}

Theorem~\ref{thm:X100structure} follows directly from this corollary.  Note that the statement in
Theorem~\ref{thm:X100structure} that the intersections of $\barW[d^2]$ and $\X[d^2](\system{2})$ are
transverse if $d>3$ follows from the coordinates above together with the fact that $\Gamma_1(d)$ is
torsion-free if $d>3$.

\paragraph{A section of $\mathcal{O}(-1)$.}

Recall that we normalized each projective class $(Y, [\eta])\in\proj\Omega\Def^0_2( \system{2})$ so that each
$\eta(\alpha_1)=1$.  We can regard this as defining a section $s$ of the canonical line bundle $\mathcal{O}(-1)$ over
$\proj\Omega\Def^0_2( \system{2})$.

\begin{prop}
  \label{prop:area1}
  The norm of this section $s$ of $\mathcal{O}(-1)$ is given by,
  \begin{equation}
    \label{eq:area1}
    h(s, s) = -\frac{1}{2\pi}\log|z| + w x.
  \end{equation}
\end{prop}

\begin{proof}
  For any genus two Riemann surface $X$ with a symplectic basis $\{\alpha_i, \beta_i\}$ of $H_1(X; \zed)$, it
  is well known that for any $\omega\in\Omega(X)$, 
  \begin{equation}
    \label{eq:area}
    \int_X |\omega|^2 =  \Im(\overline{\omega(\alpha_1)} \eta(\beta_1) + \overline{\omega(\alpha_2)} \eta(\beta_2)).
  \end{equation}
  Equation \eqref{eq:area1} follows from this formula using the fact that $w$ and $x$ are real on $\barX[d^2]$
  in these coordinates.
\end{proof}

\subsection{The stratum $X_D(T_{1,0}, T_{1,0}^1)$}
\label{subsec:onenonholnode}

We now study the stratum $\X(\system{2}, \system[1]{2})$ of Abelian differentials in $\barX$ with one nonseparating
holomorphic node.  Since this stratum is empty if $D$ is not square, we will restrict to the case $D=d^2$.

Recall that the $d$-torsion locus $\proj\Omega\moduli^0(\system{2}, \system[1]{2})(d)$ is the locus of $(X,
[\omega])\in\proj\Omega\moduli^0(\system{2}, \system[1]{2})$ such that, if we regard $(X, [\omega])$ as a
genus one differential with two marked points, then the marked points differ by exactly $d$-torsion on $X$.
The goal of this section is to prove the following theorem.

\begin{theorem}
  \label{thm:X010structure}
  $\barX[d^2]$ is an orbifold around the stratum $\X[d^2](\system{2}, \system[1]{2})$, which is a
  suborbifold equal to
  $$\proj\Omega\moduli^0(\system{2}, \system[1]{2})(d)\isom\half/\Gamma_1(d).$$
\end{theorem}

\paragraph{Vanishing of periods.}

We now prove that for Abelian differentials $(X, [\omega])\in\X[d^2]$
close to some $(X_0, [\omega_0])\in\barX[d^2]$ which has a holomorphic node,  periods along curves of $(X,
[\omega])$ which are close to a holomorphic node of $(X_0, [\omega_0])$ must vanish.

\begin{prop}
  \label{prop:vanishingperiod}
  Let $S$ be a curve system on $\Sigma_2$, let $\pi\colon\proj\Omega\Def_2(S)\to\proj\Omega\barmoduli$ be the
  natural projection, and let $(X_0, [\omega_0])\in \proj\Omega\Def_2( S)\cap\pi^{-1}(\barX[d^2])$ have a
  nonseparating holomorphic node represented by a homology class $\alpha\in H_1(\Sigma_2; \zed)$.  Then for
  $(X, [\omega])\in\proj\Omega\Def_2(S)$ sufficiently close to $(X_0, [\omega_0])$, we must have $\omega(\alpha)=0$.
\end{prop}

\begin{proof}
  Choose a representative $\omega_0$ of the projective class $[\omega_0]$, and suppose there is a sequence
  $(X_n, \omega_n)\to (X_0, \omega_0)$ with $(X_n, \omega_n)\in \pi^{-1}(\X[d^2])$ and $\omega_n(\alpha)\neq 0$
  for all $n$.  There are degree $d$ branched covers $f_n\colon(X_n, \omega_n)\to(E_n, \eta_n)$ over elliptic
  curves.  Taking a subsequence, we can assume that $(E_n, \eta_n)\to (E, \eta)$, with $E$ either a cylinder
  or an elliptic curve and that $f_n$ converges to a branched cover $f\colon (X_0, \omega_0)\to(E,
  \eta)$ of degree at most $d$.  We must then have $\eta\neq 0$.

  The injectivity radius $I(E_n, \eta_n)$ of $(E_n, \eta_n)$ is bounded by the absolute value of any nonzero
  period of $\eta$.  Since $\eta((f_n)_*\alpha) = \omega_n(\alpha)$, and $\omega_n(\alpha)\to 0$, this means
  that $I(E_n, \eta_n)\to 0$.  This is a contradiction because injectivity radius is continuous, and $(E,
  \eta)$ has nonzero injectivity radius.
\end{proof}

\paragraph{Local coordinates.}

Choose some $(X, [\omega])\in\proj\Omega\moduli^0(\system{2}, \system[1]{2})$.  By choosing a marking
$f\colon\Sigma_2\to X$, we can regard $(X, [\omega])$ as a point in
$\proj\Omega\teich^0_2(\system{2},\system[1]{2})$.  Let $\{\alpha_i, \beta_i\}_{i=1}^2$ be a symplectic
basis of $H_1(\Sigma_2;\zed)$ such that $\alpha_1$ represents the curve $\system[1]{2}$, as in \S\ref{subsec:stratum3}. 

Let $U\subset\proj\Omega\Def^0_2( \system{2}, \system[1]{2})$ be a neighborhood of $(X, [\omega])$ small
enough that the natural map $U/\Aut(X, [\omega])\to\proj\Omega\barmoduli$ is an isomorphism onto its image, and
$\eta(\alpha_2)\neq 0$ for every $(Y, [\eta])\in U$.  By normalizing each $(Y, [\eta])\in U$ so that $\eta(\alpha_2)=1$, we can identify $U$
with an open set $U_1$ in the hypersurface in $\Omega\Def^0_2( \system{2}, \system[1]{2})$ of those
Abelian differentials $(Z, \zeta)$ such that $\zeta(\alpha_2)=1$.  The coordinates $(w, x, y, z)$ from
Proposition~\ref{prop:coordinates2} then restrict to local coordinates $(x, y, z)$ on $U$.

\paragraph{$\X[d^2](\system{2}, \system[1]{2})$ in local coordinates.}

The set $U$ is not a neighborhood of $(X, [\omega])$ in $\proj\Omega\Def_2^0(\system{2})$, but rather it is a
neighborhood of $(X, [\omega])$ in the hypersurface $\proj\Omega\Def^0_2( \system{2}, \system[1]{2})$.
Proposition~\ref{prop:vanishingperiod} implies that in a sufficiently small neighborhood $V$ of $(X, [\omega])$ in
$\proj\Omega\Def_2^0(\system{2})$,
$$V\cap\X[d^2]\subset U.$$
This means that we can restrict to this hypersurface and still understand the local structure of $\barX[d^2]$
around $(X, [\omega])$.

Let $\psi\colon U\to \proj\Omega\teich^0_2( \system{2}, \system[1]{2})$ be the projection (defined in the proof of
Proposition~\ref{prop:coordinates2}) which splits each $(X, [\omega])\in U$ along a union of saddle
connections $I\cup J(I)$ which is homotopic to $\alpha_1$.  In our local coordinates on $U$, this is the map
$\phi(x, y, z) = (x, y, 0)$.

Let $\proj\Omega\teich^0_2( \system{2}, \system[1]{2})(d)$ be the inverse image under the natural projection of the $d$-torsion locus
$\proj\Omega\moduli^0(\system{2}, \system[1]{2})(d)$.

\begin{theorem}
  \label{thm:splittingpreservestoruscover}
  A $(Y,\eta)\in U\setminus \proj\Omega\teich^0_2( \system{2}, \system[1]{2})$ is in $\X[d^2]$ if and
  only if $\psi(Y, \eta)\in \proj\Omega\teich^0_2(\system{2}, \system[1]{2})(d)$.
\end{theorem}

\begin{proof}
  Let $(Y, \eta)\in U$, with saddle connections $I_1$ and $I_2$, which we split along to form an elliptic
  curve $(E, \nu)$.  Let $p$ and $q$ be the two zeros of $\eta$.  Let $J_1$ and $J_2$ be the two segments
  which are the images of the saddle connections $I_i$ in $E$.  The $J_i$ start at points $p_i$ and end at
  points $q_i$ such that upon taking the connected sum along $J_1$ and $J_2$, the $p_i$ are identified to form
  the zero $p$ and the $q_i$ are identified to form the zero $q$ (see Figure~\ref{fig:splitting}).

  \begin{figure}[htbp]
    \centering
    \input{splitting.pstex_t}
    \caption{Split along $I_1\cup I_2$ and then reglue.}
    \label{fig:splitting}
  \end{figure}
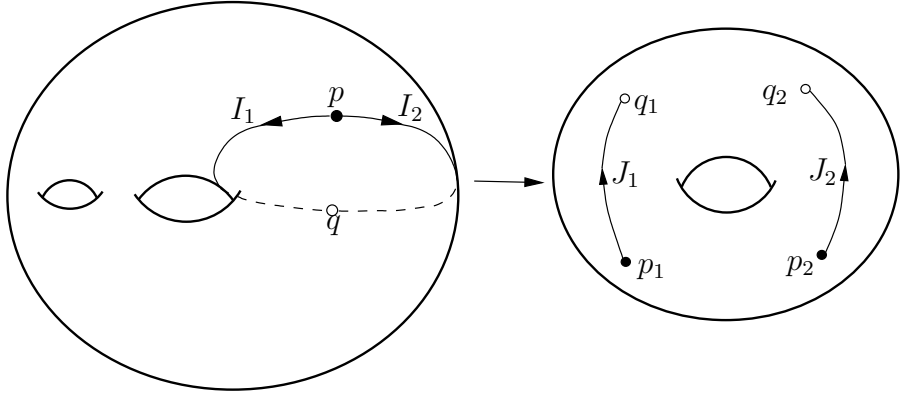

  We need to show that $(Y, \eta)$ is a degree $d$ elliptic differential if and only if the points $p_i$
  differ by exactly $d$-torsion in the group law on $E$.

  First, suppose that $(Y, \eta)$ is a degree $d$ elliptic differential.  Then $(Y, \eta)$ is branched
  over an elliptic curve by some  $g\colon(Y, \eta)\to(F, \xi)$. The saddle connections $I_i$ must have as
  their image the same segment $K$ in $F$ because they have the same direction and length, and they start at
  the same zero $p$.  It follows that when we cut $Y$ along $I_1$ and $I_2$, points which are then glued
  together to form $E$ map to the same point of $F$ under $g$.  Thus $g$ defines an isogeny $g'\colon E\to F$
  of the same degree, sending $p_1$ and $p_2$ to the same point.  This implies that $p_1$ and $p_2$ differ by
  $d$-torsion on $E$.
  
  Conversely, if $p$ and $q$ differ by $d$-torsion, then there is a degree $d$-isogeny $g\colon (E, \nu)\to(F,
  \xi)$; it must send points of $J_1$ and $J_2$ which are glued together to form $Y$ to the same point of $Y$,
  so we get a branched cover $(Y, \eta)\to (F, \xi)$ which realizes $\eta$ as a degree $d$ elliptic
  differential.
\end{proof}

\begin{cor}
  \label{cor:oneholnodecoords}
  The stratum $\X[d^2](\system{2}, \system[1]{2})$ is exactly $\proj\Omega\moduli^0(\system{2},
  \system[1]{2})(d)$.  If we have $(X, [\omega])\in \X[d^2](\system{2}, \system[1]{2})$, then in the $(x, y,
  z)$ coordinates defined above on the neighborhood $U$ of $(X, [\omega])$ in $\proj\Omega\Def^0_2(
  \system{2}, \system[1]{2})$, the variety $\barX[d^2]$ is cut out by the equation,
  \begin{equation}
    \label{eq:dtorsion}
    dy=ax+b,
  \end{equation}
  for some relatively prime integers $a$ and $b$.

  In these coordinates, the stratum $\X[d^2](\system{2}, \system[1]{2})$ is cut out by the additional equation
  $z=0$, and leaves of the foliation $\A[d^2]$ are determined by either of the equivalent equations, $x={\rm
    const}$ or $y={\rm const}$. 
\end{cor}

\begin{proof}
  It follows immediately from Theorem~\ref{thm:splittingpreservestoruscover} that $(x,y,z)\in U$ is in $\X[d^2]$
  if and only if $(x,y,0)\in\proj\Omega\teich^0_2( \system{2}, \system[1]{2})(d)$  and $z\neq 0$.  Thus
  $$U\cap\barX[d^2]=\{(x, y, z)\in U : (x, y, 0)\in\proj\Omega\teich^0_2( \system{2}, \system[1]{2})(d)\},$$
  and
  $\X[d^2](\system{2}, \system[1]{2}) = \proj\Omega\moduli^0(\system{2}, \system[1]{2})(d)$ as claimed.
  
  It remains to show that if $(X, [\omega])\in\proj\Omega\teich^0_2( \system{2}, \system[1]{2})(d)$, then
  near $(X, [\omega])$, the locus $\proj\Omega\teich^0_2( \system{2}, \system[1]{2})(d)$ is cut out by
  the equation \eqref{eq:dtorsion}.  This is because for some $(Y, \eta)\in \proj\Omega\teich^0_2(
  \system{2}, \system[1]{2})$, a relative period joining the marked points is given by $y=\eta(\beta_1)$, so
  the marked points differ by $d$-torsion if and only if $d\eta(\beta_1) = a \eta(\beta_2) + b
  \eta(\alpha_2)$ for some relatively prime integers $a$ and $b$, which is equivalent to \eqref{eq:dtorsion}.
  
  The equation for $\X[d^2](\system{2}, \system[1]{2})$ is clear from the definition of the coordinates, and
  the equations for the leaves of $\A[d^2]$ hold because either of these equations together with
  \eqref{eq:dtorsion} implies that the absolute periods are all constant.
\end{proof} 

Theorem~\ref{thm:X010structure} follows immediately from this corollary.

\subsection{The stratum $X_D(T_{2,0}, T_{2,0}^2)$}
\label{subsec:onenonpolandonenonhol}

We now study the stratum $\X(\system{3}, \system[2]{3})$ consisting of limits of eigenforms with one
nonseparating polar node and one nonseparating holomorphic node.  The arguments in this section are
straightforward combinations of those in the previous two sections, so we will omit some of the proofs.  This
stratum is empty when $D$ is not square, so we will restrict to the case $D=d^2$.

Recall that in \S\ref{subsec:stratum4} we defined the finite set $\proj\Omega\moduli^0(\system{3},
\system[2]{3})(d)$, consisting of the stable Abelian differentials $f_{q/d}\in \proj\Omega\moduli^0(\system{3},
\system[2]{3})$, where $q$ is an integer relatively prime to $d$.  The following is the main theorem of this
section.

\begin{theorem}
  \label{thm:X110classification}
  The stratum $\X[d^2](\system{3}, \system[2]{3})$ is the locus $\proj\Omega\moduli^0(\system{3}, \system[2]{3})(d)$.
  The variety $\barX[d^2]$ is nonsingular at these points, and each 
  $f_{q/d}\in\X[d^2](\system{3}, \system[2]{3})$ is a transverse intersection point of the closures of the strata
  $\X[d^2](\system{2})$ and $\X[d^2](\system{2}, \system[1]{2})$.
\end{theorem}

\paragraph{Local coordinates.}

Choose some $(X, [\omega])\in\proj\Omega\moduli^0(\system{3}, \system[2]{3})$.  By choosing a marking
$f\colon\Sigma_2\to X$, we can regard $(X, [\omega])$ as a point in
$\proj\Omega\teich^0_2(\system{3},\system[2]{3})$.  Let $\{\alpha_i, \beta_i\}_{i=1}^2$ be a symplectic
basis of $H_1(\Sigma_2;\zed)$ such that $\alpha_i$ represents the curve $\system[i]{3}$ in $\system{3}$ as in
\S\ref{subsec:stratum4}.

Let $U\subset\proj\Omega\Def^0_2( \system{3}, \system[2]{3})$ be a neighborhood of $(X, [\omega])$ small
enough that the natural map $U/\Aut(X, [\omega])\to\proj\Omega\barmoduli$ is an isomorphism onto its image.  By
normalizing each $(Y, [\eta])\in U$ so that $\eta(\alpha_1)=1$, we can identify $U$ with an open set $U_1$ in
the hypersurface in $\Omega\Def^0_2( \system{3}, \system[2]{3})$ consisting of those Abelian differentials $(Z,
\zeta)$ such that $\zeta(\alpha_1)=1$.  The coordinates $(w, x, y, z)$ from
Proposition~\ref{prop:coordinates3} then restrict to coordinates $(x, y, z)$ on $U$.

\paragraph{$\X[d^2](\system{3}, \system[2]{3})$ in local coordinates.}

Let $\psi\colon U\to\proj\Omega\teich^0_2(\system{3}, \system[2]{3})$ be the projection which sends an
Abelian differential $(Y, \eta)$ to the one obtained by unplumbing the cylinder on $(Y, \eta)$ homotopic to
the curve $\system[1]{3}$ and splitting along the pair of saddle connections $I \cup J(I)$ homotopic to
$\system[2]{2}$.  In local coordinates, $\psi(x, y, z) = (x, 0, 0)$.  Let $\proj\Omega\teich^0_2(
\system{3}, \system[2]{3})(d)$ be the inverse image under the natural projection to
$\proj\Omega\moduli^0(\system{3}, \system[2]{3})(d)$.

\begin{theorem}
  \label{thm:X110projection}
  If the neighborhood $U$ of $(X, [\omega])$ is sufficiently small,
  then a nonsingular $(Y, \eta)\in U$ is in $\X[d^2]$ if and only if $\phi(Y,
  \eta)\in\proj\Omega\teich^0_2(\system{3}, \system[2]{3})(d)$.
\end{theorem}

The proof of this theorem is a straightforward combination of the proofs of
Theorems~\ref{thm:ellipticiffcylindrical} and \ref{thm:splittingpreservestoruscover}, so we will omit the
proof.  The idea is that $f_{q/d}$ is a primitive degree $d$ branched cover of a cylinder, and given a $(Y, \eta)$
such that $\psi(Y, \eta)=f_{q/d}$, we can use this cover to exhibit $(Y, \eta)$ as a degree $d$ torus
cover.

\begin{cor}
  \label{cor:oneholonepolnodecoords}
  The stratum $\X[d^2](\system{3}, \system[2]{3})$ is equal to $\proj\Omega\moduli^0(\system{3},
  \system[2]{3})(d)$.  If $(X, [\omega])\in \X[d^2](\system{3}, \system[2]{3})$, then in a neighborhood $V$ of
  $(X, [\omega])$ in $\proj\Omega\Def^0_2(\system{3})$,
  $$V\cap \barX[d^2]\subset \proj\Omega\Def^0_2(\system{3}, \system[2]{3}).$$
  In the coordinates
  $(x, y, z)$ on a neighborhood $U$ of $(X, [\omega])$ in $\proj\Omega\Def^0_2( \system{3},
  \system[2]{3})$, the variety $\barX[d^2]$ is cut out by the equation,
  \begin{equation}
    \label{eq:qoverd}
    x = \frac{q}{d},
  \end{equation}
  for some integer $q$ relatively prime to $d$.
  
  The intersection of the closure of $\X[d^2](\system{2})$ with $U$ is cut out by the
  additional equation $y=0$, and the intersection of the closure of $\X[d^2](\system{2},
  \system[1]{2})$ with $U$ is cut out by the additional equation $z=0$.  The leaves of the foliation $\A[d^2]$ of
  $\X[d^2]$ are given by $y= {\rm const}$.
\end{cor}

\begin{proof}
  Proposition~\ref{prop:vanishingperiod} implies directly that in a neighborhood $V$ of
  $(X, [\omega])$ in $\proj\Omega\Def^0_2(\system{3})$,
  $$V\cap \barX[d^2]\subset \proj\Omega\Def^0_2(\system{3}, \system[2]{3}).$$
  It follows immediately from Theorem~\ref{thm:X110projection} that $(x,y,z)\in U$ is in $\X[d^2]$
  if and only if $(x,0,0)\in\proj\Omega\teich^0_2(\system{3}, \system[2]{3})(d)$ and $y, z\neq 0$.  Thus
  $$U\cap\barX[d^2]=\{(x, y, z)\in U : (x, 0, 0)\in\proj\Omega\teich^0_2(\system{3}, \system[2]{3})(d)\},$$
  and
  $\X[d^2](\system{3}, \system[2]{3}) = \proj\Omega\moduli^0(\system{3}, \system[2]{3})(d)$ as claimed.
  
  Suppose $(X, \omega)\in\proj\Omega\teich^0_2(\system{3}, \system[2]{3})(d)$.  The coordinate
  $x=\omega(\beta_2)$ measures a relative period on $(X, \omega)$ joining the two marked points which are
  identified to form a holomorphic node.  Since $\omega(\alpha_1)=1$, we have $x=\omega(\beta_2)=q/d$ for
  some integer $q$ relatively prime to $d$ if and only if the marked points differ by $d$-torsion.  Thus $\barX[d^2]$
  is cut out by \eqref{eq:qoverd} as claimed.
  
  The intersection of the closure of $\X[d^2](\system{2})$ with $U$ is cut out by the
  additional equation $y=0$ because $y(Y, \eta)=0$ exactly when $(Y, \eta)$ has a polar node, and the
  intersection of the closure of $\X[d^2](\system{2}, \system[1]{2})$ with $U$ is cut out by the
  additional equation $z=0$ because $z(Y, \eta)=0$ exactly when $(Y, \eta)$ has a holomorphic node.

  The equation $y={\rm const}$ defines the foliation $\A[d^2]$ because this makes the period along $\beta_1$
  constant, and we have already seen that the other periods are locally constant along $\X[d^2]$ with our
  normalization for $\eta$.
\end{proof}

Theorem~\ref{thm:X110classification} follows directly from this corollary. 

\paragraph{A section of $\mathcal{O}(-1)$.}

Recall that we normalized each projective class $(Y, [\eta])\in\proj\Omega\Def^0_2( \system{3},
\system[2]{3})$ so that each $\eta(\alpha_1)=1$.  We can regard this as defining a section $s$ of the
canonical line bundle $\mathcal{O}(-1)$ over $\proj\Omega\Def^0_2( \system{3}, \system[2]{3})$.

The following follows directly from \eqref{eq:area} together with the definition of our coordinates.
\begin{prop}
  \label{prop:area2}
  The norm of this section $s$ of $\mathcal{O}(-1)$ is given by,
  \begin{equation}
    \label{eq:area2}
    h(s, s) = -\frac{1}{2\pi}\log|z|.
  \end{equation}
\end{prop}

\subsection{The strata $X_D(T_{2,0}; 1,1 )$, $X_D(T_{2,0};2)$, and $X_D(T_{2,1})$}
\label{subsec:twononseppolarnodes}

We now turn to the strata $\X(\system{3}; 1,1 )$, $\X(\system{3};2)$, and $\X(\system{4})$ which
together consist of those limiting eigenforms which have exactly two nonseparating polar nodes, possibly with 
a separating node. 

Recall the definitions from \S\ref{subsec:stratum5} of the curves $C_\lambda$ and the points $p_\lambda$ and $w_\lambda$ in
$\proj\Omega\barmoduli$.  The following is the classification of points in the strata $\X(\system{3}; 1,1 )$,
$\X(\system{3};2)$, and $\X(\system{4})$, which we will prove over the course of this subsection.

\begin{theorem}
  \label{thm:X200classification}
  The union of the three strata $\X(\system{3}; 1,1 )$, $\X(\system{3};2)$, and $\X(\system{4})$ is equal to the
  union,
  \begin{equation}
    \label{eq:firstunion}
    \bigcup_\lambda C_\lambda,
  \end{equation}
  where the union is over all $\lambda$ such that $\lambda=\lambda(P)$ for some nondegenerate $\Y$-prototype $P$.
  The stratum $\X(\system{3};2)$ is the finite set,
  $$\bigcup_{\lambda}w_\lambda,$$
  where the union is over all $\lambda$ such that $\lambda=\lambda(P)$ for
  some  $\W$-prototype $P$.  The stratum $\X(\system{4})$ is the finite set,
  $$\bigcup_{\lambda} p_\lambda,$$
  where the union is over all $\lambda$ such that $\lambda=\lambda(P)$ for some $\P$-prototype $P$.
  The stratum $\X(\system{3}; 1, 1)$ is the complement of these finite sets in \eqref{eq:firstunion}.
\end{theorem}

We will also see that in general $\barX$ is singular along the curves $C_\lambda$, and we will give equations
for $\barX$ in local coordinates around $C_\lambda$.

\paragraph{Local coordinates.}

Choose an $(X, [\omega])$ in $\proj\Omega\moduli^0(\system{3}; 1,1)$, $\proj\Omega\moduli^0(\system{3}; 2)$ or
$\proj\Omega\moduli^0(\system{4})$.  By choosing a marking $f\colon\Sigma_2\to X$, we can regard $(X, [\omega])$
as a point in $\proj\Omega\teich^0_2(\system{3}; 1, 1)$, $\proj\Omega\teich^0_2(\system{3}; 2)$, or
$\proj\Omega\teich^0_2(\system{4})$.

Our first goal is to show that if $(X, [\omega])\in\barX$, then $(X, [\omega])$ lies in the union
\eqref{eq:firstunion}.  If $(X, [\omega])\in\barX$, then $\omega$ is an eigenform for real multiplication of
$\ord$ on $\Jac(X)$ by Theorem~\ref{cor:rmclosed}, so assume that it is such an eigenform.  It follows that
the ratio of the residues of the two poles of $\omega$ is real.  Thus we can choose a representative
$\omega$ of the projective class $[\omega]$ and choose homology classes $\alpha_1$ and $\alpha_2$ in
$H_1(\Sigma_2, \zed)$ such that the following properties hold:
\begin{itemize}
\item The $\alpha_i$ represent the two nonseparating nodes of $X$.
\item $\omega(\alpha_1)=1$.
\item $\omega(\alpha_2)=\lambda$ with $\lambda\geq1$.
\end{itemize}
If $\lambda>1$, then this choice of the $\alpha_i$ and $\omega$ is uniquely determined up to the hyperelliptic
involution by these properties.  If $\lambda=1$, then there are two possible choices for the $\alpha_i$ and
$\omega$ up to the hyperelliptic involution, and we choose one arbitrarily (the other one is then obtained by
swapping the $\alpha_i$).

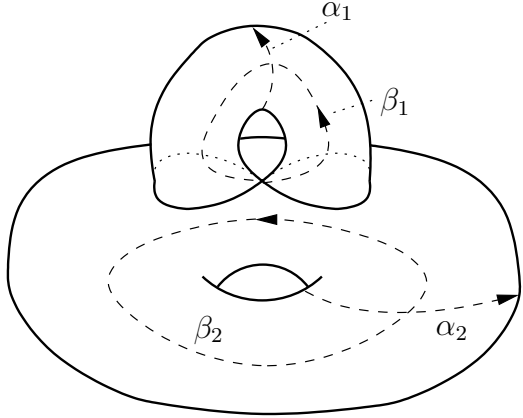
\begin{figure}[htbp]
  \centering
  \input{symplecticbasis.pstex_t}
  \caption{Symplectic basis for $H_1(\Sigma_2;\zed)$.}
  \label{fig:symplecticbasis}
\end{figure}

Extend $\{\alpha_1, \alpha_2\}$ to a basis $\{\alpha_i, \beta_i\}_{i=1}^2$ of $H_1(\Sigma_2; \zed)$ as follows:
\begin{itemize}
\item If $(X, \omega)\in\proj\Omega\moduli^0(\system{3};1,1)$, let $\beta_i$ be an arbitrary
  pair which is dual to the $\alpha_i$ with respect to the intersection pairing.
\item If $(X, \omega)\in\proj\Omega\moduli^0(\system{4})$, let $\beta_i$ be a pair dual to the $\alpha_i$ such
  that each $\beta_i$ is represented by a simple closed curve disjoint from the separating curve
  $\system[3]{4}$.  This determines each $\beta_i$ up to adding a multiple of $\alpha_i$.
\item If $(X, \omega)\in\proj\Omega\moduli^0(\system{3};2)$, then the horizontal foliation of $(X, \omega)$ consists of
  two horizontal cylinders separated by a ``figure-eight''; let $F\subset \Sigma_2$ be the inverse image of
  this figure-eight with $p$ the singular point of $F$.  In $H_1(\Sigma_2; \zed)$, let the $\beta_i$ be a pair which
  is dual to the $\alpha_i$ such that $\beta_1$ is represented by a simple closed curve which passes through
  $p$ and is disjoint from $F\setminus p$, and $\beta_2$ is represented by a simple closed curve which is
  disjoint from $F$  (see Figure~\ref{fig:symplecticbasis}).  These conditions
  determine each  $\beta_i$ uniquely up to adding a multiple of $\alpha_i$.  
\end{itemize}

Let $U$ be a small neighborhood of $(X, \omega)$ in $\proj\Omega\Def^0_2(\system{3})$ or
$\proj\Omega\Def^0_2(\system{4})$ such that the natural map $U/\Aut(X, \omega)\to\proj\Omega\barmoduli$ is an
isomorphism onto its image.  Also choose $U$ small enough so that the conclusion of
Theorem~\ref{cor:rmnearboundary} holds.  That is, if $(Y, \eta)\in U$ is an eigenform for
real multiplication by $\ord$, then the induced action of $\ord$ on $H_1(\Sigma_2;\zed)$ preserves the subgroup
$S$ spanned by the $\alpha_i$, giving $S$ the structure of an $\ord$-module.

By normalizing each projective class $(Y, [\eta])\in U$ so that $\eta(\alpha_1)=1$, we can identify $U$ with
an open set $U'$ in the hypersurface in $\Omega\Def^0_2( \system{3})$ or $\Omega\Def^0_2( \system{4})$ of
those Abelian differentials $(Z, \zeta)$ such that $\zeta(\alpha_1)=1$.  For any $\mu\in\reals$, let
$$U_\mu = \{(Y, \eta)\in U : \eta(\alpha_1)=1,\: \text{and }\eta(\alpha_2)=\mu\},$$
a hypersurface in $U$.  We
then have $(X, \omega)\in U_\lambda$.

\begin{prop}
  \label{prop:X200classification}
  If the neighborhood $U$ of $(X, \omega)$ is taken sufficiently small, then for each $(Y, \eta)\in
  U\cap\X$, the $\ord$-submodule $S$ of $H_1(\Sigma_2;\zed)$ is a quasi-invertible $\ord$-module with $\{\alpha_1,
  \alpha_2\}$ an admissible basis; furthermore, $U\cap\X\subset U_\lambda$.

  It follows that $\lambda = \lambda(P)$ for some nondegenerate $\Y$-prototype $P$.
\end{prop}

\begin{proof}
  For an Abelian differential $(Y, \eta)\in U$ which is sufficiently close to $(X, \omega)$, the classes
  $\alpha_i$ will be homologous to core curves of very tall horizontal cylinders on $Y$ by Theorem
  \ref{thm:longcylinders}.  Since $\alpha_i\cdot\beta_i=1$, a curve representing $\beta_i$ passes vertically
  through the cylinder, so the cylinder makes a large positive contribution to $\Im \eta(\beta_i)$.  Thus, if
  $(Y, \eta)\in U$ with $U$ sufficiently small, we must have $\Im \eta(\beta_i)>0$.

  The intersection pairing is an unimodular pairing between the $\ord$-modules $S$ and
  $H_1(\Sigma_2; \zed)/S$, and the $\alpha_i$ and $\beta_i$ represent bases of these modules dual with respect
  to this pairing.  The period maps associated to $\eta$ and $\Im \eta$ are nonzero $\iota_1$-linear maps of
  these modules to $\reals$ which send the $\alpha_i$ and the $\beta_i$ respectively to positive reals.  It then follows
  from Theorem \ref{thm:criterionforadmissible} that the $\alpha_i$ form an admissible basis of $S$. 

  By Proposition~\ref{prop:admissibleclassification}, for any $(Y, \eta)\in U\cap\X$ there are integers $a$, $b$
  and $c$ as in the conclusion of that proposition with $a\eta(\alpha_2)^2 + b\eta(\alpha_2)+c=0$.  Since
  there are only finitely many such integers, and $\eta(\alpha_2)\to\lambda$ as $(Y, \eta)\to(X,
  \omega)$, we must have $\eta(\alpha_2) = \lambda$ for $(Y, \eta)$ sufficiently close to $(X,
  \omega)$.

  Finally, $\lambda = \lambda(P)$ for any $\Y$-prototype $(a, b, c, \bar{q})$ with $a, b$, and $c$ as above.
\end{proof}

\begin{cor}
  \label{cor:onlylimits}
  The strata $\X(\system{3}; 1,1 )$, $\X(\system{3};2)$, and $\X(\system{4})$ are contained in the union
  \eqref{eq:firstunion}.
\end{cor}

In the coordinates $(v,w,x,y,z)$ on $\Omega\Def^0_2( \system{3})$ or $\Omega\Def^0_2( \system{4})$ from
\S\ref{subsec:stratum5}, the subspace $U_\lambda$ is cut out by the equations,
\begin{equation*}
  v =1 \text{ and } w = \lambda,
\end{equation*}
so the $(x, y, z)$ restrict to give local coordinates on $U_\lambda$, which send $(Y, \eta)\in U$ to
\begin{align*}
   x &=
  \begin{cases}
    \int_I\eta & \text{if $(X, \omega)\in \Omega\teich^0_2( \system{3}; 1, 1)$;}\\
    \left(\int_I\eta\right)^{2/3} & \text{if $(X, \omega)\in \Omega\teich^0_2( \system{3}; 2)$;}\\
    \left(\int_I\eta\right)^{2} & \text{if $(X, \omega)\in \Omega\teich^0_2( \system{4})$,}
  \end{cases}\\
   y &= e^{2 \pi i \eta(\beta_1)} \\
   z &= e^{2 \pi i \eta(\beta_2)/\lambda} ,
\end{align*}
where $I$ is a saddle connection joining distinct zeros as in \S\ref{subsec:stratum5}.

\paragraph{Prototypes.}

Now assume that our $(X, [\omega])\in C_\lambda$, where $\lambda=\lambda(P)$ for some nondegenerate
$\Y$-prototype $P$ with a small neighborhood $U$ of $(X, [\omega ])$ in the Dehn space chosen as above.
We will assign a prototype $P(Y, \eta)$ to each Abelian differential
$(Y, \eta)\in U  \cap \X$.  For now, fix such an $(Y, \eta)$.

Recall that we identify two terminal prototypes  if they differ by the involution \eqref{eq:involutionone}.
We temporarily don't want to make this identification, so we will call a \emph{fine prototype} a prototype
which is defined exactly as in \S\ref{subsec:prototypes}, except without this identification.

The marking $\Sigma_2\to Y$ is well-defined up to Dehn twist around the curves of $\system{3}$ or
$\system{4}$.  Choose a particular marking.  This allows us to consider the symplectic basis $\{\alpha_i,
\beta_i\}$ of $H_1(\Sigma_2; \zed)$ as a basis of $H_1(Y;\zed)$ as well.

By Propositions~\ref{prop:admissibleclassification} and \ref{prop:X200classification}, there is a unique $\mu$
in $K_D$ so that $\mu\cdot\alpha_1=\alpha_2$, using the real multiplication on $\Jac(Y)$.  Since $\eta$ is an
eigenform, we know that $\mu^{(1)}=\lambda$ because
$$\lambda = \eta(\alpha_2) = \eta(\mu\cdot\alpha_1) = \mu^{(1)} \eta(\alpha_1) = \mu^{(1)}.$$
If $D$ is not square, this means that $\mu$ doesn't depend on $(Y, \eta)\in U$ because $\mu$ is determined by its
first embedding.  If $D$ is square, however, there can be Abelian differentials close to $(X, \omega)$
with different $\mu^{(2)}$.

Let $a$, $b$, and
$c$ be the integers given by Proposition~\ref{prop:admissibleclassification}.  Since $\ord=\zed[a\mu]$, the
algebraic integer $a\mu$
acts on $H_1(Y; \zed)$ by an integer matrix.
 
\begin{prop}
  The matrix $T$ of the action of $a\mu$ on $H_1(Y; \zed)$ in the basis $\{\alpha_1, \alpha_2, \beta_1, \beta_2\}$
  is
   \begin{equation}
    \label{eq:matrixofamu}
    T=
    \begin{pmatrix}
      0 & -c & \phantom{-}0 & \phantom{-}q \\
      a & -b & -q & \phantom{-}0 \\
      0 & \phantom{-}0 & \phantom{-}0 & \phantom{-}a \\
      0 & \phantom{-}0 & -c & -b
    \end{pmatrix}
  \end{equation}
  with $a$, $b$, and $c$ as above, and for some integer $q$.
\end{prop}

\begin{proof}
  The matrix $T$ is self-adjoint with respect to the intersection pairing on $H_1(Y; \zed)$, which is equivalent to
  $$T^t J = J T,$$
  where
  $$J=
  \begin{pmatrix}
    \phantom{-}0 & I\\
    -I & 0
  \end{pmatrix}.
  $$
  This easily implies that
  $$T=
  \begin{pmatrix}
    A & Q \\
    0 & A^t
  \end{pmatrix}
  $$
  with $Q=-Q^t$.

  To calculate $A$, we know that
  $$a\mu\cdot\alpha_1 = a \alpha_2,$$
  and
  $$a\mu\alpha_2 = a \mu^2 \alpha_1 = (-b\mu-c)\alpha_1 = -c \alpha_1 - b \alpha_2,$$
  so
  $$A=
  \begin{pmatrix}
    0 & -c\\
    a & -b
  \end{pmatrix}.
  $$
\end{proof}

Now let $P(Y, \eta)=(a, b, c, \bar{q})$, where $q$ is as in the above Proposition, and $\bar{q}$ is the
residue class of $q$, taken modulo $\gcd(a, b, c)$.  Similarly, if $(X, \omega)$ has a double zero or a
separating node, let $P'(Y, \eta)=(a, b, c, \bar{q})$, where $\bar{q}$ is the residue class of $q$, taken
modulo $\gcd(a, c)$. Note that the definitions of $P(Y, \eta)$ and $P'(Y, \eta)$ involved some choices: there
was a choice of which $\beta_i$ to take as a basis dual to the $\alpha_i$ in $H_1(\Sigma_2; \zed)$, and there
was a choice of a marking $\Sigma_2\to Y$ to take up to Dehn twist around $\system{3}$ or $\system{4}$
which we used to transport the basis of $H_1(\Sigma_2; \zed)$ to $H_1(Y; \zed)$.

\begin{theorem}
  \label{thm:prototypewelldefined}
  The quadruple $P(Y, \eta)$ is a fine $\Y$-prototype which does not depend on the choice of the $\beta_i$ or
  the marking $\Sigma_2\to Y$.  If $(X, \omega)$ has a double zero or a separating node, then $P'(Y, \eta)$ is
  a fine $\W$ or $\P$-prototype which is also independent of these choices.
\end{theorem}

\begin{proof}
  We first show that $P(Y, \eta)$ satisfies the required properties of a prototype.  It follows directly
  from Proposition~\ref{prop:admissibleclassification} that $b^2-4ac=D$ and $a>0$.  We know $N^{K_D}_\ratls(\mu)<0$
  because $\{\alpha_1, \alpha_2\}$ is an admissible basis of the $\ord$-module $S=\langle\alpha_1, \alpha_2\rangle$, so $c<0$.  Since
  $\zed[a\mu]=\ord$ by Proposition~\ref{prop:admissibleclassification}, the matrix in \eqref{eq:matrixofamu} must be
  primitive because the action of $\ord$ on $\Jac(Y)$ is proper; thus $\gcd(a, b, c, \bar{q})=1$.  Finally,
  $a+b+c\leq 1$ because $\mu^{(1)}\geq1$.
  
  Now we will show that $P(Y, \eta)$ is independent of the choices made.  Replacing the $\beta_i$ with a new
  pair which are dual to the $\alpha_i$ or changing the marking $\Sigma_2\to Y$ by Dehn twists around
  $\system{3}$ or $\system{4}$ produces a new pair $\beta_i'$ of the form,
  \begin{align}
    \label{eq:betaiprime}
    & {\beta}_1'= {\beta}_1 + r  {\alpha}_1 + s  {\alpha}_2\\
    & {\beta}_2'= {\beta}_2 + s  {\alpha}_1 + t  {\alpha}_2 \notag
  \end{align}
  for integers $r$, $s$, and $t$.
  
  The matrix sending the basis $\{ {\alpha}_1,  {\alpha}_2,  {\beta}_1,  {\beta}_2\}$ to
  $\{ {\alpha}_1,  {\alpha}_2,  {\beta}_1',  {\beta}_2'\}$ is
  $$S=
  \begin{pmatrix}
    1 & 0 & r & s\\
    0 & 1 & s & t\\
    0 & 0 & 1 & 0\\
    0 & 0 & 0 & 1
  \end{pmatrix},
  $$

  If $T$ is the matrix of the endomorphism of $H_1(Y; \zed)$ given by multiplication by $a \lambda$ in
  the old basis. and  $T'$ is the matrix of the same endomorphism in the new basis,  then $T'$ is
  $$T'=S \cdot T\cdot S^{-1},$$
  which  is
  $$\begin{pmatrix}
    0 & -c & \phantom{-}0 & \phantom{-}q'\\
    a & -b & -q' & \phantom{-}0\\
    0 & \phantom{-}0 & \phantom{-}0 & \phantom{-}a\\
    0 & \phantom{-}0 & -c & -b
  \end{pmatrix},$$
  where
  \begin{equation}
    \label{eq:qprime}
    q'= q + a r - b s + c t.
  \end{equation}
  It follows that $q'\equiv q \mod{\gcd(a, b, c)},$ so the prototype is
  independent of the choices.
  
  Now suppose $(X, \omega)$ does have a separating node or a double zero.  We must now show that $q$ is
  well-defined modulo $\gcd(a,c)$.  In this case, the choice of $\beta_i$ is uniquely determined up to adding
  a multiple of $\alpha_i$ because we imposed additional requirements on the $\beta_i$.  That means that
  different choices give a new pair $\{\beta_i'\}$ of the form \eqref{eq:betaiprime} with $s=0$.  Then from
  \eqref{eq:qprime}, we see that $q'\equiv q \mod \gcd(a, c)$.
  
  Finally, if $(X, \omega)$ has a double zero, then $\lambda>1$ by
  Proposition~\ref{prop:doublezerosinclambda}, so $a+b+c<0$, and $P'(Y, \eta)$ is then a $\W$-prototype.
\end{proof}

\begin{remark}
  If $\lambda=1$, then there were also two possible ways to choose the $\alpha_i$.  In this case, one can
  check that swapping the $\alpha_i$ has the effect of applying the involution \eqref{eq:involutionone} to the
  prototypes $P(Y, \eta)$ and $P'(Y, \eta)$.  Thus they are independent of this choice as ordinary prototypes
  but not as fine prototypes.
\end{remark}

\paragraph{Equations for $\barX$ in local coordinates.}

Given $(X, [\omega])\in C_\lambda$ and a small neighborhood $U$ of $(X, [\omega])$ in the Dehn space as above,
we now give equations for $\barX$ in $U$ and show that $(X, [\omega])$ is indeed in $\barX$.
Since the prototype $P(Y, \eta)$ is locally constant on $U\cap\X$,
we can use these prototypes to label branches of $U\cap\X$.  Given a fine $\Y$-prototype $P$, let
$V_P\subset U$ be
the closure in $U$ of
$$\{(Y, \eta)\in U\cap\X : P(Y, \eta)=P\},$$
which in fact lies in $U_\lambda$ by Proposition~\ref{prop:X200classification}.
If $(X, \omega)$ has a double zero or a separating node,  and $P$ is a fine $\W$ or $\P$-prototype (respectively), then let
$V_P'$ the analogous locus where $P'(Y, \eta)=P$.  In this case, given a fine $\W$-prototype (or $\P$-prototype) $P$,
$$V_P= \bigcup_Q V_Q',$$
where the union is over all fine $\W$-prototypes (or $\P$-prototypes) which map to $P$.

Let
$$(a', b', c', \bar{q}') = (a, b, c, \bar{q})/\gcd(a, b, c),$$
and
$$(a'', b'', c'', \bar{q}'') = (a, b, c, \bar{q})/\gcd(a, c)$$
(here we interpret $\bar{q}'$ and $\bar{q}''$ as elements of $\ratls/\zed$).

Recall that we have coordinates $(x, y, z)$ on $U_\lambda$ defined above.

\begin{theorem}
  \label{thm:X200coordinates}
  The locus  $V_P$ is cut out by the equation,
  \begin{equation}
    \label{eq:coordinates200}
    y^{a'} = e^{-2 \pi i q'} z^{-c'},
  \end{equation}
  and $C_\lambda$ is cut out by the equations $y=z=0$.
  
  If $(X, \omega)$ has a double zero or a separating node, then $V_P'$ is cut out by the equation,
  \begin{equation}
    \label{eq:coordinates201}
    y^{a''} = e^{-2 \pi i q''} z^{-c''},
  \end{equation}
  and $\barW\cap V_P$ or respectively $\barP\cap V_P$ is cut out by the additional equation $x=0$.  The leaves
  of the foliation $\A$ of $\X$ are given in these coordinates by either of the equivalent equations $y = {\rm
    const}$ or $z = {\rm const}$.
\end{theorem}

\begin{proof}
  The proofs of \eqref{eq:coordinates200} and \eqref{eq:coordinates201} essentially the same, so we will
  only give the proof of \eqref{eq:coordinates200}.

  Choose some $(Y, \eta)\in V_P$ with a symplectic basis $\{\alpha_i, \beta_i\}$ and $\mu\in K_D$ as
  before.  Equation \eqref{eq:matrixofamu} implies
  \begin{align}
    a\mu\cdot\beta_1 &= -q \alpha_2 - c \beta_2 \label{eq:relation1}\\
    a\mu\cdot\beta_2 &= q\alpha_1 + a \beta_1 - b\beta_2.\label{eq:relation2}
  \end{align}
  Since $\eta$ is an eigenform, integrating $\eta$ over both sides gives
  \begin{align}
    a\lambda\eta(\beta_1) &= -q \lambda - c \eta(\beta_2) \label{eq:eigenform1}\\
    a\lambda\eta(\beta_2) &= q + a\eta(\beta_1) - b\eta(\beta_2)\label{eq:eigenform2}
  \end{align}
  because $\mu^{(1)}=\lambda$.
  Dividing \eqref{eq:eigenform1} by $\lambda \gcd(a, b, c)$ and exponentiating yields \eqref{eq:coordinates200}.
  
  Now consider $(Y, \eta)$ which lies in $ U_\lambda\setminus V(y,z)$ and in the locus defined by
  \eqref{eq:coordinates200} (using the notation $V(\cdot)$ for the zero locus of the enclosed functions).  We
  must show that $(Y, \eta)$ is an eigenform for real multiplication by $\ord$ with prototype $P=(a, b, c,
  \bar{q})$.

  Taking a logarithm of \eqref{eq:coordinates200}, we get
  $$a\lambda\eta(\beta_1)\equiv -\lambda\bar{q} - c\eta(\beta_2) \mod \lambda\gcd(a, b, c).$$
  We can then choose some $q\equiv\bar{q}\mod\gcd(a, b, c)$ so that \eqref{eq:eigenform1} holds.  It then
  follows automatically from the relation $a\mu^2+b\mu+c=0$ that \eqref{eq:eigenform2} holds.  We can then
  define real multiplication of $\ord$ on $\Jac(Y)$, by defining $a\mu$ to act via the matrix $T$ in
  \eqref{eq:matrixofamu}.  Equations \eqref{eq:eigenform1} and \eqref{eq:eigenform2} then imply $(Y, \eta)$
  is actually an eigenform for real multiplication with prototype $P$.

  The equations for $\barW$ and $\barP$ are obvious from the definition of the coordinates.  The leaves of
  $\A$ are as claimed because with our normalization the periods $\eta(\alpha_i)$ are constant, and either
  $y={\rm const}$ or $z={\rm const}$ implies that both of the $\eta(\beta_i)$ are constant by \eqref{eq:coordinates200}.
\end{proof}

This also completes the proof of Theorem~\ref{thm:X200classification} because we explicitly constructed
eigenforms near any point of $C_{\lambda(P)}$ for any nondegenerate $\Y$-prototype $P$ in
Theorem~\ref{thm:X200coordinates} and showed that these are the only possible limits in
Corollary~\ref{cor:onlylimits}.

\paragraph{Branches of $\barX$ through $C_\lambda$.}

Let $p=(X, [\omega])\in C_\lambda$, and let $q\in\proj\Omega \Def_2(S)$ (for the appropriate curve system
$S$), with $\pi(q)=p$.  We have seen that around $q$, the variety $\pi^{-1}(\barX)$ is the union of the $V_P$
over all fine $\Y$-prototypes $P$ such that $\lambda(P)=\lambda$.

To understand the local structure of $\barX$ around $p$, we need only to understand how $\Aut(X, [\omega])$
acts on this picture.  If $\lambda>1$, then $\Aut(X, [\omega])$ consists only of the hyperelliptic involution
which acts trivially on $\proj \Omega \Def_2(S)$.  Thus $\pi$ is a local isomorphism around $q$.  Note also
that there is no distinction between fine and ordinary prototypes when $\lambda>1$.

When $\lambda=1$, we have $\Aut(X, [\omega])\isom \zed/2 \oplus \zed/2$.  It is generated by the hyperelliptic
involution together with an involution $t$ which interchanges the two nonseparating nodes of $X$ and preserves
$\omega$.  This involution $t$ acts on $\proj\Omega\Def_2(S)$ by $t(x, y, z) = (x, z, y)$ in the coordinates
around $q$ defined above, so $t$ identifies $V_P$ with $V_{i(P)}$, where $i$ is the involution of
\eqref{eq:involutionone}.  Thus $V_P$ and $V_{i(P)}$ have the same image in $\proj\Omega\barmoduli$.  We will abuse
notation and denote the image of $V_P$ in $\proj \Omega\barmoduli$ by $V_P$ as well, where now $P$ is an
ordinary $\Y$-prototype.  This makes sense because it defines the same object if $P$ is represented by
$i(P)$.  We will think of this $V_P$ as a union of branches of $\barX$ through $p$.  To summarize, whether or
not $\lambda=1$, the germ of $\barX$ around $P$ is the union of the $V_P$, where $P$ ranges over all ordinary
$\Y$-prototypes such that $\lambda(P)=\lambda$.

If $p=(X, [\omega])$ has a double zero or a separating node, then in the same way -- just replacing $V_P$ with
$V_P'$ defined above -- we define a branch $V_P'$
of $\barX$ through $p$, where now $P$ is a $\W$ or $\P$-prototype (respectively).  This gives a finer
classification of branches then we obtain for generic $p\in C_\lambda$

To summarize this discussion, we have the following Corollary of Theorem~\ref{thm:X200coordinates}.

\begin{cor}
  \label{cor:clambdabranches}
  For any $p\in C_\lambda$, the germ of $\barX$ through $p$ is the union of germs,
  \begin{equation*}
    \bigcup_{\substack{
        P\in\Yprot\\
        \lambda(P)=\lambda}}
    V_P.
  \end{equation*}
  Each germ $V_P$ consists of $\mult(P)$ branches of $\barX$ through $P$.

  If $p=w_\lambda$ or $p=p_\lambda$, then alternatively the germ of $\barX$ through $p$ is the union,
  \begin{equation*}
    \bigcup_{\substack{
        P\in\Wprot\\
        \lambda(P)=\lambda}}
    V_P',
  \end{equation*} or respectively,
  \begin{equation*} \bigcup_{\substack{ P\in\Pprot\\ \lambda(P)=\lambda}} V_P'.
  \end{equation*} In these cases, each of these germs is actually an irreducible branch of $\barX$ through
$p$.
\end{cor}

\begin{proof} The only parts that don't follow directly from Theorem~\ref{thm:X200coordinates} and the
  discussion above are the
statements that $V_P$ consists of $\mult(P)$ branches and that $V_P'$ is irreducible.  This will
follow directly from Proposition~\ref{prop:easynormalization} because
  $$\gcd(a', c') = \frac{\gcd(a,c)}{\gcd(a, b, c)}=\mult(P),$$
  and $\gcd(a'', c'')=1$.
\end{proof}

\paragraph{A section of $\mathcal{O}(-1)$.}

We normalized each projective class $(Y, [\eta])\in U_\lambda$ so that each $\eta(\alpha_1)=1$.  We can regard
this as defining a section $s$ of the canonical line bundle $\mathcal{O}(-1)$ over $U_\lambda$.

The following follows directly from \eqref{eq:area} together with the definition of our coordinates.
\begin{prop}
  \label{prop:area3} The norm of this section $s$ of $\mathcal{O}(-1)$ over $U_\lambda$ is given by,
  \begin{equation}
   \label{eq:area3} h(s, s) = -\frac{1}{2\pi}\log|y| - \frac{\lambda^2}{2\pi}\log|z|.
  \end{equation}
\end{prop}

\subsection{The strata $X_D(T_{3,0})$ and $X_D(T_{3,0}, T_{3,0}^3)$}
\label{subsec:threenonsepnodes}

We now study Abelian differentials in $\barX$ with three nonseparating nodes.  These consist of the union of the
stratum $\X(\system{5})$, where all of the nodes are polar, and the stratum $\X(\system{5},\system[3]{5})$,
where one of the nodes is holomorphic.  We will prove the following classification of such differentials.

\begin{theorem}
  \label{thm:X300classification}
  The stratum $\X(\system{5})$ is the finite union of points,
  $$\X(\system{5}) = \bigcup_P c_{\lambda(P)},$$
  where the union is over all nonterminal $\Y$ prototypes $P$.

  The stratum $\X(\system{3},\system[3]{5})$ is the finite union of points,
  $$\X(\system{5}, \system[3]{5}) = \bigcup_P c_{\lambda(P)},$$
  where the union is over all terminal $\Y$-prototypes.
  
  If $P$ is neither degenerate nor terminal, then $c_{\lambda(P)}$ lies in the intersection of the closures of
  $C_{\lambda(P)}$ and $C_{\lambda(P^+)}$.  If $P$ is a degenerate prototype, then $c_P$ lies in the
  intersection of the closures of $\X(\system{2})$ and $C_{\lambda(P^+)}$.  If $P$ is a terminal prototype,
  then $c_{\lambda(P)}$ lies in the intersection of the closures of $C_{\lambda(P)}$ and $\X(\system{2},\system[1]{2})$.
  None of the other curves in $\bdry\X$ pass through $c_{\lambda(P)}$.
\end{theorem}

\paragraph{Local coordinates.}

Fix an $(X, [\omega])\in \proj\Omega\moduli^0(\system{5})$ or $\proj\Omega\moduli^0(\system{5}, \system[3]{5})$, and
choose a marking $\Sigma_2\to X$ so that we can regard $(X, [\omega])$
as lying in
$\proj\Omega\teich^0_2( \system{5})$ or $\proj\Omega\teich^0_2( \system{5}, \system[3]{5})$.

Our first goal is to show that if $(X, [\omega])\in\barX$, then $(X, [\omega])$ is one of the points
$c_{\lambda(P)}$ as in the statement of Theorem~\ref{thm:X300classification}.  If $(X, [\omega])\in\barX$,
then $\omega$ is an eigenform for real multiplication of $\ord$ on $\Jac(X)$ by Theorem~\ref{cor:rmclosed}, so
assume that it is such an eigenform.  It follows that the ratios of the residues of any two poles of $\omega$
are real, so we can normalize the projective class $[\omega]$ so that all of these residues are real.

Suppose $\omega$ has no holomorphic node.  In this case it also has two simple zeros; let $\{p_1, p_2\}\in
\Sigma_2$ be the inverse image of the zeros.  Let $\alpha_i \in H_1(\Sigma_2 \setminus \{p, q\})$ be homology
classes representing the three curves $\system[i]{5}$ such that $\omega(\alpha_i)\geq 0$, and choose
$\gamma_i\in H_1(\Sigma_2, \{p_1, p_2\})$ so that $\alpha_i\cdot\gamma_j=\delta_{ij}$ as in
\S\ref{subsec:stratum6}.  Then the classes,
\begin{align}
  \label{eq:betagamma}
  \beta_1 &= \gamma_1 - \gamma_3 \\
  \beta_2 &= \gamma_2 + \gamma_3, \notag
\end{align}
in $H_1(\Sigma_2; \zed)$ are dual to $\alpha_1$ and $\alpha_2$.

If $\omega$ does have a holomorphic node, then let $\alpha_i \in H_1(\Sigma_2)$ be homology classes
representing the three curves $\system[i]{5}$ such that $\omega(\alpha_i)\geq 0$, and let $\beta_i\in
H_1(\Sigma_2; \zed)$ be dual to $\alpha_1$ and $\alpha_2$.

Choose a small neighborhood $U$ of $(X, [\omega])$ in $\proj\Omega\Def^0_2(\system{5})$ or
$\proj\Omega\Def^0_2(\system{5}, \system[3]{5})$ such that the natural map $U/\Aut(X,
[\omega])\to\proj\Omega\barmoduli$ is an isomorphism onto its image and such that for each $(Y, [\eta])\in U\cap \pi^{-1}(\X)$, the
subgroup $S=\langle\alpha_1, \alpha_2\rangle$ of $H_1(\Sigma_2; \zed)$ is preserved by the real
multiplication, which is possible by Theorem~\ref{cor:rmnearboundary}.  

For the rest of this section, we will adopt the convention that each projective class $(Y, [\eta])\in U$ is
normalized so that $\eta(\alpha_1)=1$, so we can regard $U$ as an open subset of a hypersurface in
$\Omega\Def^0_2(\system{5})$ or $\Omega\Def^0_2(\system{5}, \system[3]{5})$.

\begin{prop}
  \label{prop:threenonsepclassification}
  If the neighborhood $U$ of $(X, [\omega])$ is taken sufficiently small, then for each $(Y, \eta)\in U\cap\pi^{-1}(\X)$,
  the $\ord$-submodule $S$ of $H_1(\Sigma_2, \zed)$ is a quasi-invertible $\ord$-module with $\{\alpha_1,
  \alpha_2, \alpha_3\}$ an admissible triple; furthermore, for $(Y, \eta)\in U$, the period maps $P_\eta$ and
  $P_{\omega}\colon S\to \reals$ are the same.
\end{prop}

\begin{proof}
  Since the curves of $\system{5}$ bound a pair of pants, we get the relation
  \begin{equation}
    \label{eq:pantsrelation}
    \pm\alpha_1\pm\alpha_2\pm\alpha_3.
  \end{equation}
  To see that $\{\alpha_1, \alpha_2, \alpha_3\}$ is an admissible triple, it is enough to show that
  $\{\alpha_1, \alpha_2\}$ is an admissible basis.  Just as in the proof of
  Proposition~\ref{prop:X200classification}, it is enough to show that $\Im\eta(\beta_i)>0$ for $(Y,
  \eta)$ sufficiently close to $(X, \omega)$.  This follows  from
  Theorem~\ref{thm:longcylinders}, using the fact that the $\beta_i$ pass through tall cylinders, and each
  cylinder that $\beta_i$ passes through contributes positively to $\Im\eta(\beta_i)$ because
  $\alpha_j\cdot\beta_i\geq0$ for each $j$.

  Since there are only finitely many admissible triples up to isomorphism by
  Proposition~\ref{prop:tripleclassification}, the period $\eta(\alpha_2)$ must be constant on $U$ if $U$ is
  sufficiently small, and $\eta(\alpha_3)$ must be constant by \eqref{eq:pantsrelation}.  Thus the period maps
  $P_\omega$ and $P_\eta$ are the same.
\end{proof}

Now since $\{\alpha_1, \alpha_2, \alpha_3\}$ is an admissible triple, we can reorder them and replace $\omega$
by a constant multiple so that
$$\omega(\alpha_1)=1,\, \omega(\alpha_2)=\lambda,\, \text{and } \omega(\alpha_3)=\lambda-1,$$
with
$\lambda\geq1$, and $\lambda=\lambda(P)$ for some $\Y$-prototype $P$ by
Proposition~\ref{prop:tripleclassification}.  If $\omega$ has no holomorphic node, then there is a unique such
choice of the $\alpha_i$ and $\omega$ (up to the hyperelliptic involution).  Otherwise, there are two such
choices because we can swap $\alpha_1$ and $\alpha_2$.  Also when $(X, [\omega])$ has no holomorphic nodes, reorder and
possibly change the signs of the $\gamma_i$ so that we still have $\alpha_i\cdot\gamma_j=\delta_{ij}$, and
redefine the $\beta_i$ as in \eqref{eq:betagamma} with the new $\gamma_i$.  When $(X, [\omega])$ has a
holomorphic node, reorder the $\beta_i$ so that they are still dual to the $\alpha_i$.

Thus we see that if $(X, [\omega])$ has three nonseparating nodes and $(X, [\omega])\in\barX$, then $(X,
[\omega]) = c_{\lambda(P)}$ for some $\Y$-prototype $P$ as stated in Theorem~\ref{thm:X300classification}.
Also, Proposition~\ref{prop:threenonsepclassification} implies that if $U$ is sufficiently small and if we
define for $\mu\in\reals$,
$$U_\mu = \{(Y, \eta)\in U : \eta(\alpha_1)=1\: \text{and } \eta(\alpha_2) = \mu \},$$
then $U\cap\pi^{-1}(\barX)\subset
U_\lambda$, where $\lambda=\omega(\alpha_2)$, by Proposition~\ref{prop:threenonsepclassification}.

Recall that in \S\ref{subsec:stratum6} we defined local coordinates $(v, w, x, y, z)$ on
$\Omega\Def^0_2(\system{5})$.  When $\lambda>1$, in these local coordinates, we can identify
$U_\lambda$ with the subspace defined by $v=1$ and $w=\lambda$.  Then on $U_\lambda$, the $(x, y, z)$ become
\begin{align*}
  x &= e^{2\pi i\eta(\gamma_1)} \\
  y &= e^{2\pi i\eta(\gamma_2)/\lambda} \\
  z &= e^{2\pi i \eta(\gamma_3)/(\lambda-1)},
\end{align*}
and the $(x, y, z)$ are local coordinates on $U_\lambda$.\label{page:coordinates30}

In \S\ref{subsec:stratum7} we defined local coordinates $(w, x, y, z)$ on $\Omega\Def^0_2(\system{5},
\system[3]{5})$.  In these local coordinates, we can identify $U_1$ with the subspace defined by $w=1$.  Then
on $U_1$, the $(x, y, z)$ become
\begin{align*}
  x &= e^{2\pi i\eta(\beta_1)} \\
  y &= e^{2\pi i\eta(\beta_2)} \\
  z &= \left(\int_I\eta\right)^2,
\end{align*}
and the $(x, y, z)$ are local coordinates on $U_1$.

\paragraph{Prototypes.}

Now assume that our $(X, [\omega]) = c_\lambda$, where $\lambda=\lambda(P)$ for some 
$\Y$-prototype $P$ with a small neighborhood $U$ of $(X, [\omega ])$ in the Dehn space chosen as above.
We now assign a fine $\Y$-prototype $P(Y, [\eta])$ to each $(Y, [\eta])\in U\cap\pi^{-1}(\X)$ as in
\S\ref{subsec:twononseppolarnodes}.

Let $(Y, [\eta])\in U\cap \pi^{-1}(\X)$.  We have a marking $\Sigma_2\to Y$ which is defined up to Dehn twist around
the curve system $\system{5}$. Choose a marking $\Sigma_2\to Y$ and use this to transport the symplectic basis
$\{\alpha_1, \alpha_2, \beta_1, \beta_2\}$ to $Y$.
  
There is some $\mu\in K_D$ with $\mu^{(1)}=\lambda$ such that $\mu\cdot\alpha_1 = \alpha_2$ and
$$a\mu^2 + b\mu + c=0$$
for integers $a$, $b$, and $c$ satisfying the properties of
Proposition~\ref{prop:tripleclassification}.  The generator $a\mu$ of $\ord$ defines by its action on
$\Jac(Y)$ an endomorphism $T$ of $H_1(Y;\zed)$ which is given by the matrix \eqref{eq:matrixofamu} for some integer
$q$.  Define the prototype $P(Y, \eta)=(a, b, c, \bar{q})$, where $\bar{q}$ is the reduction of $q$ modulo
$\gcd(a, b, c)$.  By the same argument as in Theorem~\ref{thm:prototypewelldefined}, $P(Y, [\eta])$ is a
well-defined fine $\Y$-prototype.

\paragraph{Local coordinates for $\barX$.}

Given $(X, [\omega]) = c_{\lambda}$ and a small neighborhood $U$ of $(X, [\omega])$ in the Dehn space as above,
we now give equations for $\pi^{-1}(\barX)$ in $U$ and show that $(X, [\omega])$ is indeed in $\barX$.
For each fine $\Y$-prototype $P$, let $\tilde{V}_P\subset U$ be the closure of
$$\{(Y, \eta)\in U\cap\pi^{-1}(\X) : P(Y, \eta) = P\},$$
which actually lies in $U_\lambda$.  The following theorem gives explicit equations for $\tilde{V}_P$ in the
coordinates $(x, y, z)$ on $U_\lambda$.

\begin{theorem}
  \label{thm:X300coordinates}
  The locus $\tilde{V}_P\subset U_{\lambda(P)}$ is cut out by the equation,
  \begin{equation}
    \label{eq:X300equation}
    x^{a'} = e^{-2 \pi i q'} y^{-c'} z^{-a'-b'-c'}.
  \end{equation}
  
  If $P$ is a nondegenerate and nonterminal $\Y$-prototype, then the closure $\overline{C}_{\lambda(P)}$ is
  cut out by the equations $x = y = 0$, and $\overline{C}_{\lambda(P^+)}$ is cut out by the equations
  $x=z=0$.
  
  If $P$ is a degenerate $\Y$-prototype, then
  the additional equation $y = 0$ cuts out $U\cap\X(\system{2})$, and the equations $x=z=0$ cuts out
  $\overline{C}_{\lambda(P^+)}$.

  If $P$ is a terminal $\Y$-prototype, then 
  the additional equation $z= 0$ cuts out $U\cap\X(\system{2}, \system[1]{2})$, and the equation $x=y=0$ cuts out
  $\overline{C}_{\lambda(P)}$.

  The leaves of the foliation $\A$ of $\X$ are given in these local coordinates by $y^\lambda z^{\lambda-1} =
  {\rm const}$.
\end{theorem}

\begin{proof}
  If the prototype $P$ is terminal, then \eqref{eq:X300equation} reduces to \eqref{eq:coordinates200}, and the
  proof that $\tilde{V}_P$ is given by this equation is the same as in Theorem~\ref{thm:X200coordinates}.
  
  Now assume that $P$ is not terminal.  Using the same symplectic basis $\{\alpha_i, \beta_i\}$ and $\mu\in
  K_D$ such that $\mu\cdot\alpha_1=\alpha_2$ as above, equation \eqref{eq:matrixofamu} again implies
  \eqref{eq:relation1} and \eqref{eq:relation2}.  Integrating $\eta$ over both sides yields
  \eqref{eq:eigenform1} and \eqref{eq:eigenform2}, and substituting in \eqref{eq:betagamma} gives
  \begin{equation}
    \label{eq:eigenform3}
    a \eta(\gamma_1) + \frac{c}{\lambda} \eta(\gamma_2) + \frac{a+b+c}{\lambda-1}\eta(\gamma_3) = -q 
  \end{equation}
  Dividing by $\gcd(a, b, c)$ and exponentiating yields \eqref{eq:X300equation}.

  Conversely if $(Y, \eta)$ is in $U_\lambda$ minus the coordinate axes and is in the locus defined
  by \eqref{eq:X300equation}, then \eqref{eq:eigenform3} holds modulo $\gcd(a, b, c)$, and by changing
  $q$ modulo $\gcd(a, b, c)$, we can ensure that it holds exactly.  It then follows that \eqref{eq:eigenform1}
  and \eqref{eq:eigenform2} hold, and so $(Y, \eta)$ is an eigenform.

  The equations for the $C_\lambda$ and the other strata are straightforward.
  
  The periods $\eta(\alpha_1)$ and $\eta(\alpha_2)$ are fixed because of our normalization for $\eta$, and the periods
  $\eta(\beta_1)$ and $\eta(\beta_2)$ are related by the real multiplication, so if one is constant, then the
  other is.  Thus the leaves of the foliation $\A$ in these coordinates are given by $\eta(\beta_2)={\rm
    const}$, which is equivalent to $y^\lambda z^{\lambda-1} = {\rm const}$.
\end{proof}

\begin{proof}[Proof of Theorem~\ref{thm:X300classification}]
  We saw as a consequence of Proposition~\ref{prop:threenonsepclassification} that each stable Abelian
  differential in $\barX$ with three nonseparating nodes is one of the $c_{\lambda(P)}$ for some
  $\Y$-prototype $P$, and it follows from Theorem~\ref{thm:X300coordinates} that these $c_{\lambda(P)}$ are in
  fact in $\barX$.  The other statements also follow directly from Theorem~\ref{thm:X300coordinates}.
\end{proof}

\paragraph{Branches of $\barX$ through $c_\lambda$.}

The varieties cut out by \eqref{eq:X300equation} have a single branch through the origin, as we will see in
Proposition~\ref{prop:hardnormalization}, so we obtain a bijective correspondence between branches of
$\pi^{-1}(\barX)$ through $\pi^{-1}(c_\lambda)$ and fine $\Y$-prototypes $P$ such that $\lambda(P)=\lambda$
given by $P\mapsto \tilde{V}_P$.  Just as in \S\ref{subsec:twononseppolarnodes}, the action of $\Aut(X,
[\omega])$ identifies branches corresponding to the same ordinary $\Y$-prototype, so we obtain the first
statement of the following theorem.

\begin{theorem}
  \label{thm:branchintersection}
  The germ of $\barX$ through $c_\lambda$ is the following union of irreducible branches:
   \begin{equation*}
    \bigcup_{\substack{
        P\in\Yprot\\
        \lambda(P)=\lambda}}
    \tilde{V}_P.
  \end{equation*}
  
  If $P$ is a nondegenerate $\Y$-prototype such that $\lambda(P)=\lambda$, then the branch $\tilde{V}_P$
  through $c_\lambda$ intersects $C_\lambda$ in the germ $V_P$ of Corollary~\ref{cor:clambdabranches}.  If $P$
  is a nonterminal $\Y$-prototype such that $\lambda(P)=\lambda$, then the branch $\tilde{V}_P$ through
  $c_\lambda$ intersects $C_{\lambda(P^+)}$ in the germ $V_{P^+}$.
\end{theorem}

\begin{proof}  
  Suppose that $P=(a, b, c, \bar{q})$ is nonterminal.  Let $(Y, \eta)\in\tilde{V}_P$ be a nonsingular Abelian
  differential close enough to a point $p\in C_{\lambda(P^+)}$ that it is contained in some well-defined
  subgerm $V_{P'}$ of $\barX$ through $p$.  We need to show that $P'=P^+$.
  
  We will continue to use the notation that we used to define the local coordinates around $c_\lambda$, so we
  have classes $\alpha_1$, $\alpha_2$, and $\alpha_3\in H_1(Y; \zed)$ representing the nodes of $c_\lambda$ such that
  $\eta(\alpha_1)=1$, $\eta(\alpha_2)=\lambda$, and $\eta(\alpha_3)=\lambda-1$; we have $\beta_1$,
  and $\beta_2\in H_1(Y; \zed)$ such that $\{\alpha_1, \alpha_2, \beta_1, \beta_2\}$ is a symplectic basis of
  $H_1(Y; \zed)$; and we have a $\mu\in K_D$ such that $\mu\cdot\alpha_1=\alpha_2$ and $a\mu^2 + b\mu + c=0$.

  In order to calculate the prototype $P'$, we must choose a symplectic basis $\{{\alpha}'_1, {\alpha}'_2, {\beta}'_1,
  {\beta}'_2\}$ of $H_1(Y; \zed)$ such that the ${\alpha}'_i$ are homologous to cylinders on $Y$ representing
  the nodes of $p\in C_{\lambda(P^+)}$ and such that if we set
  $$\tilde{\lambda} = \frac{\eta({\alpha}'_2)}{\eta({\alpha}'_1)},$$
  then $\tilde{\lambda}\geq1$ and $\tilde{\lambda}'<0$.  Given such a basis, let $\tilde{\mu}\in K_D$ be
  such that $\tilde{\mu}\cdot{\alpha}'_1 = {\alpha}'_2$.  Let $a'$, $b'$, and $c'$ be the unique integers such
  that
  $$a'\tilde{\mu}^2 + b'\tilde{\mu} + c' = 0,$$
  $a'>0$, and $(b')^2 -4a'c'=D$.  For $x\in K_D$, set $T_x$ and
  ${T}'_x$ to be the matrices of the action of $x$ on $H_1(Y; \zed)$ in the bases $\{\alpha_i, \beta_i\}$ and
  $\{{\alpha}'_i, {\beta}'_i\}$ respectively.  Then the prototype $P'$ is $P'=(a', b', c',
  \bar{q}')$, where $q'$ is the entry in the upper right corner of $T_{a'\tilde{\mu}}'$.

  There are two cases to consider, depending on whether or not $\lambda\geq2$.  First suppose that
  $\lambda\geq2$, or equivalently $4a + 2b + c\leq 0$.  Define the new symplectic basis of $H_1(Y; \zed)$,
  \begin{align*}
    \alpha_1' &= \alpha_1  & \beta_1' &= \beta_1 + \beta_2 \\
    \alpha_2' &= \alpha_3=\alpha_2 - \alpha_1 & \beta_2' &=  \beta_2.
  \end{align*}
  Then $\tilde{\mu} = \mu-1$, and $\tilde{\lambda} = \lambda-1$, which satisfies the required properties.  We
  have
  $$(a', b', c') = (a, 2a+b, a+b+c).$$
  The matrix $T_{a \mu}$ is given by \eqref{eq:matrixofamu}, and an easy calculation shows that the new matrix
  is
  $$
   T_{a'\tilde{\mu}}'=
    \begin{pmatrix}
      0 & -a-b-c & 0 & q \\
      a & -2a-b & -q & 0 \\
      0 & 0 & 0 & a \\
      0 & 0 & -a-b-c & -2a-b
    \end{pmatrix},
    $$
    so $P'=P^+$ is as claimed.

    Now consider the case where $\lambda <2$, or equivalently $4a + 2b + c >0$.  In this case, define a new
    symplectic basis of $H_1(Y; \zed)$ by,
    \begin{align*}
    \alpha_1' &= \alpha_2 - \alpha_1 & \beta_1' &= \beta_2 \\
    \alpha_2' &= \alpha_1 & \beta_2' &= \beta_1 + \beta_2.
  \end{align*}
  Then
  $$\tilde{\mu} = \frac{1}{\mu-1} = \frac{-a \mu - a - b}{a + b + c},$$
  and $\tilde{\lambda} = 1/(\lambda-1)$ satisfies the required properties that $\tilde{\lambda}\geq 1$ and
  $\tilde{\lambda}' <0$, using the fact that $\lambda>1$ because $P$ is nonterminal.  We have
  $$(a', b', c') = (-a-b-c, -2a-b, -a),$$
  and an easy calculation shows that the new matrix is,
  $$T_{a'\tilde{\mu}}' = T_{a\mu + a + b}' =
  \begin{pmatrix}
    0 & a & 0 & q \\
    -a-b-c & 2a + b & -q & 0 \\
    0 & 0 & 0 & -a -b -c \\
    0 & 0 & a & 2a +b
  \end{pmatrix}.
  $$
  Thus $P'=P^+$ as claimed.
  
  Calculating the prototype of the intersection of $\tilde{V}_P$ with $C_\lambda$ is much easier. Here we can
  use the old basis $\{\alpha_i, \beta_i\}$ of $H_1(Y; \zed)$ to calculate the prototype, and thus it is just $P$.
\end{proof}

\paragraph{A section of $\mathcal{O}(-1)$.}

We normalized each projective class $(Y, \eta)\in U_\lambda$ so that each
$\eta(\alpha_1)=1$.  We can regard this as defining a section $s$ of the canonical line bundle $\mathcal{O}(-1)$ over
$U_\lambda$.

The following follows directly from \eqref{eq:area} together with the definition of our coordinates.
\begin{prop}
  \label{prop:area4}
  The norm of this section $s$ of $\mathcal{O}(-1)$ over $U_\lambda$ is given by,
  \begin{equation}
    \label{eq:area4}
    h(s, s) = -\frac{1}{2\pi}\log|x| - \frac{\lambda^2}{2\pi}\log|y| - \frac{(\lambda-1)^2}{2\pi}\log|z|.
  \end{equation}
\end{prop}


%% file: splitting.pstex_t
\begin{picture}(0,0)%
\includegraphics{splitting.pstex}%
\end{picture}%
\setlength{\unitlength}{3947sp}%
\begingroup\makeatletter\ifx\SetFigFont\undefined%
\gdef\SetFigFont#1#2#3#4#5{%
  \reset@font\fontsize{#1}{#2pt}%
  \fontfamily{#3}\fontseries{#4}\fontshape{#5}%
  \selectfont}%
\fi\endgroup%
\begin{picture}(5626,2464)(2646,-5346)
\put(4677,-3519){\makebox(0,0)[lb]{\smash{{\SetFigFont{12}{14.4}{\rmdefault}{\mddefault}{\updefault}{\color[rgb]{0,0,0}$p$}%
}}}}
\put(4060,-3628){\makebox(0,0)[lb]{\smash{{\SetFigFont{12}{14.4}{\rmdefault}{\mddefault}{\updefault}{\color[rgb]{0,0,0}$I_1$}%
}}}}
\put(5109,-3614){\makebox(0,0)[lb]{\smash{{\SetFigFont{12}{14.4}{\rmdefault}{\mddefault}{\updefault}{\color[rgb]{0,0,0}$I_2$}%
}}}}
\put(4663,-4322){\makebox(0,0)[lb]{\smash{{\SetFigFont{12}{14.4}{\rmdefault}{\mddefault}{\updefault}{\color[rgb]{0,0,0}$q$}%
}}}}
\put(6448,-4029){\makebox(0,0)[lb]{\smash{{\SetFigFont{12}{14.4}{\rmdefault}{\mddefault}{\updefault}{\color[rgb]{0,0,0}$J_1$}%
}}}}
\put(7687,-4013){\makebox(0,0)[lb]{\smash{{\SetFigFont{12}{14.4}{\rmdefault}{\mddefault}{\updefault}{\color[rgb]{0,0,0}$J_2$}%
}}}}
\put(6619,-4597){\makebox(0,0)[lb]{\smash{{\SetFigFont{12}{14.4}{\rmdefault}{\mddefault}{\updefault}{\color[rgb]{0,0,0}$p_1$}%
}}}}
\put(7556,-4581){\makebox(0,0)[lb]{\smash{{\SetFigFont{12}{14.4}{\rmdefault}{\mddefault}{\updefault}{\color[rgb]{0,0,0}$p_2$}%
}}}}
\put(7401,-3506){\makebox(0,0)[lb]{\smash{{\SetFigFont{12}{14.4}{\rmdefault}{\mddefault}{\updefault}{\color[rgb]{0,0,0}$q_2$}%
}}}}
\put(6595,-3566){\makebox(0,0)[lb]{\smash{{\SetFigFont{12}{14.4}{\rmdefault}{\mddefault}{\updefault}{\color[rgb]{0,0,0}$q_1$}%
}}}}
\end{picture}%

%% file: symplecticbasis.pstex_t
\begin{picture}(0,0)%
\includegraphics{symplecticbasis.pstex}%
\end{picture}%
\setlength{\unitlength}{3947sp}%
\begingroup\makeatletter\ifx\SetFigFont\undefined%
\gdef\SetFigFont#1#2#3#4#5{%
  \reset@font\fontsize{#1}{#2pt}%
  \fontfamily{#3}\fontseries{#4}\fontshape{#5}%
  \selectfont}%
\fi\endgroup%
\begin{picture}(3260,2682)(4457,-5147)
\put(5645,-4645){\makebox(0,0)[lb]{\smash{{\SetFigFont{12}{14.4}{\rmdefault}{\mddefault}{\updefault}{\color[rgb]{0,0,0}$\beta_2$}%
}}}}
\put(7167,-4627){\makebox(0,0)[lb]{\smash{{\SetFigFont{12}{14.4}{\rmdefault}{\mddefault}{\updefault}{\color[rgb]{0,0,0}$\alpha_2$}%
}}}}
\put(6448,-2621){\makebox(0,0)[lb]{\smash{{\SetFigFont{12}{14.4}{\rmdefault}{\mddefault}{\updefault}{\color[rgb]{0,0,0}$\alpha_1$}%
}}}}
\put(6817,-3215){\makebox(0,0)[lb]{\smash{{\SetFigFont{12}{14.4}{\rmdefault}{\mddefault}{\updefault}{\color[rgb]{0,0,0}$\beta_1$}%
}}}}
\end{picture}%

%% file: geometric.tex
\section{Geometric compactification of $X_D$}
\label{sec:geometric}

We have compactified the Hilbert modular surface $\X$ by taking its closure $\barX$ in
$\proj\Omega\barmoduli$; however, this compactification is not suitable for out purposes.  One problem is that
$\barX$ has non-normal singularities, and another is that several cusps of the curves $\W$ and $\P$ can come
together to one point in $\barX$.  It turns out that taking the normalization of $\Y$ solves both of these
problems and produces a useful compactification which has only quotient singularities.

\begin{definition}
  The \emph{geometric compactification} $\Y$ of the Hilbert modular surface $\X$ is the normalization of
  $\barX$.
\end{definition}

This defines $\Y$ as a variety.  We will also give $\Y$ the structure of a complex orbifold and of a
projective variety in \S \ref{subsec:curvesinY}.

In this section, we will study in detail the geometry of $\Y$.  In \S \ref{subsec:curvesinY}, we will study
curves in $\Y\setminus\X$, and their intersections with the curves $\barW$ and $\barP$.  We show in
\S\ref{subsec:Ymaps} that $\Y$ maps to the Baily-Borel compactification $\bX$ by a map which is the identity
on $\X$, and we draw conclusions about the cohomology of $\Y$.  In \S\ref{subsec:involutionofY}, we study the
extension of involution $\tau$ of $\X$ to $\Y$.  We will assume that the reader is familiar with the operation
of normalization.  For a summary of the results about normalization which we use in this section, see \S\ref{sec:normal}.

\subsection{Geometry of  $Y_D$}
\label{subsec:curvesinY}

The goal of this subsection is to classify the curves of $\Y\setminus\X$ and to understand the ways in which
these curves intersect each other and the curves $\barW, \barP\subset \Y$.  This is basically a
straightforward translation of the results of \S\ref{sec:limitsofeigenforms} on the local structure of $\barX$
to the setting of the normalization.

\paragraph{Local normalization of $\barX$.}

Normalization is a local operation: constructing the normalization of $\barX$ just amounts to replacing each
small open set with its normalization and then gluing together these normalizations.  We saw in
\S\ref{sec:limitsofeigenforms} that the singular points of $\barX$ are locally modeled on the varieties
$V(X^p-Y^q)$ and $V(X^p - Y^qZ^r)$ in $\cx^3$, so we need to understand the normalizations of these
varieties.  

Let
$\theta_m = e^{2\pi i/m}$.  The proof of the following proposition is easy and will be left to the reader.
\begin{prop}
  \label{prop:easynormalization}
  The normalization of the variety $V=V(X^p-Y^q)\subset\cx^3$ is $\gcd(p,q)$ copies of $\cx^2$ which map to $V$ by
  $$f_r(x,y)=(x^q, \theta_{\gcd(p,q)}^r x^p, y).$$
  It follows that $V$ has $\gcd(p,q)$ branches through each singular point $(0,0,Z)$.
\end{prop}

Given $p$, $q$, and $r\in\nats$ relatively prime, let
\begin{align}
  \label{eq:msn}
  m&=\frac{p}{\gcd(p,q)\gcd(p, r)} \\
  s&\equiv\left(\frac{r}{\gcd(p, r)}\right)^{-1} \mod \frac{p}{\gcd(p, r)}\notag\\
  n &= \frac{-q}{\gcd(p, q)} s,\notag
\end{align}
and define $Q(p,q,r)$ to be the normal analytic space $\cx^2 /G$,
where $G$ is the order $m$ cyclic group generated by the transformation,
$$(x, y)\mapsto (\theta_m x, \theta_m^n y).$$

\begin{lemma}
  \label{lem:kernelcyclic}
  If $p$ $q$ and $r$ are relatively prime integers, then the map
  $$f\colon\zed/\left(\frac{p}{\gcd(p, q)}\right)\times\zed/\left(\frac{p}{\gcd(p, r)}\right)\to\zed/p$$
  defined by $f(\alpha, \beta)=q\alpha+r\beta$ is surjective, and its kernel is a cyclic group of order $m$
  which is generated by $(\gcd(p, r), -q s)$, where $m$ and $s$ are as defined in \eqref{eq:msn}.
\end{lemma}

\begin{proof}
  It follows from the fact that $p$, $q$, and $r$ are relatively prime that $f$ is surjective, so the kernel
  $K$ has order $m$.  The element $(\gcd(p, r), -q s)$ clearly belongs to $K$ and has order at least $m$
  because $\gcd(p, r)$ has order $m$ in
  $$\zed/\left(\frac{p}{\gcd(p, q)}\right).$$
  Thus it is a generator of $K$.  
\end{proof}

\begin{prop}
  \label{prop:hardnormalization}
  The normalization of the variety $V=V(X^p-Y^qZ^r)\subset\cx^3$ is $p\colon Q(p,q,r)\to V$, where $p$ is
  induced by the map $\tilde{p}\colon \cx^2\to V$ defined by
  $$\tilde{p}(x, y) = 
  \begin{pmatrix}
    x^{q/\gcd(p, q)} y^{r/\gcd(q,r)}\\
    x^{p/gcd(p, q)} \\
    y^{p/\gcd(p,r)}
  \end{pmatrix}.
  $$
\end{prop} 

\begin{proof}
  For $i=1$ or $2$, let $(x_i, y_i)$ be a point in $\cx^2$ minus its coordinate axes.  We claim that
  $$\tilde{p}(x_1, y_1) = \tilde{p}(x_2, y_2),$$
  if and only if the $(x_i, y_i)$ are related by the action of
  $G$.  It would then follow that $p$ is the normalization of $V$ because $p$ would be well-defined,
  finite-to-one, and biholomorphic on the complement of the coordinate axes

  It is easy to show that $\tilde{p}(x_1, y_1) = \tilde{p}(x_2, y_2)$ if the $(x_i, y_i)$ are related by $G$.
  Conversely, suppose the $(x_i, y_i)$ have the same image.  Then we must have
  $$(x_2, y_2) = \left( e^{2\pi i\alpha\gcd(p, q)/p} x_1, e^{2\pi i\beta\gcd(p, r)/p} x_2\right)$$
  for some pair,
  $$(\alpha, \beta)\in \zed/\left(\frac{p}{\gcd(p, q)}\right)\times\zed/\left(\frac{p}{\gcd(p, r)}\right),$$
  such that
  $$q \alpha + r \beta \equiv 0 \mod p.$$
  By Lemma~\ref{lem:kernelcyclic} there is some $k\in\zed$ such that
  \begin{align*}
    \alpha&\equiv k\gcd(p,r) \mod\frac{p}{\gcd(p, q)} \\
    \beta&\equiv -k q s \mod \frac{p}{\gcd(p, r)},
  \end{align*}
  and it follows that the $(x_i, y_i)$ are related by an element of $G$.
\end{proof}

Since the normalizations of Propositions~\ref{prop:easynormalization} and \ref{prop:hardnormalization} are
either nonsingular or quotient singularities, we can regard $\Y$ as a complex orbifold.

\paragraph{One nonseparating polar node.}

We now begin to study the curves which make up $\Y \setminus \X$.  Define $\Sone[D]\subset\Y$ by
$$\Sone[D] = \pi^{-1}(\X(\system{2})).$$
By Proposition~\ref{prop:noonenode}, $\Sone[D]$ is empty unless $D$
is square.  By Corollary~\ref{cor:X100coordinates}, points of $\Sone$ correspond precisely to degree $d$
cylinder covering differentials; $\Sone$ is a suborbifold of $\Y[d^2]$ by Theorem~\ref{thm:X100structure}, and
$$\Sone\isom\half/\Gamma_1(d).$$
If $d>3$, then $\Sone$ is actually a submanifold of $\Y[d^2]$ because
$\Gamma_1(d)$ is torsion-free.

There are finitely many intersection points of $\barW[d^2]$ with $\Sone$, and these intersections are
transverse if $d>3$ by Theorem~\ref{thm:X100structure}.  These intersection points are exactly the
one-cylinder cusps of $\W[d^2]$.

\paragraph{One nonseparating holomorphic node.}

Define $\Stwo[D]\subset\Y$ by
$$\Stwo[D] = \pi^{-1}(\X(\system{2}, \system[1]{2})).$$
By Proposition~\ref{prop:noonenode}, $\Stwo[D]$ is
empty unless $D$ is square.  $\Stwo$ is a suborbifold of $\Y[d^2]$ isomorphic to $\half/\Gamma_1(d)$ by
Theorem~\ref{thm:X010structure}, and $\Stwo$ is a submanifold when $d>3$.  The points of $\Stwo$ correspond to
genus one differentials with two marked points which differ by exactly $d$-torsion identified to form a node.

$\SLtwoR$ acts on $\Stwo$ with discrete stabilizer, so we can regard $\Stwo$ as a degenerate Teichm\"uller curve.  Just
as for $\W[d^2]$, cusps of $\Stwo$ correspond to periodic directions on a given $(X, \omega)\in\Stwo$, and we
can divide the cusps of $\Stwo$ into one and two-cylinder cusps just as for $\W[d^2]$.

Given $q\in\nats$ such that $0<q\leq d/2$ and $\gcd(q,d)=1$, let $s_{q/d}\in\Y[d^2]$ be the inverse image of
$f_{q/d}\in\barX[d^2]$.  By Theorem~\ref{thm:X110classification}, these are transverse intersection points of
$\barSone$ and $\barStwo$, and these are all of their intersection points, because no other points in our
classification of points in $\barX[d^2]$ lie in this intersection.  These points $s_{q/d}$ are also
exactly the one-cylinder cusps of $\Stwo$.

\paragraph{Two or more nodes.}

For each curve $C_\lambda\subset\barX$, we associated a nondegenerate $\Y$-prototype $P$ to each branch of
$\barX$ through each point $p\in C_\lambda$. This means that we can associate a nondegenerate $\Y$-prototype
to each point in $\pi^{-1}(p)$ by Theorem~\ref{thm:pointsarebranches}.

For any nondegenerate $\Y$-prototype $P$, let $C_P^0\subset \Y$ be set of points in $\pi^{-1}(C_{\lambda(P)})$
associated to $P$.  The projection $C_P^0\to C_{\lambda(P)}$ is surjective and $\mult(P)$-to-one by
Corollary~\ref{cor:clambdabranches}.  The union of the curves $C_P$ parameterizes those points in $\Y$
representing differentials with two nonseparating nodes and possibly a separating node.

If $\lambda>1$, then $C_\lambda$ contains exactly one point $w_\lambda \in \barW$ which represents an Abelian
differential with a double zero by Proposition~\ref{prop:doublezerosinclambda}.  We associated a $\W$-prototype to each
branch of $\barX$ through $w_\lambda$, and we saw in Corollary~\ref{cor:clambdabranches} that $\W$-prototypes
$P$ such that $\lambda(P)=\lambda$ correspond bijectively to such branches.  Given any $\W$-prototype $P$, let
$w_P\in \Y$ be the point in $\pi^{-1}(w_{\lambda(P)})$ corresponding to the branch associated to $P$.  The
point $w_P$ is a two-cylinder cusp of $\W$, and these are all of the two-cylinder cusps.  This correspondence
between $\W$-prototypes and two-cylinder cusps of $\W$ is equivalent to the one described in
\cite{mcmullenspin}; more precisely, the point $w_P$ described here is the cusp $w_P$ of Theorem~\ref{thm:cuspstoprots}.

Similarly, for any $\P$-prototype $P$, define $p_P\in \Y$ to be the unique point in $\pi^{-1}(p_{\lambda(P)})$
coming from the branch associated to $P$.  These points $p_P$ are all of the cusps of $\P$.

We know that $\W$ and $\P$ are disjoint since they are $\SLtwoR$ orbits.  Their intersections with
$\Y\setminus\X$ are also disjoint because we have accounted for all of these intersection points above,
so it follows that $\barP\cap\barW=\emptyset$ in $\Y$.

Suppose $P$ is a $\W$ or a $\P$-prototype that maps to the $\Y$-prototype $Q$.  Then the point $w_P$ or $p_P$
is contained in $C_Q^0$.  This means that if $Q$ is a nonterminal prototype, then $C_Q$ contains $\mult(P)$
points of $\barW$ and the same number of points of $\barP$.  If $Q$ is a terminal prototype, then $C_Q$
contains a single point of $\barP$ (because $\mult(Q)=1$) and no point of $\barW$ because no $\W$-prototype
maps to a terminal $\Y$-prototype.

\begin{prop}
  The locus $C_P^0$ is a nonsingular curve in $\Y$. 
  The restriction $\pi\colon C_P^0\to C_{\lambda(P)}$ is an unramified cover of degree $\mult(P)$, and $C_P^0$
  is a twice-punctured sphere.

  Furthermore, $\barW$ and $\barP$ are  smooth suborbifolds of $\Y$ which intersect $C_P^0$ transversely.
\end{prop}

\begin{proof}
  Let $p\in C_P^0$, and $q=\pi(p)$.  By Theorem~\ref{thm:X200coordinates}, the branch $V$ of $\barX$ through $q$
  corresponding to $p$ is isomorphic to the germ at the origin of the variety in $\cx^3$ cut out by
  $$x^{a''} = \theta y^{-c''}$$
  for some root of unity $\theta$.  In these coordinates, $C_{\lambda(P)}^0$ is cut out by $x=y=0$.  If
  $q=w_{\lambda(P)}$ or $q=P_{\lambda(P)}$, then $\barW$ or $\barP$ is cut out by $z=0$.
  
  In coordinates $(u, v)$ on a neighborhood of $p$, the map
  $\pi\colon\Y\to\X$ is given by
  $$\pi(u, v) = (u^{-c''}, \theta'u^{a''}, v),$$
  for some root of unity $\theta'$.
  In these coordinates, $C_P$ is cut out by $u=0$, and if $p\in\barW$ or $\barP$, then $\barW$ or $\barP$ is
  cut out by $v=0$.  The smoothness and transversality statements follow immediately.

  Since every point in $C_{\lambda(P)}^0$ has $\mult(P)$ preimages in $C_P$, the projection is unramified of
  that degree.  Thus $C_P^0$ is a twice-punctured sphere because it is an unramified cover of one.

  Finally, $\barW$ and $\barP$ are suborbifolds because the coordinates above show that $\barW$ and $\barP$ are
  nonsingular where they cross the $C_P^0$.  These are the only points where $\barP$ and $\barW$ meet
  $\Y\setminus\X$ except for the intersections of $\barW$ with $\Sone[D]$ in the case when $D$ is square,
  and we have already seen that $\barW$ is nonsingular there in Theorem~\ref{thm:X100structure}.
\end{proof}

\begin{remark}
  The curves $C_P^0$ contain no orbifold points of $\Y$.  This is clear when $\lambda(P)>1$ because the image
  $C_{\lambda(P)}$ of $C_P^0$ in $\proj\Omega\barmoduli$ contains no orbifold points.  When $\lambda(P)=1$,
  the image $C_1$ is an orbifold locus of order two, so \emph{a priori} $C_P$ could be also.  This doesn't
  happen because for any $(X, [\omega])\subset C_1$, the group $\Aut(X, [\omega])$ fixes none of the branches
  of the inverse image of $\barX$ through the inverse image of $(X, [\omega])$ in $\proj\Omega\Def_2(S)$.
\end{remark}

We now study the closure of $C_P^0$, which we call $C_P$.  By Proposition~\ref{prop:Clambda}, $C_P\setminus
C_P^0$ consists of points representing differentials $c_\lambda$ with three nonseparating nodes, defined in
\S\ref{subsec:stratum5}.

By Theorems~\ref{thm:X300classification} and \ref{thm:branchintersection}, for every $\Y$-prototype $P$, there
is a point $c_{\lambda(P)}\in\barX$ and a branch of $\barX$ through $c_{\lambda(P)}$ labeled by $P$.  Let
$c_P\in\Y$ be the point corresponding to the branch labeled with $P$.  This gives a bijective correspondence
between the points of $\Y$ representing Abelian differentials with three nonseparating nodes and
$\Y$-prototypes.  Unless $P$ is terminal, all three nodes of $c_P$ are polar.

Recall that for a prototype $P=(a, b, c, \bar{q})$, we set
$$(a', b', c') = (a, b, c)/\gcd(a, b, c).$$
The following proposition describes the structure of $\Y$ in a neighborhood of $c_P$.
\begin{prop}
  \label{prop:aroundcP}
  A neighborhood of the point $c_P$, where $P=(a, b, c, \bar{q})$ is a nondegenerate, nonterminal $\Y$-prototype,
  is isomorphic to a neighborhood $U/G$ of the origin in the quotient singularity $Q(a', -c', -a'-b'-c')$.  In 
  coordinates $(u, v)$ on $U$, we have 
  \begin{align*}
    V(uv) &= (\Y\setminus\X)\cap U \\
    V(u) &= C_{P}\cap U\\
    V(v) &= C_{P^+}\cap U.
  \end{align*}

  If $P$ is a degenerate prototype, then $c_P$ is a smooth point of $\Y$.  In a neighborhood of $c_P$, the
  complement $\Y\setminus \X$ is contained in $\barSone[D] \cup C_{P^+}$, and these curves meet transversely at
  $c_P$.
  
  If $P$ is a terminal prototype, then $c_P$ is also a smooth point of $\Y$.  In a neighborhood of $c_P$, the
  complement $\Y\setminus \X$ is contained in $\barStwo[D] \cup C_{P}$, and these curves meet transversely at
  $c_P$.
\end{prop}

\begin{proof}
  These statements all follow directly from Theorems~\ref{thm:X300coordinates} and
  \ref{thm:branchintersection}, using Proposition~\ref{prop:hardnormalization} to translate these theorems
  into local coordinates on $\Y$.
\end{proof}

It follows from this proposition that the order of the orbifold point $c_P$ is
$$m_P = \frac{a'}{\gcd(a', c')\gcd(a', b'+c')} = \frac{a}{\gcd(a, c)\gcd(a, b+c)\gcd(a, b, c)}.$$

\begin{prop}
  The curve $C_P$ is a connected rational curve which meets the points $c_P$ and $c_{P^-}$ and no other point
  $c_Q$.  If $C_P$ is given the structure of an orbifold with orbifold points of order $m_P$ at $c_P$ and
  $m_{P^-}$ at $c_{P^-}$, then $C_P$ is a suborbifold of $\Y$.
\end{prop}

\begin{proof}
  The complement $C_P\setminus C_P^0$ is contained in the union of the points $c_Q$, and by
  Proposition~\ref{prop:aroundcP} it consists of exactly the points $c_P$ and $c_{P^-}$.  It follows from the
  coordinates in Proposition~\ref{prop:aroundcP} that $C_P$ is a suborbifold at these two points.

  The curve $C_P$ is rational because it is the closure of a twice-punctured sphere.
\end{proof}

When $D$ is not square, the operations $P^+$ and $P^-$ are defined for all $P$, so the curves $C_P$ are
divided into finitely many closed chains of rational curves.  When $D=d^2$, the curves $C_P$ are divided into
finitely many chains of rational curves joining $\barSone$ to $\barStwo$.

To summarize, we have established the following:
\begin{theorem}
  \label{thm:Ysummary}
  $\Y$ has the following properties:
  \begin{enumerate}
  \item $\Y$ is a compact, complex orbifold.  Its orbifold points are located at the elliptic points of $\X$,
    the elliptic points of $\Si[d^2]\isom \half/\Gamma_1(d)$, and the points $c_P$ for which $m_P>1$.
  \item The curves $C_P$, $\barW$, $\barP$, and $\barSi[D]$ are all suborbifolds of $\Y$.
  \item The curves $\barStwo[D]$, $\barW$, and $\barP$ are pairwise disjoint, and $\barSone[D]$ is disjoint
    from $\barP$.
  \item The curves $\barSone[D]$, $\barStwo[D]$, $\barW$, and $\barP$ meet the curves $C_P$ transversely.  These
    intersections correspond to two cylinder cusps of $\Stwo[D]$ and $\W$.  If $D=d^2$ with $d>3$, then
    $\barW[d^2]$ intersects $\Sone$ transversely.
  \item The point $w_P\in\barW$ lies on the curve $C_Q$, where $Q$ is the $\Y$-prototype associated to $P$,
    and $w_P$ lies on no other curve $C_{q'}$. This gives a bijective correspondence between the intersection
    points of $\barW$ with $C_Q$ and the set of $\W$-prototypes associated to the $\Y$-prototype $Q$.
    Similarly, the intersection points of $\barP$ with $C_P$ correspond bijectively in the same way to the set
    of $\P$-prototypes associated to the $\Y$-prototype $P$.

    If $P$ is neither initial nor terminal, then $C_P$ meets $\W$ and $\P$ in $\mult(P)$ points.

    If $P$ is terminal, then $C_P$ meets $\barStwo$ and $\barP$ once each and is disjoint from $\barW$.

    If $P$ is initial, then $C_P$ meets $\barSone$, $\barP$, and $\barW$ once each.
  \item The curves $\barSone$ and $\barStwo$ meet transversely  in the $N$ points $s_{q/d}$, where
  $$
  N=
  \begin{cases}
    1, & \text{if } d=2;\\
    \frac{1}{2}\phi(d), & \text{if } d>2
  \end{cases}
  $$
  (here $\phi$ is the Euler phi-function).
  \end{enumerate}
\end{theorem}

\begin{cor}
  If $D$ is not square, then $\W$ and $\P$ have the same number of cusps.  The number of
  cusps of $\P[d^2]$ is equal to the number of two-cylinder cusps of $\W[d^2]$ plus the number of two-cylinder
  cusps of $\Stwo$.
\end{cor}

\paragraph{Example: $\Y[17]$.} We will now illustrate these results in some specific cases.  Most of the
prototypes which arise will have $\bar{q}=0$, and we will abbreviate those prototypes by omitting $q$.

There are five $\Y[17]$ prototypes:
$$(1,1,-4),\quad (2, -3, -1),\quad (2, -1, -2),\quad (1, -3, -2),\quad (1, -1, -4).$$
They form one orbit
under the operation $P\mapsto P^+$, with each prototype being sent to the next on the list.  We have
$\mult(2,-1,-2)=2$, so $C_{(2,-1,-2)}$ meets $\barW[17]$ in the points $w_{(2, -1,-2, 0)}$ and $w_{(2, -1, 2,
  1)}$ and $\barP[17]$ in the points $p_{(2, -1,-2, 0)}$ and $p_{(2, -1, 2, 1)}$.  The other prototypes have
multiplicity one, so the other curves $C_P$ meet $\barW[17]$ and $\barP[17]$ once each.  Since $m_P=1$ for
each $P$, none of the points $c_P$ are orbifold points, and the $C_P$ meet each other transversely.

This configuration is shown in Figure~\ref{fig:Y17}.  The curves $C_P$ form the pentagon and are marked by
their prototypes.  The curves representing $\barW[17]$ and $\barP[17]$ in this figure accurately represent the connected
components of these curves and their intersections with the $C_P$.

\begin{figure}[htbp]
  \centering
  \input{Y_17.pstex_t}
  \caption{$Y_{17}$}
  \label{fig:Y17}
\end{figure}

\paragraph{Example: $\Y[12]$.}

There are three $\Y[12]$ prototypes:
$$(2, -2, -1), \quad (1, -2, -2), \quad (1, 0, -3).$$
They form one orbit under the operation $P\mapsto P^+$.  We have $m_{(2, -2, -1)}=2$, so the point $c_{(2, -2,
  -1)}$ is an orbifold point of order two in the intersection of $C_{(1, -2, -2)}$ and $C_{(2, -2, -1)}$.  The
other two points $c_P$ are nonsingular.  Each of these prototypes has multiplicity one, so $\W[12]$ and
$\P[12]$ each have three cusps and intersect each $C_P$ once.

\paragraph{Example: $\Y[25]$.}

There are nine $\Y[25]$-prototypes, which are divided into two orbits under $P\mapsto P^+$:
$$(1, -5, 0), \quad (1, -3, -4), \quad (1, -1, -6), \quad (1, 1, -6), \quad (1, 3, -4);$$
and
$$(2, -5, 0), \quad (2, -1, -3), \quad (2, -3, -2), \quad (2, 1, -3).$$
The prototypes after the first one on
each list are nondegenerate and correspond to seven curves $C_P$ in $\Y[25]$.  Each of these prototypes has
$m_P=1$, so none of the $c_P$ are orbifold points.  The curve $\W[25]$ has six two-cylinder cusps; it meets
$C_{(2, -3, -2)}$ twice and the other four curves $C_P$ corresponding to nonterminal prototypes once each.
$\W[25]$ also has two one-cylinder cusps (see \S\ref{sec:eulerWe}, \cite{mcmullenspin}, or
\cite{lelievreroyer} for formulas for the number of one-cylinder cusps) corresponding to two intersections
with $\Sone[25]$.  The curve $\barP[25]$ also meets $C_{(2, -3, -2)}$ twice and each of the other $C_P$ once
each.

Each of the two-cylinder cusps of $\W[25]$ can be represented by a square-tiled surface with five squares
having periodic horizontal direction.   Figure~\ref{fig:Y25} depicts the curves $C_P$ and $\Si[25]$ in $\Y[25]$.
Next to each $C_P$ are diagrams representing square tiled surfaces associated to each cusp of $\W[25]$ which
meets $C_P$.

\begin{figure}[htbp]
  \centering
  \input{Y_25.pstex_t}
  \caption{$Y_{25}$}
  \label{fig:Y25}
\end{figure}

\paragraph{Projective structure of $\Y$.}

So far, we have given $\Y$ the structure of an algebraic variety and a complex orbifold.   We can also give
$\Y$ the structure of a projective variety:

\begin{prop}
  $\Y$ is a projective variety.
\end{prop}

\begin{proof}
  Let $Z\subset\barmoduli$ be the closure of the image of the natural map $\X\to\barmoduli$.  The image of
  $\pi\colon\Y\to\barmoduli$ is contained in $Z$, so we have a finite surjective morphism $\Y\to Z$.  This
  means that $\Y$ is the normalization of $Z$ in $K(\Y)=K(\X)$.  By Theorem~\ref{thm:normalizationprojective},
  $\Y$ is a projective variety because $\barmoduli$ is projective.
\end{proof}

\subsection{$Y_D$ maps to $\widehat{X}_D$}
\label{subsec:Ymaps}

Recall that we introduced the Baily-Borel compactification $\bX$ of $\X$ in
\S\ref{sec:abelianvarieties}.  If $D$ is not square, $\bX$ consists of $\X$ together with a finite set $C$ of
cusps.  If $D=d^2$, then $\bX[d^2]$ consists of $\X$ together with the curves $\Rone$ and $\Rtwo$ isomorphic to
$\half/\Gamma_1(d)$ and the finite set $C$ of cusps.  Let $i\colon\X\to\Y$ be the natural inclusion.

\begin{theorem}
  \label{thm:YtoX}
  There is a unique morphism $p\colon \Y\to\bX$ such that $p\circ i = {\rm id}_{\X}$.  This map $p$ has
  the following properties:
  \begin{itemize}
  \item $p$ maps each curve $C_P$ onto a cusp of $\bX$, and this induces a bijection between the connected
    components of $\bigcup_P C_P$ and the set of cusps $C\subset \bX$.
  \item The restriction,
    $$p\colon \Y\setminus \bigcup_P C_P \to \X\setminus C$$
    is an isomorphism.
  \item $p$ maps $\Si$ isomorphically onto $\Ri$.
  \end{itemize}
\end{theorem}

\begin{proof}
  Recall from \S\ref{sec:abelianvarieties} that there is the natural morphism $j\colon\bX\to\satsiegelmod$
  with image $Z$ such that $j\colon\bX\to Z$ is the normalization of $Z$ in $K(\X)$.

  The varieties $\Y$ and $\bX$ both contain $\X$ as an open dense set, so we can identify both of their
  function fields with $K(\X)$.  The morphism $q\colon\Y\to\satsiegelmod$ which is the composition of the maps
  $$\Y\to\barmoduli\to\satsiegelmod$$
  has image $Z$.  By the universal property of normalization, Theorem~\ref{thm:universalproperty}, there is a
  map $p\colon\Y\to\bX$ such that $p^*$ is the identity on $K(\X)$ and $j\circ p = q$.  Since $p^*$ is the
  identity on $K(\X$), we must have $p\circ i= {\rm id}_{\X}$.

  The locus of points in $\Y$ representing  stable Riemann surfaces whose Jacobian is $(\cx^*)^2$ is the union
  $\bigcup C_P$.   Thus
  $q^{-1}(\siegelmod[0])=\bigcup C_P$.  We have $j^{-1}(\siegelmod[0])=C$, so $p^{-1}(C)=\bigcup C_P$.  By
  Theorem~\ref{thm:zariskismain}, Zariski's Main Theorem, the fiber of $p$ over each point in $C$ is
  connected, so $p$ induces a bijection between the components of $\bigcup C_P$ and the points of $C$.

  For each point $t\in\bX\setminus C$, the image $j(t)$ represents either an Abelian surface in $\siegelmod$
  or an elliptic curve in $\siegelmod[1]$  In the either case, there are only finitely many
  points in $\Y$ which map to $j(t)$, so the fiber $p^{-1}(t)$ is finite.  This fiber is also connected by
  Zariski's Main Theorem, so it is a single point.  This means that the restriction,
  $$p\colon \Y\setminus \bigcup_P C_P \to \X\setminus C,$$
  is a bijection, and it is thus an isomorphism because these are normal varieties.
  
  The locus of points $t$ in $\bX[d^2]$ such that $j(t)\in\siegelmod[1]$ is the union of $\Rone$ and $\Rtwo$,
  and the locus of points in $\Y[d^2]$ whose Jacobians are extensions of elliptic curves is the union of
  $\Sone$ and $\Stwo$.  This means that $p$ must take the curves $\Si$ to $\Ri$, and the restriction of $p$ to
  the $\Si$ is an isomorphism onto the $\Ri$ because $\Si$ is disjoint from each $C_P$.
\end{proof}

\paragraph{Cohomology of $\Y$.}

We can use Theorem~\ref{thm:YtoX} to get information about the cohomology of $\Y$.  In $H^2(\Y; \ratls)$, let
$B$ be the subspace generated by the fundamental classes of the rational curves $C_P$.

\begin{theorem}
  \label{thm:negativedefinite}
  The intersection form on $B$ is negative definite.
\end{theorem}

\begin{proof}
  Consider the map $q\colon \Y'\to\Y$ obtained by resolving the singularities at the points $c_P$ of $\Y$.
  Define $r\colon \Y'\to\bX$ by $r=p\circ q$, and let $E= r^{-1}(C)$.  The locus $E$ is the union of
  irreducible curves $E_i$.  By possibly taking some blow-ups of $\Y'$, we can assume that the curves $E_i$
  are nonsingular and transverse (actually, performing these blow-ups is not necessary; see the remark below).
   
  Let $B'\subset H^2(\Y';\ratls)$ be the subspace spanned by the fundamental classes of the $E_i$.  By
  Theorem~\ref{thm:resolutionnegativedefinite}, the intersection form on $B'$ is negative definite.
  
  The map $q^*\colon B\to B'$ preserves the intersection forms on $B$ and $B'$ because $q$ is degree one.  We
  claim that $q^*\colon B\to B'$ is injective.  To see this, define a linear map $r\colon B'\to B$ as follows.
  If $E_i\subset E$, and $q(E_i)$ is one of the curves $C_P$, let $r([E_i])$ be the fundamental class of that
  curve.  Otherwise, $q(E_i)$ is a point; in this case let $r([E_i])=0$.  We have $r\circ
  q^*={\rm id}_B$ on $B$, so $q^*$ is injective on $B$.

  Thus $q^*$ embeds $B$ in $B'$, preserving the intersection forms.  Since the intersection form on $B'$ is
  negative definite, the intersection form on $B$ must be as well.
\end{proof}

\begin{remark}
  We can actually describe the variety $\Y'$ above very explicitly.  The exceptional fiber of a resolution of
  a cyclic quotient singularity consists of a chain of rational curves $\{E_i\}_{i=1}^n$ with $E_i$
  having one transverse intersection with $E_{i+1}$ if $i \leq i \leq n-1$.  This is described in
  \cite{vandergeer88}.
  
  Thus resolving each $c_P\in\Y$ which is singular replaces each chain of curves $C_P$ with a longer chain of
  rational curves with transverse crossings.
\end{remark}

Since $\Y$ is an orbifold, Poincar\'e duality holds for rational cohomology, as Satake showed in
\cite{satake56}.  This together with Theorem~\ref{thm:negativedefinite} means that we can define the
orthogonal projection $\pi_B\colon H^2(\Y; \ratls)\to B$.

\begin{cor}
  \label{cor:boundaryrelation}
  If $D$ is not square, then
  \begin{equation*}
    \pi_B[\barP] = \pi_B[\barW].
  \end{equation*}
  as classes in $H^2(\Y; \ratls)$.  If $D=d^2$, then
  \begin{equation*}
    \pi_B[\barP[d^2]] = \pi_B[\barW[d^2]] + \pi_B[\barStwo[d^2]]
  \end{equation*}
\end{cor}

\begin{proof}
  For each curve $C_P$, the number of intersection points of $C_P$ with $\barP$ is equal to the number of
  intersection points of $C_P$ with $\barW$ and $\barStwo[D]$ by Part~4 of Theorem~\ref{thm:Ysummary} (using
  the convention that $\Stwo[D]$ is empty if $D$ is not square).  Since these intersection points are
  transverse, this implies
  $$[C_P]\cdot[\barP] = [C_P]\cdot([\barW] + [\barStwo[D]]).$$
  The desired equations then follow because the intersection matrix of the $C_P$ is nondegenerate.  
\end{proof}

\subsection{Involution of $Y_D$}
\label{subsec:involutionofY}

We now study the extension of the involution $\tau$ of $\X$ to $\Y$.   We then use this involution to deduce
information about the fundamental classes of the curves $\barWe$.

\begin{lemma}
  \label{lem:normalizationofgraph}
  Let $\Gamma\subset\X\times\barmoduli$ be the graph of the natural map $\X\to\barmoduli$, and let
  $\overline{\Gamma}$ be the closure of $\Gamma$ in $\bX\times\barmoduli$.  Then there is a morphism
  $q\colon\Y\to\overline{\Gamma}$ which is the normalization of $\overline{\Gamma}$.
\end{lemma}

\begin{proof}
  We have a map $q = p\times\pi\colon\Y\to\bX\times\barmoduli$ by Theorem~\ref{thm:YtoX}, and the image of $q$
  is in $\overline{\Gamma}$.  To show that $q\colon\Y\to\overline{\Gamma}$ is the normalization of
  $\overline{\Gamma}$, we need only to show that $q$ is finite and birational.  The map $q$ is birational
  because $p\colon\Y\to\bX$ is birational, as is the projection $\overline{\Gamma}\to\bX$.  The map $q$ is
  finite because the map $\pi\colon\Y\to\barmoduli$ is finite.
\end{proof}

\begin{theorem}
  The involution $\tau$ of $\X$ extends to an involution $\tau$ of $\Y$.  This involution satisfies
  \begin{equation}
    \label{eq:tauonSone}
    \tau(\Sone[D])=\Stwo[D],
  \end{equation}
  and
  \begin{equation}
    \label{eq:tauonCP}
    \tau(C_P) = C_{t(P)},
  \end{equation}
  where $t$ is the involution on the set of nondegenerate $\Y$-prototypes defined in \S\ref{subsec:prototypes}.
\end{theorem}

\begin{proof}
  Let $\overline{\Gamma}\subset \bX\times\barmoduli$ as in Lemma~\ref{lem:normalizationofgraph}.  The
  involution,
  $$\tilde{\tau} = \tau\times{\rm id}_{\barmoduli}\colon\bX\times\barmoduli\to\bX\times\barmoduli,$$
  restricts to an involution $\tilde{\tau}$ of $\overline{\Gamma}$.  By the universal property,
  Theorem~\ref{thm:universalproperty}, the composition $\tilde{\tau}\circ
  q\colon\Y\to\overline{\Gamma}$ lifts to a map $\tau\colon\Y\to\Y$ which extends $\tau$ on $\X$.  Since
  $\tau^2$ is the identity on the open, dense subset $\X$ of $\Y$, the map $\tau$ is an involution of $\Y$.

  The involution $\tau$ of $\bX$ sends $\Rone[D]$ to $\Rtwo[D]$.  Since the map $p\colon\Y\to\X$ commutes with
  the involutions, we must have $\tau(\Sone[D]) = \Stwo[D]$ by Theorem~\ref{thm:YtoX}.

  We are regarding $\X\subset\proj\Omega\barmoduli$ as the set of stable Abelian differentials $(X, [\omega])$, where
  $X$ has a compact Jacobian, and $\Jac(X)$ has real multiplication $\rho\colon\ord\to\End(\Jac(X))$ with
  $\omega$ as an $\iota_1$-eigenform defined up to constant multiple.  In these terms, $\tau(X, [\omega] )$ is
  the pair $(X, [\omega'])$, where $\omega'$ is an $\iota_2$-eigenform for the same real
  multiplication on $\Jac(X)$ ($\omega'$ is then an $\iota_1$-eigenform for the Galois conjugate real
  multiplication $\rho'$).

  Now suppose that $(X, [\omega])\in\X$ is close to a point in the curve $C_P^0$, and $(X, [\omega'])=\tau(X, [\omega])$ is close
  to a point in $C_{P'}^0$.  We must show that
  $P'=t(P)$.

  Recall the definition of the prototype $P$ from \S\ref{subsec:twononseppolarnodes}.  On $(X, \omega)$, there
  are two tall cylinders $C_1$ and $C_2$.  Let $\{\alpha_i, \beta_i\}_{i=1}^2$ be a symplectic basis of
  $H_1(X;\zed)$ such that $\alpha_i$ represents the cylinder $C_i$ and
  $\omega(\alpha_2)/\omega(\alpha_1)\geq1$.  Let $\mu\in K_D$ be such that $\mu\cdot\alpha_1=\alpha_2$,
  and let $\psi(t)=at^2+bt+c$ be as in Proposition~\ref{prop:admissibleclassification}, in particular $\psi(\mu)=0$.
  For $\lambda\in K_D$, let $T_\lambda$ be the matrix of the action of $\lambda$ on $H_1(X;\ratls)$ in the
  symplectic basis $\{\alpha_i, \beta_i\}$.  Then $T_{a\mu}$ is as in \eqref{eq:matrixofamu}.  The prototype
  is then $P=(a, b, c, \bar{q})$, where $a$, $b$, and $c$ are as above, and $\bar{q}$ is the reduction of the
  upper right entry of $T_{a\mu}$ modulo $\gcd(a, b, c)$.

  There are two cases to consider.  First suppose that $\mu^{(2)}\leq -1$, which is equivalent to
  $a-b+c\leq0$.  Let $\{\alpha_i', \beta_i'\}_{i=1}^2$ be the symplectic basis,
  $$\{\alpha_1', \alpha_2', \beta_1', \beta_2'\} = \{\alpha_1, -\alpha_2, \beta_1, -\beta_2\},$$
  of $H_1(X;\zed)$.  Let $\tilde{\mu}=-\mu$, and normalize $[\omega']$ so that $\omega'(\alpha_1')=1$.  Then
  $\tilde{\mu}\cdot\alpha_1'=\alpha_2'$, and
  $$\omega'(\alpha_2') = \omega'(\tilde{\mu}\cdot\alpha_1') =
  \tilde{\mu}^{(2)}\omega'(\alpha_1')=-\mu^{(2)}\geq1.$$
  This means that we can use the symplectic basis
  $\{\alpha_i', \beta_i'\}$ to calculate $P'$.  If $\psi'(t)=a't^2+b't+c'$ is the polynomial from
  Proposition~\ref{prop:admissibleclassification} with $\psi'(\tilde{\mu})=0$, then
  $$\psi'(t) = \psi(-t) = a t^2 - bt +c,$$
  so $(a', b', c')= (a, -b, c)$.  For $\lambda\in K_D$, let $T_\lambda'$ be the matrix $T_\lambda$ in
  the new basis $\{\alpha_i, \beta_i\}$.  Then a simple calculation shows,
  \begin{equation*}
     T_{a'\tilde{\mu}}' 
    =  T_{-a\mu}' 
    = 
    \begin{pmatrix}
      0 & -c & \phantom{-}0 & \phantom{-}q \\
      a & \phantom{-}b & -q & \phantom{-}0 \\
      0 & \phantom{-}0 & \phantom{-}0 & \phantom{-}a \\
      0 & \phantom{-}0 & -c & \phantom{-}b
    \end{pmatrix}.
  \end{equation*}
  Thus $P'= (a, -b, c, \bar{q}) = t(P)$.

  Now suppose that $-1<\mu^{(2)}<0$, which is equivalent to
  $a-b+c>0$.  Let $\{\alpha_i', \beta_i'\}_{i=1}^2$ be the symplectic basis,
  $$\{\alpha_1', \alpha_2', \beta_1', \beta_2'\} = \{\alpha_2, -\alpha_1, \beta_2, -\beta_1\},$$
  of $H_1(X;\zed)$,  let $\tilde{\mu}=-1/\mu$, and normalize $\omega'$ so that $\omega'(\alpha_1')=1$.  Then
  $\tilde{\mu}\cdot\alpha_1'=\alpha_2'$, and
  $$\omega'(\alpha_2') =-\mu^{(2)}>1$$
  as before.  If $\psi'(t)=a't^2+b't+c'$ is the polynomial from
  Proposition~\ref{prop:admissibleclassification}  with $\psi'(\tilde{\mu})=0$, then
  $$\psi'(t) = -t^2\psi(-1/t) = -c t^2 + bt -a,$$
  so $(a', b', c')= (-c, b, -a)$.  For $\lambda\in K_D$, let $T_\lambda'$ be the matrix $T_\lambda$ in
  the new basis $\{\alpha_i, \beta_i\}$.  Then,
  \begin{equation*}
    T_{a' \tilde{\mu}} = T_{c/\mu} = ac(T_{a\mu})^{-1} 
    =
    \begin{pmatrix}
      -b & c &  \phantom{-}0 & \phantom{-}q \\
      -a & 0 & -q & \phantom{-}0 \\
      \phantom{-}0 &  0 &  -b & -a  \\
      \phantom{-}0 &  0 & \phantom{-}c & \phantom{-}0
    \end{pmatrix},
  \end{equation*}
  and in the new basis,
  \begin{equation*}
    T_{a'\tilde{\mu}}' =
    \begin{pmatrix}
      0 &  \phantom{-}a & \phantom{-}0 &  \phantom{-}q \\
      -c & -b & -q &  \phantom{-}0\\
      0 &  \phantom{-}0 & \phantom{-}0 & -c \\
      0 &  \phantom{-}0 & \phantom{-}a & -b
    \end{pmatrix}.
  \end{equation*}
  Thus $P'= (-c, b, -a, \bar{q}) = t(P)$.
\end{proof}

It is not true that $\tau(\Wzero) = \Wone$, but we will prove that something like this
is true on the level of cusps.

\begin{theorem}
  \label{thm:cuspsbyspin2}
  Suppose $D\equiv 1 \pmod 8$.  If $D$ is not square, then for any nondegenerate $\Y$-prototype,
  \begin{equation}
    \label{eq:cuspsbyspin1}
    \#(\barWone \cap C_P) = \#(\barWzero \cap \tau(C_P)).
  \end{equation}
  If $D=d^2$, then
  \begin{equation}
    \label{eq:cuspsbyspin2}
    \#(\barWone[d^2] \cap C_P) + \#(\barStwo\cap C_P) = \#(\barWzero[d^2] \cap \tau(C_P)).
  \end{equation}
\end{theorem}

\begin{proof}
  Suppose $P$ is a nondegenerate, nonterminal, noninitial $\Y$-prototype.  Recall that the cusps of
  $\W$ lying on $C_P$, where $P=(a, b, c, \bar{q})$, are the points $w_Q$, where $Q=(a, b, c, \bar{r})$ is a
  $\W$-prototype with $\bar{r}\in\zed/\gcd(a, c)$, and $\bar{r}\equiv q \pmod {\gcd(a, b, c)}$.  The cusp $w_Q$
  lies on $\W^{\epsilon(Q)}$, where $\epsilon(Q)$ is given by \eqref{eq:epsilonofP}.
  
  There are two cases to consider, depending on whether or not $a-b+c\leq0$.  First suppose $a-b+c\leq0$.  In
  this case,
  $$\epsilon(t(Q)) \equiv \epsilon(Q) + b \equiv \epsilon(Q) + 1 \pmod 2$$
  for any $\W$-prototype $Q$ associated to $P$, using the fact that $b$ is odd because $b^2-4ac=D\equiv 1\pmod
  8$.  It follows that $C_P$ has as many cusps of $\Wone$ as $C_{t(P)}$ has cusps of $\Wzero$.
  
  Now suppose that $a-b+c>0$.  Here there are four cases to consider, depending on the parity of $a$ and $c$.
  First suppose $a\equiv 1 \pmod 2$, and $c\equiv 0 \pmod 2$.  In this case, for any $\W$-prototype $Q$
  associated to $P$, we have $\epsilon(Q)\equiv 0\pmod 2$ and $\epsilon(t(Q))\equiv 1 \pmod 2$.  Thus, every
  cusp of $\W$ on $C_P$ is in $\Wzero$, and $C_{t(P)}$ has the same number of cusps, all in $\Wone$.

  Suppose $a\equiv 0 \pmod 2$, and $c\equiv 1 \pmod 2$.  In this case, for any $\W$-prototype $Q$
  associated to $P$, we have $\epsilon(Q)\equiv 1\pmod 2$ and $\epsilon(t(Q))\equiv 0 \pmod 2$.  Thus, every
  cusp of $\W$ on $C_P$ is in $\Wone$, and $C_{t(P)}$ has the same number of cusps, all in $\Wzero$.
  
  Suppose $a\equiv c\equiv 0\pmod 2$.  In this case, for any $\W$-prototype $Q=(a, b, c,
  \bar{r})$ associated to $P$, we have $\epsilon(Q)\equiv \epsilon(t(Q)) \bar{r}\pmod 2$.  Since $\gcd(a, b,
  c)\equiv 1 \pmod 2$, this means that exactly half of the cusps of $\W$ on $C_P$ are on $\Wone$, and the same
  is true for $C_{t(P)}$.

  The last case, $a\equiv c \equiv 1 \pmod 2$ doesn't occur because $b^2-4ac\equiv 1\pmod 8$.  Thus we have
  shown \eqref{eq:cuspsbyspin1}.
  
  It remains to prove \eqref{eq:cuspsbyspin2} when $P$ is an initial or terminal prototype.  We claim that if
  $P$ is an initial $\Y[d^2]$-prototype, then $C_P$ intersects $\barWzero[d^2]$ once and is disjoint from
  $\barWone[d^2]$.  Let $P=(a, b, c, \bar{q})$ be an initial prototype.  Since $a-b+c=0$,
  $$\gcd(a, c) = \gcd(a, b, c),$$
  so $\mult(P)=1$, and $\barW[d^2]$ intersects $C_P$ once.  We have
  \begin{align*}
    &b = a+c \\
    &\Rightarrow (a+c)^2-4ac=d^2 \\
    &\Rightarrow a-c = d \\
    &\Rightarrow a \equiv c + 1 \pmod 2,
  \end{align*}
  and furthermore $f=d$.  It follows that $\epsilon(P)=0$, thus the intersection point of $C_P$ with
  $\barW[d^2]$ is in $\barWzero[d^2]$.

  We saw in Theorem~\ref{thm:Ysummary} that if $P$ is a terminal $\Y$-prototype, then $C_P$ intersects
  $\barStwo$ once and is disjoint from $\barW[d^2]$.  Thus if $P$ is a terminal prototype, then both sides of
  \eqref{eq:cuspsbyspin2} are one, and if $P$ is an initial prototype, then both sides of \eqref{eq:cuspsbyspin2}
  are zero.
\end{proof} 

\begin{cor}
  \label{cor:tauoncusps}
  If $D$ is not square, then
  \begin{equation}
    \label{eq:tauoncusps}
    \tau_* \pi_B[\barWzero] = \pi_B [\barWone]
  \end{equation}
  in $H^2(\Y; \ratls)$.  If $D=d^2$, then
  \begin{equation}
    \label{eq:tauoncuspssquare}
    \tau_* \pi_B[\barWzero[d^2]] = \pi_B[\barWone[d^2]] + \pi_B[\barStwo]
  \end{equation}
  in $H^2(\Y[d^2];\ratls)$.
\end{cor}

\begin{proof}
  Since the intersection pairing on $B\subset H^2(\Y; \ratls)$ is nondegenerate by
  Theorem~\ref{thm:negativedefinite}, it suffices to show that
  \begin{equation*}
    [C_P]\cdot(\tau_* \pi_B[\barWzero]) = [C_P]\cdot\pi_B [\barWone],
  \end{equation*}
  and
  \begin{equation*}
    [C_P]\cdot(\tau_* \pi_B[\barWzero[d^2]]) = [C_P]\cdot(\pi_B[\barWone[d^2]] + \pi_B[\barStwo])
  \end{equation*}
  for each $C_P$.  Since the intersections of $\barW$ and $\barStwo$ with each $C_P$ are transverse, this
  follows directly from Theorem~\ref{thm:cuspsbyspin2}.
\end{proof}


%% file: Y_17.pstex_t
\begin{picture}(0,0)%
\includegraphics{Y_17.pstex}%
\end{picture}%
\setlength{\unitlength}{3947sp}%
\begingroup\makeatletter\ifx\SetFigFont\undefined%
\gdef\SetFigFont#1#2#3#4#5{%
  \reset@font\fontsize{#1}{#2pt}%
  \fontfamily{#3}\fontseries{#4}\fontshape{#5}%
  \selectfont}%
\fi\endgroup%
\begin{picture}(5734,7439)(3301,-7419)
\put(8176,-2086){\makebox(0,0)[lb]{\smash{\SetFigFont{12}{14.4}{\rmdefault}{\mddefault}{\updefault}{\color[rgb]{0,0,0}$P_{17}$}%
}}}
\put(8476,-2536){\makebox(0,0)[lb]{\smash{\SetFigFont{12}{14.4}{\rmdefault}{\mddefault}{\updefault}{\color[rgb]{0,0,0}$W_{17}$}%
}}}
\put(3301,-2011){\makebox(0,0)[lb]{\smash{\SetFigFont{12}{14.4}{\rmdefault}{\mddefault}{\updefault}{\color[rgb]{0,0,0}$Y_{17}\setminus X_{17}$}%
}}}
\put(8101,-4111){\makebox(0,0)[lb]{\smash{\SetFigFont{12}{14.4}{\rmdefault}{\mddefault}{\updefault}{\color[rgb]{0,0,0}$(1, -3, -2)$}%
}}}
\put(5026,-2911){\makebox(0,0)[lb]{\smash{\SetFigFont{12}{14.4}{\rmdefault}{\mddefault}{\updefault}{\color[rgb]{0,0,0}$(2, -3, -1)$}%
}}}
\put(3901,-3586){\makebox(0,0)[lb]{\smash{\SetFigFont{12}{14.4}{\rmdefault}{\mddefault}{\updefault}{\color[rgb]{0,0,0}$(1, 1, -4)$}%
}}}
\put(7801,-3061){\makebox(0,0)[lb]{\smash{\SetFigFont{12}{14.4}{\rmdefault}{\mddefault}{\updefault}{\color[rgb]{0,0,0}$(2, -1, 2)$}%
}}}
\put(6151,-5161){\makebox(0,0)[lb]{\smash{\SetFigFont{12}{14.4}{\rmdefault}{\mddefault}{\updefault}{\color[rgb]{0,0,0}$(1, -1, -4)$}%
}}}
\end{picture}

%% file: Y_25.pstex_t
\begin{picture}(0,0)%
\includegraphics{Y_25.pstex}%
\end{picture}%
\setlength{\unitlength}{3947sp}%
\begingroup\makeatletter\ifx\SetFigFont\undefined%
\gdef\SetFigFont#1#2#3#4#5{%
  \reset@font\fontsize{#1}{#2pt}%
  \fontfamily{#3}\fontseries{#4}\fontshape{#5}%
  \selectfont}%
\fi\endgroup%
\begin{picture}(5349,8387)(2914,-8651)
\put(3376,-2236){\makebox(0,0)[lb]{\smash{\SetFigFont{12}{14.4}{\rmdefault}{\mddefault}{\updefault}{\color[rgb]{0,0,0}$(1, -3, -4)$}%
}}}
\put(4651,-2236){\makebox(0,0)[lb]{\smash{\SetFigFont{12}{14.4}{\rmdefault}{\mddefault}{\updefault}{\color[rgb]{0,0,0}$(1, -1, -6)$}%
}}}
\put(5926,-2161){\makebox(0,0)[lb]{\smash{\SetFigFont{12}{14.4}{\rmdefault}{\mddefault}{\updefault}{\color[rgb]{0,0,0}$(1, 1, -6)$}%
}}}
\put(7126,-2161){\makebox(0,0)[lb]{\smash{\SetFigFont{12}{14.4}{\rmdefault}{\mddefault}{\updefault}{\color[rgb]{0,0,0}$(1, 3, -4)$}%
}}}
\put(3451,-5161){\makebox(0,0)[lb]{\smash{\SetFigFont{12}{14.4}{\rmdefault}{\mddefault}{\updefault}{\color[rgb]{0,0,0}$(2, -1, -3)$}%
}}}
\put(5401,-5161){\makebox(0,0)[lb]{\smash{\SetFigFont{12}{14.4}{\rmdefault}{\mddefault}{\updefault}{\color[rgb]{0,0,0}$(2, -3, -2)$}%
}}}
\put(6976,-5161){\makebox(0,0)[lb]{\smash{\SetFigFont{12}{14.4}{\rmdefault}{\mddefault}{\updefault}{\color[rgb]{0,0,0}$(2, 1, -3)$}%
}}}
\put(3076,-4261){\makebox(0,0)[lb]{\smash{\SetFigFont{12}{14.4}{\rmdefault}{\mddefault}{\updefault}{\color[rgb]{0,0,0}$S_{25}^1$}%
}}}
\put(7726,-4186){\makebox(0,0)[lb]{\smash{\SetFigFont{12}{14.4}{\rmdefault}{\mddefault}{\updefault}{\color[rgb]{0,0,0}$S_{25}^2$}%
}}}
\end{picture}

%% file: bundles.tex
\section{Line bundles over $Y_D$}
\label{sec:bundles}

In this section, we will define extensions of some line bundles over $\X$ to $\Y$, and we will calculate their
Chern classes.  In \S\ref{subsec:orbifolvectorbundles}, we discuss background material on vector bundles over orbifolds.  In
\S\ref{subsec:goodmetrics}, we recall Mumford's  notion of a good Hermitian metric, which allows one to
express the Chern classes of a vector bundle in terms of curvatures of singular Hermitian metrics.  In \S\ref{subsec:extensionofbundle}, we
discuss the extension of $\Omega\X$ to $\Y$, and we discuss the extension of $T^*\A$ to $\Y$ in \S
\ref{subsec:extensionoffoliation}.

\subsection{Orbifold vector bundles}
\label{subsec:orbifolvectorbundles}

Most of the theory of vector bundles over manifolds works for vector bundles over orbifolds.  Here following
\cite{chenruan} we discuss the facts that we need and refer the reader to the Appendix to \cite{chenruan} for
details.

Let $X$ be a complex orbifold with an atlas $\{U_\alpha/ G_\alpha\to X\}$, where $U_\alpha\subset\cx^n$, and
$G_\alpha$ is a finite automorphism group of $U_\alpha$ (we allow $G_\alpha$ to have elements which act
trivially on $U_\alpha$, and $X$ is said to be \emph{reduced} if no such element exists).
A \emph{rank $r$
  holomorphic orbifold bundle over $X$} is an orbifold $E$ with a map $\pi\colon E\to X$ such that for each
$\alpha$,
$$\pi^{-1}(U_\alpha)\isom (U_\alpha\times \cx^r)/G_\alpha,$$
where $G_\alpha$ acts on $U_\alpha\times\cx^r$ by
$$g\cdot(z, v) = (g\cdot z, \rho(z,g)(v))$$
for some $\rho\colon U_\alpha\times G_\alpha\to{\rm GL}_r\cx$ which is a homomorphism if we fix $z$ and is
holomorphic if we fix $g$.

Nonreduced orbifold structures occur naturally.  For example, $\moduli$ has a nonreduced orbifold structure, where the local group $G$ at a point $X$ is
$\Aut(X)$.  This is nonreduced because the hyperelliptic involution acts trivially.  The associated reduced
orbifold structure has local group $G=\Aut(X)/J$ at $X$.  With this definition of an orbifold vector bundle,
$\Omega\moduli$ is only an orbifold vector bundle if we take the nonreduced orbifold structure on $\moduli$.
When necessary, we will implicitly use this orbifold structure.
Similarly, $\X$ has a nonreduced orbifold structure
with local group $G=\Aut(A, \rho)$ at $(A, \rho)$.

We want to be able to pull back bundles along maps between orbifolds.  The operation of pulling back a bundle
is much more complicated in the category of orbifolds then in the category of topological spaces.  In
particular, given a map $f\colon X\to Y$ of orbifolds, it is not always possible to pull back a bundle over
$Y$ by $f$; however, this is possible if $f$ is what is called a \emph{good} map.  We will informally discuss
pullbacks of orbifold bundles without actually defining a good map.

Given a map $f\colon X\to Y$ of orbifolds and charts $V/H\subset Y$ and $U/G\subset X$ such that
$f(U/G)\subset V/H$, there is a lift $\tilde{f}\colon U\to V$ such that for each $g\in G$ there is some $h\in
H$ such that
$$\tilde{f}(g\cdot z) = h\cdot\tilde{f}(z).$$
Note that there may not be a group homomorphism $\sigma\colon G\to H$ such that
$$\tilde{f}(g\cdot z) = \sigma(g)\tilde{f}(z).$$
If $f$ is good in the sense of \cite{chenruan}, then such a homomorphism does exist.  The actual definition of
a good map is some complicated global condition which locally amounts to the existence of such a homomorphism;
we will not get into this here.

In the notation of the previous paragraph, suppose that $f\colon U/G\to V/H$ is good and we have over $V/H$ a
vector bundle $E=(V\times\cx^r)/H \to V/H$, with $H$ acting on $V\times\cx^r$ by $\rho\colon V\times H\to{\rm
  GL}_r\cx$.  We define $f^* E$ over $U$ by
\begin{equation}
  \label{eq:localpullback}
  f^*E = (U\times\cx^r)/G,
\end{equation}
where $G$ acts on $U\times\cx^r$ by $g\cdot(z, v) = (g\cdot z, \tau(z, g)(v))$, with
$$\tau(z, g) = \rho(\tilde{f}(z), \sigma(g)).$$

Globally, the pullback of a bundle $E\to Y$ along a map $f\colon X\to Y$ is an orbifold vector bundle $f^*E$
which locally satisfies \eqref{eq:localpullback}.

Given an orbifold $X$, write $X_{\rm reg}$ for the set of $p$ points of $X$ such that for a chart $U/G\ni p$,
no element of $G$ fixes $p$.  The following theorem follows from Lemmas~4.4.3 and 4.4.11 of \cite{chenruan}.

\begin{theorem}
  \label{thm:pullbackdefined}
  Let $f\colon X\to Y$ be a map of orbifolds such that $f^{-1}(Y_{\rm reg})$ is an open, dense, and connected
  subset of $X$.  Then $f$ is a good map of orbifolds, and for any orbifold vector bundle $E\to Y$, there is a
  well-defined pullback bundle $f^*E\to X$.  The pullback bundle satisfies the naturality property,
  $$c_1(f^*E) =f^*c_1(E).$$
\end{theorem}

\begin{remark}
  Chern classes for orbifold vector bundles are defined using the Chern-Weil construction, which associates a
  de Rham cohomology class to a metric on an orbifold bundle.  This theory is worked out in \cite{chenruan}.
\end{remark}

\subsection{Good metrics}
\label{subsec:goodmetrics}

It is well known that the Chern classes of a complex vector bundle can be given in terms of the curvature of a
Hermitian metric.  In \cite{mumford77}, Mumford showed that for a Hermitian metric with certain mild
singularities, called a \emph{good metric}, the Chern classes of the vector bundle can still be given in terms
of the curvature of the metric.  Mumford's results are in terms of nonsingular complex projective varieties,
but his results and proofs work just as well in the setting of complex orbifolds.  In this section, following
\cite{mumford77}, we will summarize what we need about good metrics, translating to the setting of orbifolds.

Let $X\subset \overline{X}$ be complex $n$-dimensional orbifolds with $D=\overline{X}-X$ a divisor and
$\overline{X}$ compact.  Suppose that we can cover $D$ with coordinate charts $\Delta^n/G$, where
$\Delta^n\subset\cx^n$ is a polydisk, and $G$ is a finite group of automorphisms of $\Delta^n$ such that:
\begin{itemize}
\item Each transformation $g\in G$ is of the form,
$$g(z_1, \ldots, z_n) = (\theta_1 z_1, \ldots, \theta_n z_n),$$
for some roots of unity $\theta_i$.
\item $D\cap \Delta$ is a union of the coordinate axes $z_1=0, \ldots, z_r=0$ for $1\leq r \leq n$.
\end{itemize}

In such a neighborhood $\Delta^n/G$, we have $\Delta^n\cap X=(\Delta^*)^r\times\Delta^{n-r}$.  We give $\Delta^n\cap
X$ a metric by putting the Poincar\'e metric,
$$ds^2=\frac{|dz|^2}{|z|^2(\log|z|)^2},$$
on the $\Delta^*$ factors and putting the Euclidean metric $|dz|^2$ on the $\Delta$ factors.  Call the product
metric $\omega$.  This metric is $G$-invariant, so it is a metric on $\Delta^n/G$ in the orbifold sense.

\begin{definition}
  A $p$-form $\eta$ on $X$ has \emph{Poincar\'e growth} if there is a cover of $\overline{X}\setminus X$ by
  polycylinders $\Delta^n_\alpha/G_\alpha$ as above such that in each $\Delta^n_\alpha$ we have some constant
  $C_\alpha$ such that for any $p$ vectors $t_i$ at any point $x\in \Delta^n_\alpha\cap X$,
  \begin{equation}
    \label{eq:poincaregrowth}
    |\eta(t_1, \ldots, t_p)|^2 < C_\alpha \omega_\alpha(t_1, t_1) \cdots \omega_\alpha(t_p, t_p),
  \end{equation}
  where $\omega_\alpha$ is the metric on $\Delta^n_\alpha/G_\alpha$ defined above.
\end{definition}

In this paper, given a differential form $\eta$, we will write $\langle\eta\rangle$ for the current defined by
$\eta$, and $[\eta]$ for the de Rham cohomology class defined by $\eta$.

\begin{prop}
  \label{prop:poincaregrowth}
  A $p$-form $\eta$ with Poincar\'e growth satisfies,
  $$\int_X |\eta\wedge\zeta|<\infty$$
  for any smooth form $\zeta$ on $\overline{X}$ of complementary
  dimension.  Thus, $\eta$ defines a $p$-current $\langle\eta\rangle$.
\end{prop}

\begin{definition}
A $p$-form $\eta$ in $X$ is {\it good} if $\eta$ and $d\eta$ both have Poincar\'e growth.
\end{definition}

\begin{prop}
  \label{prop:dgood}
  If $\eta$ is a good $p$-form, then   
  $$d\langle\eta\rangle=\langle d\eta\rangle.$$
\end{prop}

Note that it if $\eta$ is a closed good form, it is not necessarily true that $\eta=d\zeta$ for a
good form $\zeta$; however, this will be the case if $\eta$ is the Chern form of a {\it good} Hermitian
metric.

Let $\overline{L}$ be a holomorphic orbifold line bundle over $\overline{X}$ with $L$ the restriction to $X$,
and let $H$ be a Hermitian metric on $E$.

\begin{definition}
  $H$ is a good metric if for every point $x\in D$, and every polycylindrical neighborhood
  $\Delta^n/G$ as described above with $D$ given by
  $$\prod_{i=1}^kz_i,$$
  and for every holomorphic section $e$ of $\overline{E}$ over $\Delta^n/G$, setting $h=H(e, e)$,
  we have
  \begin{itemize}
  \item
    $$|h|, |h|^{-1}<C\left(\sum_{i=1}^k \log|z_i|\right)^{2 m}$$
    for some $C>0$ and $m>1$.
  \item
    The 1-forms $(\partial h)/h$ are good on $\Delta^n/G \cap X$.
  \end{itemize}
\end{definition}

\begin{theorem}
  \label{thm:goodmetric}
  If $H$ is a good metric, then the Chern form, $c_1(E, H)$ is good and the current $\langle c_1(E, H)\rangle$ represents the
  Chern class $c_1(\overline{E})$ in $H^2(\overline{X}; \ratls)$.  Furthermore,
  $$c_1(E, H)- d\eta$$
  is a smooth 2-form on $\overline{X}$ for some good 1-form $\eta$.
\end{theorem}

The last statement follows from the proof of Theorem 1.4 in \cite{mumford77}.

We now show that Chern forms of good metrics behave as they should with respect to the cup product pairing and
integrating over curves.

\begin{prop}
  \label{prop:goodpairingone}
  Suppose $\dim_\cx \overline{X}=2$.  If $\omega$ and $\eta$ are Chern forms of two good metrics, then
  \begin{equation}
    \label{eq:pairing}
    \langle [\omega], [\eta]\rangle=\int_X\omega\wedge\eta.
  \end{equation}
\end{prop}

\begin{proof}
  By Theorem~\ref{thm:goodmetric}, we can write
  \begin{align*}
    \omega &= \omega' + d\omega'' \\
    \eta &= \eta' + d \eta'',
  \end{align*}
  where $\omega'$ and $\eta'$ are smooth, closed $2$-forms on $\overline{X}$, and $\omega''$ and $\eta''$ are
  good $1$-forms.  We have,
  \begin{align}
    \langle [\omega], [\eta]\rangle &= \int_{\overline{X}} \omega'\wedge \eta' \notag \\
    &= \int_X (\omega\wedge\eta - \omega \wedge d\eta'' - d \omega'' \wedge \eta + d\omega'' \wedge d \eta'')
    \notag \\
    &= \int_X (\omega\wedge\eta + d(\omega\wedge\eta'') - d(\omega''\wedge\eta) + d(\omega''\wedge d\eta'')). \label{eq:above}
  \end{align}

  We claim that $\int_X d\alpha=0$ for any good $3$-form $\alpha$.  The last three terms in \eqref{eq:above}
  are of this form, so \eqref{eq:pairing} follows from this claim.  By Proposition~\ref{prop:dgood},
  \begin{equation}
    \label{eq:intdalphazero}
    \int_X d\alpha = \langle d\alpha\rangle(1) = d \langle\alpha\rangle(1) = - \langle\alpha\rangle(d 1) = 0,
  \end{equation}
  so the claim follows.
\end{proof}

\begin{prop}
  \label{prop:goodpairingtwo}
  Suppose $C\subset\overline{X}$ is a curve with no irreducible component contained in $D$, and suppose
  $\omega$ is the Chern form of a good Hermitian metric. Then
  $$\langle [C], [\omega]\rangle=\int_C \omega.$$
\end{prop}

\begin{proof}[Sketch of proof]
  Let $\omega = \omega' + d \omega''$ with $\omega'$ a smooth 2-form on $\overline{X}$ and $\omega''$ a good
  1-form.  We have
  \begin{equation*}
    \langle [C], [\omega]\rangle = \int_C \omega' = \int_C \omega - \int_C d\omega''.
  \end{equation*}
  We need to show that
  \begin{equation}
    \label{eq:integralfinite}
    \int_C |\omega| < \infty
  \end{equation}
  and
  \begin{equation}
    \label{eq:integralzero}
    \int_C d\omega'' = 0.
  \end{equation}

  Let $C_0 = C\setminus C\cap D.$  We can speak of Poincar\'e growth for forms on $C_0$ using the
  compactification of $C_0$ obtained by adding a point to each cusp of $C_0$.  It follows from the Schwartz
  Lemma that any good $p$-form on $X$ restricts to a good $p$-form on $C_0$, so $\omega|_{C_0}$ and $\omega''|_{C_0}$
  are good.  Then \eqref{eq:integralfinite} follows from Proposition~\ref{prop:poincaregrowth}, and
  \eqref{eq:integralzero} follows from \eqref{eq:intdalphazero}.
\end{proof}

\subsection{Extension of the bundle $\Omega X_D$ to $Y_D$}
\label{subsec:extensionofbundle}

\paragraph{Extension of $\Omega\X$.}

Recall that we have the map $\pi\colon\Y\to\proj\Omega\barmoduli$ from Theorem~\ref{thm:YtoX}.  Let
$$\pi_1\colon\proj\Omega\barmoduli\to\barmoduli$$
be the natural projection, and let
$$\pi_2 = \pi_1\circ\pi\colon\Y\to\barmoduli.$$
Let $\Omega^0\Y = \pi_2^*(\Omega\barmoduli)$.  It does not follow directly from
Theorem~\ref{thm:pullbackdefined} that this pullback is well defined because we must take the nonreduced
orbifold structure on $\barmoduli$ for $\Omega\barmoduli$ to be a bundle, and then
$\pi_2^{-1}(\barmoduli)_{\rm reg}= \Y$.  Instead, first pull back the square $(\Omega\barmoduli)^{\otimes 2}$.
This is a bundle over $\barmoduli$ with the usual reduced orbifold structure, so
$\pi_2^{-1}({\barmoduli})_{\rm reg}$ consists of the elliptic points of $\X$ together with the elliptic points
of $\Si$ and the curves $C_P$.  This is an open, dense, and
connected subset of $\Y$, thus the pullback of the square is well defined by
Theorem~\ref{thm:pullbackdefined}.  Then define $\Omega^0\Y$ to be the quotient of this pullback by $\pm 1$.

Given $p\in \Y$ with $\pi_2(p)=X\in\barmoduli$, the fiber of $\Omega^0\Y$ over $p$ is isomorphic to
$\Omega(X)/G_p$ for some finite group $G_p$.  If $p$ is not an orbifold point of $\Y$, then $G=\{\pm 1\}$.  If
$p$ is an elliptic point of $\X$ or $\Si[D]$, then $G_p = \Aut(X)$.  If $p=c_P$ with $P$ a nondegenerate,
nonterminal $\Y$-prototype, then \textit{a priori} $G_p$ could be bigger than $\{\pm 1\}$ because $c_P$ is a
singular point of $\Y$.  In fact, $\pi(c_P)$ is a nonsingular point of $\proj\Omega\barmoduli$ because the residues
at the nodes of $\pi(c_P)$ are all different, so by \eqref{eq:localpullback}, $G_p \isom\{\pm 1\}$.  Thus we
have shown that the fiber of $\Omega^0\Y$ over $p$ is $\Omega(X)/\{\pm 1\}$ unless $p$ is an elliptic point of
$\X$ or $\Si[D]$.

Over $\proj\Omega\barmoduli$, there is the canonical orbifold line bundle $\mathcal{O}(-1)\to\proj\Omega\barmoduli$ whose
fiber over  $(X, [\omega])$ is the subspace of $\Omega(X)$ spanned by $\omega$.  Define
\begin{gather*}
  \Omega^1\Y = \pi^* \mathcal{O}(-1)\\
  \Omega^2\Y = \tau^*\Omega^1\Y.
\end{gather*}
Since $\mathcal{O}(1)\subset \pi_1^*\Omega\barmoduli$ as a sub line bundle, $\Omega^1\Y$ is a sub line bundle
of $\Omega^0\Y$.  The involution $\tau$ lifts to an involution of $\Omega^0\Y$, so we can also regard
$\Omega^2\Y$ as a sub line bundle of $\Omega^0\Y$.  Given $p \in \Y$ with $\pi_2(p) = X$, the Jacobian $\Jac(X)$
comes with a choice of real multiplication by $\ord$, and we can regard the fiber of $\Omega^i\Y$ over $p$ as
$\Omega^i(X)$, the space of $\iota_i$-eigenforms.  Since
$$\Omega(X) = \Omega^1(X)\oplus\Omega^2(X),$$
we obtain
\begin{equation}
  \label{eq:directsum}
  \Omega^0\Y = \Omega^1\Y\oplus\Omega^2\Y.
\end{equation}

Define
$$Q^i\Y = (\Omega^i\Y)^2,$$
which parameterizes quadratic differentials which are squares of $\iota_i$-eigenforms.  In what follows, we
will sometimes abbreviate $\Omega^1\Y$ and $Q^1\Y$ by $\Omega\Y$ and $Q\Y$.

\begin{prop}
  \label{prop:Qrestrictiontrivial}
  The restriction of $Q\Y$ to each curve $C_P\subset \Y$ is trivial.
\end{prop}

\begin{proof}[Sketch of proof]
  Define a global nonzero section of $Q\Y$ over $C_P$ as follows.  If $X\in C_P^0$, then let $q_X\in Q(X)$ be
  the unique quadratic differential which is the square of the Abelian differential $\omega_X\in\Omega(X)$
  which has residue $1$ at both nonseparating nodes of $X$.  Otherwise, $X$ has three nonseparating nodes.  In
  this case, define $q_X$ in the same way, using the two nodes of $X$ which are limits of nodes of surfaces in
  $C_P^0$.
\end{proof}

\paragraph{Chern class of $\Omega\Y$.}

In \S\ref{subsec:hilbertmodular}, we gave $\Omega\X$ a Hermitian metric $h_\Omega$.  In terms of Riemann
surfaces, the metric is given on the fiber over $X$ by
$$h_\Omega(\omega, \omega) = \int_X|\omega|^2.$$
The induced metric $h_Q$ on $Q\X$ is given in the fiber over $X$ by
$$h_Q(q, q) = \left(\int_X |q|\right)^2.$$
On $\Y$, these metrics become singular along the curves $C_P$ and $\barSone[D]$ because stable Abelian
differentials representing points on these curves have infinite area.

\begin{theorem}
  \label{thm:hqgood}
  The metric $h_Q$ on $Q\Y$ is a good metric.
\end{theorem}

\begin{lemma}
  \label{lem:goodcalculation}
  Let $H$ be a singular Hermitian metric on the trivial bundle $L=\Delta^n\times\cx$ over the polydisk
  $\Delta^n\in\cx^n$, and let $s$ be a holomorphic, nonzero section with
  $$H(s, s) = \left(\sum_{i=1}^r c_i \log|z_i| + c\right)^2,$$
  for some constants $c_i$ and $c$ with $c_i<0$.  Then $H$ is a good metric on $L$.
\end{lemma}

\begin{proof}
  Let $h = H(s, s)$.  Clearly,
  $$h, h^{-1} < O\left(\sum_{i=1}^r \log|z_i|\right)^2,$$
  so it remains to show that $\alpha=\partial h/h$ is a good 1-form.

  We have,
  $$\alpha = \frac{\sum_i\frac{c_idz_i}{z_i}}{\sum_i c_i\log|z_i|+c},$$
  so
  \begin{align*}
    \left|\alpha\left(\frac{\partial}{\partial z_j}\right)\right| &= \frac{c_j}{|z_j|\left|\sum_i c_i\log|z_i|
        + c\right|} \\
    & \leq \frac{1}{|z_j|\log|z_j|} \\
    & = \omega\left(\frac{\partial}{\partial z_j}, \frac{\partial}{\partial z_j}\right)^{1/2}.
  \end{align*}
  Thus $\alpha$ has Poincar\'e growth.  We have,
  $$d\alpha = \frac{1}{2}\frac{\sum_{ij}c_i c_j \frac{dz_i \wedge d\bar{z}_j}{z_i\bar{z}_j}}{\left(\sum c_i
      \log|z_i|+ c\right)^2},$$
  so
  \begin{align*}
    \left|d\alpha\left(\frac{\partial}{\partial z_k}, \frac{\partial}{\partial \bar{z}_l}\right)\right| &= \frac{1}{2}\frac{c_k c_l}{|z_k
      z_l|\left(\sum_i c_i \log|z_i| + c\right)^2} \\
    & \leq\frac{1}{4}\frac{1}{|z_k|\log|z_k|}\frac{1}{|z_l|\log|z_l|} \\
    &=\frac{1}{4}\omega\left(\frac{\partial}{\partial z_k}, \frac{\partial}{\partial
        z_k}\right)^{1/2}\omega\left(\frac{\partial}{\partial \bar{z}_l}, \frac{\partial}{\partial
        \bar{z}_l}\right)^{1/2}.
 \end{align*}
 Thus, $d\alpha$ also has Poincar\'e growth, so $\alpha$ is a good $1$-form, and $H$ is a good metric.
\end{proof}

\begin{proof}[Proof of Theorem~\ref{thm:hqgood}]
  We need to give holomorphic sections of $Q\Y$ around every point of $\barSone[D]$ and $C_P$ and show that
  the norms of these sections satisfy the required bounds.
  
  If $p\in\barSone$, then we saw in \S\ref{subsec:onenonpolarnode} and \S\ref{subsec:onenonpolandonenonhol}
  that $\barX[d^2]$ is normal around $p$, so we can work in $\barX[d^2]$.  We gave a section $s$ of
  $\mathcal{O}(-1)$ around $p$ whose norm is given by \eqref{eq:area1} or \eqref{eq:area2}.  The norm of $s^2$
  as a section of $\mathcal{O}(-2)$ is the square of the norm of $s$, and $s^2$ can be regarded as a section of $Q\Y$
  around $p$.  Thus $h_Q$ is good around each $p\in\barSone[d^2]$ by Lemma~\ref{lem:goodcalculation}.

  If $p\in C_P^0$, then we can consider $\pi(p)$ as a point in $\proj\Omega\Def^0_2(\system{3})$ or
  $\proj\Omega\Def^0_2(\system{4})$.  As we saw in \S\ref{subsec:twononseppolarnodes}, $\pi(p)$ is
  contained in a hypersurface $U_\lambda$ with coordinates $(x, y, z)$ on $U_\lambda$.  The variety
  $\barX$ is contained in $U_\lambda$ near $\pi(p)$ and is cut out by the equation,
  \eqref{eq:coordinates200}.  In a neighborhood of $p$ in $\Y$, there are coordinates $(u, v)$ with
  \begin{equation*}
    x = u,\:
    y = v^q, \: \text{and}\:
    z = \theta v^r,
  \end{equation*}
  for some $p, q\in \nats$ and root of unity $\theta$ by Proposition~\ref{prop:easynormalization}.  We defined
  a section $s$ of $\mathcal{O}(-1)$ over $U_\lambda$ with norm given by \eqref{eq:area3}.  The pullback of
  $s^2$ to $\Y$ gives a section $t$ of $Q\Y$ with norm,
  $$h_Q(t, t) = (c_1 \log|v| + c)^2,$$
  for some constants $c_1$ and $c$ with $c_1<0$.  Thus by Lemma~\ref{lem:goodcalculation}, $h_Q$ is good
  around every $p\in C_P^0$.

  Now suppose $p=c_P$ for some prototype $P$.  We can consider $\pi(p)$ as a point in $\proj\Omega\Def^0_2(\system{5})$ or
  $\proj\Omega\Def^0_2(\system[3]{5})$.  As we saw in \S\ref{subsec:threenonsepnodes}, $\pi(p)$ is
  contained in a hypersurface $U_\lambda$ with are coordinates $(x, y, z)$ on $U_\lambda$.  The variety
  $\barX$ is contained in $U_\lambda$ near $\pi(p)$ and is cut out by the equation,
  \eqref{eq:X300equation}.  In a neighborhood of $p$ in $\Y$, there are coordinates $(u, v)$ with
  \begin{equation*}
    x = u^p v^q, \:
    y = u^r, \:\text{and}\:
    z = v^s,
  \end{equation*}
  for some $p, q, r, s\in \nats$ by Proposition~\ref{prop:hardnormalization}.  We defined a section $s$ of
  $\mathcal{O}(-1)$ over $U_\lambda$ with norm given by \eqref{eq:area4}.  The pullback of $s^2$ to $\Y$ gives
  a section $t$ of $Q\Y$ with norm,
  $$h_Q(t, t) = (c_1 \log|u| + c_2\log|v|)^2,$$
  for some negative constants $c_1$ and $c_2$.  Thus by Lemma~\ref{lem:goodcalculation}, $h_Q$ is good
  around every $c_P$.
\end{proof}

\begin{cor}
  \label{cor:chernQ}
  The first Chern class of $Q^i\Y$ is
  \begin{equation}
    \label{eq:chernQ}
    c_1(Q^i\Y) = [\omega_i].
  \end{equation}
\end{cor}

\begin{proof}
  Since $h_Q$ is a good metric, $c_1(Q\Y, h_Q)$ represents $c_1(Q\Y)$ by Theorem~\ref{thm:goodmetric}.  We
  showed in Proposition~\ref{prop:chHQ} that $c_1(Q\Y, h_Q) = \omega_1$.  This shows that $c_1(Q^1\Y) =
  [\omega_1]$.  Also,
  \begin{equation*}
    c_1(Q^2\Y) = c_1(\tau^* Q^1\Y) = \tau^* c_1(Q^1\Y) = [\omega_2].
  \end{equation*}
\end{proof}

\begin{cor}
  \label{cor:pairingzero}
  For $i=1,2$, and for any nondegenerate $\Y$-prototype,
  $$[\omega_i]\cdot [C_P] =0.$$
\end{cor}

\begin{proof}
  This pairing is equal to the degree of $Q\Y$ restricted to $C_P$.  Since this restriction is trivial by
  Proposition~\ref{prop:Qrestrictiontrivial}, the degree is zero.
\end{proof}

\subsection{Extension of the foliation ${\cal A}_D$ to $Y_D$}
\label{subsec:extensionoffoliation}

\begin{prop}
  \label{prop:extensionofA}
  The foliation $\A$ of $\X$ extends to a foliation (which we will continue to call $\A$) of
  $$\Y\setminus\bigcup_P c_P,$$
  where the union is over all nonterminal $\Y$-prototypes.  The curves $C_P^0$ and $\Sone[D]$ are leaves of $\A$,
  and $\A$ is transverse to $\barStwo[D]$.
\end{prop}

\begin{proof}
  It follows directly from the equations for $\A$ given in Corollaries~\ref{cor:X100coordinates}, and
  \ref{cor:oneholnodecoords} that $\A$ extends over $\Sone[D]$ and $\Stwo[D]$ to a foliation which contains
  $\Sone[D]$ as a leaf and is transverse to $\Stwo[D]$.  The equations for $\A$ in
  Corollary~\ref{cor:oneholonepolnodecoords} show that $\A$ extends over the intersection points of $\Sone[D]$
  and $\Stwo[D]$ and is transverse to $\Stwo[D]$ there.  For $p\in C_P^0$, equations for $\A$ near
  $\pi(p)\in\barX$ are given in Theorem~\ref{thm:X200coordinates}.  There are local coordinates $(u, v)$ in a
  neighborhood of $p$ in $\Y$ such that $x=u$, $y=v^p$, and $z=\theta v^q$ for some $p, q\in\nats$ and root of
  unity $q$.  In these coordinates, $C_P$ is cut out by $v=0$, and $\A$ is given by $v=\const$.  Thus $\A$
  extends over $C_P^0$ to a foliation which contains $C_P^0$ as a leaf.
\end{proof}

We will now study $T^*\A$, the cotangent bundle to the leaves of $\A$.  Since every orbifold line bundle on a
complex orbifold minus a subvariety of codimension at least two extends to an orbifold line bundle over the
entire orbifold, we can regard $T^*\A$ as an orbifold bundle over all of $\Y$, even though the foliation is
singular at the points $c_P$.  The holomorphic sections of $\A$ over an open set $U$ are exactly the
holomorphic sections over $U$ minus the singular points of $\A$.

\begin{prop}
  \label{prop:singularfibers}
  The only singular fibers of $T^*\A$ are over the elliptic points of $\X$ and the elliptic points of $\Si$.
  For each nonterminal $\Y$-prototype $P$, the restriction of $T^*\A$ to the curve $C_P$ is trivial.  For each
  terminal $\Y$-prototype, the restriction of $T^*\A$ to the curve $C_P$ is isomorphic to $T^*C_P(c_{P^-})$.
\end{prop}

\begin{proof}
  Besides the elliptic points, the only points over which $T^*\A$ could have singular fibers are the singular
  points $c_P$ of $\A$.  Every section of an orbifold line bundle vanishes on a singular fiber, so to show
  that the fiber over $c_P$ is nonsingular, it suffices to show that $T^*\A$ has a nonzero section over
  $c_P$.

  Let $(x, y, z)$ be the coordinates on  in a neighborhood of $\pi(c_P)$ in the hypersurface $U_\lambda\subset
  \proj\Omega\Def_2(\system{5})$ as on p.~\pageref{page:coordinates30}.  By Theorem~\ref{thm:X300coordinates} and
  Proposition~\ref{prop:hardnormalization}, a neighborhood of $c_P$ in $\Y$ is of the form $\Delta^2/G$, where
  $\Delta^2$ is a polydisk with coordinates $(u, v)$ such that
  $$x= u^p v^q, \quad y=u^r, \quad z = v^s$$
  for some $p,q, r, s\in\nats$, and $G$ is a cyclic group whose action on $\Delta^2$ is generated by
  $$(u, v) = (\theta u, \theta^n v)$$
  for some root of unity $\theta$.  By the equation,
  $$y^\lambda z^{\lambda-1},$$
  for the foliation $\A$ from
  Theorem~\ref{thm:X300coordinates}, $\A$ is given on $\Delta^2$ by the equation
  $$u^a v^b = {\rm const}$$
  for some $a, b\geq0$.  If $P$ is nonterminal, then $\lambda>1$, and $a, b>0$.  Otherwise $\lambda=1$, and
  $b=0$, so the foliation is nonsingular at $c_P$.  Define a vector field $X$ on $\Delta^2$ by
  $$X = b u \frac{\partial}{\partial u} - a v \frac{\partial}{\partial v}.$$
  $X$ is tangent to $\A$ and is invariant under $G$, so it gives a nonzero section of $T^*\A$ in a
  neighborhood of $c_P$.  Thus the fiber of $T^*\A$ over $c_P$ is not singular.

  The restriction of $X$ to $C_P$ is a vector field on $C_P$ in a neighborhood of $c_P$ which vanishes at
  $c_P$.  That means that a vector field $Y$ on $C_P$ defines a nonzero section of $T^*\A|_{C_P}$ if and only
  if $Y$ vanishes at $c_P$.  If $P$ is nonterminal, then $C_P$ has two points $c_P$ and $c_{P^-}$ which pass
  through singular points of $\A$.  A holomorphic vector field on $C_P$ which has a zero at each of these
  points and no other zeros determines a nonzero holomorphic section of $T^*\A|_{C_P}$.  Thus this restriction
  is trivial.

  If $P$ is terminal, then $c_P$ is a nonsingular point of $\A$.  Thus,  the only singular
  point of $\A$ which $C_P$ intersects is $c_{P^-}$, so
  $$T^*\A|_{C_P}\isom T^*C_P(c_{P^-}).$$
\end{proof}

\begin{cor}
  \label{cor:chernpairing}
  The pairings of $c_1(Y^*\A)$ with the fundamental classes $[C_P]$ are
  \begin{equation*}
    c_1(T^*\A)\cdot[C_P] =
    \begin{cases}
      0 & \text{if $P$ is not terminal;}\\
      -1 & \text{if $P$ is terminal.}
    \end{cases}
  \end{equation*}
\end{cor}

\begin{proof}
  The pairing $c_1(T^*\A)\cdot[C_P]$ is equal to the degree of the restriction of $T^*\A$ to $C_P$.  If $P$ is
  nonterminal, then this restriction is trivial, so the pairing is zero.  If $P$ is terminal, then this
  restriction has degree $-1$ by Proposition~\ref{prop:singularfibers}.
\end{proof}

\paragraph{Chern class of $T^*\A$ when $D$ is not square.}

We now calculate the first Chern class of $T^*\A$ by relating it to the bundle $Q^2\Y$.

\begin{lemma}
  \label{lem:sectionmeromorphic}
  If $D$ is not square, and $\mathcal{L}\to\Y$ is an orbifold line bundle with a nonzero section
  $s\in\Gamma(\X, \mathcal{L})$ defined over $\X$, then $S$ is a meromorphic section of $\mathcal{L}$ over $\Y$.
\end{lemma}

\begin{proof}
  Let $\mathcal{O}(\mathcal{L})$ be the sheaf of sections of $\mathcal{L}$.  By \cite{serregaga},
  $\mathcal{O}(\mathcal{L})$ is a coherent algebraic sheaf on $\Y$. We claim that $s$ is defined algebraically
  over $X$.  It would follow from this because every algebraic section of a line bundle over a dense,
  Zariski-open subset is meromorphic by \cite[Lemma~II.5.3]{hartshorne}.
  
  The direct image $\pi_*\mathcal{O}(\mathcal{L})$ is a coherent algebraic sheaf on $\bX$.  Given an open set
  $U\subset\bX$ and $t\in\Gamma(U, \pi_*\mathcal{O}(\mathcal{L}))$, let $r=t/s$, a holomorphic function on
  $U\setminus C$, where $C=\bX\setminus\X$ is the set of cusps of $\X$.  By Koecher's principle (see
  \cite{vandergeer88}), $r$ extends uniquely to a holomorphic function $\tilde{r}$ on $U$.  This defines an
  injective analytic map of sheaves,
  $$i\colon\pi_*\mathcal{O}(\mathcal{L})\to\mathcal{O}_{\bX},$$
  with $i(t) = \tilde{r}$.  By \cite{serregaga}, $i$ is actually defined algebraically.  Since $i(s) = 1$, and
  $1$ is algebraic, it follows that $s$ is an algebraic section of $\mathcal{L}$ over $\X$.
\end{proof}

\begin{theorem}
  When $D$ is not square,
  \begin{equation*}
    Q^2\Y\isom T^*\A.
  \end{equation*}
\end{theorem}

\begin{proof}
  Define
  $$\mathcal{L} = Q^2\Y\otimes(T^*\A)^{-1}.$$
  By Propositions~\ref{prop:chHQ} and \ref{prop:c1Tstar},
  $$Q^2\Y|_{\X} \isom L_2 \isom T^*\A|_{\X},$$
  so there is a nonzero, holomorphic section $s$ of $\mathcal{L}$ over $\X$, which is a meromorphic section of
  $\mathcal{L}$ over $\Y$ by Lemma~\ref{lem:sectionmeromorphic}.

  The divisor of $s$ is
  $$(s) = \sum e_P C_P,$$
  so
  \begin{equation*}
    c_1(\mathcal{L}) = \sum e_P [C_P].
  \end{equation*}
  We need to show that the $e_P$ are all zero, for then $s$ would be nonzero and holomorphic, so $\mathcal{L}$
  would be trivial.  By Propositions~\ref{prop:Qrestrictiontrivial} and \ref{prop:singularfibers}, the
  restriction of $\mathcal{L}$ to each $C_P$ is trivial, so
  $$[C_Q]\cdot\sum e_P[C_P] = [C_Q]\cdot c_1(\mathcal{L}) = 0,$$
  for each $\Y$-prototype $Q$.  It then follows from Theorem~\ref{thm:negativedefinite} that $e_P=0$ for each $P$.
\end{proof}

\begin{cor}
  \label{cor:chernTstarA}
  If $D$ is not square, then the first Chern class of $T^*\A$ is
  \begin{equation}
    \label{eq:chernTstarA}
    c_1(T^*\A) = [\omega_2].
  \end{equation}
\end{cor}

\paragraph{Chern class of $T^*\A[d^2]$.}

The foliation $\A[d^2]$ of $\X[d^2]$ extends to a foliation $\bA$ of the orbifold $\bX[d^2]$.  Give $T\bA$ the
metric $\rho$ induced by the hyperbolic metric along the leaves of $\A$.  This metric is singular along
$\barRtwo$.

\begin{theorem}
  \label{thm:rhogood}
  The metric $\rho$ is a good metric for $T\bA(-\barRtwo)$.
\end{theorem}

\begin{proof}
  Let $G\subset{\rm SL}_2\ord[d^2]$ be the cyclic subgroup generated by
  $$\left(
    \begin{pmatrix}
      1 & 0 \\
      0 & 1
    \end{pmatrix},
    \begin{pmatrix}
      1 & d \\
      0 & 1
    \end{pmatrix}
  \right).$$
  We then have a map,
  $$\half\times\half/G\isom \half\times\Delta^*\to\X[d^2]$$
  (where $\Delta^*$ is the punctured unit disk), which extends to an unramified map of orbifolds,
  $$p\colon\half\times\Delta\to \bX[d^2]$$
  sending $\half\times\{0\}$ to $\Rtwo$, and sending leaves of the vertical foliation by  disks to
  leaves of $\bA$.  The vector field,
  $$X = z_2 \frac{\partial}{\partial z_2},$$
  on $\half\times\Delta^*$ locally defines a nonzero, holomorphic section of $T\bA(-\barRtwo)$ around points
  in $\Rtwo$.  The norm of this section is
  $$\rho(X, X) = (\log|z_2|)^{-2},$$
  so by Lemma~\ref{lem:goodcalculation}, $\rho$ is good near points of $\Rtwo$ (note that the exponent of the
  metric in Lemma~\ref{lem:goodcalculation} is irrelevant because changing the exponent only changes
  $\partial h/h$ by a constant).

  The proof that $\rho$ is good near the cusps of $\bX$ is the same, except these points are covered by
  $\Delta\times\Delta$, rather then $\half\times\Delta$.
\end{proof}

\begin{cor}
  \label{cor:chernTstarAdsquared}
  The first Chern class of $T^*\A[d^2]$ is
  \begin{equation}
    \label{eq:chernTstarAdsquared}
    c_1(T^*\A[d^2]) = [\omega_2] - [\barStwo]
  \end{equation} 
\end{cor}

\begin{proof}
  By Theorem~\ref{thm:rhogood},
  $$c_1(T^*\bA) = [\omega_2] - [\barRtwo].$$
  Since $\pi^*(T^*\bA)$ is isomorphic to $T^*\A[d^2]$ over $\X[d^2]$,
  $$\pi^*(T^*\bA) = T^*\A[d^2]\left(\sum_P e_P C_P\right),$$
  where the sum is over all nondegenerate $\Y$-prototypes $P$.  Thus
  \begin{equation}
    \label{eq:c1withcp}
    \sum_P e_P[C_P] = [\omega_2] - [\barSone] -  c_1(T^*\A[d^2])
  \end{equation}
  because $\pi^*[\barRtwo] = [\barStwo]$.

  The pairing of $[C_P]$ with the right hand side of \eqref{eq:c1withcp} is trivial for all $P$ by
  Corollary~\ref{cor:chernpairing} together with the fact that from Theorem~\ref{thm:Ysummary} that $C_P$ is
  disjoint from $\barStwo[d^2]$ if $P$ is nonterminal, and $C_P$ intersects $\barStwo[d^2]$ in one transverse
  intersection if $P$ is terminal.  Theorem~\ref{thm:negativedefinite} then implies that the $e_P$ are all
  zero.  Thus $c_1(T^*\A[d^2])$ is as claimed.
\end{proof}


%% file: euler.tex
\section{Euler characteristic of $W_D$}

\label{sec:modular}

In this section, we calculate  $\chi(\W)$.  We  construct a meromorphic section of a line bundle over
$\Y$ which vanishes along $\barW$ and has simple poles along $\barP$ and $\barStwo[D]$.  This allows us to
relate $\chi(\W)$ to $\chi(\P)$, $\chi(\X)$, and $\chi(\Stwo[D])$.

In \cite{mcmullenhilbert}, McMullen defined a quadratic differential on the leaves of $\A$ in $\X$, and used
this to construct a transverse measure for the foliation of $\X$ by $\SLtwoR$ orbits.  This quadratic
differential was also studied in \cite{schmoll05} in the case when $D$ is square.  Here we recall the
construction of this quadratic differential from \cite{mcmullenhilbert} and define it on all of $\Y$.

\begin{theorem}
  \label{thm:quadraticdifferential}
  On each leaf $L$ of $\A$, there is locally a quadratic differential $q$ which has a simple zero on each point of
  $\barW\cap L$, has a simple pole on each point of $(\barP\cup\barStwo[D])\cap L$, and is elsewhere
  holomorphic and nonzero.
\end{theorem}

\begin{proof}
  For each $z\in L$, let $X_z$ be the corresponding Riemann surface with real multiplication.  Choose a
  basepoint $z_0\in L$, and choose some eigenform $\omega_{z_0}\in \Omega^1(X_{z_0})$.  There is a unique
  section $z\mapsto \omega_z$ of $\Omega\Y$ over $L$ such that the absolute periods of the $\omega_z$ are
  locally constant.  When $L\subset \X$, this follows from Proposition~\ref{prop:periodsconstant}.  If
  $L=C_P$, this is clear because the periods of $\omega_z$ are determined by the residues at the nonseparating
  nodes of $X_z$, and along $C_P$ these residues have constant ratio, and if $L=\Sone$, this is also clear
  because if $(X_z, \omega_z)$ is a cylinder covering differential, then the periods of $\omega_z$ are all
  rational multiples of the period around the node of $X_z$, so if this period is constant along $L$, then all
  of the periods are constant.

  If $z\in L\setminus(\barW\cup\barP\cup\barStwo[D])$, then $(X_z, \omega_z)$ has two distinct simple zeros.  Let
  \begin{equation}
    \label{eq:fdef}
    f(z) = \int_\gamma \omega_z,
  \end{equation}
  where $\gamma$ is some path joining the zeros of $\omega_z$.  This $f(z)$
  is a multivalued holomorphic function on $L\setminus(\barW\cup\barP\cup\barStwo[D])$ because the value of
  $f(z)$ depends on a choice of an oriented path joining the zeros of $\omega_z$.  Define
  $$q=(\partial f)^2.$$
  We claim that $q$ is a well-defined quadratic differential on $L$.  To see this, suppose we replace an
  oriented path $\gamma$ joining the zeros $p$ and $q$ of $\omega_z$ with a new path $\gamma'$ joining $p$ to $q$.  Since
  $\gamma-\gamma'$ is a closed path, this changes $f$ by an absolute period of $\omega_z$.  Since the absolute
  periods of $\omega_z$ are constant along $L$, replacing $\gamma$ with $\gamma'$ does not change $\partial
  f$.  If we replace $\gamma$ with $-\gamma$, the same path with the opposite orientation, this changes $f$ to
  $-f$, which does not affect $q$ because of the exponent.  Thus $q$ is well defined.
  
  We now identify the zeros and poles of $q$.  Since the absolute and relative periods give a system of local
  coordinates on the strata in $\Omega\barmoduli$, the relative periods give local coordinates on $L$ because
  the absolute periods are constant, so $q$ is holomorphic and nonzero on the complement of
  $\barP$, $\barW$, and $\barStwo[D]$.

  Suppose $z\in L\cap\barW$.  For $z\in\cx$, let $I(z)\in\cx$ be the segment joining $0$ to $z$.  For
  $w\in\Delta_\epsilon$, where $\Delta_\epsilon$ is some small $\epsilon$-ball around $0$ in $\cx$, let
  $$\phi(w) = (X_w, \omega_w) \#_{I(w^{3/2})},$$
  (this is the operation of splitting a double zero defined in \S\ref{subsec:flatgeometry}).  As in
  \S\ref{subsec:stratum1}, $\phi$ is a conformal isomorphism onto some neighborhood of $z$ in $\cx$.  In these
  coordinate, we have $f(w) = w^{3/2}$, so
  $$q = \frac{9}{4} w\, dw^2.$$
  Thus $q$ has a simple zero at $z$.

  Now suppose $z\in L\cap \barP$.  Then
  $$(X_z, \omega_z) = (X_1, \omega_1) \#(X_2, \omega_2),$$
  the one-point connected sum of two genus one differentials.  Define a conformal mapping $\phi\colon\Delta_\epsilon\to L$ by
  $$\phi(w) = (X_1, \omega_1) \#_{I(w^{1/2})}(X_2, \omega_2),$$
  taking a connected sum along $I(w^{1/2})$.
  Similarly, if $z \in L\cap\barStwo[D]$, then we can regard $(X_z, \omega_z)$ as an genus one differential
  with two points identified to a node.  Define a conformal mapping $\phi\colon\Delta_\epsilon\to L$ by
  $$\phi(w) = (X_z, \omega_z) \#_{I(w^{1/2})},$$
  taking a
  self-connected sum along $I(w^{1/2})$ as described in
  \S\ref{subsec:flatgeometry}.  In either case, we have in these
  coordinates $f(w)= w^{1/2}$, so
  $$q = \frac{1}{2} w^{-1} dw^2.$$
  Thus $q$ has a simple pole at $z$.
\end{proof}

This construction locally defines meromorphic sections of $(T^*\A)^2$, but does not give a global section.
The problem is that the definition of the quadratic differential $q$ on a leaf $L$ depended on a choice of
$\omega_{z_0}\in\Omega(X_{z_0})$ for some basepoint $z_0\in L$.  There is no obvious way to choose these
quadratic differentials coherently to get a global section.  To get a global section of a bundle, we must
twist by some power of $\Omega\Y$.

\begin{theorem}
  \label{thm:meromorphicsection}
  There is a meromorphic section of the line bundle
  $$\mathcal{L} = (Q\Y)^* \otimes (T^*\A)^2$$
  on $\Y$ which has a simple zero along $\barW$, has a simple pole along $\barP$ and $\barStwo[D]$ (which is
  empty if $D$ is nonsquare), and is elsewhere nonzero and finite.
\end{theorem}

\begin{proof}
  The construction in the proof of Theorem~\ref{thm:quadraticdifferential} defined for each point,
  $$z\in \Y\setminus\bigcup c_P,$$
  together with a choice of $\omega\in\Omega^1(X_z)$ a quadratic differential $q$ on the leaf of $\A$ through
  $z$.  Thus, we have a map $h\colon \Omega^1(X_z)\to (T^*\A)^2|_z$.  With $f$ as in \eqref{eq:fdef},
  if we replace $\omega\in \Omega^1(X_z)$ with $a\omega$, then $f$ becomes $af$, and $q$ becomes $a^2 q$.
  Thus $h$ satisfies,
  $$h(a \omega) = a^2 h(\omega).$$
  We can thus regard $h$ as a linear map,
  $$\Omega^1(X_z)^{\otimes 2} \to (T^*\A)|_z.$$
  This defines a meromorphic section $s$ of $\mathcal{L}$ over
  $\Y\setminus\bigcup c_P$ which has the same zeros and poles as the quadratic differentials on the leaves of
  $\A$.  This section extends holomorphically over the points $c_P$ to give a section over all of $\Y$ because
  any holomorphic section of a line bundle defined on a neighborhood of a normal point extends over that
  point.
\end{proof}

\begin{cor}
  \label{cor:funWD}
  If $D$ is not square, then the fundamental class $[\barW]$ of $\barW$ in $H^2(\Y;\ratls)$ is given by
  \begin{equation}
    \label{eq:funWDone}
    [\barW] =[\barP] - [\omega_1] + 2 [\omega_2]. 
  \end{equation}
  The fundamental class $[\barW[d^2]]$ of $\barW[d^2]$ in
  $H^2(\Y[d^2]; \ratls)$ is given by
  \begin{equation}
    \label{eq:funWdtwoone}
    [\barW[d^2]] =  [\barP[d^2]] - [\barStwo] - [\omega_1] + 2[\omega_2].
  \end{equation}
\end{cor}

\begin{proof}
  By \cite[p.~141]{griffithsharris}, for any line bundle $L$ over a compact, complex manifold $M$ with a
  meromorphic section $s$ of $L$ having divisor $D$,
  $$c_1(L) = [D],$$
  This is proved for complex manifolds, but this is still true and the  proof works just
  as well for orbifolds.

  In our situation, the section $s$ of $\mathcal{L}$ from Theorem~\ref{thm:meromorphicsection} implies
  $$[\barW] - [\barP]- [\barStwo[D]] = c_1(\mathcal{L}) = -c_1(Q\Y) + 2 c_1(T^*\A).$$
  This, together with \eqref{eq:chernQ}, \eqref{eq:chernTstarA}, and \eqref{eq:chernTstarAdsquared} yield the
  desired formulas.
\end{proof}

\begin{cor}
  \label{cor:eulercharacteristic}
  If $D\neq1$ is a fundamental discriminant, then
  \begin{align}
    \label{eq:corWD}
    \chi(\W[f^2D]) &= \chi(\P[f^2D])-2\chi(\X[f^2D]) \\ 
    \label{eq:corWDa}
    &= -9 \zeta_{K_D}(-1) f^3\sum_{r|f}\kron{D}{r}\frac{\mu(r)}{r^2}.
  \end{align}
   If $D=d^2$, then
   \begin{align}
     \label{eq:corWd2}
     \chi(W_{d^2}) &= \chi(\P[d^2]) - \chi(\Stwo) -2 \chi(\X[d^2]) \\
     \label{eq:corWd2a}
     &= -\frac{1}{16} d^2(d-2)\sum_{r|d}\frac{\mu(r)}{r^2}.
  \end{align}
\end{cor}

\begin{proof}
  If we pair $-[\omega_1]$ with both sides of \eqref{eq:funWDone} and \eqref{eq:funWdtwoone}, then by
  Propositions~\ref{prop:goodpairingone} and \ref{prop:goodpairingtwo}, we get
  \begin{align}
    -\int_{\W}\omega_1 &= - \int_{\P}\omega_1 - 2\int_{\X}\omega_1\wedge\omega_2, \quad \text{and} \label{eq:intone}\\
    -\int_{\W[d^2]}\omega_1 &= - \int_{\P[d^2]}\omega_1  -
    2\int_{\X[d^2]}\omega_1\wedge\omega_2 + \int_{\Stwo[d^2]}\omega_1. \label{eq:inttwo}
  \end{align}
  Since $\W$, $\P$, and $\Stwo$ are transverse to $\A$, the form $-\omega_1$ restricts to the Chern form of
  the hyperbolic metric on these Riemann surfaces.  Thus by the Gauss-Bonnet theorem we get \eqref{eq:corWD}
  and \eqref{eq:corWd2}.  We calculated $\chi(\X)$ and $\chi(\P)$ in Theorems~\ref{thm:chiX} and \ref{thm:chiP}.  We
  will calculate $\chi(\Stwo[d^2])$ in the following proposition.  Putting all of this together yields
  \eqref{eq:corWDa} and \eqref{eq:corWd2a}.
\end{proof}

\begin{prop}
  We have,
  \begin{equation}
    \label{eq:chiStwo4}
    \chi(\Si[4]) = -\frac{1}{2},
  \end{equation}
  and
  \begin{equation}
    \label{eq:chiStwo}
    \chi(\Si) = - \frac{1}{12} d^2\sum_{r|d}\frac{\mu(r)}{r^2}
  \end{equation}
  when $d>2$.
\end{prop}

\begin{proof}
  We know that
  $$\Si \isom \half/\Gamma_1(d).$$
  It follows from \cite[Theorem~4.2.5]{miyake} that,
  \begin{equation*}
    \chi(\half/\Gamma_1(d)) =
    \begin{cases}
      -\frac{1}{2} & \text{if $d=2$;}\\
      -\frac{1}{12}\phi(d) d \prod_{p|d}\left(1+\frac{1}{p}\right) & \text{when $d>2$,}
    \end{cases}
  \end{equation*}
  where $\phi$ is the Euler $\phi$-function, and the product is over all primes dividing $d$.  By
  \cite[Proposition~2.2.5]{ir},
  $$\phi(d) = d \prod_{p|d}\left(1-\frac{1}{p}\right),$$
  so when $d>2$,
  $$\chi(\half/\Gamma_1(d)) = -\frac{1}{12} d^2\prod_{p|d}\left(1-\frac{1}{p^2}\right) = - \frac{1}{12}
  d^2\sum_{r|d}\frac{\mu(r)}{r^2},$$
  as claimed.
\end{proof}


%% file: fundamental.tex
\section{Fundamental class of $\overline{W}_D$}

\label{sec:fundamental}

In this section, we will calculate the fundamental class of $\barW$.  By Corollary~\ref{cor:funWD}, we just need
to know the fundamental classes of $\barSi[D]$ and $\barP$.

\begin{theorem}
  \label{thm:fundamentalSi}
  The fundamental class of $\barSi$ in $H^2(\Y[d^2];\ratls)$ is
  \begin{equation}
    \label{eq:fundamentalSdsquared}
    [\barSi] = \frac{6}{d}[\omega_i] + \pi_B[\barSi]
  \end{equation}
\end{theorem}

\begin{proof}
  From Theorem~\ref{thm:YtoX}, we have a map $p\colon\Y[d^2]\to\bX[d^2]$ which collapses the curves $C_P$ to
  cusps of $\bX[d^2]$.  Let
  $$q\colon\bX[d^2]\to\overline{\half/\SLtwoZ} = \barmoduli[1,1]$$
  be the map induced by the projection of $\half\times\half$, the universal cover of $\X[d^2]$, onto its
  second factor.  Let $r = q\circ p$, and let $\infty\in\barmoduli[1,1]$ by the single point added to
  $\half/\SLtwoZ$.  Then
  $$r^{-1}(\infty) = \Sone \cup \bigcup_P C_P.$$
  Let $f$ be a holomorphic function defined on a neighborhood of $\infty$ in $\barmoduli[1,1]$ which has a simple
  zero at $\infty$.  We claim that $f\circ r$ vanishes to order $d$ along $\Sone$.
  
  It is enough to show that $f\circ q$ vanishes to order $d$ along $\Rone[d^2]$.  Let $G\subset {\rm SL_2}
  \ord[d^2]$ be the cyclic group generated by
  $$\left(
    \begin{pmatrix}
      1 & 0 \\
      0 & 1
    \end{pmatrix},
    \begin{pmatrix}
      1 & d \\
      0 & 1
    \end{pmatrix}
  \right),
  $$
  and let $N\subset\SLtwoZ$ be the subgroup of upper-triangular matrices.  We have the following commutative diagram:
  $$\xymatrix{
    \half\times\half/G \ar[r] \ar[d] & \half/N \ar[d] \\
    \half\times\Delta \ar[r]_t \ar[d]_i & \Delta \ar[d]^j \\
    \bX[d^2] \ar[r]_q & \barmoduli[1,1]}$$
  where $\Delta\subset\cx$ is the unit disk.  The two topmost
  vertical maps are isomorphisms onto $\half\times\Delta^*$ and $\Delta^*$ respectively, where $\Delta^*$ is
  the punctured disk.  The maps $i$ and $j$ are unramified maps of orbifolds satisfying $i^{-1}(\Rone) =
  \half\times\{0\}$ and $j^{-1}(\infty) = \{0\}$.  The map $t$ is given by $t(z, w) = w^d$, and it follows that
  $f\circ q$ vanishes to order $d$ along $\Rone$ as claimed because
    $f\circ j \circ t$ vanishes to order $d$ along $\half\times\{0\}$.

  It follows from this claim that
  \begin{equation}
    \label{eq:funky1}
    r^*[\infty] = d[\barSone] + \sum_Pe_P[C_P],
  \end{equation}
  where $e_P$ is the order of vanishing of $f\circ r$ along $C_P$.  Let $\eta$ be the 2-form on
  $\half/\SLtwoZ$ induced by
  $$\frac{1}{2\pi}\frac{dx\wedge dy}{y^2}$$
  on $\half$.  Since $\chi(\half/\SLtwoZ)=-1/6$,
  $$\int_{\barmoduli[1,1]}\eta=\frac{1}{6},$$
  by the Gauss-Bonnet Theorem, so the cohomology class defined
  by the closed current $\langle\eta\rangle$ on $\barmoduli[1,1]$ satisfies
  \begin{equation*}
    6[\eta] = [\infty]
  \end{equation*}
  in $H^2(\barmoduli[1,1]; \ratls)$.  We have
  $$r^*\eta = \omega_1,$$
  so
  \begin{equation}
    \label{eq:funky2}
    r^*[\infty] = 6[\omega_1].
  \end{equation}
  Equations \eqref{eq:funky1} and \eqref{eq:funky2} imply \eqref{eq:fundamentalSdsquared} for $i=1$.  The case
  $i=2$ follows from the same argument or by applying $\tau^*$ to \eqref{eq:fundamentalSdsquared}.
\end{proof}

\begin{theorem}
  \label{thm:fundamentalP}
  If $D$ is not square, then the fundamental class of $\barP$ in $H^2(\Y;\ratls)$ is
  \begin{equation}
    \label{eq:fundamentalPD}
    [\barP]= \frac{5}{2}([\omega_1] + [\omega_2]) + \pi_B[\barP].
  \end{equation}
  The fundamental class of $\barP[d^2]$ in $H^2(\Y[d^2];\ratls)$ is
  \begin{equation}
    \label{eq:fundamentalPdsquared}
    [\barP[d^2]] = \left(\frac{5}{2} - \frac{3}{d}\right)([\omega_1] + [\omega_2])+ \pi_B [\barP[d^2]]
  \end{equation}
\end{theorem}

\begin{proof}
  In $\barmoduli$, let $\Delta_0$ be the divisor which is the closure of the locus of stable Riemann surfaces
  with one nonseparating node, and let $\Delta_1$  be the divisor which is the closure of the locus of stable Riemann surfaces
  with one separating node.  Let $\delta_i = [\Delta_i]$, the fundamental class of $\Delta_i$.  Define
  $$\lambda_1 = c_1(\Omega\barmoduli).$$
  These cohomology classes satisfy the well-known relation,
  \begin{equation}
    \label{eq:divisorrelation}
    \delta_1 = 5\lambda_1 - \frac{1}{2}\delta_0,
  \end{equation}
  proved in \cite{mumford83}.
  
  Since $\pi_2^*(\Omega\barmoduli) = \Omega^1\Y\oplus\Omega^2\Y$, where $\pi_2\colon\Y\to\barmoduli$ is the
  natural map, we have
  \begin{equation}
    \label{eq:pullbacklambda}
    \pi_2^*(\lambda_1) = \frac{1}{2}([\omega_1]+[\omega_2])
  \end{equation}
  by Corollary~\ref{cor:chernQ}.

  We claim that
  \begin{equation}
    \label{eq:pullbackdelta1}
    \pi_2^*(\delta_1) = [\barP].
  \end{equation}
  Since $\pi_2^{-1}(\Delta_1) = \barP$, it suffices to show that $\Delta_1$ is generically transverse to
  $\pi_2(\Y)$.  In $\siegelmod[2]$, the divisor $\Delta_1$ corresponds to the surface $\X[1]$ parameterizing products of
  elliptic curves.  In $\siegelmod[2]$, the intersection of $\X[1]$
  with $\X$ is transverse because these are both linear subspaces of $\siegelmod$, so if their intersection
  was not transverse, one would be contained in the other.  Equation \eqref{eq:pullbackdelta1} follows.

  We now claim that
  \begin{equation}
    \label{eq:pullbackdelta0}
    \pi_2^*(\delta_0) = [\barSone[D]] + [\barStwo[D]] + \sum_P e_P [C_P]
  \end{equation}
  for some integers $e_P$.  Since
  $$\pi_2^{-1}(\Delta_0) = \bigcup_{i=1}^2\barSi[D] \cup \bigcup_P C_P,$$
  $\pi_2\circ\tau = \pi_2$, and $\tau^*\Sone[D] = \Stwo[D]$, it
  suffices to show that $\Delta_0$ meets $\pi_2(\Y)$ transversely
  along $\pi_2(\Sone[D])$.
  
  Let $p\in\Sone$, and let $(X, \omega)\in\Omega\barmoduli$ be a corresponding eigenform.  We claim that
  $\pi_2(\Y)$ meets $\Delta_0$ transversely at $\pi_2(p)$.  Let $\{\alpha_i, \beta_i\}$ be a symplectic basis for
  $H_1(X;\zed)$ as in \S\ref{subsec:stratum2}.  Then we get coordinates $(v, w, x, y, z)$ on a neighborhood
  $U$ of $(X, \omega)$ in $\Omega\barmoduli$ as in \S\ref{subsec:stratum2}.  The subspace $H$ of $U$ defined
  by the equations,
  $$v=\omega(\alpha_1) \quad w=\omega(\alpha_2),$$
  maps locally biholomorphically to $\barmoduli$, so the coordinates $(x, y, z)$ on $H$ induce coordinates on
  $\barmoduli$ on a neighborhood $W$ of $X$.  In these coordinates,
  $$\Delta_0 = V(z).$$
  By Corollary~\ref{cor:X100coordinates},
  $$\pi_2(\Y)\cap W = V(x - \omega(\beta_2)).$$
  Thus $\pi_2(\Y)$ and $\Delta_0$ meet transversely at $p$ as claimed, and \eqref{eq:pullbackdelta0} follows.

  Equations \eqref{eq:fundamentalPD} and \eqref{eq:fundamentalPdsquared} then follow from
  \eqref{eq:divisorrelation}, \eqref{eq:pullbacklambda}, \eqref{eq:pullbackdelta1}, and \eqref{eq:pullbackdelta0}.
\end{proof}

\begin{remark}
  I am grateful to Gerard van der Geer for providing the idea of the proof of this theorem.
\end{remark}

\begin{remark}
  We can also use \eqref{eq:fundamentalPD} and \eqref{eq:fundamentalPdsquared} to get a new proof of
  Theorem~\ref{thm:chiP}.  Since this proof doesn't use Siegel's formula, Theorem~\ref{thm:siegelformula},
  this together with the previous proof of Theorem~\ref{thm:chiP} can be used to give a proof of Siegel's
  formula.  
\end{remark}

\begin{cor}
  \label{cor:fundamentalW}
  If $D$ is not square, then the fundamental class of $\barW$ in $H^2(\Y;\ratls)$ is given by
  \begin{equation}
    \label{eq:fundamentalW}
    [\barW] = \frac{3}{2} [\omega_1] + \frac{9}{2}[\omega_2] + \pi_B[\barW].
  \end{equation}
  The fundamental class of $\barW[d^2]$ in $H^2(\Y[d^2];\ratls)$ is given by
  \begin{equation}
    \label{eq:fundamentalWdsquared}
    [\barW[d^2]] = \frac{3}{2}\left(1 - \frac{2}{d}\right) [\omega_1] +
    \frac{9}{2}\left(1-\frac{2}{d}\right)[\omega_2] + \pi_B[\barW[d^2]].
  \end{equation}
\end{cor}

\begin{proof}
  This follows from plugging the formulas from Theorems~\ref{thm:fundamentalSi} and \ref{thm:fundamentalP}
  into the formulas from Corollary~\ref{cor:funWD}.
\end{proof}


%% file: normalbundles.tex
\section{Normal Bundles}
\label{sec:normalbundles}

We now study the normal bundles of the curves $\barW$, $\barP$, and $\barStwo[D]$ with the goal of calculating
the self-intersection numbers of these curves.  For any curve $C\subset \Y$, we will write $N(C)$ for its normal bundle.

\begin{prop}
  \label{prop:normalequalstangent}
  For any connected component $C$ of $\barW$, $\barP$, or $\barStwo[D]$,
  \begin{equation}
    \label{eq:normalequalstangent}
    N(C) \isom T\A|_C
  \end{equation}
  as holomorphic line bundles over $C$.
\end{prop}

\begin{proof}
  We claim that the foliation $\A$ is transverse to $\barW$, $\barP$, and $\barStwo[D]$.  Equation
  \eqref{eq:normalequalstangent} follows directly from this claim.
  
  The curves $\W$ and $\P$ are transverse to $\A$ because the inverse images of these curves in
  $\half\times\half$ are the unions of graphs of holomorphic functions $\half\to\half$.  The curve $\Stwo[D]$
  is also transverse to $\A$ by Proposition~\ref{prop:extensionofA}.  The closures of these curves are then
  transverse to $\A$ because they intersect $\Y\setminus\X$ in the curves $C_P$, these intersections are
  transverse by Theorem~\ref{thm:Ysummary}, and the curves $C_P$ are leaves of $\A$ by
  Proposition~\ref{prop:extensionofA}.
\end{proof}

\begin{theorem}
  \label{thm:selfintersections}
  For any connected component $C$ of $\P$ or $\Stwo[D]$,
  \begin{equation*}
    [\overline{C}]^2 = \chi(C).
  \end{equation*}
  For any component $C$ of $\W$,
  \begin{equation*}
    [\overline{C}]^2 = \frac{1}{3}\chi(C).
  \end{equation*}
\end{theorem}

\begin{proof}
  Let $C$ be a connected component of $\P$, $\Stwo[D]$, or $\W$, and choose a tubular neighborhood
  $U\subset\Y$ of $\overline{C}$ which is small enough that each point in $U\setminus \overline{C}$ represents a stable Abelian
  differential which has a unique shortest saddle connection joining distinct zeros.  We wish to define a map,
  \begin{equation*}
    \Phi\colon U\to(Q^1 \overline{C})^*,
  \end{equation*}
  where $Q^1 \overline{C}$ is the restriction of $Q^1\Y$ to $\overline{C}$.
  
  Consider $p\in U$ representing the projective class of a stable Abelian differential $(X, [\omega])$.  Let
  $I\subset X$ be the unique shortest saddle connection connecting distinct zeros, and let $(Y, \eta)$ be the
  stable Abelian differential obtained by collapsing $I$ as in \S\ref{subsec:flatgeometry}.  The projective
  class $(Y, [\eta])$ then represents a point of $\overline{C}$.  Since $\omega$ and $\eta$ are both
  eigenforms for real multiplication, we have an isomorphism,
  \begin{equation*}
    T\colon \Omega^1(Y)\to\Omega^1(X),
  \end{equation*}
  defined by $T(\eta) = \omega$.  Define $S\in Q^1(Y)^* = \Omega^1(Y)^{-2}$ by
  \begin{equation*}
    S(\nu) = \left(\int_I T(\sqrt{\nu})\right)^2,
  \end{equation*}
  where the integral along $I$ is with respect to some choice of orientation of $I$.  It doesn't matter which
  orientation we take for $I$ or which square root of $\nu$ we take, so $S$ is well-defined.  Now define
  $\Phi$ by
  \begin{equation*}
    \Phi(p) = ((Y, [\eta]), S).
  \end{equation*}
  For $p\in \overline{C}$, we define $\Phi(p) = (p, 0)$.  Note that $\Phi$ takes leaves of $\A$ to fibers of $(Q^1\overline{C})^*$.
  
  Suppose that $C$ is a connected component of $\P$.  We claim that in this case, $\Phi$
  is injective.  To see this, let $q\in(Q^1\overline{C})^*$ be represented by $(Y, [\eta])\in \proj\Omega\barmoduli$ and
  $S\in Q^1(Y)^*$.  This $(Y, \eta)$ is the one point union of two genus one differentials or cylinders:
  $$(Y, \eta) = (Y_1, \eta_1) \# (Y_2, \eta_2).$$
  Normalize $\eta$ so that $S(\eta^2)=1$.  If $\Phi(p) = q$,
  then $p$ is represented by the connected sum
  $$(X, \omega) = (Y_1, \eta_1) \#_I (Y_2, \eta_2),$$
  where $I\subset\cx$ is the segment joining $0$ to $1$.  Thus $p$ is determined uniquely by $q$, so $\Phi$ is
  injective as claimed.  When $C=\Stwo[D]$, then $\Phi$ is injective by the same argument, using a
  self-connected sum in place of the connected sum operation above.
  
  Recall that in \S\ref{subsec:flatgeometry}, we defined the operation of splitting a double zero, which is
  inverse to the operation of collapsing a saddle connection, and associates to a sufficiently small segment
  $I\subset(Y, \eta)$ starting at the zero of $\eta$ the Abelian differential $(Y, \eta)\#_I$.
  
  Suppose now that $C$ is a connected component of $\W$.  We claim that in this case, $\Phi$ is branched of
  order three along $\overline{C}$.  Again, let $q\in(Q^1\overline{C})^*$ be represented by $(Y, \eta)\in \proj\Omega\barmoduli$
  and $S\in Q^1(Y)^*$, and normalize $\eta$ so that $S(\eta^2)=1$.  This $(Y, \eta)$ is a stable Abelian
  differential with a double zero $z$.  If $\Phi(p) = q$, then $p$ is represented by
  $$(X, \omega) = (Y, \eta)\#_I$$
  for some oriented segment $I$ starting at $z$ such that
  \begin{equation}
    \label{eq:integral1}
    \int_I\eta=1.
  \end{equation}
  There are at most three such segments because there are three positively oriented horizontal directions at
  the zero of $\eta$; therefore, $q$ has at most three preimages.  If $S$ is small, then there is an embedded
  ball around $z$ with large radius.  This means that there are three embedded segments starting at $z$
  satisfying \eqref{eq:integral1}, and we can split along each of these segments.  Thus any point in $Q^1 \overline{C}$
  sufficiently close to the zero section has exactly three preimages under $\Phi$, and the claim follows.

  Now if $C$ is a connected component of $\P$ or $\Stwo[D]$, then we have seen that $U$ is homeomorphic
  to a neighborhood of the zero section in $Q^1 \overline{C}$.  Since $U$ is a tubular neighborhood of $\overline{C}$, $U$ is also
  homeomorphic to a neighborhood of the zero section in $N(\overline{C})$.  Therefore,
  $$[\overline{C}]^2 = \deg N(\overline{C}) = \deg Q^1 \overline{C}.$$
  Since $c_1(Q^1\Y) = [\omega_1]$ by Corollary~\ref{cor:chernQ},
  $$[\overline{C}]^2 = \deg Q^1 \overline{C} = -\int_C\omega_1 = \chi(C),$$
  as claimed.
  
  Now suppose $C$ is a connected component of $\W$.  Regard $U$ as a neighborhood of the zero section in
  $N(\overline{C})$.  By the above claim, we can choose $U$ so that $\Phi\colon U\to (Q^1 \overline{C})^*$ preserves fibers and is
  exactly three-to-one onto its image.  For any line bundle $B\to \overline{C}$, let $\tau_B\in H^2(B, B\setminus
  \overline{C};\reals)$ be its Thom class.  The Thom class $\tau_{(Q^1\overline{C})^*}$ is represented by a $2$-form which is
  supported in $\Phi(U)$ which satisfies
  $$\int_F \tau_{(Q^1\overline{C})^*}=1$$
  for each fiber $F$ of $(Q^1\overline{C})^*$.  For each fiber $F$ of $N(\overline{C})$, we have
  $$\int_F \Phi^* \tau_{(Q^1\overline{C})^*}=3,$$
  thus
  $$\Phi^*\tau_{(Q^1\overline{C})^*} = 3 \tau_{N(\overline{C})}.$$
  It follows that
  \begin{align*}
    [\overline{C}]^2 &= \tau_{N(\overline{C})}\cdot[\overline{C}] \\
    &= \frac{1}{3}\tau_{(Q^1 \overline{C})^*}\cdot[\overline{C}] \\
    &=\frac{1}{3}\deg(Q^1\overline{C})^* \\
    &= \frac{1}{3}\chi(C),
  \end{align*}
  as claimed.
\end{proof}

\begin{cor}
  \label{cor:intomega2}
  For any connected component $C$ of $\W$, we have
  \begin{equation*}
    \int_C\omega_2 = \frac{1}{3}\int_C\omega_1.
  \end{equation*}
\end{cor}

\begin{proof}
  By Proposition~\ref{prop:normalequalstangent}, we have
  $$[C]^2 = \deg T\A|_C = c_1(T\A)\cdot[C] = -\int_C\omega_2.$$
  By Theorem~\ref{thm:selfintersections},
  $$[C]^2 = -\frac{1}{3}\int_C\omega_1,$$
  and the claim follows.
\end{proof}


%% file: euler2.tex
\section{Euler characteristic of $W_D^\epsilon$}
\label{sec:eulerWe}

\paragraph{Cohomology of $\Y$.}

We saw in Theorem~\ref{thm:negativedefinite} that the intersection pairing on $H^2(\Y; \ratls)$ is negative
definite on the subspace $B$ generated by the fundamental classes of the curves $C_P$.  By
Corollary~\ref{cor:pairingzero}, the subspace $\langle[\omega_1], [\omega_2]\rangle\subset H^2(\Y;\ratls)$ is
orthogonal to $B$.  Thus, if we let $J\subset H^2(\Y;\ratls)$ be the orthogonal complement to
$B\oplus\langle[\omega_1],[\omega_2]\rangle$, then we have the orthogonal direct sum,
\begin{equation}
  \label{eq:directsum2}
  H^2(\Y; \ratls) = B \oplus \langle[\omega_1], [\omega_2]\rangle \oplus J.
\end{equation}
Since
$$B\oplus\langle[\omega_1], [\omega_2]\rangle\subset H^{1,1}(\Y; \ratls),$$
$J$ contains all of $H^{2,0}(\Y; \ratls)$ and $H^{0,2}(\Y; \ratls)$.

\paragraph{$\chi(\We)$ when $D$ is not square.}

We now calculate $\chi(\We)$ when $D$ is not square.  Until further notice, we will assume
that $D$ is not square.

According to Corollary~\ref{cor:boundaryrelation},
\begin{equation*}
  \pi_B[\barW] = \pi_B[\barP].
\end{equation*}
Let $B_D = \pi_B[\barP]$, and let $B_D^\epsilon = \pi_B[\barWe]$ for $\epsilon=1, 2$.  Since
$$[\barW] = [\barWzero] + [\barWone],$$
we have
\begin{equation*}
  B_D^0 + B_D^1 = B_D.
\end{equation*}

\begin{lemma}
  For any nonsquare $D$, we have
  \begin{equation}
    \label{eq:BDsquared}
    (B_D)^2 = -15 \chi(\X).
  \end{equation}
\end{lemma}

\begin{proof}
  Since $\W$ and $\P$ are disjoint,
  $$[\barW]\cdot[\barP]=0.$$
  Equation \eqref{eq:BDsquared} follows directly from this together with the equations \eqref{eq:fundamentalPD}
  and \eqref{eq:fundamentalW} for these fundamental classes.
\end{proof}

\begin{theorem}
  \label{thm:fundamentalWe}
  If $D$ is not square, then the fundamental class of $\barWe$ in $H^2(\Y;\ratls)$ is
  \begin{equation}
    \label{eq:fundamentalWe}
    [\barWe] = \frac{3}{4}[\omega_1] + \frac{9}{4} [\omega_2] + B_D^\epsilon + j
  \end{equation}
  for some $j\in J$.
\end{theorem}

\begin{proof}
  Since $\Wzero$ and $\P$ are disjoint, we have
  \begin{equation}
    \label{eq:1}
    [\barWzero]\cdot[\P] = 0,
  \end{equation}
  and by \eqref{eq:tauoncusps}, we have
  \begin{equation}
    \label{eq:4}
    (B_D^0)^2 = (B_D^1)^2.
  \end{equation}

  By \eqref{eq:directsum2} and Corollary~\ref{cor:intomega2}, the fundamental classes of the $\barWe$ are of
  the form,
  \begin{align*}
    [\barWzero] &= a[\omega_1] + 3a[\omega_2] + B_D^0 + j\\
    [\barWone] &= \left(\frac{3}{2} - a\right)[\omega_1] + \left(\frac{9}{2} - 3a\right)[\omega_2] + B_D^1 - j,
  \end{align*}
  for some $a\in\ratls$ and $j\in J$.  In terms of $a$, \eqref{eq:1} becomes
  \begin{equation}
    \label{eq:5}
    10a\chi(\X) + (B_D^0)^2 + B_D^0\cdot B_D^1 = 0.
  \end{equation}
  
  From \eqref{eq:BDsquared} and \eqref{eq:4}, we obtain
  \begin{equation}
    \label{eq:9}
    (B_D^0)^2 + B_D^0\cdot B_D^1 = \frac{1}{2}(B_D^0 + B_D^1)^2 = \frac{1}{2}(B_D)^2=-\frac{15}{2}\chi(\X).
  \end{equation}
  Combining \eqref{eq:5} and \eqref{eq:9} yields
  \begin{equation*}
        a=\frac{3}{4},
  \end{equation*}
  as desired.
\end{proof}

\begin{remark}
  It seems likely that $j=0$, but we don't know how to prove this.
\end{remark}

\begin{cor}
  If $D$ is not square, then
  \begin{equation*}
    \chi(\Wzero) = \chi(\Wone). 
  \end{equation*}
\end{cor}

\paragraph{$\chi(\We)$ when $D$ is square.}

We now turn to the calculation of $\chi(\We[d^2])$.  The idea is the same as the proof of
Theorem~\ref{thm:fundamentalWe}, but the calculation is more complicated because of the presence of the curves
$\Si$.  We will restrict to the case $d>2$ because $\W[4]=\emptyset$.  We start by calculating the
intersections of various classes in $H^2(\Y[d^2]; \ratls)$.

\begin{lemma}
  For any $d>2$, we have the following intersection numbers:
  {\allowdisplaybreaks
    \begin{align}
      \label{eq:omega1dotomega2}
      [\omega_1]\cdot[\omega_2] &= \frac{1}{72}d^3\Sum \\
      \label{eq:SidotP}
      \pi_B[\barSi]\cdot\pi_B[\barP[d^2]] &= \left(-\frac{5}{24} d^2 + \frac{1}{4}d\right)\Sum \\
      \label{eq:SidotSi}
      (\pi_B[\barSi])^2 &= -\frac{1}{12} d^2 \Sum \\
      \label{eq:S1dotS2}
      \pi_B[\barSone]\cdot\pi_B[\barStwo] &= -\frac{1}{2}d\Sum + \frac{1}{2}\phi(d) \\
      \label{eq:S2dotW}
      \pi_B[\barStwo]\cdot\pi_B[\barW[d^2]]&=\left(-\frac{1}{8} d^2 + \frac{1}{4}d\right)\Sum \\
      \label{eq:WdotW}
      (\pi_B[\barW[d^2]])^2 &= \left(-\frac{5}{24}d^3 + \frac{19}{24}d^2 -
        \frac{3}{4}d\right)\Sum\\
      \label{eq:PdotP}
      (\pi_B[\barP[d^2]])^2 &= \left(-\frac{5}{24}d^3 + \frac{11}{24}d^2 -
        \frac{1}{4}d\right)\Sum\\
      \label{eq:S1dotW}
      \pi_B[\barSone]\cdot\pi_B[\barW[d^2]] &= \left(-\frac{5}{24}d^2 + \frac{3}{4} d\right) \Sum -
      \frac{1}{2}\phi(d) 
    \end{align}}
\end{lemma}

\begin{proof}
  Equation \eqref{eq:omega1dotomega2} is $\chi(\X[d^2])$, which is given in Theorem~\ref{thm:chiX}.

  By Theorem~\ref{thm:Ysummary}, $\barP[d^2]\cap\barSi=\emptyset$. Thus $[\barP[d^2]]\cdot[\barSi]=0$, from
  which \eqref{eq:SidotP} follows.

  By Theorem~\ref{thm:selfintersections} and \eqref{eq:chiStwo},
  $$[\barSi]^2 = \chi(\Si) = -\frac{1}{12}d^2\Sum,$$
  from which \eqref{eq:SidotSi} follows.

  By Theorem~\ref{thm:Ysummary},
  $$[\barSone]\cdot[\barStwo] = -\frac{1}{2}\phi(d).$$
  Equation \eqref{eq:S1dotS2} follows.

  By Theorem~\ref{thm:Ysummary},
  $$[\barStwo]\cdot[\barW[d^2]]=0.$$
  Equation~\eqref{eq:S2dotW} follows.

  By Theorem~\ref{thm:selfintersections} and \eqref{eq:corWd2},
  $$[\barW[d^2]]^2 = \frac{1}{3}\chi(\W[d^2]) = -\frac{1}{48}d^2(d-2)\Sum,$$
  from which \eqref{eq:WdotW} follows.  Equation \eqref{eq:PdotP} is proved similarly.
  
  By Corollary~\ref{cor:boundaryrelation},
  $$\pi_B[\barSone]\cdot\pi_B[\barW[d^2]] = \pi_B[\barSone]\cdot\pi_B[\barP[d^2]] -
  \pi_B[\barSone]\cdot\pi_B[\barStwo].$$
  Then \eqref{eq:S1dotW} follows from \eqref{eq:SidotP} and \eqref{eq:S1dotS2}.

\end{proof}

\begin{theorem}
  For any $d>2$, the fundamental class of $\barWe[d^2]$ in $H^2(\Y[d^2];\ratls)$ is given by
  \begin{align*}
    [\barWzero[d^2]] &= \frac{3}{4}\left( 1 - \frac{1}{d}\right)[\omega_1] + \frac{9}{4} \left( 1 -
      \frac{1}{d}\right)[\omega_2] + \pi_B[\barWzero[d^2]] +j\\
    [\barWone[d^2]] &= \frac{3}{4}\left( 1 - \frac{3}{d}\right)[\omega_1] + \frac{9}{4} \left( 1 -
      \frac{3}{d}\right)[\omega_2] + \pi_B[\barWone[d^2]] -j.
  \end{align*}
  for some $j\in J$.
\end{theorem}

\begin{proof}
  By \eqref{eq:directsum2}, Corollary~\ref{cor:fundamentalW}, and Corollary~\ref{cor:intomega2}, the
  fundamental classes of the $\barWe[d^2]$ are given by
  \begin{align*}
    [\barWzero[d^2]] &= a[\omega_1] + 3a[\omega_2] + \pi_B[\barWzero[d^2]] + j \\
    [\barWone[d^2]] &= \left(\frac{3}{2}-\frac{3}{d}-a\right)[\omega_1] +
    \left(\frac{9}{2}-\frac{9}{d}-3a\right)[\omega_2] + \pi_B[\barWone[d^2]] - j
  \end{align*}
  for some $a\in\ratls$ and $j\in J$.

  From Corollary~\ref{cor:tauoncusps}, we have
  \begin{align*}
    (\pi_B[\barWzero[d^2]])^2 &= ( \pi_B[\barWone[d^2]] + \pi_B[\barStwo[d^2]])^2 \\
    &= (\pi_B[\barW[d^2]] - \pi_B[\barWzero[d^2]] + \pi_B[\barStwo])^2.
  \end{align*}
  Using \eqref{eq:SidotSi}, \eqref{eq:S2dotW}, and \eqref{eq:WdotW}, this simplifies to 
  \begin{multline}
    \label{eq:7a}
    \pi_B[\barW[d^2]]\cdot\pi_B[\barWzero[d^2]] + \pi_B[\barWzero[d^2]]\cdot\pi_B[\barStwo] \\ =
    \left(-\frac{5}{48}d^3 + \frac{11}{48}d^2 - \frac{1}{8}d\right)\Sum.
  \end{multline}
  From Corollary~\ref{cor:boundaryrelation}, we have
  \begin{equation}
    \label{eq:6a}
     \pi_B[\W[d^2]] + \pi_B[\Stwo[d^2]] = \pi_B[\P[d^2]].  
  \end{equation}
  Multiplying \eqref{eq:6a} by $\pi_B[\barWzero[d^2]]$ and subtracting the result from \eqref{eq:7a}, we
  obtain
  \begin{equation}
    \label{eq:8a}
    \pi_B[\barWzero[d^2]]\cdot\pi_B[\barP[d^2]] = \left(-\frac{5}{48}d^3 + \frac{11}{48}d^2 - \frac{1}{8}d\right)\Sum.
  \end{equation}
  Since $\barWzero[d^2]$ and $\barP[d^2]$ are disjoint, we have
  \begin{equation*}
    [\barWzero[d^2]]\cdot[\barP[d^2]] = 0.
  \end{equation*}
  Expanding, this becomes
  \begin{equation*}
    \left(\frac{5}{32}d^3 - \frac{1}{6} d^2\right)a\Sum + \pi_B[\barWzero[d^2]]\cdot\pi_B[\barP[d^2]] = 0,
  \end{equation*}
  which with \eqref{eq:8a} yields
  \begin{equation*}
    a = \frac{3}{4}\left(1-\frac{1}{d}\right).
  \end{equation*}
\end{proof}

\begin{cor}
  For any $d>2$ with $d\equiv 1\pmod 2$,
  \begin{align*}
    \chi(\Wzero[d^2])&=-\frac{1}{32} d^2 (d-1)\sum_{r | d}\frac{\mu(r)}{r^2}\\
     \chi(\Wone[d^2])&=-\frac{1}{32} d^2 (d-3)\sum_{r | d}\frac{\mu(r)}{r^2}.\\
  \end{align*}
\end{cor}

\paragraph{Once cylinder cusps.}

As an application of the calculation of $[\barWe[d^2]]$, we give formulas for the number of one-cylinder cusps of $\We[d^2]$.
These formulas were established independently by Leli\`evre and Royer in \cite{lelievreroyer}.

\begin{theorem}
  For any $d>3$, the number of one-cylinder cusps of $\W[d^2]$ is
  \begin{equation}
    \label{eq:onecyl1}
    \frac{1}{6}d^2\Sum - \frac{1}{2}\phi(d);
  \end{equation}
  the number of one-cylinder cusps of $\Wzero[d^2]$ is
  \begin{equation}
    \label{eq:onecyl2}
    \frac{1}{24}d^2\Sum;
  \end{equation}
  and the number of one-cylinder cusps of $\Wone[d^2]$ is
  \begin{equation}
    \label{eq:onecyl3}
    \frac{1}{8}d^2\Sum - \frac{1}{2}\phi(d).
  \end{equation}
\end{theorem}

\begin{proof}
  The one-cylinder cusps of $\W[d^2]$ are the points of the intersection,
  $$\Sone\cap\barW[d^2].$$
  When $d>3$ this intersection is transverse by Theorem~\ref{thm:Ysummary}, so
  $[\barW[d^2]]\cdot[\barSone]$ is equal to the number of one-cylinder cusps of $\W[d^2]$.  Similarly,
  $[\barWe[d^2]]\cdot[\barSone]$ is equal to the number of one-cylinder cusps of $\We[d^2]$.

  Using \eqref{eq:S1dotW}, we obtain
  \begin{equation*}
    [\barW[d^2]]\cdot[\barSone] = \frac{1}{6}d^2\Sum - \frac{1}{2}\phi(d),
  \end{equation*}
  which implies \eqref{eq:onecyl1}.

  Now let's calculate $[\barSone]\cdot[\barWzero[d^2]]$.  By Corollary~\ref{cor:tauoncusps}, we have
  \begin{align}
    \notag
    &\pi_B[\barSone]\cdot\pi_B[\barWzero[d^2]] \\
    \notag
    &= \tau^*\left(\pi_B[\barSone]\right)\cdot\tau^*\left(\pi_B[\barWzero[d^2]]\right) \\
    \notag
    &=\pi_B[\barStwo]\cdot\pi_B[\barWone[d^2]] + (\pi_B[\barStwo])^2 \\
    \label{eq:onecyl4}
    &= \pi_B[\barStwo]\cdot\pi_B[\barW[d^2]] + \pi_B[\barStwo]\cdot\pi_B[\barWzero[d^2]] + (\pi_B[\barStwo])^2
  \end{align}
  Since $\barWzero[d^2]$ and $\barStwo$ are disjoint, we have $[\barWzero[d^2]]\cdot[\barStwo]=0$, and it
  follows that
  \begin{equation}
    \label{eq:onecyl5}
    \pi_B[\barStwo]\cdot\pi_B[\barWzero[d^2]] = \left(-\frac{1}{16}d^2 + \frac{1}{16}d\right)\Sum.
  \end{equation}
  Substituting \eqref{eq:SidotSi}, \eqref{eq:S2dotW}, and \eqref{eq:onecyl5} into \eqref{eq:onecyl4}, we
  obtain
  \begin{equation}
    \label{eq:onecyl6}
    \pi_B[\barSone]\cdot\pi_B[\barWzero[d^2]] = \left(-\frac{7}{48}d^2 + \frac{3}{16}d\right)\Sum.
  \end{equation}
  Equation \eqref{eq:onecyl6} yields
  \begin{equation*}
    [\barSone]\cdot[\barWzero[d^2]]=\frac{1}{24}d^2\Sum,
  \end{equation*}
  which implies \eqref{eq:onecyl2}.  Equation \eqref{eq:onecyl3} follows from \eqref{eq:onecyl1} and \eqref{eq:onecyl2}.
\end{proof}


%% file: siegelveech.tex
\section{Siegel-Veech Constants}
\label{sec:siegelveech}

As an application of our results, we record the Siegel-Veech constants
counting cylinders on translation surfaces on the Teichm\"uller curves
$\W$ and $\We$.  This is basically a matter of plugging in our results
into known formulas for these constants.

Given a translation surface $(X, \omega)$, let
\begin{equation*}
  N((X, \omega), L) = \#\{\text{maximal cylinders of length at most 
  $L$ on $(X, \omega)$}\}.
\end{equation*}
If $(X, \omega)$ lies on a Teichm\"uller curve, then Veech
\cite{veech89} showed
\begin{equation}
  \label{eq:svdef}
  N((X, \omega),L)\sim  \frac{c}{\pi\Area(X, \omega)} L^2.
\end{equation}
The constant $c$, known as a Siegel-Veech constant,
only depends on the Teichm\"uller curve on which $(X, \omega)$ lies.
Let $c_D$ be and $c_D^\epsilon$ be the Siegel-Veech constants associated to
$\W$ and $\We$ respectively.

Given a $\W$-prototype $P=(a, b, c, \bar{q})$, define
\begin{equation*}
  v(P) = \frac{c}{\gcd(a, c)}\left( 1 -
  \frac{a}{c}\lambda^2\right)\left(1 + \frac{1}{\lambda^2}\right),
\end{equation*}
where $\lambda = \lambda(P)$, the positive root of $ax^2+bx+c=0$.

\begin{theorem}
  If $D$ is not square, then
  \begin{equation*}
    c_D = \frac{\sum_{P\in\Wprot} v(P)}{-2 \chi(\W)}, \quad \text{and}
    \quad  c_D^\epsilon = \frac{\sum_{P\in\Wprot^\epsilon} v(P)}{-2 \chi(\W^\epsilon)},
  \end{equation*}
  where $\Wprot^\epsilon$ is the set of $\W$-prototypes of
  spin invariant $\epsilon$.
\end{theorem}

\begin{proof}[Sketch of proof]
  To fix notation, assume $\W$ is connected as the proof is the same
  otherwise.  Each cusp of $\W$ corresponds to a $\W$-prototype $P$,
  and we associated to this cusp on p.~\pageref{page:cusp} a surface $(X_P,
  \omega_P)$ on $\W$ with a decomposition into two cylinders $C_1$ and
  $C_2$ (say $C_1$ is the short cylinder).  The subgroup of the Veech
  group of $(X_P, \omega_P)$ (the stabilizer of this surface in $\SLtwoR$) which preserves the horizontal direction
  is generated by
  \begin{equation*}
    g =
    \begin{pmatrix}
      1 & t \\
      0 & 1
    \end{pmatrix},
  \end{equation*}
  where
  \begin{equation*}
    t  =
    \frac{c}{\gcd(a, c)}.
  \end{equation*}
  Let $i(C_i)$ be the order of the Dehn twist which $g$ induces on
  $C_i$.  We have,
  \begin{align*}
    i(C_1) &= -\frac{c}{\gcd(a, c)}, \quad\text{and}\\
    i(C_2) &=\phantom{-} \frac{a}{\gcd(a, c)}.
  \end{align*}
  We have
  \begin{align*}
    v(P) &= \Area(X_P, \omega_P) \sum_i \frac{i(C_i)}{\Area(C_i)}.
  \end{align*}
  It follows from Theorem~6.5 of \cite{gutkinjudge} that for any
  $(X, \omega)\in\Omega_1\W$, 
  \begin{equation*}
    N((X, \omega), L) \sim \frac{\sum_{P\in\Wprot}v(P)}{\Area \W} L^2,
  \end{equation*}
  as desired.
\end{proof}

We list the Siegel-Veech constants for $D<100$ in Table~\ref{tab:sv},
using the convention $c_D = c_D^0$ if $D\equiv 1\pmod 8$ ($c_D^1$ is
the Galois conjugate of $c_D^0$ as we prove below).  From numerical
calculations, as $D\to\infty$, the constants appear to converge to
$10$, which by \cite{emz} is the Siegel-Veech constant for counting cylinders on a generic
$(X, \omega)\in \Omega\moduli(2)$.  It would be
interesting to find a closed formula for $c_D$.

\paragraph{Arithmetic of Siegel-Veech constants.} According to
\cite{gutkinjudge}, the Siegel-Veech constant $c$ in \eqref{eq:svdef}
lies in the trace field of the Veech group of $(X, \omega)$.  We get more precise
information in the case of the $\W$:

\begin{theorem}
  Suppose $D$ is not square.  If $D\not\equiv 1 \pmod 8$, then $c_D\in \ratls$.
  Otherwise, $c_D^0$ and $c_D^1$ are Galois conjugate elements of $\ratls(\sqrt{D})$.
\end{theorem}

\begin{proof}
  First assume $D\not\equiv 1 \pmod 8$.  We have the involution $t$ on
  the set of $\W$-prototypes, defined in \S\ref{sec:prototypes}.
  Actually, it was defined on $\Y$-prototypes, but the definition
  works just as well in this case.  This involution satisfies
  \begin{equation*}
    v(t(P)) = v(P)',
  \end{equation*}
  it follows that $c_D'=c_D$ as desired.

  Now suppose $D\equiv 1\mod 8$.  If it were true that
  $\epsilon(t(P)) = \epsilon(P)+1$ -- where $\epsilon$ is the spin
  invariant -- then we would be done.  This is not true, but we can
  modify $t$ so that it is.  Define a bijection on the set of $\W$-prototypes,
  \begin{equation*}
    s(P) =
    \begin{cases}
      P, & \text{if $a-b+c<0$ or $a\not\equiv b \pmod 2$}; \\
      P + (0,0,0, \gcd(a. b. c)), & \text{otherwise}.
    \end{cases}
  \end{equation*}
  Then $t' = s\circ t$ is also a bijection.  One can check that it
  satisfies
  \begin{align*}
    \epsilon(t'(P)) &= \epsilon(P)+1, \quad\text{and} \\
    v(t'(P)) &= v(P)',
  \end{align*}
  from which the second claim follows.
\end{proof}

\paragraph{Applications to billiards.}

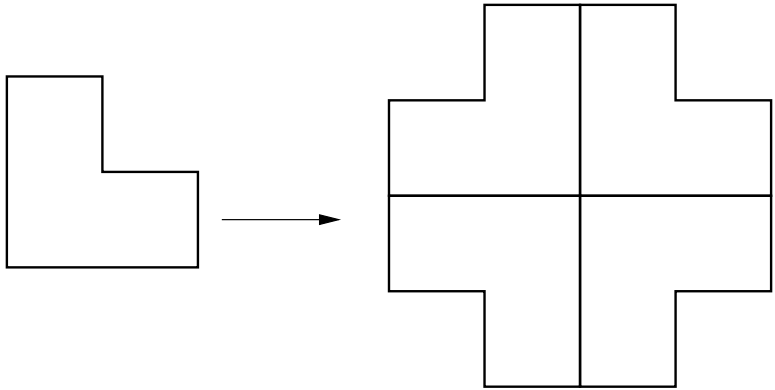
\begin{figure}[htbp]
  \centering
  \input{cross.pstex_t}
  \caption{Unfolding a L-shaped polygon}
  \label{fig:cross}
\end{figure}

Given a L-shaped polygon $P$ (or more generally a rational angled
polygon), there is a construction called \emph{unfolding} which
produced a translation surface from $P$.  To construct the unfolding
of $P$, join four copies of $P$ to form a cross and glue opposite
sides as in Figure~\ref{fig:cross}.  This yields a genus two translation
surface  $(X, \omega)$ with a single double zero.  Billiards paths on
$P$ unfold to closed geodesics on $(X, \omega)$, and we obtain a one-to
one correspondence between closed paths of length $L$ on $P$ and on the
unfolding.  Recall that we defined in \S\ref{sec:intro} a L-shaped
polygon $P(D)$ for each real quadratic discriminant $D$, and we
defined $N(P(D), L)$ to be the number of families of closed billiards
paths of length at least $L$.

\begin{theorem}
  If $D$ is not square, then
  \begin{equation*}
    N(P(D), L) \sim c(D)\frac{\pi}{4 \Area(P(D))} L^2,
  \end{equation*}
  where $c(D)= c_D$, if $D\not\equiv 1\pmod 8$, and $c(D) = c_D^{(1+f)/2}$
  if $D\equiv 1 \pmod 8$, where $f$ is the conductor of $D$.
\end{theorem}

\begin{proof}
  This follows directly from \eqref{eq:svdef} applied to the unfolding
  of $P(D)$.  The factor of $4$ is because the unfolding has $4$ times
  the area.  To see that the unfolding has discriminant $D$ and spin
  invariant $\epsilon=(1+f)/2$, note that
  the  unfolding is $\GLtwoR$-equivalent to the surface $(X_P,
  \omega_P)$ associated to the $\W$-prototype,
  \begin{equation*}
    P =
    \begin{cases}
      \left(1, -1, \frac{1-D}{4}, 0\right); &\text{if $D$ is odd};\\
      \left(1, 0, -\frac{D}{4}, 0\right); &\text{if $D$ is even},
    \end{cases}
  \end{equation*}
  defined on p.~\pageref{page:cusp}.  The spin invariant is then given
  by Theorem~\ref{thm:cuspsbyspin}.
\end{proof}

\input{svtable}


%% file: cross.pstex_t
\begin{picture}(0,0)%
\includegraphics{cross.pstex}%
\end{picture}%
\setlength{\unitlength}{3947sp}%
\begingroup\makeatletter\ifx\SetFigFont\undefined%
\gdef\SetFigFont#1#2#3#4#5{%
  \reset@font\fontsize{#1}{#2pt}%
  \fontfamily{#3}\fontseries{#4}\fontshape{#5}%
  \selectfont}%
\fi\endgroup%
\begin{picture}(4844,2444)(4779,-4733)
\end{picture}%

%% file: svtable.tex
\begin{table}[p] \centering \begin{tabular}{|cc|cc|} \hline$D$ & $c_D$ & $D$ & $c_D$ \\ \hline$5$
 & $\displaystyle \frac{\vphantom{\int}25
}{\vphantom{\int}3
}$ & $56$
 & $\displaystyle \frac{\vphantom{\int}1796
}{\vphantom{\int}195
}$ \\ $8$
 & $\displaystyle \frac{\vphantom{\int}28
}{\vphantom{\int}3
}$ & $57$
 & $\displaystyle \frac{23693
}{2352
}+ \frac{29
}{784
} \sqrt{57
}$ \\ $12$
 & $\displaystyle \frac{\vphantom{\int}26
}{\vphantom{\int}3
}$ & $60$
 & $\displaystyle \frac{\vphantom{\int}2158
}{\vphantom{\int}231
}$ \\ $13$
 & $\displaystyle \frac{\vphantom{\int}91
}{\vphantom{\int}9
}$ & $61$
 & $\displaystyle \frac{\vphantom{\int}194651
}{\vphantom{\int}19305
}$ \\ $17$
 & $\displaystyle \frac{221
}{24
}+ \frac{1
}{8
} \sqrt{17
}$ & $65$
 & $\displaystyle \frac{52429
}{5376
}+ \frac{113
}{1792
} \sqrt{65
}$ \\ $20$
 & $\displaystyle \frac{\vphantom{\int}31
}{\vphantom{\int}3
}$ & $68$
 & $\displaystyle \frac{\vphantom{\int}413
}{\vphantom{\int}39
}$ \\ $21$
 & $\displaystyle \frac{\vphantom{\int}133
}{\vphantom{\int}15
}$ & $69$
 & $\displaystyle \frac{\vphantom{\int}26611
}{\vphantom{\int}2805
}$ \\ $24$
 & $\displaystyle \frac{\vphantom{\int}148
}{\vphantom{\int}15
}$ & $72$
 & $\displaystyle \frac{\vphantom{\int}18868
}{\vphantom{\int}1785
}$ \\ $28$
 & $\displaystyle \frac{\vphantom{\int}82
}{\vphantom{\int}9
}$ & $73$
 & $\displaystyle \frac{3285
}{352
}+ \frac{23
}{864
} \sqrt{73
}$ \\ $29$
 & $\displaystyle \frac{\vphantom{\int}377
}{\vphantom{\int}35
}$ & $76$
 & $\displaystyle \frac{\vphantom{\int}2822
}{\vphantom{\int}285
}$ \\ $32$
 & $\displaystyle \frac{\vphantom{\int}190
}{\vphantom{\int}21
}$ & $77$
 & $\displaystyle \frac{\vphantom{\int}116699
}{\vphantom{\int}12597
}$ \\ $33$
 & $\displaystyle \frac{473
}{48
}- \frac{11
}{144
} \sqrt{33
}$ & $80$
 & $\displaystyle \frac{\vphantom{\int}12631
}{\vphantom{\int}1254
}$ \\ $37$
 & $\displaystyle \frac{\vphantom{\int}9139
}{\vphantom{\int}945
}$ & $84$
 & $\displaystyle \frac{\vphantom{\int}487
}{\vphantom{\int}51
}$ \\ $40$
 & $\displaystyle \frac{\vphantom{\int}1924
}{\vphantom{\int}189
}$ & $85$
 & $\displaystyle \frac{\vphantom{\int}336821
}{\vphantom{\int}32319
}$ \\ $41$
 & $\displaystyle \frac{8897
}{960
}- \frac{23
}{320
} \sqrt{41
}$ & $88$
 & $\displaystyle \frac{\vphantom{\int}182236
}{\vphantom{\int}18837
}$ \\ $44$
 & $\displaystyle \frac{\vphantom{\int}7682
}{\vphantom{\int}735
}$ & $89$
 & $\displaystyle \frac{702833
}{68640
}- \frac{831
}{22880
} \sqrt{89
}$ \\ $45$
 & $\displaystyle \frac{\vphantom{\int}299
}{\vphantom{\int}33
}$ & $92$
 & $\displaystyle \frac{\vphantom{\int}204178
}{\vphantom{\int}21945
}$ \\ $48$
 & $\displaystyle \frac{\vphantom{\int}325
}{\vphantom{\int}33
}$ & $93$
 & $\displaystyle \frac{\vphantom{\int}2823449
}{\vphantom{\int}270963
}$ \\ $52$
 & $\displaystyle \frac{\vphantom{\int}1283
}{\vphantom{\int}135
}$ & $96$
 & $\displaystyle \frac{\vphantom{\int}3194
}{\vphantom{\int}345
}$ \\ $53$
 & $\displaystyle \frac{\vphantom{\int}228695
}{\vphantom{\int}21021
}$ & $97$
 & $\displaystyle \frac{44329
}{4488
}- \frac{1145
}{40392
} \sqrt{97
}$ \\  \hline   \end{tabular}   \caption{Siegel-Veech constants for $\W$}   \label{tab:sv} \end{table}

%% file: lyapunov.tex
\section{Lyapunov Exponents}
\label{sec:lyapunov}

Consider the rank $2g$ bundle $\mathcal{H}(\reals) \to
\Omega_1\moduli[g]$ whose fiber over the surface $(X, \omega)$ is
$H_1(X; \reals)$.  The action of the diagonal subgroup of $\SLtwoR$
lifts to a linear action on $\mathcal{H}_1(\reals)$.  This action is
called the Kontsevich-Zorich cocycle.  The $2g$ Lyapunov exponents of
the Kontsevich-Zorich cocycle with respect to a finite, ergodic,
$\SLtwoR$-invariant measure $\mu$ are of the form,
$$1 = \lambda_1(\mu) > \dots > \lambda_g(\mu) > 0 > -\lambda_g(\mu) >
\dots -\lambda_1(\mu)= -1.$$
The goal of this section is to prove the following unpublished theorem
of Kontsevich and Zorich as an application of our previous results.

\begin{theorem}
  \label{thm:kz2}
  If $\mu$ is any finite, ergodic, $\SLtwoR$-invariant measure on
  $\Omega_1\moduli$, then
  $$\lambda_2(\mu)=
  \begin{cases}
    1/3, & \text{if $\mu$ is supported on } \Omega_1\moduli(2);\\
    1/2, & \text{if $\mu$ is supported on } \Omega_1\moduli(1,1).
  \end{cases}
  $$
\end{theorem}

\begin{remark}
  When $\mu$ is period measure on $\E$, with $D$ nonsquare, the proof
  will require Theorems~\ref{thm:Fextends} and
  \ref{thm:Fintersectionzero} below.  These results will also be used
  in a future paper, and we will defer the proofs until then.
\end{remark}

The proof of Theorem~\ref{thm:kz2} is based on a formula due to
Kontsevich \cite{kontsevich} for the sums of Lyapunov exponents in
terms of some integrals over $\proj\Omega\moduli[g]$.  Kontsevich states his
formula for the Lyapunov exponents on the entire strata
$\Omega_1\moduli[g]({\bf n})$, but his formula is equally valid for
all ergodic, $\SLtwoR$-invariant probability measures on $\Omega_1\moduli$ and we will
state it for these.

Let $E$ be the bundle over $\proj \Omega\moduli$ obtained by pulling
back $\Omega\moduli$ by $\proj\Omega\moduli\to\moduli$, and let $L$ be
the canonical sub-line-bundle whose fiber over $(X, [\omega])$ is $\cx
\omega$.  Give $E$ the Hodge metric,
$$h(\omega, \eta) = \int_X \omega\wedge \bar{\eta},$$
and define two-forms on $\proj\Omega\moduli,$
\begin{equation*}
  \gamma_1 = c_1(E, h) \quad \text{and} \quad
  \gamma_2 = c_1(L, h).
\end{equation*}

For any $\SLtwoR$-invariant measure $\mu$ on $\Omega_1\moduli[g]$, let
$\pi_*\mu$ be the pushforward to $\proj\Omega\moduli[g]$, and let
$\proj\mu$ be the disintegration of $\pi_*\mu$ with respect to the
foliation $\mathcal{F}$ of $\proj\Omega\moduli[g]$ by images of
$\SLtwoR$-orbits.  That is, $\proj\mu$ is the unique transverse
invariant measure such that the product of $\proj \mu$ with the
hyperbolic area measure on the leaves of $\mathcal{F}$ is $\pi_*\mu$.
See \cite{bourbaki} for a discussion of disintegration.

\begin{theorem}[Kontsevich]
  For any finite ergodic, $\SLtwoR$-invariant measure $\mu$ on
  $\Omega_1\moduli$,
  \begin{equation}
    \label{eq:kontsevichformula}
    \lambda_1(\mu) + \lambda_2(\mu) = \frac{\int_{\proj\mu} \gamma_1}{\int_{\proj\mu} \gamma_2}.
  \end{equation}
\end{theorem}

\paragraph{Uniform distribution.}

In order to evaluate \eqref{eq:kontsevichformula}, we will need the
following weak uniform distribution result for Teichm\"uller curves.
Given a space $X$ with an $\SLtwoR$-action, let $M(X)$ be the space
of finite $\SLtwoR$-invariant measures on $X$ with the weak${}^*$
topology.  Let $C({\bf n})\subset M(\Omega_1\moduli[g]({\bf n}))$ be
the convex cone spanned by the measures defined by  hyperbolic area
on those Teichm\"uller curves in $\Omega_1\moduli[g]({\bf n})$ which are
generated by square-tiled surfaces.   Similarly, let  $C_{d^2}\subset
M(\Omega_1\E[d^2](1,1))$ be the convex cone spanned by those Teichm\"uller
curves in $\Omega_1\E[d^2](1,1)$ generated
by square-tiled surfaces.

\begin{theorem}
  \label{thm:uniformdistribution}
  We have 
  \begin{align*}
    \mu'({\bf n}) & \in C({\bf n}), \quad \text{and}\\
    \mu'_{d^2} &\in C_{d^2},
  \end{align*}
  where the measures on the left are the period measures on these
  strata defined in \S\ref{sec:abelian}.
\end{theorem}

Let $\eta_m({\bf n})$ and $\eta_m[d]$ be the measures on
$\Omega\moduli[g]({\bf n})$ and $\Omega\E[d^2](1,1)$ defined by
putting $\delta$-masses of equal weight on the square-tiled surfaces
with at most $m$ squares, normalized to have the same volumes as
$\mu'({\bf n})$ and $\mu'_{d^2}$ respectively.  Let $\eta'_m({\bf n})$
and $\eta'_m[d]$ be the projections to $\Omega_1\moduli[g]({\bf n})$
and $\Omega_1\E[d^2](1,1)$.  The following lemma is well-known.

\begin{lemma}
  \label{lem:uniformdistribution}
  We have
  \begin{align*}
    \lim_{m\to\infty} \eta'_m({\bf n}) &= \mu'({\bf n}),\\
    \lim_{m\to\infty} \eta'_m[d] &= \mu'_{d^2}
  \end{align*}
  in the weak${}^*$ topology on measures.
\end{lemma}

\begin{proof}
  We will prove the first statement, the other having the same proof.
  Let $\mu''({\bf n})$ be the restriction of $\mu'({\bf n})$ to
  $\Omega_{\leq 1}\moduli[g]({\bf n})$, and let
  \begin{equation*}
    \eta_m''({\bf n}) = \left(t_{1/\sqrt{m}}\right)_*\eta_m({\bf n}),
  \end{equation*}
  where $t_r$ is the multiplication-by-$r$ map on $\Omega\moduli[g]$.
  Since  $\mu''({\bf n})$ and $\eta_m''({\bf n})$ project to $\mu'({\bf
    n})$ and $\eta_m'({\bf n})$ respectively, it is enough to show that
  \begin{equation*}
    \lim_{m\to\infty} \eta''_m({\bf n}) = \mu''({\bf n}).
  \end{equation*}
  This is easy to see in the period coordinates on
  $\Omega\moduli[g]({\bf n})$ (described in \S\ref{sec:abelian}).  In
  period coordinates, the square-tiled surfaces are exactly the points
  on the integer lattice, so $\eta''_m({\bf n})$ consists of
  $\delta$-masses of uniform weight on points of norm less than one
  on the $\zed/\sqrt{m}$ lattice.  Since these lattices
  have mesh approaching zero, the measures converge to the uniform
  measure which is just $\mu''({\bf n})$ in period coordinates.
\end{proof}

\begin{proof}[Proof of Theorem~\ref{thm:uniformdistribution}]
  The set of square-tiled surfaces is $\SLtwoZ$-invariant, so the
  measure $\eta_m'({\bf n})$ is also $\SLtwoZ$-invariant.   Thus we
  can define
  $$\nu_m({\bf n}) = \int\limits_{\SLtwoR/\SLtwoZ} \eta_m'({\bf n})\,d\rho \in
  M(\Omega_1\moduli[g]({\bf n})),$$
  where $\rho$ is Haar measure on $\SLtwoR$.
  The measure $\nu_m({\bf n})$ is $\SLtwoR$-invariant and supported
  on finitely many Teichm\"uller curves generated by square-tiled
  surfaces, so $\nu_m({\bf n})$ lies in the cone $C({\bf n})$.

  Let $\Delta\subset \SLtwoR$ be a fundamental domain for $\SLtwoZ$.
  Given a compactly supported, continuous function $f$ on
  $\Omega_1\moduli[g]({\bf n})$, we have,
  \begin{align*}
    \lim_{m\to\infty} \int\limits_{\Omega_1\moduli[g]({\bf n})} f
    \,d\nu_m({\bf n})&=
    \lim_{m\to\infty}\int\limits_{A\in\Delta}\int\limits_{\Omega_1\moduli[g]({\bf
        n})} f\circ A \,d\eta_m'({\bf n}) \,d\rho\\
    &=\int\limits_{A\in\Delta}\left(\lim_{m\to\infty}\int\limits_{\Omega_1\moduli[g]({\bf
          n})} f\circ A \, d\eta_m'({\bf n})\right)d\rho\\
    &=\int\limits_{A\in\Delta}\int\limits_{\Omega_1\moduli[g]({\bf
        n})}f\circ A\, d\mu'({\bf n})\,d\rho\\
    &= \int\limits_{A\in\Delta}\int\limits_{\Omega_1\moduli[g]({\bf
        n})}f \,d\mu'({\bf n})\,d\rho\\
    &= \vol(\SLtwoR/\SLtwoZ)\int\limits_{\Omega_1\moduli[g]({\bf n})}
    f \,d\mu'({\bf n}),
  \end{align*}
  where the second equality follows from the Dominated Convergence
  Theorem; the third equality follows from
  Lemma~\ref{lem:uniformdistribution}; and the fourth equality uses
  the $\SLtwoR$-invariance of $\mu'({\bf n})$.  It follows that
  $\mu'({\bf n})\in C({\bf n})$.  The proof for $\mu'_{d^2}$ is the same.
\end{proof}

\paragraph{The Siegel-Veech transform.}

Given any translation surface $(X, \omega)$ with an oriented saddle connection $I$, we associate the complex
number $v(I)=\int_I\omega$.  Let $V(X, \omega)$ be the collection of vectors in $\cx$ associated with saddle
connections on $(X, \omega)$, together with the multiplicities with which they appear.

For any continuous, compactly supported function $f\colon\cx\to\reals$, the Siegel-Veech transform
$\hat{f}\colon\Omega_1\moduli[g]({\bf n})\to\reals$ is defined by
\begin{equation*}
  \hat{f}(X, \omega) = \sum_{v\in V(X, \omega)} f(v).
\end{equation*}

\begin{theorem}[\cite{veech98}]
  Given any ergodic, $\SLtwoR$-invariant measure $\mu$ on $\Omega_1\moduli[g]({\bf n})$, we have
  \begin{equation}
    \label{eq:svtransform}
    \frac{1}{\vol \mu} \int_{\Omega_1\moduli[g]({\bf n})} \hat{f} \,d\mu= c_{\rm sc}(\mu) \int_{\reals^2} f.
  \end{equation}
  The Siegel-Veech constant $c_{\rm sc}(\mu)$ depends only on the
measure $\mu$.
\end{theorem}

If $\mu$ is supported on a Teichm\"uller curve $C$, we will also write $c_{\rm sc}(C)$ for $c_{\rm sc}(\mu)$.
It follows from \cite{veech98} that for any $(X, \omega)$ on a Teichm\"uller curve $C$,
\begin{equation*}
  N_{\rm sc}((X, \omega), L)\sim c_{\rm sc}(C) \frac{\pi}{\Area(X, \omega)} L^2,
\end{equation*}
where $N_{\rm sc}$ is the counting function for saddle connections.

\paragraph{Bounded Siegel-Veech constants.}

We now show that Siegel-Veech constants are uniformly bounded over Teichm\"uller curves of a given genus.

\begin{theorem}
  \label{thm:svbounded}
  There is a uniform bound,
  \begin{equation*}
    c_{\rm sc}(C) \leq D,
  \end{equation*}
  where $C$ ranges over all Teichm\"uller curves in $\Omega\moduli[g]$, and $D$ only depends on the genus.
\end{theorem}

This follows easily from the following results.

\begin{theorem}[{\cite{masur90}, \cite{eskinmasur}}]
  \label{thm:masurbound}
  Given any compact subset $K\subset\Omega\moduli[g]({\bf n})$, there is a constant $c(K)$ such that for any $(X,
  \omega)\in K$,
  $$N_{\rm sc}((X, \omega), L) < c(K) L^2.$$
\end{theorem}

\begin{theorem}
  \label{thm:recurrence}
  Given any connected component $S$ of a stratum $\Omega_1\moduli[g]({\bf n})$, there is a compact subset
  $K\subset S$ which intersects the $\SLtwoR$ orbit of each point in $S$.
\end{theorem}

Theorem~\ref{thm:recurrence} is a corollary of the main result in Athreya's thesis \cite{athreya06}.

\begin{proof}[Proof of Theorem~\ref{thm:svbounded}]
  For each connected component $S$ of a stratum, choose a $K$ as in Theorem~\ref{thm:recurrence}.  Then for every
  Teichm\"uller curve $C\subset S$, we have $c_{\rm sc}(C)< c(K)$ by Theorem~\ref{thm:masurbound}.  Since
  there are only finitely many connected components of strata in each genus, we obtain a bound depending only
  on the genus.
\end{proof}

Let $K_\epsilon\subset \Omega_1\moduli[g]({\bf n})$ be the locus of translation surfaces which have no saddle
connection of length less than $\epsilon$, and let $\tilde{K}_\epsilon$ be the complement of $K_\epsilon$.

\begin{cor}
  \label{cor:cuspvolbound}
  For any Teichm\"uller curve $C\subset \Omega_1\moduli[g]({\bf n})$, we have
  \begin{equation*}
    \frac{\vol(C\cap \tilde{K}_\epsilon)}{\vol(C)}< D\epsilon^2,
  \end{equation*}
  for some constant $D$ depending only on the genus $g$.
\end{cor}

\begin{proof}
  Let $\chi_\epsilon$ be the indicator function of the disk of radius $\epsilon$.  We have
  $$\hat{\chi}_\epsilon\geq\chi_{\tilde{K}_\epsilon}.$$
  Even though $\chi_\epsilon$ is not continuous, we can apply \eqref{eq:svtransform} by the monotone
  convergence theorem.  We have
  \begin{equation*}
    \frac{\vol(C\cap \tilde{K}_\epsilon)}{\vol(C)} \leq \frac{1}{\vol(C)}\int_C \hat{X}_\epsilon = c_{\rm
      sc}(C) \pi \epsilon^2 \leq D \pi \epsilon^2,
  \end{equation*}
  where $D$ is the constant of Theorem~\ref{thm:svbounded}.
\end{proof}

\paragraph{Extension of $\F$ to $\Y$.}

We equipped the foliation $\F$ of $\X$ by $\SLtwoR$-orbits with a transverse
invariant measure $\proj \mu'_{D}$, so we can integrate 2-forms over $\F$
to obtain a closed current $\langle\F\rangle$ on $\X$.  The foliation $\F$ does not
extend to a foliation of $\Y$.  In fact, $\F$ has isolated
singularities at the points $c_P$ as well as along the one-dimensional
loci in the curves $C_P$ consisting of forms whose horizontal
foliation has two cylinders.  Nevertheless, the current defined by
$\F$ extends to a closed current on $\Y$.

\begin{theorem}
  \label{thm:Fextends}
  The foliation $\F$ defines a closed current $\langle\F\rangle$ on
  $\Y$, defined by
  $$\langle\F\rangle(\omega) = \int_{\F}\omega$$
  for each $C^\infty$ 2-form $\omega$ on $\Y$.
\end{theorem}

We will defer the proofs of this theorem and the next to a future
paper.  The idea of the proof of Theorem~\ref{thm:Fextends} is to
compare $\int_{\F}\omega$ with $\vol(\mu'_{D})$ using the Schwartz
lemma, and then to apply Theorem~\ref{thm:volfinite} that these
volumes are finite.

Let $[\F]$ be the cohomology class in $H^2(\Y; \reals)$ defined by
$\langle\F\rangle$.

\begin{theorem}
  \label{thm:Fintersectionzero}
  For any component $C$ of $\barW$, $\barP$, or $\barStwo[D]$, we have
  \begin{equation*}
    [C]\cdot[\F] = 0.
  \end{equation*}
\end{theorem}

The conclusion of this theorem is true for any closed leaves of a
measured foliation which are not atoms of the transverse measure.  The
curves $\W$, $\P$. and $\Stwo[D]$ are leaves of $\F$, but their cusps
pass through singular points of the extension of $\F$ to $\Y$.  The proof of
Theorem~\ref{thm:Fintersectionzero} amounts to  estimating
integrals of smooth forms around these singularities.

\paragraph{Proof of Theorem~\ref{thm:kz2}.}

Consider the pullbacks of the bundles $E$ and $L$ over $\proj\Omega\moduli$ by
$\pi\colon\X\to\proj\Omega\moduli$.  We have
\begin{align*}
  \pi^*E &= \Omega^1\X\oplus \Omega^2\X, \quad \text{and} \\
  \pi^*L &= \Omega^1\X.
\end{align*}
Proposition~\ref{prop:chHQ} implies
\begin{align*}
  \pi^*\gamma_1 &= \omega_1 + \omega_2, \quad \text{and}\\
  \pi^*\gamma_2 &= \omega_1.
\end{align*}
Let $\mu$ be an ergodic $\SLtwoR$-invariant measure on $\Omega_1\moduli$ whose support is not one of the
strata $\Omega_1\moduli(2)$ or $\Omega_1\moduli(1,1)$.  By McMullen's classification of ergodic, invariant
measures, each such $\mu$ is supported on one of the eigenform loci $\Omega_1\E$.  Then $\proj\mu$ is a
transverse, invariant measure to the foliation $\F$ of $\X$, and
\eqref{eq:kontsevichformula} becomes
\begin{equation}
  \label{eq:lambda2mu}
  \lambda_2(\mu) = \frac{\int_{\proj\mu}\omega_2}{\int_{\proj\mu}\omega_1}.
\end{equation}

If the ergodic, invariant measure $\mu$ is supported on a connected component $C$ of $\W$, then $\lambda_2(\mu) = 1/3$ by
\eqref{eq:lambda2mu} and Corollary~\ref{cor:intomega2}.

Now suppose $\mu$ is supported on a Teichm\"uller curve $C\subset\proj\Omega\moduli(1,1)$ (so $C$ is either
the decagon curve $D_{10}$ or a component of $\W[d^2][n]$).  By Corollary~\ref{cor:funWD}, we have
\begin{equation}
  \label{eq:wpsrelation}
  [\barW] - [\barP] + [\barStwo[D]] = -[\omega_1] + 2[\omega_2].
\end{equation}
Since $\overline{C}$ is disjoint from each of the three curves on the left, pairing $[\overline{C}]$ with both
sides of \eqref{eq:wpsrelation} yields $\lambda_2(\mu)= 1/2$.  If $\mu$ is the period measure $\mu'_D$ on $\E$, then
$$\int_{\proj\mu'_D}\omega_i = [\F](\omega_i).$$
Then we obtain $\lambda_2(\mu'_D)=1/2$ by pairing $[\F]$ with both sides of \eqref{eq:wpsrelation} and applying
Theorem~\ref{thm:Fintersectionzero}.

It remains to deal with the measures $\mu'({\bf n})$ on the strata $\Omega_1\moduli({\bf n})$.  Given a point
$p\in\proj\Omega\moduli$, let $L$ be the leaf of $\mathcal{F}$ through $p$.  Define
\begin{equation*}
  f(p) = \frac{\gamma_1|_L}{\gamma_2|_L}.
\end{equation*}
This is a well-defined real number because it is the ratio of top-degree forms.

\begin{lemma}
  As a function on $\proj\Omega\moduli$,  we have
  \begin{equation*}
    1 \leq f \leq 2.
  \end{equation*}
\end{lemma}

\begin{proof}
  Since $f$ is continuous, it suffices to prove this bound for every Teichm\"uller curve $C\subset
  \proj\Omega\moduli$.  The curve $C$ lies on some Hilbert modular surface $\X$.  Let
  $$\tilde{C}\subset \half\times\half$$
  be a connected component of the inverse image of $C$ in the universal cover of $\X$, and let $\pi_i\colon
  \tilde{C} \to \half$ be the two projections.  Since $C$ is transverse to the vertical foliation $\A$ of
  $\X$, the first projection $\pi_1$ is a conformal isomorphism of $\tilde{C}$ with $\half$.  By the Schwartz
  lemma, the projection $\pi_2$ is a contraction.  Since $\omega_i$ is the pullback of the hyperbolic area
  form on $\half$ by $\omega_i$, we have
  \begin{equation*}
    0\leq \frac{\omega_2|_C}{\omega_1|_C}\leq 1.
  \end{equation*}
  The claim follows, since
  \begin{equation*}
    f|_C = \frac{\omega_1|_C + \omega_2|_C}{\omega_1|_C}.
  \end{equation*}
\end{proof}

We need to show that
\begin{equation*}
  \frac{1}{\vol\mu'({\bf n})} \int_{\Omega_1\moduli({\bf n})} f \, d\mu'({\bf n}) =
  \begin{cases}
    1/3; & \text{if ${\bf n} = (2)$};\\
    1/2; & \text{if ${\bf n} = (1,1)$},
  \end{cases}
\end{equation*}
where we regard $f$ as a function on $\Omega_1\moduli({\bf n})$ by pulling it back from the projectivization.
Let $\{\nu_m\}$ be a sequence of measures supported on Teichm\"uller curves in $\Omega_1\moduli({\bf n})$ and
converging to $\mu'({\bf n})$.  Since the average of $f$ over any Teichm\"uller curve $C$ is either $1/2$ or
$1/3$, depending on the stratum in which $C$ lies, it is enough to show that
\begin{equation*}
  \frac{1}{\vol \nu_m}\int f \, d\nu_m \to \frac{1}{\vol\mu'({\bf n})} \int f \, d\mu'({\bf n}).
\end{equation*}

Recall that we defined $K_\epsilon$ to be the locus of translation surfaces with no saddle connection shorter
than $\epsilon$.  Since $K_\epsilon$ is compact  (as is shown in
\cite{kms}), $\vol \Omega\moduli({\bf n}) <\infty$,  and the $K_\epsilon$ exhaust
$\Omega_1\moduli({\bf n})$ as $\epsilon\to 0$, we can choose for any $\delta>0$ an  $\epsilon$ small enough
that
\begin{equation}
  \label{eq:bound1}
  \frac{\mu'({\bf n})(\tilde{K}_\epsilon)}{\vol\mu'({\bf n})}< \frac{\delta}{2}.
\end{equation}
By Corollary~\ref{cor:cuspvolbound}, we can also choose $\epsilon$ small enough that
\begin{equation}
  \label{eq:bound2}
  \frac{\nu_m({\bf n})(\tilde{K}_\epsilon)}{\vol\nu_m({\bf n})}< \frac{\delta}{2}.
\end{equation}

Let $g_\epsilon$ be a continuous, compactly supported function on $\Omega_1\moduli({\bf n})$ such that $g\equiv 1$ on $K_\epsilon$, and
$0\leq g \leq 1$.  We have,
\begin{align*}
  & \lim_{m\to\infty} \left| \frac{1}{\vol \nu_m}\int f\,d\nu_m - \frac{1}{\vol\mu'({\bf n})}\int f \,
    d\mu'(\bf n) \right| \\
  & \leq \lim_{m\to\infty} \left|\frac{1}{\vol \nu_m}\int g_\epsilon f\,d\nu_m - \frac{1}{\vol\mu'({\bf n})}\int
    g_\epsilon f \,
    d\mu'(\bf n) \right| \\  &  \hfill \hspace{2.5in} +  \frac{\mu'({\bf
    n})(\tilde{K}_\epsilon)}{\vol\mu'({\bf n})}  + \frac{\nu_m({\bf
      n})(\tilde{K}_\epsilon)}{\vol\nu_\epsilon({\bf n})} \\
  & < \delta,&
\end{align*}
since the limit of the first term is zero by Theorem~\ref{thm:uniformdistribution}, and the other two terms
are bounded by \eqref{eq:bound1} and \eqref{eq:bound2}.

This completes the proof of Theorem~\ref{thm:kz2}.  Note that the last part of the proof applies to the
$\Omega_1\E[d^2]$, so we can prove Theorem~\ref{thm:kz2} for these spaces without appealing to
Theorems~\ref{thm:Fextends} and \ref{thm:Fintersectionzero}.


%% file: normalization.tex
\section{Normal  varieties}
\label{sec:normal}

In this section, we record standard facts about normality for algebraic varieties and analytic spaces which
we use in this paper.  We will consider all of our algebraic varieties to be over $\cx$.

\paragraph{Normal algebraic varieties.}

A point $p$ on an algebraic variety $X$ is said to be \emph{normal} if the local ring $\mathcal{O}_p$ of $X$
at $p$ is an integral domain which is integrally closed in its field of fractions.  A variety $X$ is \emph{normal} if
it is normal at each of its points.

A \emph{normalization} of a variety $X$ is a normal variety $Y$ together with a finite surjective morphism
$\pi\colon Y\to X$.  More generally (following \cite{mumfordred}), if $X$ is an irreducible variety, and if $L$
is a finite algebraic extension of $K(X)$, then a \emph{normalization of $X$ in $L$} is a normal variety $Y$ with
function field $K(Y)=L$, together with a finite surjective morphism $\pi\colon Y\to X$ such that $\pi^*\colon
K(X)\to K(Y)=L$ is the given inclusion of $K(X)$ in $L$.  If $L=K(X)$, this is just the usual normalization of
$X$.

\begin{theorem}[{\cite[Theorem III:8.3]{mumfordred}}]
  \label{thm:normalizationexists}
  For any irreducible variety $X$ and finite algebraic extension $L$ of $K(X)$, there is a normalization of
  $X$ in $L$, and any two such normalizations are equivalent.
\end{theorem}

Normalization is also characterized by a universal property.  In the case, where $L = K(X)$, this is an
exercise in \cite[p.91]{hartshorne}.

\begin{theorem}
  \label{thm:universalproperty}
  The normalization $\pi\colon Y \to X$ of $X$ in $L=K(Y)$ has the following universal property.  Let $Z$ be a
  normal variety with a dominant morphism $q\colon Z\to X$, and let $r\colon K(Y)\to K(Z)$ be an isomorphism
  such that $r\circ \pi^* = q^*$ as inclusions $K(X)\to K(Z)$.  Then there is a unique morphism $s\colon Z\to Y$
  such that $s^* = r$ and $\pi\circ s =q$.
\end{theorem}

We will also have use for the following theorems:

\begin{theorem}[{\cite[III:8.4]{mumfordred}}]
  \label{thm:normalizationprojective}
  The normalization of a projective variety in a finite algebraic extension is projective.
\end{theorem}

\begin{theorem}[Zariski's Main Theorem, \cite{hartshorne}]
  \label{thm:zariskismain}
  Let $f\colon X\to Y$ be a birational morphism of projective varieties, and assume that $Y$ is normal.  Then
  for every $y\in Y$, the fiber $f^{-1}(y)$ is connected.
\end{theorem}

\begin{theorem}[{Zariski, \cite[Theorem~VIII.32]{zs}}]
  \label{thm:completionintegrallyclosed}
  If $X$ is a normal variety, $p\in X$, and $\widehat{\mathcal{O}}_p$ is the completion of $\mathcal{O}_p$, then
  $\widehat{\mathcal{O}}_p$ is an integrally closed integral domain.  
\end{theorem}

\paragraph{Normal analytic spaces.}

A complex analytic space is a ringed space which is locally modeled on analytic subvarieties of $\cx^n$.  We
start by recalling some basic notions about complex analytic spaces.  See \cite{gunningrossi} and
\cite{gunning70} for details.

Consider the set of pairs $(U, X)$, where $U\subset\cx^n$ is a neighborhood of some point $p$, and $X\subset U$ is
an analytic subvariety of $U$.  We consider two such pairs to be equivalent if the varieties agree on some
common neighborhood of $p$.  An equivalence class of such pairs is called a \emph{germ on an analytic variety
  at $p$}.  A germ is considered to be irreducible if it contains no proper subgerms.

\begin{theorem}[{\cite[p.11]{gunning70}}]
  Any germ of an analytic subvariety at $p$ can be written uniquely as the union of finitely many irreducible
  germs of analytic varieties at $p$.
\end{theorem}

Given a point $p\in X$, an analytic subvariety of an open subset of $\cx^n$, a \emph{branch} of $X$ through
$p$ is an irreducible subgerm of the germ of $X$ at $p$.  Informally, a branch of $X$ through $p$ is a
connected component of $U\cap (X\setminus X_{\rm sing})$, where $U$ is a small and sufficiently regular
neighborhood of $p$, and $X_{\rm sing}$ is the singular set of $X$.  Given a point $p$ in an analytic space
$X$, let $\mathcal{O}_p$ be the local ring of germs of holomorphic functions on $X$ at $p$.  Just as for
algebraic varieties, we say that $p$ is a \emph{normal point} of $X$ if $\mathcal{O}_p$ is an integrally
closed integral domain.  We say that $X$ is \emph{normal} if each of its points is normal.

There is also a geometric characterization of normality.  Following \cite{whitney}, we say that a function $f$
on an open subset $V$ of an analytic variety $X$ is \emph{weakly holomorphic} if there is a subvariety
$V'\subset V$ with the following properties:
\begin{itemize}
\item $V'$ is nowhere dense in $V$, and $V_{\rm sing}\subset V'$.
\item $f$ is holomorphic in $V\setminus V'$.
\item $f$ is locally bounded in $V$.
\end{itemize}
It follows from \cite[Theorem~4.10I]{whitney} that an analytic space $X$ is weakly holomorphic at $p$ if and
only if every germ of a weakly holomorphic function at $p$ is in fact holomorphic.

A \emph{normalization} of an analytic space $X$ is a
normal analytic variety $Y$ together with a holomorphic map $p\colon Y\to X$ such that
\begin{itemize}
\item $p\colon Y\to X$ is proper and has finite fibers.
\item If $X_{\rm sing}$ is the singular set of $X$, and $A=p^{-1}(X_{\rm sing})$, then $Y\setminus A$ is dense
  in $Y$, and $p\colon Y\setminus A\to X$ is biholomorphic onto $X\setminus X_{\rm sing}$.
\end{itemize}

\begin{theorem}[{\cite[p.38]{laufer}}]
  \label{thm:analyticunique}
  Every analytic space has a unique normalization.  
\end{theorem}

\begin{theorem}[{\cite[p.37]{laufer}}]
  \label{thm:pointsarebranches}
  If $\pi\colon Y\to X$ is the normalization of $X$, then for each $p\in X$, the map $\pi$ induces a bijective
  correspondence between the points of $\pi^{-1}(p)$ and the branches of $X$ through $P$.
\end{theorem}

\begin{theorem}
  \label{thm:algebraicnormal}
  If $X$ is an algebraic variety, and $\pi\colon Y\to X$ is the normalization of $X$ as an algebraic variety
  then it is also the normalization of $X$ as an analytic space.  
\end{theorem}

\begin{proof}
  By the definition of the normalization of an analytic space, we need only to show that $Y$ is normal as an
  analytic space.  For this proof, write $\mathcal{O}_p^{\rm an}$ for the local ring of germs of holomorphic
  functions on $X$ at $p$, to differentiate from the local ring of $\mathcal{O}_p$ of algebraic functions at
  $p$.  Let $\widehat{\mathcal{O}}_p$ be the completion of $\mathcal{O}_p$.
  
  Let $f$ in the quotient field $K(\mathcal{O}_p^{\rm an})$ be integral over $\mathcal{O}_p^{\rm an}$.  Since
  $\mathcal{O}_p^{\rm an}\subset \widehat{\mathcal{O}}_p$, it follows that $f$ is integral over
  $\widehat{\mathcal{O}}_p$.  Since $\widehat{\mathcal{O}}_p$ is integrally closed by
  Theorem~\ref{thm:completionintegrallyclosed}, it follows that $f\in\widehat{\mathcal{O}}_p$.  Then
  $f\in\mathcal{O}_p^{\rm an}$ because $K(\mathcal{O}_p^{\rm an})\cap \widehat{\mathcal{O}}_p = \mathcal{O}_p^{\rm
    an}$ in $K(\widehat{\mathcal{O}}_p)$.
\end{proof}

A \emph{resolution } of an analytic space $X$ is a complex manifold $Y$ together with a proper analytic map
$\pi\colon Y\to X$ such that $\pi\colon Y\setminus\pi^{-1}(X_{\rm sing})\to X\setminus X_{\rm sing}$ is biholomorphic
and such that $\pi^{-1}(X\setminus X_{\rm sing})$ is dense in $Y$.  The following was proved by Mumford in
\cite{mumford61}, by
Grauert in \cite{grauert}, and also by Laufer in \cite[Theorem~4.4]{laufer}.  

\begin{theorem}
  \label{thm:resolutionnegativedefinite}
  Let $\pi\colon Y\to X$ be a resolution of $X$, and let $p\in X$ be an isolated normal point.  Suppose that
  $A=\pi^{-1}(p)$ is the union of irreducible curves $A_i$ which are nonsingular and intersect transversely,
  with at most two branches of $\bigcup_i A_i$ passing through any point.  Then the intersection matrix
  $$(A_i\cdot A_j)$$
  is negative definite.
\end{theorem}
